%% file: Main_AbsorbingConditions.tex
\documentclass[leqno]{amsart}
\usepackage{paralist}
\usepackage{stmaryrd}
\usepackage{amsmath}
\usepackage{verbatim}
\usepackage{multirow}
\usepackage{pgfplots}
\pgfplotsset{compat=1.14}
\usepackage{tikz}
\usepackage{caption}
\usepackage{placeins}
\usepackage[linesnumbered,ruled]{algorithm2e}
\usepackage{setspace}
\usepackage[noend]{algpseudocode}
\usepackage{tabularx,ragged2e,booktabs,caption}
\usepackage{graphics} 
\usepackage{epsfig} 
\usepackage{graphicx}  \usepackage{epstopdf}
\usepackage[colorlinks=true]{hyperref}
\usepackage{bm}
\hypersetup{urlcolor=blue, citecolor=red}
\usepackage{gensymb}
\usepackage[margin=1.15in]{geometry}
\newtheorem{remark}{Remark}

\newtheorem{definition}{Definition}
\title[Self-adaptive absorbing boundary conditions] {Self-adaptive absorbing boundary conditions for quasilinear acoustic wave propagation}
\author[Markus Muhr, Vanja Nikoli\' c, and Barbara Wohlmuth]{}
\email{\href{mailto:muhr@ma.tum.de}{muhr@ma.tum.de}}
\email{\href{mailto:vanja.nikolic@ma.tum.de}{vanja.nikolic@ma.tum.de}}
\email{\href{mailto:wohlmuth@ma.tum.de}{wohlmuth@ma.tum.de}}

\thanks{$^*$Corresponding author: Markus Muhr, \href{mailto:muhr@ma.tum.de}{muhr@ma.tum.de}}
\usepackage[utf8]{inputenc}

\usepackage{cleveref}
\crefname{section}{§}{§§}
\Crefname{section}{§}{§§}

\crefformat{section}{§#2#1#3}
\begin{document} 
\maketitle
\centerline{\scshape Markus Muhr$^*$, Vanja Nikoli\' c, and Barbara Wohlmuth}

\vspace{4mm}
{\footnotesize
   \centerline{Technical University of Munich, Department of Mathematics, Chair of Numerical Mathematics}
   \centerline{Boltzmannstra\ss e 3, 85748 Garching, Germany}
}

\vspace{6mm}
\begin{abstract}
We propose a self-adaptive absorbing technique for quasilinear ultrasound waves in two- and three-dimensional computational domains. As a model for the nonlinear ultrasound propagation in thermoviscous fluids, we employ Westervelt's wave equation solved for the acoustic velocity potential. The angle of incidence of the wave is computed based on the information provided by the wave-field gradient which is readily available in the finite element framework. The absorbing boundary conditions are then updated with the angle values in real time. Numerical experiments illustrate the accuracy and efficiency of the proposed method. 
\end{abstract}
\vspace{8mm}
\section{Introduction}
Accurate simulation of nonlinear ultrasound offers a path to a better quality of many procedures in industry and medicine, from non-destructive detection of material damages~\cite{fierro2015nonlinear, jhang2009nonlinear, matlack2012evaluation} to non-invasive treatments of medical disorders~\cite{illing2005safety,kennedy2004high, kennedy2005high, wu2003preliminary, yoshizawa2009high}. When studying such procedures, there is always a region of interest: a kidney stone that will be disintegrated or a propagating fatigue crack in a component of an aircraft. The large physical space then often has to be truncated for numerical simulations. To accurately simulate ultrasound, we have to avoid spurious reflections of the wave at the boundary of the truncated domain.  \\
\indent Absorbing boundary conditions provide a simple and effective way of dealing with unwanted reflections. They were introduced by B. Engquist and E. Majda in their seminal work \cite{EnquistMajda}. Since then many approaches have been developed for the non-reflecting boundary conditions; we refer the reader to the review papers~\cite{givoli1991non,nataf2013absorbing} and the references given therein. In spite of such comprehensive research in this area, only a small portion of the results focus on nonlinear models.\\
\indent A class of semilinear wave equations and nonlinear Schr\"odinger equations was investigated in~\cite{szeftel2006absorbing, szeftel2006nonlinear}. Results for nonlinear hyperbolic systems of the form $u_t+A(u)u_x=0$ were obtained in~\cite{hedstrom1979nonreflecting}. An approach based on the operator splitting method was used in~\cite{li2011local} to derive absorbing conditions for a semilinear wave equation of the form $u_{tt}-a^2\Delta u=f(u)$. In~\cite{ShevchenkoKaltenbacher, shevchenko2015absorbing}, nonlinear ultrasound propagation was investigated in this context for the first time, and absorbing conditions were developed for the Westervelt equation in the pressure form. \\
\indent Another commonly used approach for avoiding spurious reflections is the Perfectly Matched Layer (PML) technique. Developed by J.-P. B\' erenger in~\cite{berenger1994perfectly}, this method introduces an artificial absorbing layer around the computational domain. Linear acoustic wave equations have been extensively studied in this context; see, for example,~\cite{abarbanel1999well, hein2007acoustic, hesthaven1998analysis, kaltenbacher2013modified, qi1998evaluation}. \\ 
\indent The main downside of absorbing conditions is that they are sensitive to the angle of incidence of the wave. As a rule of thumb, they perform bad if the range of incidence angles is large. The information on the incidence angles can be included in the conditions to tackle this issue; we refer to the work in~\cite{higdon1986absorbing, higdon1987numerical, keys1985absorbing}. However, these angles are not a priori known in a realistic computational setting. The idea behind the self-adaptive technique is to compute the local wave vector and then update the absorbing conditions with the angle information \emph{in real time}. This approach was applied in~\cite{xu2007adaptive} to the Schr\"odinger-type equations, where the wave number was computed via the Gabor transform. In~\cite{shevchenko2012self}, the linear wave equation was investigated in this context. It was proposed to divide the absorbing boundary into segments and compute the local incidence angle by employing the Fourier transformation only in the vicinity of the boundary. The self-adaptive approach to absorption can also be found in earlier works on acoustic scattering~\cite{jin1996application, li1994adaptive}.  \\
\indent The goal of our work is to develop an efficient self-adaptive absorbing technique for nonlinear ultrasound propagation. As a model equation, we employ a classical quasilinear acoustic model - Westervelt's equation. We first extend the results from~\cite{ShevchenkoKaltenbacher, shevchenko2015absorbing} by considering the potential form of Westervelt's equation and the non-zero angle of incidence $\theta$. In addition, we derive the absorbing conditions for two- and three-dimensional computational domains. The derivation relies on choosing an appropriate linearization of the equation around a reference solution. The absorbing conditions are then formally derived for this linearization, after which we bring back the nonlinear term.  \\
\indent To obtain the angle $\theta$ in practice, we develop a self-adaptive method that locally computes the incidence angle and updates the absorbing conditions on the fly. Unlike the self-adaptive approach taken in~\cite{shevchenko2012self}, we base the local angle computation on the gradient of the wave field.  Computation of the local propagation direction in isotropic media based on the information provided by the wave-field gradient has been investigated in~\cite{jia2009calculation, patrikeeva2013comparison, Poynting1, ye2012effective, Poynting2, zhang2011direct}. This approach is particularly suitable for our finite-element framework since the gradient information is already available at every time step in our simulations. The use of the field gradient information in the absorbing conditions was already investigated for the Helmholtz equation in~\cite{gordon2015compact}. There it was proposed to replace the normal derivatives that appear in the absorbing conditions by the derivatives in the direction of the wave propagation. In the linear regime, our approach can be understood as an extension of the gradient method in~\cite{gordon2015compact} for a time-dependent wave model.  \\
\indent We organize the rest of the paper as follows. We begin in Section~\ref{Section:Modeling} by introducing the model and setting the problem. Section~\ref{Section:Derivation} contains the derivation of absorbing conditions for a given angle of incidence of the wave. In Section~\ref{Section:NumTreatment}, we present the numerical scheme for solving the initial-boundary value problem for the Westervelt equation. Section~\ref{Section:Angle} describes the computation of the local incidence angle via the information provided by the wave field gradient. Finally, in Section~\ref{Section:NumResults}, we present numerical experiments which illustrate the accuracy of the proposed adaptive boundary conditions.
\section{Modeling and problem setting} \label{Section:Modeling}
The weakly nonlinear models in thermoviscous acoustics that are commonly used are obtained as an approximation of the compressible Navier-Stokes system. We here briefly reflect upon the derivation, which will give us a better understanding of the often-employed Westervelt equation. Mathematically rigorous justification of the classical acoustic models can be found in~\cite{kaltenbacher2017fundamental}. For a detailed insight into the acoustic field theory, we refer to~\cite{Crighton, hamilton1998nonlinear,MKaltenbacher}.\\
\indent  Propagation of waves can be described by the time, the density $\bar{\varrho}$, the pressure $\bar{u}$ and the velocity $\bar{\mathbf{v}}$, decomposed into their ambient value and the acoustic perturbation
\begin{equation*}
\begin{aligned}
& \bar{\varrho}=\varrho_0+\varrho, \\
& \bar{u}=u_0+u, \\
& \bar{\mathbf{v}}=\mathbf{v}_0+\mathbf{v};
\end{aligned}
\end{equation*}
see~\cite{MKaltenbacher}. We call $\varrho$ the acoustic density, $u$ the acoustic pressure, and $\mathbf{v}$ the acoustic particle velocity. The equations governing the wave propagation are then given by
\begin{itemize}
\item the equation of momentum conservation
$$(\varrho_0+\varrho)\mathbf{v}_t+\frac{\varrho_0}{2}\nabla (\mathbf{v} \cdot \mathbf{v})+\nabla u=\left(\frac{4\nu_V}{3}+\eta_V \right)\Delta \mathbf{v},$$
\item the equation of mass conservation
$$\varrho_t+\varrho_0 \nabla \cdot \mathbf{v}=-\varrho \nabla \cdot \mathbf{v}-\mathbf{v}\cdot \nabla \varrho, $$
\item the pressure-density relation
$$\varrho=\frac{1}{c^2}u-\frac{1}{\varrho_0 c^4}\frac{B}{2A}u^2-\frac{\kappa}{\varrho_0 c^4}\left(\frac{1}{c_\Omega}-\frac{1}{c_u}\right)u_t.$$
\end{itemize}
Above, $\eta_V$ denotes the bulk viscosity and $\nu_V$ the shear viscosity. The constant $\kappa$ stands for the adiabatic exponent, $c_u$ and $c_\Omega$ denote the specific heat capacitance at constant pressure and constant volume, respectively. The parameter of nonlinearity $B/A$ is an indicator of the nonlinearity of the medium. Finally, $c$ denotes the speed of sound in the fluid. \\
\indent This system of equations is approximated by one model, whereby every term of order two and higher in the acoustic Mach number is neglected. This approach results in
\begin{equation} \label{Kuznetsov_pressure_velocity}
\begin{cases}
\frac{1}{c^2}u_{tt}-\displaystyle \Delta u-\frac{b}{c^2}\Delta u_t=\frac{1}{\varrho_0 c^4}\frac{B}{2A}u_{tt}+\frac{\varrho_0}{c^2}\frac{\partial^2 }{\partial t^2}(\mathbf{v} \cdot \mathbf{v}), \vspace{1mm}\\
\varrho_0 \mathbf{v}_t =-\nabla u,
\end{cases}
\end{equation}
where the so-called sound diffusivity $b$ is given by $$b=\frac{1}{\varrho_0}\left(\frac{4\nu_V}{3}+\eta_V\right)+\frac{\kappa}{\varrho_0}\left(\frac{1}{c_V}-\frac{1}{c_u}\right).$$
The acoustic velocity potential $\psi$ is then introduced to obtain a scalar equation; it is related to the acoustic pressure by
\begin{align}  \label{pressure_potential_relation} 
& u \approx \varrho_0 \psi_t,
\end{align}
and to the acoustic particle velocity by
\begin{align}
& \mathbf{v}=-\nabla \psi.
\end{align}
By expressing \eqref{Kuznetsov_pressure_velocity} in terms of $\psi$, integrating with respect to time and taking the resulting constant of integration to be zero, we arrive at the Kuznetsov equation
\begin{align} \label{Kuznetsov}
\frac{1}{c^2}\psi_{tt}-\Delta \psi-\frac{b}{c^2}\Delta \psi_t=\frac{B/A}{c^4}\psi_t \psi_{tt}+\frac{2}{c^2}\nabla \psi \cdot \nabla \psi_t;
\end{align}
cf.~\cite{kuznetsov}. If the cumulative nonlinear effects dominate the local ones in the sense of
\begin{align} \label{psi_approximation}
\frac{2}{c^2}\nabla \psi \cdot \nabla \psi_t \approx \frac{2}{c^2} \psi_t \psi_{tt},
\end{align}
a simplification of \eqref{Kuznetsov}, known as the Westervelt equation~\cite{westervelt1963parametric}, is obtained
\begin{align} \label{Westervelt}
\frac{1}{c^2}\psi_{tt}-\Delta \psi-\delta \Delta \psi_t=\frac{k}{c^2}\, \psi_t \psi_{tt}.
\end{align}
Above we have introduced the notation
\begin{align}
\delta= \frac{b}{c^2}, \ k= \frac{1}{c^2}(B/A+2).
\end{align}
After numerically solving \eqref{Westervelt}, the pressure field can be obtained in a post-processing step via the relation \eqref{pressure_potential_relation}. We mention as well that the Westervelt equation in the pressure form is given by
\begin{align} \label{Westervelt_pressure}
\frac{1}{c^2}u_{tt}-\Delta u-\delta \Delta u_t=\frac{k}{\varrho c^2}\, (u u_{tt}+u_t^2).
\end{align}
Equation \eqref{Westervelt_pressure} can be obtained from \eqref{Kuznetsov_pressure_velocity} by employing the approximation \eqref{psi_approximation} which in terms of the velocity and pressure reads as $$\frac{\varrho_0}{2}\textbf{v} \cdot \textbf{v} \approx \frac{u^2}{2 \varrho_0 c^2}.$$
\newpage
\section{Absorbing conditions for the Westervelt equation}\label{Section:Derivation}
We consider the Westervelt equation \eqref{Westervelt} on a three-dimensional spatial domain $\Omega=\{(x,y, z): x<0, \, y, z \in \mathbb{R}\}$. We restrict ourselves in the problem description to the case of constant coefficients $c$, $b$, and $k$
\subsection{Linearization of the Westervelt equation}   
Following the approach from~\cite{ShevchenkoKaltenbacher} where equation \eqref{Westervelt_pressure} was considered, we derive absorbing conditions for the Westervelt equation in the potential form \eqref{Westervelt} by first transforming it into a linear equation, deriving the non-reflecting conditions for this linear model and then ``plugging'' back in the nonlinear term. \\
\indent It is clear that the linearization of the equation plays a crucial role in deriving the absorbing conditions. A possible linearization of the Westervelt equation around a reference solution $\psi^{\textup{ref}}$ is given by
\begin{align} \label{Westervelt_linearization_KS}
\displaystyle \frac{1}{c^2}\psi_{tt}-\Delta \psi-\delta \Delta \psi_t=\frac{k}{c^2}  \beta(x, y, z, t)\, \psi_{tt},
\end{align}
where $\beta=\psi_t^{\textup{ref}}$. We note that, unlike in the derivation of the Westervelt equation in Section~\ref{Section:Modeling}, here we do not split $\psi$ into a background and oscillatory part, instead we assume $\psi^{\textup{ref}}$ to be a solution of the equation. This linearization is analogous to the one employed in ~\cite{ShevchenkoKaltenbacher} for the equation in the pressure formulation \eqref{Westervelt_pressure}. \\
\indent We propose an alternative linearization. Note that the right-hand side of the Westervelt equation \eqref{Westervelt} can be rewritten as $\displaystyle \frac{1}{2}\frac{k}{c^2}(\psi_t^2)_t$. We linearize the term $\psi_t^2$ as $\psi_t^{\textup{ref}} \psi_t$ and study the following equation
\begin{align} \label{Westervelt_linearization_1/2}
\displaystyle \frac{1}{c^2}\psi_{tt}-\Delta \psi-\delta \Delta \psi_t=\frac{1}{2}\frac{k}{c^2}(\beta(x, y, z, t) \psi_{t})_t,
\end{align}
where again $\beta =\psi_t^{\textup{ref}}$. The absorbing boundary conditions based on linearizations \eqref{Westervelt_linearization_KS} and \eqref{Westervelt_linearization_1/2} are numerically compared in Section~\ref{Section:NumResults}. \\ 
\indent We first derive absorbing conditions for the linear equation \eqref{Westervelt_linearization_1/2} and a given angle of incidence. After obtaining the conditions for such a model, the coefficient $\beta$ is set back to $\psi_t$ to obtain nonlinear conditions.\\
\indent We remark that the linearization via Taylor expansion around a reference solution $\psi^{\textup{ref}}$ is not considered here since it would introduce the term $\psi_t^{\textup{ref}}\psi_{tt}^{\textup{ref}}$ into the linearized equation; we also refer to the discussion in~\cite{ShevchenkoKaltenbacher}. 
\subsection{Derivation of absorbing conditions for the angle of incidence \mathversion{bold}$\theta$} We study here the derivation of the absorbing conditions for the linearization \eqref{Westervelt_linearization_1/2}; equation \eqref{Westervelt_linearization_KS} can be treated analogously. To derive the conditions, we could use the frozen coefficient approach which first transforms the variable coefficient equation into its constant-coefficient counterpart by "freezing" its coefficients at a given point before employing the Fourier transform in the $(y,z,t)$ coordinates; cf.~\cite[Section 1]{EnquistMajda}. Although the main focus of the present work is the derivation of the zero-order (adaptive) absorbing conditions, we still follow the approach based on the pseudo-differential calculus since it allows to arrive at a general system for determining the correcting terms beyond order zero in the absorbing conditions; see system~\eqref{system_symbols} below.\\
\indent We first rewrite the linearized equation \eqref{Westervelt_linearization_1/2} in the operator form as
\begin{align} \label{operator_equation}
\mathcal{P}u=0,
\end{align}
where the operator $\mathcal{P}$ is given by
\begin{align} \label{operator_P}
\mathcal{P}=\left(\frac{1}{c^2}-\frac{1}{2}\frac{k}{c^2}\beta(x, y, z, t)\right)\, \partial^2_t-\partial_x^2-\partial_y^2-\partial_z^2-\delta \partial^2_{x} \partial_t-\delta \partial^2_{y} \partial_t-\delta \partial^2_{z} \partial_t-\frac{1}{2}\frac{k}{c^2} \beta_t(x, y, z, t) \partial_t.
\end{align}
At this point we also introduce
\begin{equation} \label{def_alpha}
\begin{aligned}
&    \alpha_0(x, y, z, t)=\sqrt{\frac{1}{c^2}-\frac{1}{2}\frac{k}{c^2} \beta(x, y, z, t)}, \\
&    \alpha_1(x, y, z, t)=\frac{1}{2}\frac{k}{c^2} \beta_t(x, y, z, t) .
\end{aligned}
\end{equation}
Note that the well-posedness results for the Westervelt equation rely on the fact that the factor $1-k \psi_t$ next to the second time derivative remains positive; see~\cite{kaltenbacher2009global,kaltenbacher2011well, kawashima1992global}. Therefore, it is reasonable to assume that the term under the square root in \eqref{def_alpha} is positive for sufficiently small data. We note that we proceed heuristically since the rigorous justification of the derivations given below would also require $C^\infty$ regularity of $\alpha_0$ and $\alpha_1$ which is not proven here.\\
\indent Absorbing boundary conditions for \eqref{operator_equation} can be derived by employing the pseudo-differential calculus and factorization of the operator $\mathcal{P}$ according to L. Nirenberg's procedure \cite[Chapter II]{Nirenberg}. We briefly summarize the procedure here for the convenience of the reader. A detailed account on the pseudo-differential operators can be found in \cite{Hoermander, Nirenberg, Wong}.  
\begin{definition}\cite{hormander1994analysis, Wong} Let the set $\mathcal{S}^m$, where $m \in \mathbb{R}$, be defined as the set of all functions $q(t, \tau) \in C^\infty(\mathbb{R}^d \times \mathbb{R}^d)$ such that for any two multi-indices $k$ and $l$, there is a positive constant $C_{kl}$ depending only on $k$ and $l$, such that 
\begin{align*}
|\partial_t^k \partial _\tau^l q(t, \tau)| \leq C_{kl}(1+|\tau|)^{m-|l|}, \quad t, \tau \in \mathbb{R}^d.
\end{align*}
$\mathcal{S}^m$ is called the space of symbols of order $m$. We set $S^{-\infty}=\displaystyle \cap_{m \in \mathbb{R}} \mathcal{S}^m$.
\end{definition}
\begin{definition}\cite[Definition 5.2]{Wong}
Let $q$ be a symbol. The pseudo-differential operator $\mathcal{Q}$ associated to $q$ is defined by
\begin{align*}
(\mathcal{Q} \varphi )(t)=(2 \pi)^{-d/2} \int_{\mathbb{R}^d} \, e^{it \cdot \tau}\, q(t, \tau)\, \mathcal{F} \varphi (\tau)\, d\tau, \quad  
\end{align*}
where $\varphi$ is a function from the Schwartz space, and  $\mathcal{F}$ denotes the Fourier transform.
\end{definition}
\subsubsection{Propagation without losses}
\begin{figure}[h!]
\begin{center}
\includegraphics[trim=0cm 0cm 0cm 0cm, clip, scale=0.5]{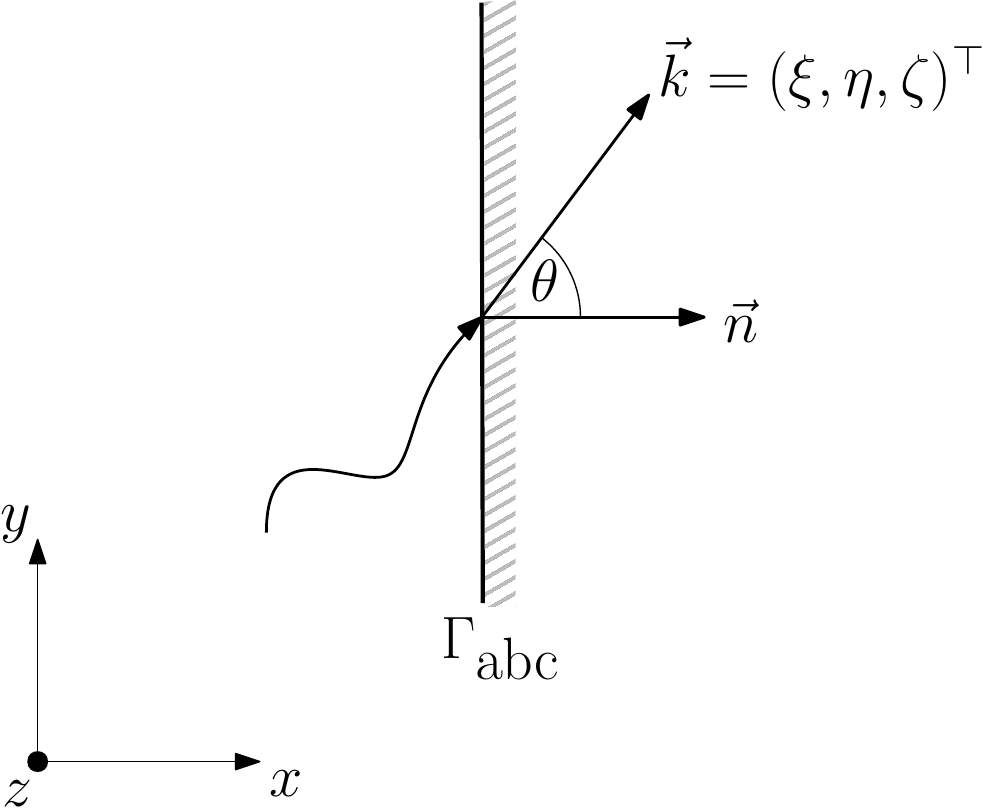}\hspace{1cm}\includegraphics[trim=0cm 0cm 0cm 0cm, clip, scale=0.4]{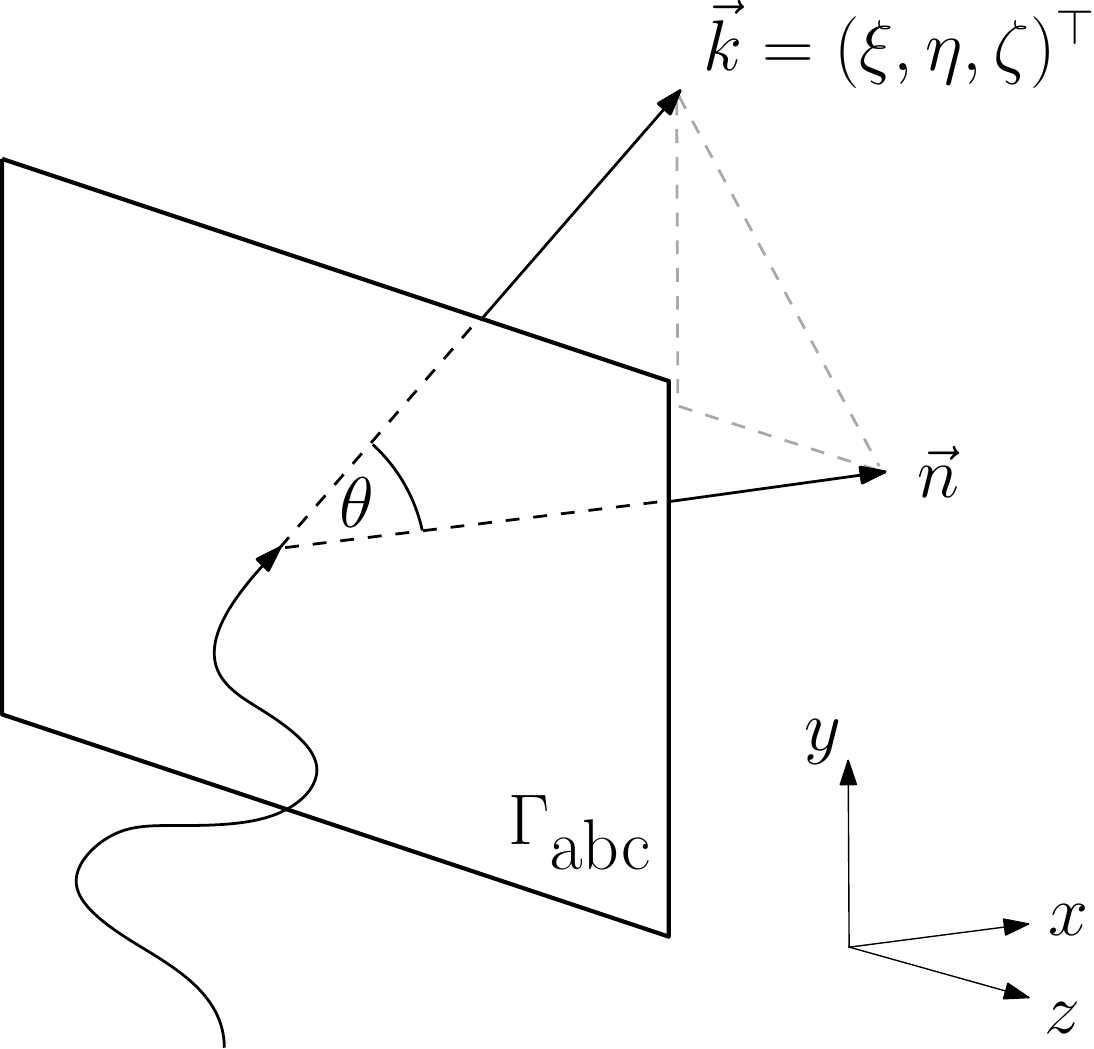}
\caption{Illustration of the interplay between the wave vector $\vec{k}$, outward normal vector $\vec{n}$, and the angle of incidence $\theta$. \label{fig:Boundary_sketch}}
\end{center}
\end{figure}

Following the approach in~\cite{ShevchenkoKaltenbacher}, we first derive the absorbing conditions with the assumption that $\delta=0$. This assumption facilitates the derivation of the conditions based on the pseudo-differential factorization. The $\delta$-term will be included as a post-processing step based on energy arguments. \\
\indent The derivation of the conditions relies on the fact that the operator $\mathcal{P}$ can be factorized into in the form 
\begin{align} \label{factorization_2D}
\mathcal{P}=-\displaystyle \left(\partial_x -\mathcal{A}\left(x, y, z, t, \partial_y, \partial_z, \partial_t\right )\right) \left(\partial_x-\mathcal{B} \left(x, y, z, t, \partial_y, \partial_z, \partial_t\right) \right)+\mathcal{R}\left(x, y, z, t, \partial_y, \partial_z, \partial_t\right);
\end{align}
see \cite[Lemma 1]{Nirenberg}. In \eqref{factorization_2D}, the operators $\mathcal{A}$ and $\mathcal{B}$ are pseudo-differential operators with the symbols $a(x,y,z,t, \eta,\zeta,\tau)$ and $b(x,y,z,t, \eta,\zeta, \tau)$, respectively, from the space $\mathcal{S}^1$. The conditions on $\mathcal{A}$ that we will develop will have the effect of associating $\mathcal{A}$ with waves that travel out of the computational domain. The pseudo-differential operator $\mathcal{R}$ is a smoothing operator with the full symbol $r(x,y,z,t, \eta,\zeta, \tau)$ that belongs to $\mathcal{S}^{-\infty}$. \\
\indent The symbols $a$ and $b$ formally admit asymptotic expansions
\begin{align*}
& a(x,y,z,t, \eta, \zeta, \tau) \sim \displaystyle \sum_{j \geq 0} a_{1-j}(x,y,z,t, \eta, \zeta, \tau), \\
& b(x,y,z,t, \eta, \zeta, \tau) \sim \displaystyle \sum_{j \geq 0} b_{1-j}(x,y,z,t, \eta, \zeta, \tau),
\end{align*}
where $a_{1-j}$ and $b_{1-j}$ denote homogeneous functions of degree $1-j$ with respect to $\tau$; see \cite[Theorem 5.10]{Wong}. We note that the symbol $a(x,y,z,t, \eta, \zeta, \tau)b(x,y,z,t, \eta, \zeta,\tau)$ of the product of $A$ and $B$ has an asymptotic expansion as well
\begin{align*}
a(x,y,z,t, \eta,\zeta,\tau)b(x,y,z,t, \eta, \zeta, \tau) \sim \displaystyle \sum_{j \geq 0, \atop {k+l+n=j, \atop k,l,n \geq 0}} \frac{(-1)^n}{n!} \partial^{n}_\tau a_{1-l}(x,y,z,t, \eta, \zeta, \tau)\partial^n_t b_{1-k}(x,y,z,t, \eta,\zeta,\tau);
\end{align*}
see \cite[Theorem 7.1]{Wong}.\\
\indent According to \cite[Theorems 1 and 2]{MajdaOsher}, absorbing boundary conditions on the boundary $x=0$ are given in the form 
\begin{align*}
\left(\partial_x-\displaystyle \mathcal{A}\left(x,y,z,t, \partial_y, \partial_z,\partial_t \right)\right) \psi \Bigl \vert_{x=0}=0.
\end{align*}
Since the symbol $a$ of $A$ has an infinite expansion, in numerical simulations this expansion is truncated after a certain number of terms. Absorbing conditions of order $k \in \mathbb{N}_0$ are then on the symbolic level given by
\begin{align} \label{abc_order_k}
\left(\partial_x-\displaystyle \sum_{j=0}^{k} a_{1-j}(0,y,z,t, \eta, \zeta, \tau)\right) \psi \Bigl \vert_{x=0}=0.
\end{align}
The higher-order absorbing conditions, although numerically more accurate, are also significantly more involved when it comes to implementation. We compute the absorbing conditions of order zero for the given angle of incidence $\theta$. Combined with the proposed self-adaptive technique, this approach allows to improve the accuracy of zero-order conditions, yet keeps them easy to implement.  \\
\indent We recall how the operator $\mathcal{P}$ was defined in \eqref{operator_P} and then develop factorization \eqref{factorization_2D} to obtain
\begin{equation} \label{eq_factored}
\begin{aligned}
    &\alpha_0^2 \, \partial^2_t-\partial_x^2-\partial_y^2-\partial_z^2-\alpha_1 \partial_t  \\
    =& \,-\partial_x^2+(\mathcal{A}+\mathcal{B})\partial_x+\mathcal{B}_x-\mathcal{A}\mathcal{B}+\mathcal{R}.
\end{aligned}
\end{equation}
By employing the asymptotic expansion of symbols $a$, $b$, and $ab$, equation \eqref{eq_factored} reduces on the symbolic level to
\begin{equation} \label{symbolic_level}
 \begin{aligned}
 &\alpha_0^2(i\tau)^2-(i \eta)^2-(i \zeta)^2-\alpha_1 (i \tau)\\
 \cong& \,\displaystyle \sum_{j \geq 0} (a_{1-j}+b_{1-j})\partial_x
 +\sum_{j \geq 0} \partial_x b_{1-j}-\sum_{{j \geq 0, \atop k+l+n=j,} \atop k,l,n \geq 0} \, \frac{(-1)^n}{n!}\partial_\tau^n a_{1-l} \partial_t^n b_{1-k},
 \end{aligned}
\end{equation}
 Above, we have denoted the dual variables to $t$, $y$, and $z$ by $\tau$, $\eta$, and $\zeta$, with the correspondence $\partial_t \leftrightarrow i \tau$, $\partial_y \leftrightarrow i \eta$, and $\partial_z \leftrightarrow i \zeta$.
Following the notation in~\cite{EnquistMajda}, $\cong$ stands for "within a smooth error" since we have dropped $\mathcal{R}$. We note that the operator $R$ can only be controlled in terms of its smoothness. However, its action on the solution is expected to be negligible for high frequencies that are present in ultrasound waves.\\
\indent To determine $a_1$ and $b_1$, we equate the symbols with the same degree of homogeneity with respect to $\tau$ and get the system 
\begin{align} \label{system_a_1_different}
 \begin{cases}
 a_1+b_1= 0 \\
 a_1b_1=-(\alpha_0^2(i\tau)^2-(i\eta)^2-(i\zeta)^2),
 \end{cases}
\end{align}
assuming that $\alpha_0^2\tau^2 \cong \eta^2+\zeta^2$. The system that determines the coefficients $\{a_{1-j}, b_{1-j}\}_{j \geq 1}$ is then given by
\begin{equation} \label{system_symbols}
\begin{aligned}
\begin{cases}
a_{1-j}+b_{1-j}&= \, 0, \quad \quad j \geq 1, \\[1mm]
-\alpha_1 (i \tau) \delta_{j1}&=\, \displaystyle -\sum_{{j \geq 1, \atop k+l+n=j,} \atop k, l, n \geq 0} \frac{(-1)^n}{n!}\partial_\tau^n a_{1-l} \partial_t^n b_{1-k}+\partial_x b_{1-j}, 
\end{cases}
\end{aligned}
\end{equation}
where $\delta$ denotes the Kronecker delta. From \eqref{system_a_1_different}, we find that
\begin{align*}
a_1=& \, -\sqrt{\alpha_0^2(i \tau)^2-(i \eta)^2-(i \zeta)^2}
\end{align*}
and
\begin{align*}
b_1= \sqrt{\alpha_0^2(i \tau)^2-(i \eta)^2-(i \zeta)^2}.
\end{align*}
Note that the sign of $a_1$ determines the propagation direction of the wave. To obtain the absorbing conditions for the given angle of incidence, we freeze the coefficient $\beta$ in \eqref{operator_P} by assuming that it is constant in space and time. The dispersion relation for \eqref{Westervelt_linearization_1/2} when $\beta$ is constant is as follows
\begin{align} \label{dispersion_relation}
    \alpha_0^2 (i \tau)^2-(i\xi)^2-(i\eta)^2-(i\zeta)^2=0.
\end{align}
The wave vector is given by $(\xi, \eta, \zeta)$. If we denote by $\theta \in [0 \, ^{\circ}, 90 \, ^{\circ}]$ the angle between the incident wave and the outer normal to the boundary, we have for $\tau > 0$
\begin{align*}
    \sin \theta =\frac{\sqrt{\eta^2+\zeta^2}}{\sqrt{\xi^2+\eta^2+\zeta^2}}=\frac{\sqrt{\eta^2+\zeta^2}}{\alpha_0 \tau};
\end{align*}
see Figure~\ref{fig:Boundary_sketch}.~\\[-5mm]

~\\
\noindent Therefore, we can express $a_1$ as 
\begin{align*}
a_1=& \,-\alpha_0(i \tau) \sqrt{1-\frac{\eta^2+\zeta^2}{\alpha_0^2 \tau^2}}= \, -\alpha_0(i \tau) \cos \theta.
\end{align*}
According to \eqref{abc_order_k}, the absorbing conditions of order zero are then given by
\begin{align} \label{ABC_Westervelt_inviscid_angle_1/2_linear}
\frac{\partial \psi}{\partial n}=-\alpha_0 \, \psi_t \cos \theta.
\end{align}
We mention that for a variable coefficient problem with jumps outside the computational domain, it is not possible to build an exact ABC based on local computations. After returning to $\beta=\psi_t$ and $\displaystyle \alpha_0=\sqrt{\frac{1}{c^2}-\frac{1}{2}\frac{k}{c^2}\psi_t}$, we obtain the absorbing boundary conditions for the inviscid Westervelt equation in the potential formulation for a given angle of incidence $\theta$:
\begin{align} \label{ABC_Westervelt_inviscid_angle_1/2}
c\frac{\partial \psi}{\partial n}+\sqrt{1-\frac{k}{2}\psi_t} \, \psi_t \, \cos \theta=0.
\end{align}
\subsubsection{Propagation with losses} We next want to incorporate the $b$ term into the conditions. This was not possible before since we needed $b=0$ to make use of the pseudo-differential factorization and the dispersion relation \eqref{dispersion_relation}. Instead, we employ a reasoning based on an energy argument.\\
\indent To this end, we test the linearized equation \eqref{Westervelt_linearization_1/2} with $\psi_t$, integrate over space and $(0,t)$, where $t \leq T$, and integrate by parts with respect to time, to arrive at the following identity:
\begin{equation} \label{energy_identity}
\begin{aligned}
& E_0[\psi](t)+\delta \int_0^t \|\nabla \psi_t\|^2_{L^2(\Omega)}\, \textup{d}s\\
=& \,E_0[\psi](0) + \int_0^t \int_{\Omega} \left(\tfrac12(\alpha^2_{0})_t+\alpha_1 \right)\psi_t^2 \, \textup{d}x \textup{d}s+ \int_0^t \int_{\partial \Omega}\left(\frac{\partial \psi}{\partial n}+\delta \frac{\partial \psi_t}{\partial n}\right)\, \psi_t \, \textup{d}x \textup{d}s,
\end{aligned}
\end{equation}
where the energy is given by
\begin{equation*}
\begin{aligned}
E_0[\psi](t)= \frac12 \left(\|\alpha_0(t) \psi_t(t)\|^2_{L^2(\Omega)}+ \|\nabla \psi(t)\|^2_{L^2(\Omega)} \right). 
\end{aligned}
\end{equation*}
This identity suggests to modify the conditions \eqref{ABC_Westervelt_inviscid_angle_1/2_linear} to include the sound diffusivity as follows
 \begin{align} \label{ABC_Westervelt_angle_1/2_linear_delta}
 	\frac{\partial \psi}{\partial n}+\delta \frac{\partial \psi_t}{\partial n}=-\alpha_0 \, \psi_t \cos \theta \quad \text{on } \ \Gamma_{\textup{abc}}.
 \end{align}
 These conditions facilitate the extraction of energy through the boundary since \eqref{energy_identity} becomes
\begin{equation} \label{energy_identity_new}
\begin{aligned}
& E_0[\psi](t)+\delta \int_0^t \|\nabla \psi_t\|^2_{L^2(\Omega)}\, \textup{d}s+\int_0^t \left\|\sqrt{\alpha_0 \cos \theta}\, \psi_t \right\|^2_{L^2(\partial \Omega)} \, \textup{d}s\\
=& \,E_0[\psi](0) + \int_0^t \int_{\Omega} (\tfrac12(\alpha^2_{0})_t+\alpha_1)\psi_t^2 \, \textup{d}x \textup{d}s,
\end{aligned}
\end{equation} 
from which by employing Gronwall's inequality it follows that
\begin{equation*}
    \begin{aligned}
  E_0[\psi](t) +\delta \int_0^t \|\nabla \psi_t\|^2_{L^2(\Omega)}\, \textup{d}s+\int_0^t \left\|\sqrt{\alpha_0 \cos \theta}\, \psi_t \right\|^2_{L^2(\partial \Omega)} \, \textup{d}s \leq C(T) E_0[\psi](0), 
    \end{aligned}
\end{equation*}
provided $\tfrac12(\alpha^2_{0})_t+\alpha_1 \in L^\infty(0,T; L^\infty(\Omega))$. \\
\indent We therefore adopt conditions \eqref{ABC_Westervelt_angle_1/2_linear_delta}. After returning to $\beta=\psi_t$ and $\displaystyle \alpha_0=\sqrt{\frac{1}{c^2}-\frac{1}{2}\frac{k}{c^2}\psi_t}$ in \eqref{ABC_Westervelt_angle_1/2_linear_delta} and recalling that $\delta=b/c^2$, we obtain the nonlinear conditions
 \begin{align} \label{ABC_Westervelt_angle_1/2}
 	c^2\frac{\partial \psi}{\partial n}+b\frac{\partial \psi_t}{\partial n}=-\, c\, \sqrt{1-\frac{k}{2}\psi_t} \, \psi_t \, \cos \theta \quad \text{on } \ \Gamma_{\textup{abc}}.
 \end{align}
We note that in realistic settings the sound diffusivity $b$ in fluids is small; see, e.g.,~\cite{MKaltenbacher}. It is also known that the presence of a large $b$ damping in the model would imply a parabolic instead of a wave-like behavior of the equation resulting in an exponential decay of the energy; cf.~\cite[Theorem 3.3]{kaltenbacher2009global}.\\
\indent Setting $k$ to zero in \eqref{ABC_Westervelt_angle_1/2} corresponds to conditions for a linear, strongly damped wave equation. If in addition $b=0$, we end up with the standard linear absorbing conditions for the angle $\theta$
\begin{align} \label{ABC_EM_angle}
c\frac{\partial \psi}{\partial n}+\, \psi_t \, \cos \theta=0 \quad \text{on } \ \Gamma_{\textup{abc}};
\end{align}
see~\cite{higdon1986absorbing, keys1985absorbing}.
\begin{remark}[One- and two-dimensional domains]
In a one-dimensional setting, system \eqref{system_a_1_different} for determining the symbols $a_1$ and $b_1$ simplifies to
\begin{align*} 
\begin{cases}
a_1+b_1= 0 \\
a_1b_1=-\alpha_0^2(i\tau)^2.
\end{cases}
\end{align*}
In a two-dimensional setting, system \eqref{system_a_1_different} simplifies to
\begin{align*} 
\begin{cases}
a_1+b_1= 0 \\
a_1b_1=-(\alpha_0^2(i\tau)^2-(i\eta)^2).
\end{cases}
\end{align*}
It is then straightforward to show that conditions \eqref{ABC_Westervelt_angle_1/2} hold in $1$D and $2$D as well, where in 1D the angle $\theta$ can be interpreted as being set to $\theta=0^{\circ}$. 
\end{remark}
\begin{remark}[A different linearization]
Employing linearization \eqref{Westervelt_linearization_KS} would result in the following absorbing conditions 
\begin{align} \label{ABC_Westervelt_angle_KS}
c^2\frac{\partial \psi}{\partial n}+b\frac{\partial \psi_t}{\partial n}=-\, c\, \sqrt{1-k\psi_t} \, \psi_t \, \cos \theta \quad \textup{on } \ \Gamma_{\textup{abc}}.
\end{align}
The performance of conditions \eqref{ABC_Westervelt_angle_1/2} and \eqref{ABC_Westervelt_angle_KS} is compared in Section~\ref{Section:NumResults}, where the proposed conditions \eqref{ABC_Westervelt_angle_1/2} significantly outperform \eqref{ABC_Westervelt_angle_KS}.
\end{remark}

\begin{remark}
In our experiments, we employ the gradient information to compute the angle of incidence $\theta$ via
\begin{align*}
    \cos \theta=\dfrac{|\nabla \psi \cdot n|}{\sqrt{\psi^2_x+\psi^2_y+\psi_z^2}}\quad \textup{on } \ \Gamma_{\textup{abc}},
\end{align*}
assuming that $\nabla \psi \neq 0$ on the absorbing boundary. The linear conditions \eqref{ABC_EM_angle} for the angle $\theta$ are then equivalent to
 \begin{align} \label{equiv_conditions}
     c|\nabla \psi|+\psi_t=0 \quad \textup{on } \ \Gamma_{\textup{abc}}.
 \end{align}
Therefore, in the linear regime, conditions \eqref{equiv_conditions} can be seen as the extension of the absorbing conditions proposed in \cite{gordon2015compact} for the Helmholtz equation to the linear time-dependent wave model. 
\end{remark}
\section{Numerical treatment}\label{Section:NumTreatment}
After deriving the absorbing boundary conditions for the potential form of the Westervelt equation \eqref{Westervelt}, we next focus on the numerical schemes used in simulations. We begin by formulating the initial-boundary value problem that has to be solved.
\subsection{The initial-boundary value problem for the Westervelt equation}
We consider the following problem for the Westervelt equation:
\begin{equation} \label{IBVP}
\begin{cases}
\begin{alignedat}{3}
\psi_{tt}-c^2\Delta \psi-b\Delta \psi_t&=\frac{1}{c^2}(B/A+2)\, \psi_t \psi_{tt}&&\qquad \textup{in } \Omega\times (0,T),  \\[1mm]
\psi &= g &&\qquad \textup{on } \Gamma_{\textup{exc}}\times (0,T),  \\[1mm]
c\frac{\partial\psi}{\partial n}+\frac{b}{c}\frac{\partial\psi_t}{\partial n} &=- \sqrt{1-\sigma k\psi_t}\, \psi_t \cos \theta(\psi) &&\qquad \textup{on } \Gamma_{\textup{abc}}\times (0,T),\\[1mm]
\frac{\partial\psi}{\partial n} &= 0 &&\qquad \textup{on } \Gamma_{\textup{N}}\times(0,T),  \\[1mm]
\psi = \psi_t &= 0 &&\qquad \textup{in } \Omega\times\lbrace 0\rbrace.
\end{alignedat}
\end{cases}
\end{equation}
The wave source is given in the form of inhomogeneous Dirichlet conditions on the excitation part of the domain boundary $\Gamma_{\textup{exc}}\subset\partial\Omega$. In our numerical tests, the excitation signal is always taken to be a modulated sine wave, that is growing over time until its maximal amplitude is reached, i.e.
\begin{equation} \label{excitation}
    g(t) = \begin{cases} \displaystyle
    (f^2/4)t^2 \mathfrak{A}\, \sin(\omega t), & t<2/f, \\[2mm]
    \hspace{1.3cm}\mathfrak{A}\, \sin(\omega t), & t\geq 2/f,
    \end{cases}
\end{equation}
where $f$ denotes the frequency, $\omega = 2\pi f$ the angular frequency, and $\mathfrak{A}$ the maximal amplitude of the signal. \\
\indent We have introduced the parameter $\sigma$ within the square root of the absorbing conditions in \eqref{IBVP}. In this way, we generalized in one formula all the absorbing conditions that we want to compare. Indeed, setting $\sigma = 0$ yields the adaptive absorbing conditions for the linear strongly damped wave equation 
\begin{align} \label{ABC_W_0}
c\frac{\partial \psi}{\partial n}+\frac{b}{c}\frac{\partial \psi_t}{\partial n}=-\psi_t \cos \theta(\psi),
\end{align}
which we denote from now on in experiments by ``ABC$^0_{\textup{W}}$ adaptive". In case also $\theta=0$ everywhere, we denote them just by $\textup{ABC}_{\textup{W}}^0$. Setting $\sigma = 1/2$ recovers our new nonlinear adaptive conditions, denoted by ``ABC$_{\textup{W}}^{1/2}$ adaptive". If the angle is always set to zero, we denote them just by ABC$_{\textup{W}}^{1/2}$. Finally, $\sigma = 1$ leads to conditions based on the second linearization \eqref{Westervelt_linearization_KS}, which are denoted in the experiments by ``ABC$_{\textup{W}}^{1}$ adaptive" and ABC$_{\textup{W}}^{1}$.  \\
\indent We start from the weak form of the problem \eqref{IBVP}. We are looking for a solution in  $$\{\psi \in C^1([0,T]; H^2(\Omega)) \cap C^2([0,T]; H^1(\Omega))\, \vert \ \psi= g \ \textup{on } \Gamma_{\textup{exc}}\times (0,T)\}$$ such that
\begin{equation*}
\begin{aligned}
&\int_{\Omega}((1-k \psi_t)\psi_{tt}v + c^2\nabla\psi\cdot\nabla v + b\nabla\psi_t\cdot\nabla v) \, \textup{d}\Omega 
+\int_{\Gamma_{\textup{abc}}}  c \, \sqrt{1-\sigma k\psi_t}\, \psi_t \, \cos \theta(\psi) \, v~\textup{d}S
=  \, 0
\end{aligned}
\end{equation*}
for all test functions in $\lbrace v \in H^1(\Omega)\,|\, v= 0 \ \textup{on } \Gamma_{\textup{exc}}\times (0,T)\rbrace$ a.e. in time, with $(\psi, \psi_t)\vert_{t=0}=(0,0)$. We assume that the problem \eqref{IBVP} is well-posed, although the rigorous proof is beyond the scope of the current work. Results on the well-posedness of the Westervelt equation in the pressure form with nonlinear absorbing conditions for the angle of incidence $\theta=0^{\circ}$ can be found in~\cite{shevchenko2015absorbing, simonett2017well}.
\subsection{Finite element discretization and time integration} We follow the standard discretization methods for nonlinear acoustics based on finite elements~\cite{fritz2018well, Hoffelner, kagawa1992finite, MKaltenbacher, muhr2017isogeometric, WalshTorres}. The finite element method is employed in space with lowest order conforming elements on simplicial meshes. \\ 
\indent The mass \textbf{M}, stiffness \textbf{K}, and damping matrix \textbf{C} as well as the nonlinearity tensor $\boldsymbol{\mathcal{T}}$ are assembled in the usual manner; ~see~\cite{MKaltenbacher, muhr2017isogeometric}. By dividing the set of degrees of freedom into the set of Dirichlet degrees $\textup{D}$ and the set of interior degrees $\textup{I}$, the semi-discrete problem reads as follows
\begin{equation} \label{ODE_system}
 \begin{aligned}
 \begin{cases}
   \, {\bf M}_{\textup{I},\textup{I}}\underline{\ddot{\psi}}_{\textup{I}}+{\bf K}_{\textup{I},\textup{I}}\underline{\psi}_{\textup{I}}+{\bf C}_{\textup{I},\textup{I}}\underline{\dot{\psi}}_{\textup{I}}-\boldsymbol{\mathcal{T}}_{\textup{I},\textup{I},\textup{I}}[\underline{\ddot{\psi}}_{\textup{I}},\underline{\dot{\psi}}_{\textup{I}},\cdot]-\boldsymbol{\mathcal{T}}_{\textup{I},\textup{D},\textup{I}}[\underline{\ddot{\psi}}_{\textup{I}},\underline{\dot{\psi}}_{\textup{D}},\cdot]\\[2mm]
  \,  -\boldsymbol{\mathcal{T}}_{\textup{D},\textup{I},\textup{I}}[\underline{\ddot{\psi}}_{\textup{D}},\underline{\dot{\psi}}_{\textup{I}},\cdot]-\underline{\bf A}_{\textup{I}}(\underline{\psi},\underline{\dot{\psi}},\theta(\psi))= F(t) \hspace{43mm} \textup{in~} (0,T),\\[2mm]
   \, \underline{\psi}=\underline{\psi_t}=0 \hspace{90mm} \textup{at~} t=0.
    \end{cases}
\end{aligned}
\end{equation}
The right-hand side of the equation is given by
$$F(t)=-{\bf M}_{\textup{I},\textup{D}}\underline{\ddot{\psi}}_{\textup{D}}-{\bf K}_{\textup{I},\textup{D}}\underline{\psi}_{\textup{D}}-{\bf C}_{\textup{I},\textup{D}}\underline{\dot{\psi}}_{\textup{D}}-\boldsymbol{\mathcal{T}}_{\textup{D},\textup{D},\textup{I}}[\underline{\ddot{\psi}}_{\textup{D}},\underline{\dot{\psi}}_{\textup{D}},\cdot].$$
The underlined quantities $\underline{\psi},\underline{\dot{\psi}}$, and $\underline{\ddot{\psi}}$ denote the coefficient vectors of $\psi, {\psi}_t$, and ${\psi}_{tt}$ resulting from the spatial finite element discretization. The compact notation with index sets \textup{D} and \textup{I} is used to extract the respective rows and columns of matrices and vectors that belong to Dirichlet and interior degrees of freedom. \\
\indent The absorbing boundary vector $\underline{\textbf{A}}$ is formally given by

\begin{equation}
\begin{aligned}
& \displaystyle \underline{\bf A}(\psi_t, \theta(\psi)) =(A_{i}(\psi_t, \theta(\psi)))_{i\in\textup{DOF}(\Gamma_{\textup{abc}})}, \\[2mm]
& \displaystyle A_i(\psi_t, \theta(\psi)) =\int_{\Gamma_{\textup{abc}}}c\sqrt{1-\sigma k\psi_t}\, \psi_t\cos \theta(\psi) \, N_i~\textup{d}S.     
\end{aligned}
 \end{equation}
Above, $N_i$ stands for the finite element ansatz function of the $i$-th global degree of freedom, while DOF($\Gamma_{\textup{abc}}$) is the set of degrees of freedom belonging to the absorbing boundary. \\
\indent The nonlinearity tensor $\boldsymbol{\mathcal{T}}$ is used to resolve the nonlinear bulk term in the weak formulation in a fixed-point iteration. The same also holds for the absorbing boundary vector which is iteratively updated with the current values of $\psi_t$, $\psi_{tt}$ as well as updates for the angle $\theta$.\\
\indent The system \eqref{ODE_system} is a nonlinear system of ordinary differential equations of second order with $|\textup{I}|$ components. It remains to solve it by using a suitable time integrator. Following~\cite{muhr2017isogeometric}, we employ the Generalized-$\alpha$ scheme in combination with the Newmark relations for time integration. Values of the Generalized-$\alpha$ parameters $(\alpha_m, \alpha_f)$ and the Newmark parameters $(\beta_{\textup{nm}}, \gamma_{\textup{nm}})$ that are used in experiments are chosen according to the stability and accuracy criteria stated in \cite{bonelli2002analyses}:
$$\alpha_{m}= \frac{2\rho_{\infty}-1}{1+\rho_{\infty}},~~ \alpha_f = \frac{\rho_{\infty}}{1+\rho_{\infty}},~~ \beta_{\textup{nm}} = \frac{1}{(1+\rho_{\infty})^2},~~ \gamma_{\textup{nm}} = \frac{1}{2}\frac{3-\rho_{\infty}}{1+\rho_{\infty}}$$  
where we take $\rho_{\infty}=1/2$; see also Table~\ref{table_param} for the resulting values.\\

\indent In comparison to the numerical solvers proposed in~\cite{fritz2018well, muhr2017isogeometric}, a new aspect of our method is the computation of the angle of incidence $\theta(\psi)$. The angle is computed once in every time step before the first assembly of the absorbing boundary vector. Details on how we compute the angle can be found in Section~\ref{Section:Angle}.

 \begin{center}
 	\renewcommand{\arraystretch}{1.6}
 	\captionof{table}{Time stepping parameters} \label{table_param} 
 	\begin{tabular}{||c | c | c ||} 
 		\cline{2-3}
 		\multicolumn{1}{c|}{}& parameter & value   \\ 
 		\cline{2-3}\hline
 		\multirow{2}{*}{Newmark-parameters} & $\beta_{\textup{nm}}$ & 4/9    \\
 		\cline{2-3}
 		& $\gamma_{\textup{nm}}$ & 5/6 \\ 
 		\hline
 		\multirow{2}{*}{Generalized-$\alpha$ parameters}& $\alpha_m$ & 0\\ 
 		\cline{2-3}
 	    & $\alpha_f$  & 1/3    \\ 
 		\hline
 		\multirow{2}{*}{Nonlinear iteration parameters}& TOL & $10^{-6}$\\ 
 		\cline{2-3}
 	    & $\kappa_{\max}$  & 100    \\ 
 		\hline
 	\end{tabular}\vspace{2mm} \\
 \end{center}

 
\section{Computation of the angle of incidence}\label{Section:Angle}
It remains to see how we can compute the angle of incidence $\theta$ that the local wave vector encloses with the outward normal to the absorbing boundary at a given point on $\Gamma_{\textup{abc}}$. To reduce the computational cost, we compute the angle $\theta$ once per time step and do not update it within the nonlinear fixed-point iteration. \\ \enlargethispage*{2cm}
\indent According to \cite{Poynting1,Poynting2}, the Poynting vector $\textbf{P}(\psi)$ of a wave field $\psi$ can be used to compute its local propagation direction. The vector $\textbf{P}(\psi)$ is given by\\
\vspace*{-3mm}
\begin{equation}\label{eq:Poynting}
\textbf{P}(\psi)=-\frac{\partial\psi}{\partial t}\nabla\psi.
\end{equation}

Since neither the sign nor the norm of the propagation vector plays a role when computing the incidence angle, we restrict ourselves to the spatial gradient alone to determine the main propagation direction of the wave. Such an approach was also taken in \cite{jia2009calculation, ye2012effective, zhang2011direct} for wave fields in an isotropic medium. In our case, this method of computing the local propagation direction works especially well since, although globally discontinuous, the element-wise gradient information is readily available at every time step in the finite element framework. We also refer to~\cite{patrikeeva2013comparison} for a further discussion on the use of the Poynting vector in angle decomposition methods. \\
\indent We conduct experiments with zero initial data and Dirichlet conditions on part of the boundary, and so most of the potential field is at rest at the beginning of the simulation. However, numerical noise of low magnitude can be present at the absorbing boundary before the wave reaches it. Such behavior could be accounted to weak ill-posedness; see~\cite{mulder1997experiments}. To tackle this issue, we implement a switch. When going over all elements adjacent to the absorbing boundary, we only compute the element-wise angle of incidence once the local wave amplitude (in terms of absolute value of the elements degrees of freedom) exceeds a certain percentage $p_1$ of a reference value; see~Algorithm \ref{alg:AngleComp}, line \ref{alg:p1_criterion}. We take the reference value to be the source amplitude of the wave. In the case that the source amplitude is not known a priori, an alternative would be to compute the maximum field amplitude in the interior of the domain and take this as a reference value. As long as the criterion is not matched, the local angle of incidence is set to 0; see~Algorithm \ref{alg:AngleComp}, line \ref{alg:set_angle_to_zero}. We note that a similar approach was taken in~\cite{shevchenko2012self}. In all our numerical experiments, we set $p_1 = 0.1$. \\
\indent Algorithm~\ref{alg:AngleComp} summarizes our method of computing the incidence angle. Within the algorithm, indices in the exponent indicate the time step.\\
~\\
\let\oldnl\nl
\newcommand{\nonl}{\renewcommand{\nl}{\let\nl\oldnl}}

\begin{algorithm}[H]
\setstretch{1.1}
\caption{Angle-computation algorithm}\label{alg:AngleComp}\vspace{1mm}
\nonl \underline{Initialization:}~\\[1mm]
\textup{Formally \textbf{set}}~$|\nabla\psi_{\textup{el}}^{(-1)}|=\infty$ \textup{~ and $\theta_{\textup{el}}^0=0\, ^{\circ}$ for all elements}~\textup{el}~\\[5mm]
\nonl \underline{\textup{In time step}~$n=1,2,...$~\textbf{do}}\\[1mm]
\SetKwBlock{forloop}{for \textup{el} $\in \lbrace \textrm{elements}\,:\,\textrm{element has an edge/face on the absorbing boundary}\rbrace$}{end}
\forloop{
\If{$\max\lbrace|\psi_i^{n}|\,:\,\textup{dof}_i~\textup{belongs to}~\textup{el}\rbrace > p_1\cdot \mathfrak{A}$\label{alg:p1_criterion}}{
\textup{enable angle computation for}~\textup{el}}
\eIf{\textup{angle computation is enabled for}~\textup{el}}
    {
    \textup{Evaluate}~$\nabla\psi(\vec{x},t^{n-1})$~\textup{within}~ \textup{el} $\rightarrow \textup{save as}~\nabla\psi_{\textup{el}}^{n-1}$~\\
    \eIf{$|\nabla\psi_{\textup{el}}^{n-1}|\leq p_2 \, \max\limits_{k<n-1}|\nabla\psi_{\textup{el}}^k|$\label{alg:p2_criterion}}{
    \textup{\textbf{set}} $\theta_{\textup{el}}^{n} = \theta_{\textup{el}}^{n-1}$ \label{alg:angle_from_last_timestep}
    }
    {
    \textup{\textbf{compute}} $\theta_{\textup{el}}^n = \arccos\left(\dfrac{|\langle\nabla\psi_{\textup{el}}^{n-1}, \vec{n}_{\textup{el}}\rangle|}{|\nabla\psi_{\textup{el}}^{n-1}|}\right)$\label{alg:compute_angle}~\\
    }
    }
    {
    \textup{\textbf{set}} $\theta_{\textup{el}}^n=0\, ^{\circ}$\label{alg:set_angle_to_zero}
    }}
\end{algorithm}
~\\
\indent Note that even once the local amplitude of the wave at a given element is large enough for the angle computation to start, unreliable angle values can be computed at points where a local wave maximum or minimum hits the boundary since the gradient is close to zero. As a remedy, we propose that gradients with the Euclidean norm below some threshold should not influence the angle of incidence. Whenever such a small gradient appears, we use the angle of the last time step; see~Algorithm \ref{alg:AngleComp}, line \ref{alg:angle_from_last_timestep}. We employ a percentage $p_2$ of a reference value for the threshold. As the reference value we take the local gradient history of the given element and pick the norm-wise maximum over the past time steps; see~Algorithm \ref{alg:AngleComp}, line \ref{alg:p2_criterion}. We compute a new angle of incidence only in cases where the threshold is surpassed, i.e., when the local gradient is sufficiently large; cf. Algorithm~\ref{alg:AngleComp}, line \ref{alg:compute_angle}. To also reduce the oscillations with respect to time in the angle distribution, in all experiments we choose a relatively high threshold of $p_2=0.5$. \\
\indent Together with the previously introduced switch, the above approach allows steering the sensitivity of the angle computation algorithm by adapting the parameters $p_1$ and $p_2$. If these parameters are close to 1, the angle is only computed for very high amplitudes and local gradients, while for most other parts of the boundary the angle remains zero. On the other hand, small values of $p_1$ and $p_2$ lead to highly sensitive angles reacting to even small perturbations in the wave field. \\

\section{Numerical results}\label{Section:NumResults}
\begin{figure}[h!]
\begin{center}
\includegraphics[scale=0.55]{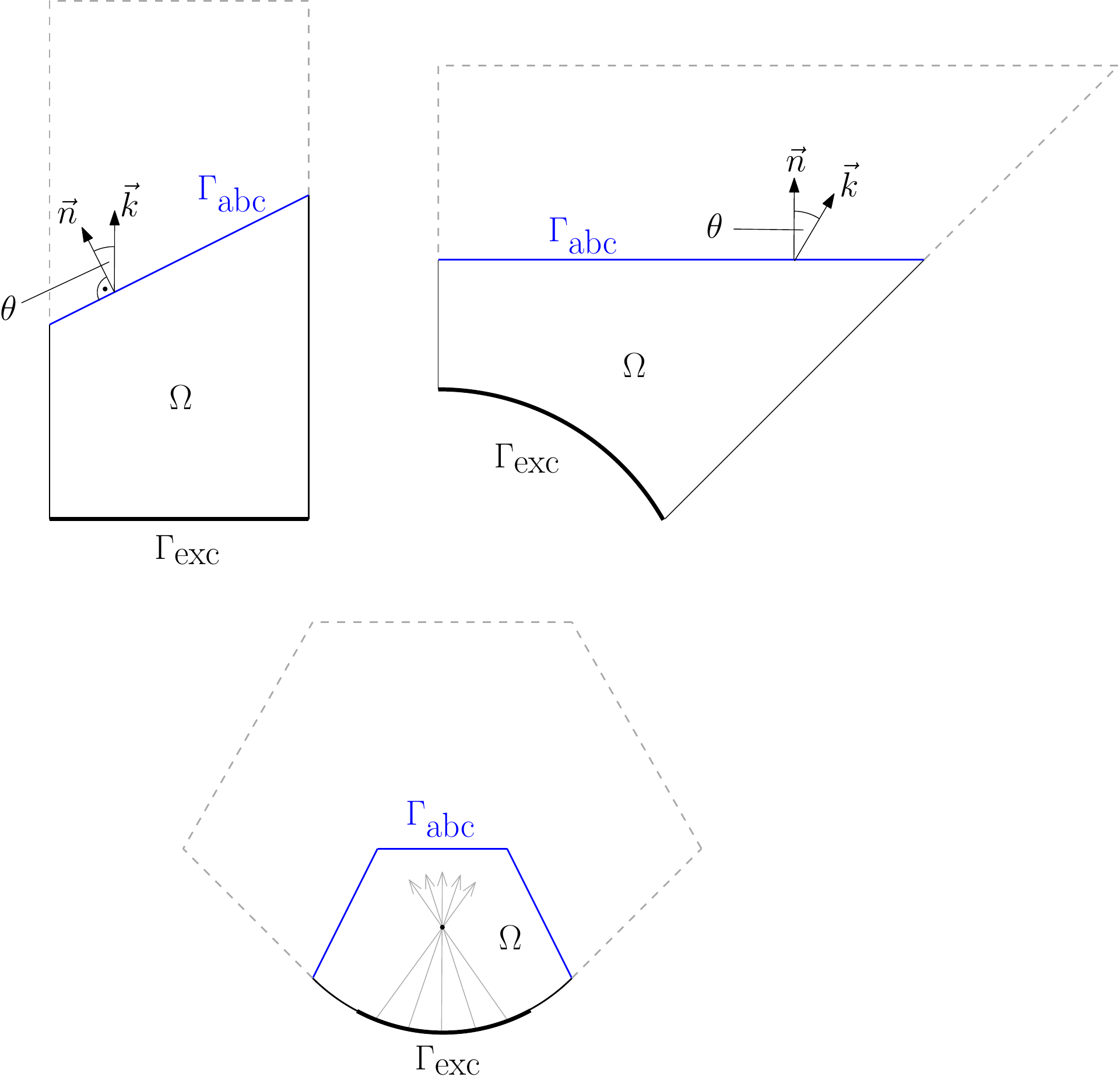}
\caption{Computational domains used in simulations. \textbf{(Top left)} Channel with an inclined absorbing boundary, \textbf{(Top right)} Octant of a ``plate with hole'' geometry, \textbf{(Bottom)} Geometry of a focusing transducer.
\label{fig:Domains}}
\end{center}
\end{figure}

We proceed with numerical simulations where we put our self-adaptive technique to the test. Computational domains $\Omega$ used in numerical experiments are sketched in Figure~\ref{fig:Domains}. The dashed lines symbolize the boundaries of the reference domain $\Omega_{\textup{ref}}$ where we compute the reference solution $\psi^{\textup{ref}}$. The reference solution is always first computed on $\Omega_{\textup{ref}}$, then restricted to the actual domain $\Omega$, and compared with the potential field obtained on $\Omega$ by employing the absorbing conditions.\\
\indent We mention again that in all numerical simulations conforming finite elements of lowest order on simplicial meshes are used. Geometry and mesh are generated by using the Gmsh software package~\cite{gmsh}.

\subsection{Domain with an inclined absorbing boundary} \label{SubSec:InclinedBoundary}
\begin{figure}[h!]
\begin{center}
\includegraphics[scale=0.2, trim=1.5cm 0cm 38cm 0cm, clip]{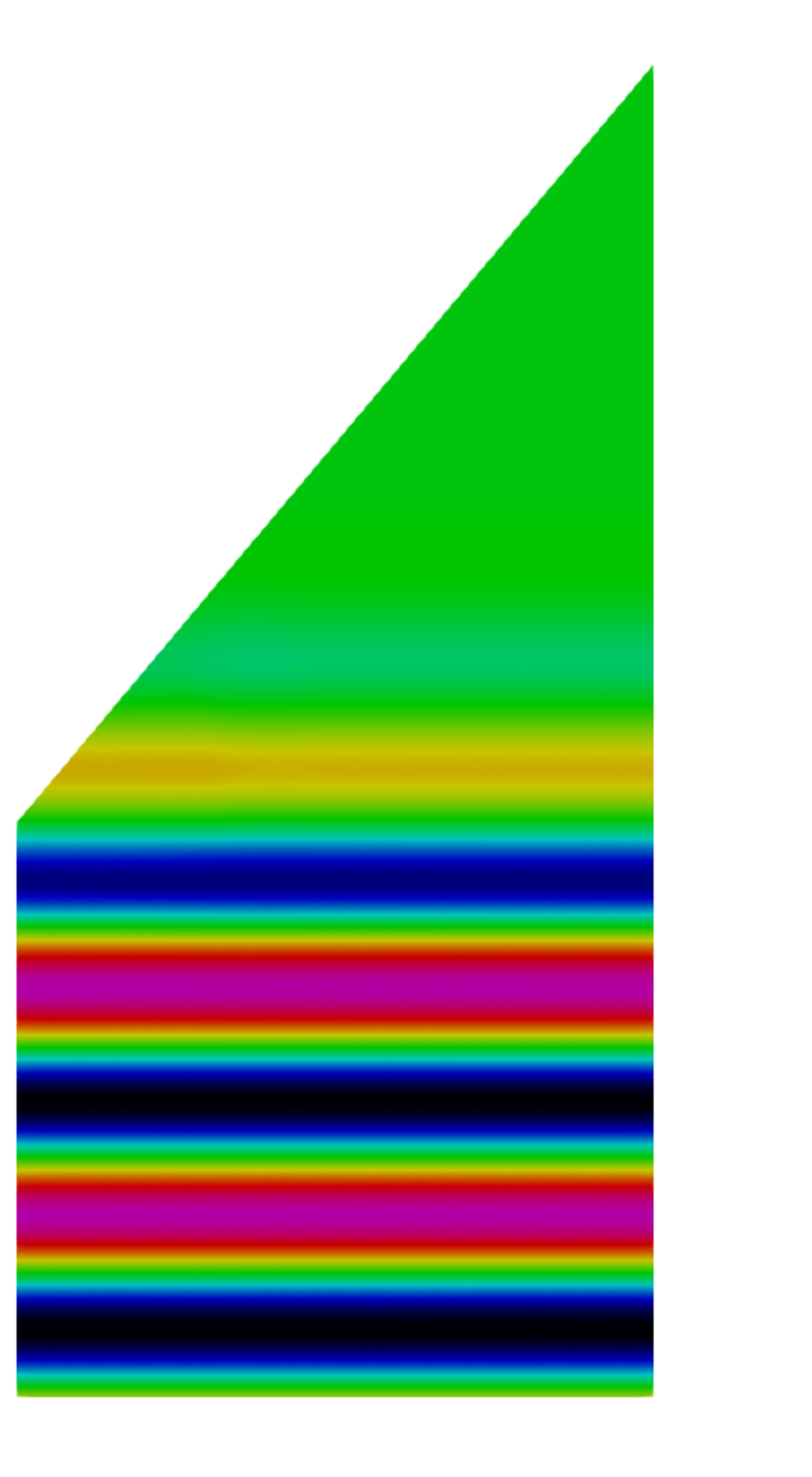}\includegraphics[scale=0.2, trim=1.5cm 0cm 38cm 0cm, clip]{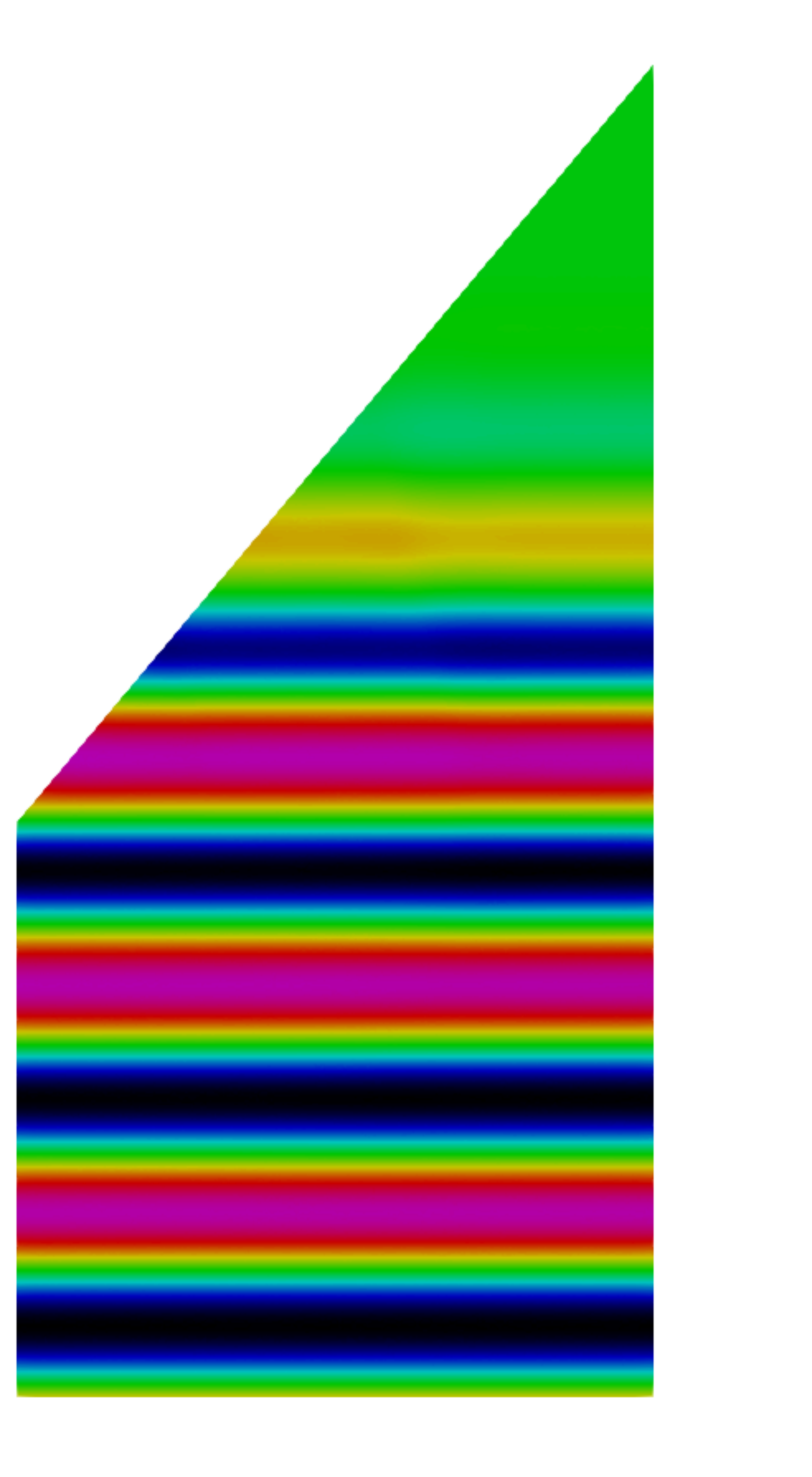}\includegraphics[scale=0.2, trim=1.5cm 0cm 38cm 0cm, clip]{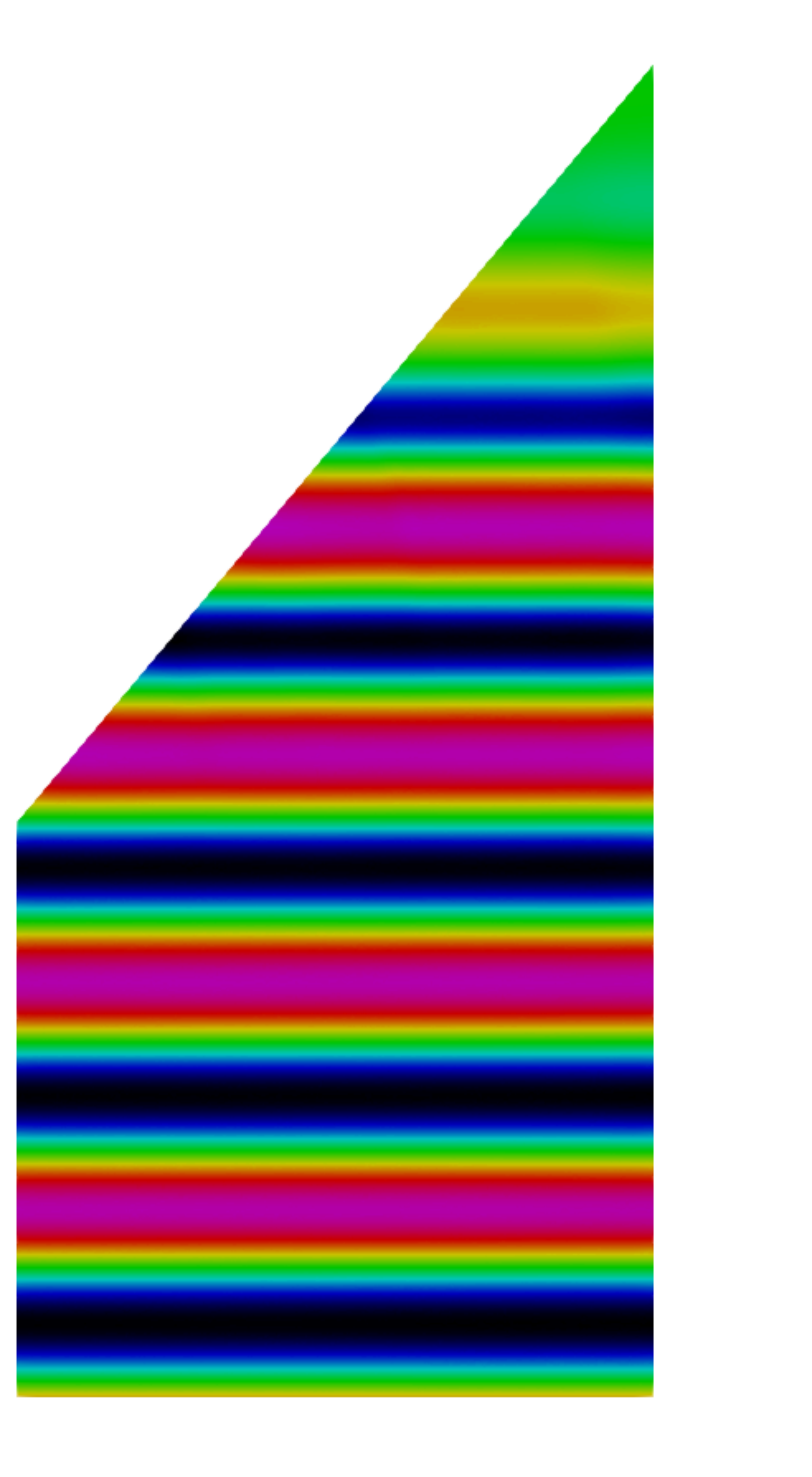}\includegraphics[scale=0.2, trim=1.5cm 0cm 38cm 0cm, clip]{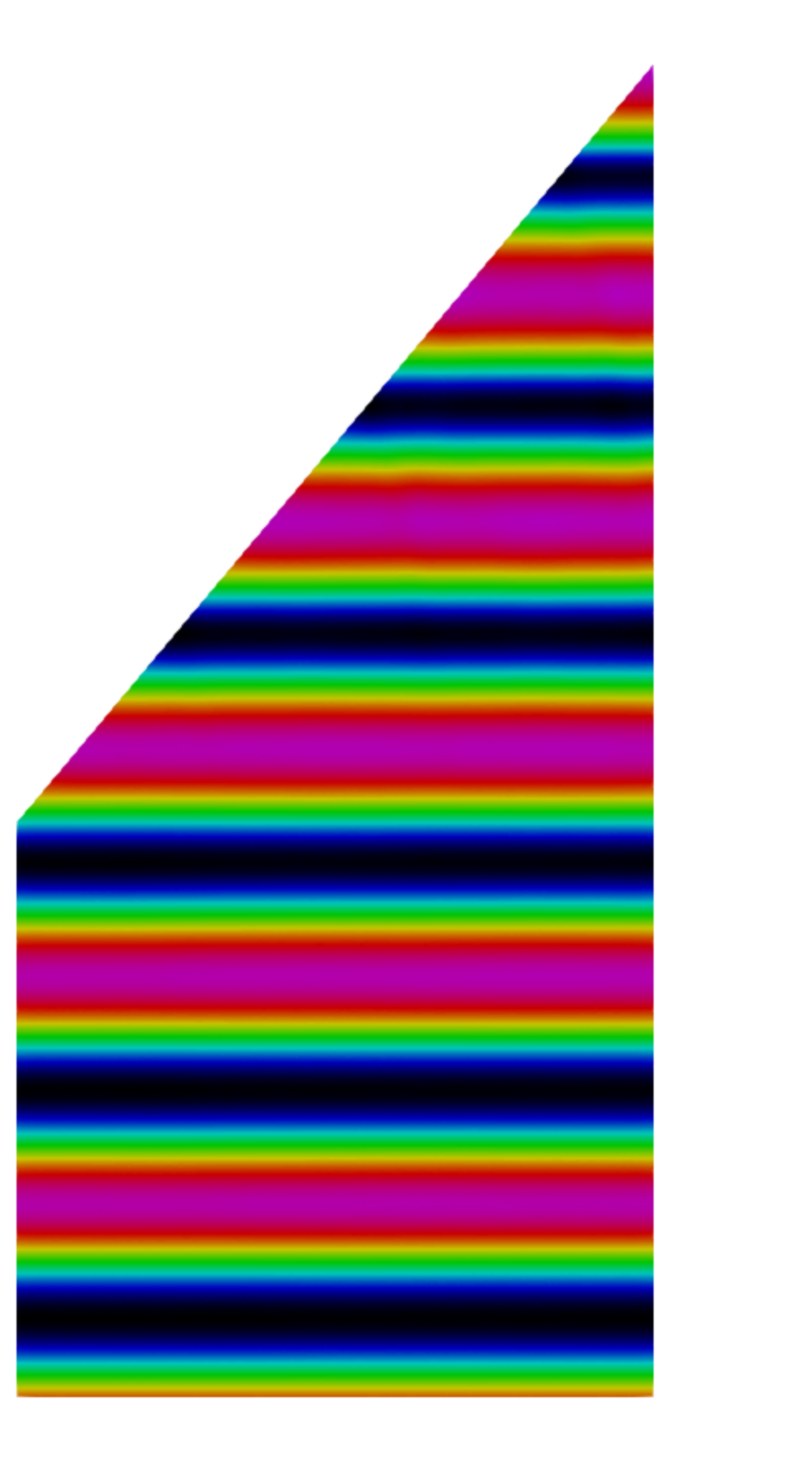}
\caption{Time-snapshots of the potential field.
\label{fig:Plane50Wave}}
\end{center}
\end{figure}

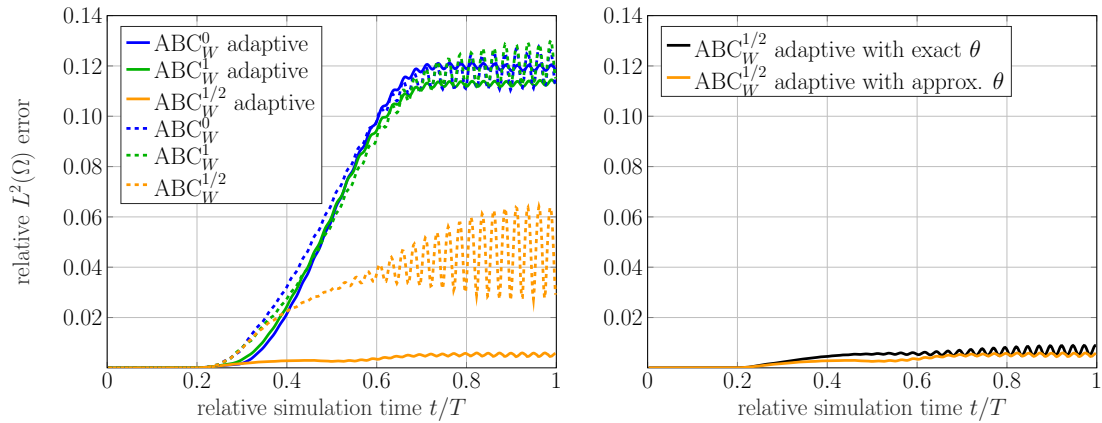
\begin{figure}[h!]
\begin{center}
\input{images/Plane20/Cost.tex}\hspace*{2mm} \input{images/Plane20/Angle_Quality.tex}
\caption{\textbf{Inclined plane boundary:} Relative $L^2(\Omega)$ error of the potential $\psi(t)$ over the simulation time with $\theta=20^{\circ}$ . \textbf{(Left)} Nonlinear vs. linear conditions with and without adaptivity, \textbf{(Right)} Performance of adaptive conditions with numerically computed vs. the exact incidence angle $\theta$.
\label{fig:Plane20Cost}}
\end{center}
\end{figure}

In our first experiment, we consider a two-dimensional channel geometry, where the upper (absorbing) boundary $\Gamma_{\textup{abc}}$ is tilted by a given angle of $\theta$; see Figure~\ref{fig:Domains}. The waves originate from the excitation boundary $\Gamma_{\textup{exc}}$ at the bottom of the rectangle and travel straight upwards. We impose homogeneous Neumann boundary conditions on the sides of the domain. The wave vector $\vec{k}$ in this setting points straight upwards. Therefore, the angle that the wave vector and the outward normal $\vec{n}$ of $\Omega$ at $\Gamma_{\textup{abc}}$ enclose is exactly $\theta$ for the absorbing boundary. \\
\indent We compute the reference solution $\psi^{\textup{ref}}$ on a larger domain without the absorbing boundary and then conduct a simulation on $\Omega$. To get an impression on how the wave propagates in the present setting, Figure~\ref{fig:Plane50Wave} shows the potential field $\psi$ at different time snapshots. Material parameters were chosen to be the ones of water, i.e. 
$$c=1500 \, \textup{m}/\textup{s}, \ b=6\cdot 10^{-9} \, \textup{m}^2/\textup{s}, \ \rho = 1000 \, \textup{kg}/\textup{m}^3, \ B/A=5;$$ see~\cite[Chapter 5]{MKaltenbacher}, while the excitation \eqref{excitation} has an amplitude of $\mathfrak{A}=0.01\,\textup{m}^2/\textup{s}^2$ and a frequency of $210\,\textup{kHz}$. The experiment was conducted for two different angles $$\theta\in\lbrace 20^{\circ},  50^{\circ}\rbrace.$$
\indent For spatial discretization, we take 13045 ($20^{\circ}$ case) and 13046 ($50^{\circ}$ case) degrees of freedom in space for the channel width of $0.02\,$m and channel length (in the middle) of $0.03\,\textup{m}$. In time, 9800 time steps were taken to cover the interval from $t_0 = 0$ until $T = 9.45\cdot 10^{-5}\,$s.\\
\indent The error plots are given in Figure~\ref{fig:Plane20Cost}. We observe that the conditions that do not take the angle of incidence into consideration significantly deteriorate when the angle increases. For the angle of incidence $\theta = 50^{\circ}$, the maximal relative error is more than $20$\,\% for the linear conditions and around $17$\,\% for the nonlinear conditions with the fixed angle $\theta=0^{\circ}$. In comparison, the self-adaptive technique when combined with the nonlinear conditions \eqref{ABC_Westervelt_angle_1/2} allows for the error to remain around $1$\,\%. Note that linearization \eqref{Westervelt_linearization_1/2} and the resulting absorbing conditions \eqref{ABC_Westervelt_angle_1/2} clearly outperform conditions \eqref{ABC_Westervelt_angle_KS}. We, therefore, proceed in the following experiments with testing only \eqref{ABC_Westervelt_angle_1/2} in combination with the self-adaptive technique. We also observe that the nonlinear conditions  \eqref{ABC_Westervelt_angle_1/2} that use an approximate angle computation via the gradient of the wave field perform similarly to the conditions that employ the exact angle.\\

\begin{figure}[h!]
\begin{center}
\includegraphics[scale=0.15, trim=1.5cm 0cm 38cm 0cm, clip]{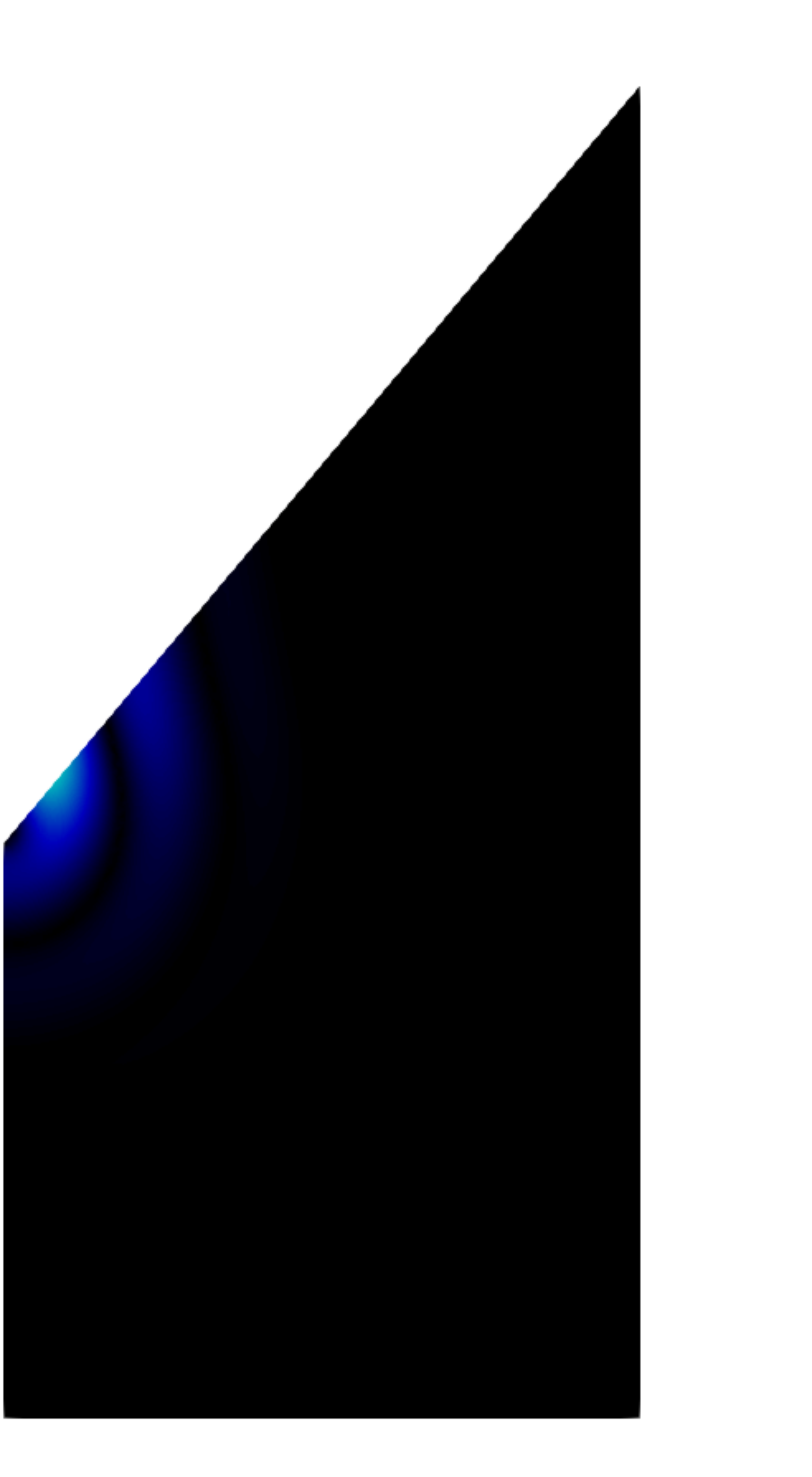}\includegraphics[scale=0.15, trim=1.5cm 0cm 38cm 0cm, clip]{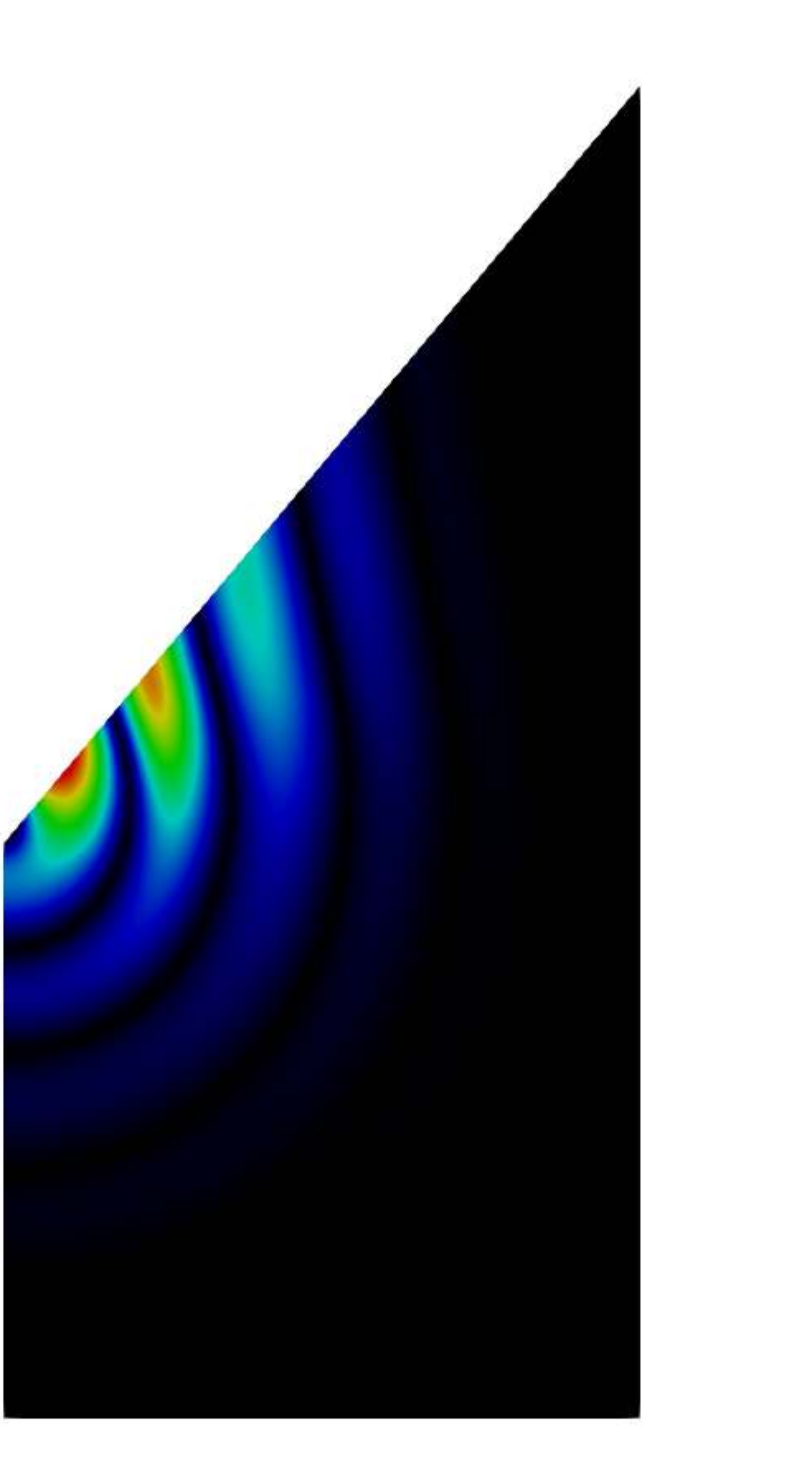}\includegraphics[scale=0.15, trim=1.5cm 0cm 38cm 0cm, clip]{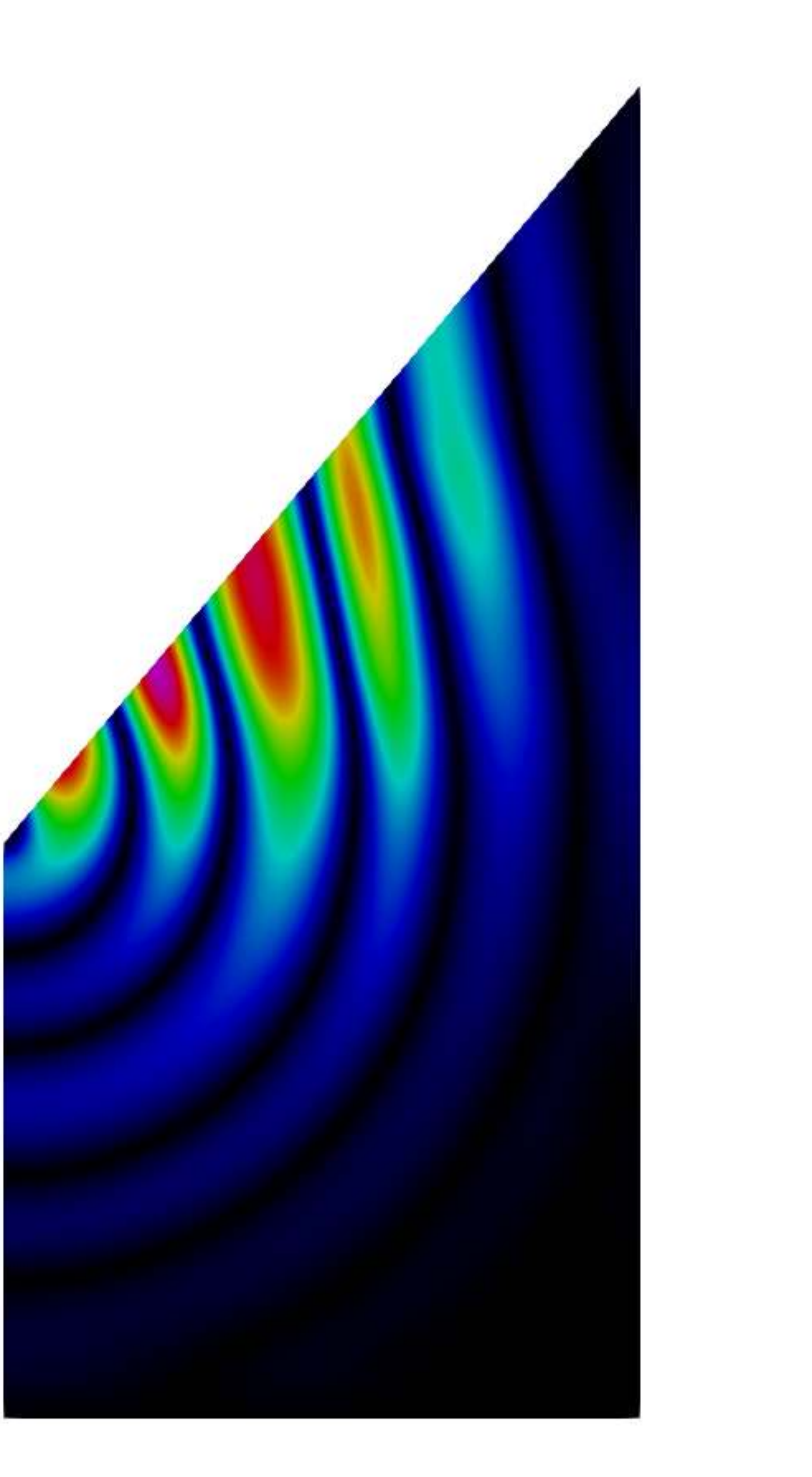}\includegraphics[scale=0.15, trim=1.5cm 0cm 38cm 0cm, clip]{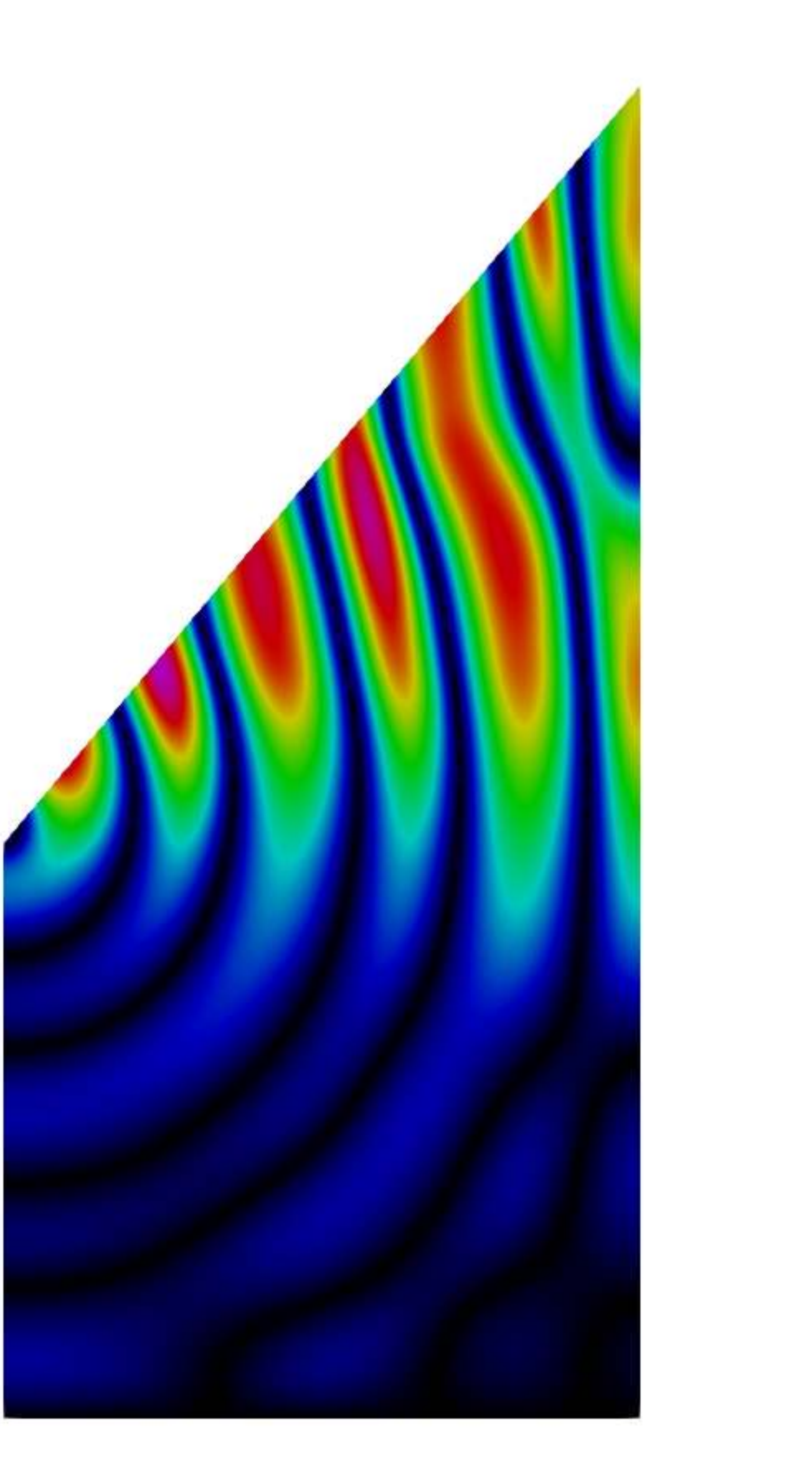}~\\[-5mm]
\includegraphics[scale=0.15, trim=1.5cm 0cm 38cm 0cm, clip]{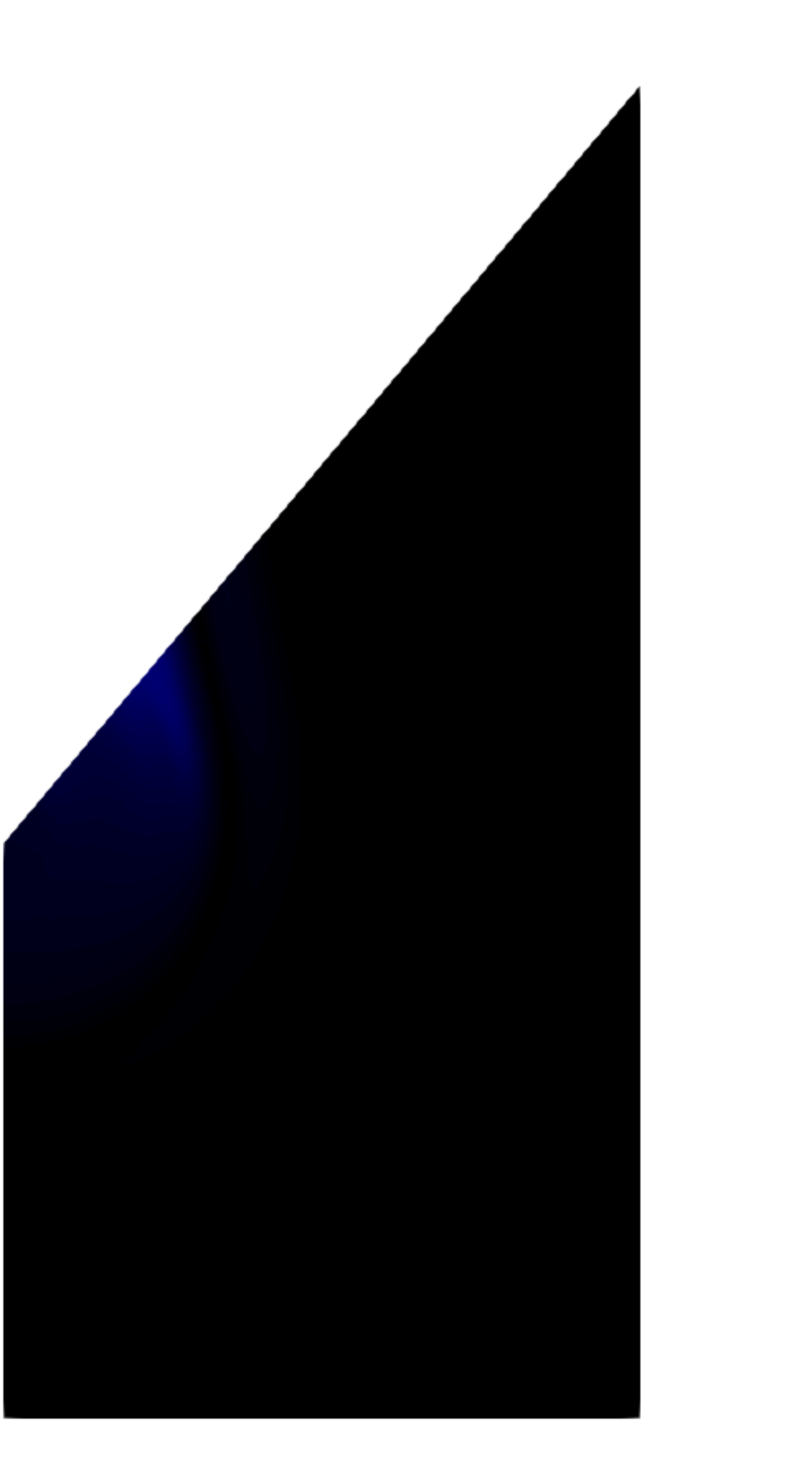}\includegraphics[scale=0.15, trim=1.5cm 0cm 38cm 0cm, clip]{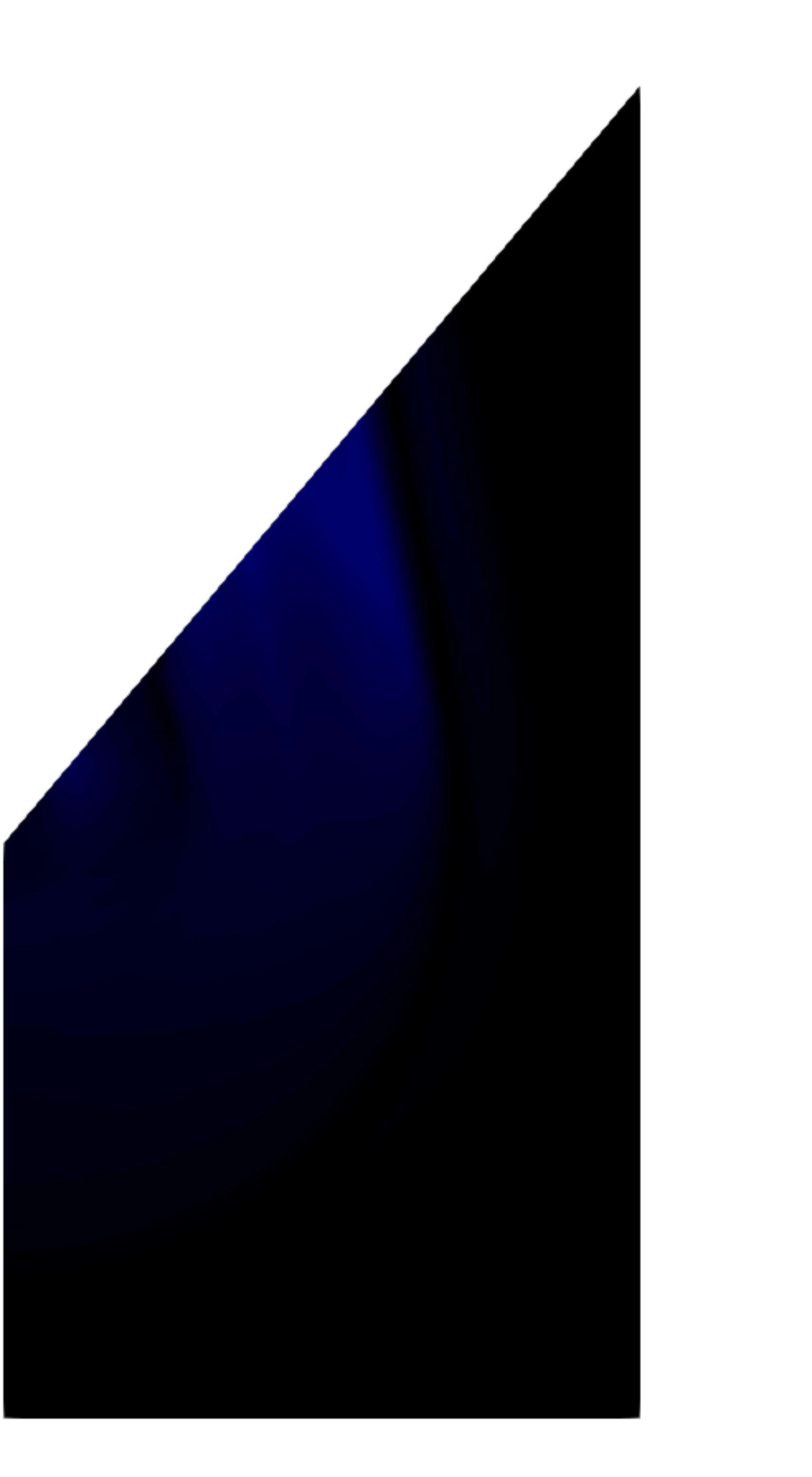}\includegraphics[scale=0.15, trim=1.5cm 0cm 38cm 0cm, clip]{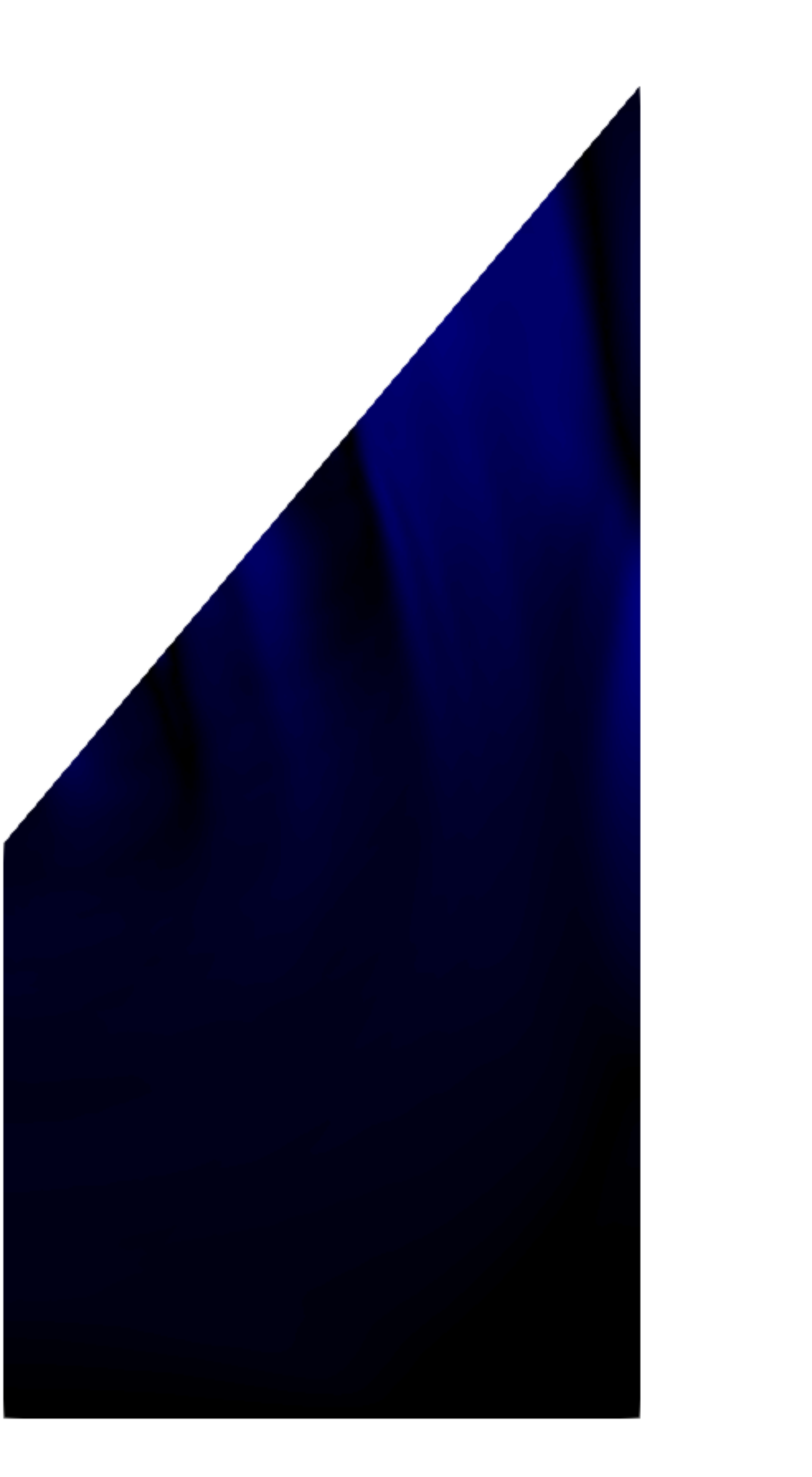}\includegraphics[scale=0.15, trim=1.5cm 0cm 38cm 0cm, clip]{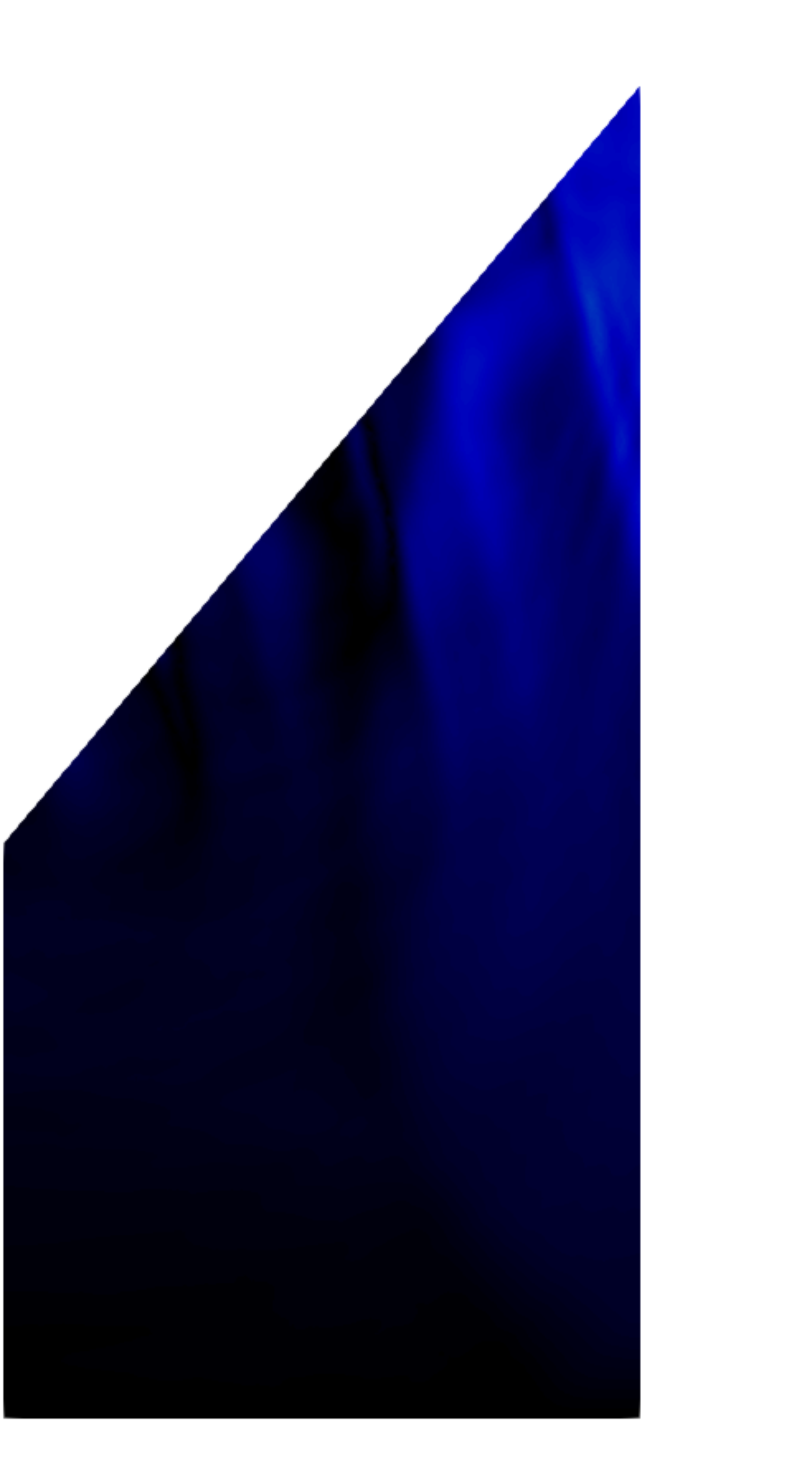}
\caption{Potential Difference $|\psi(t)-\psi^{\textup{ref}}(t)|$ plotted over time \textbf{(Horizontal)} for ABC$_W^{1/2}$ \textbf{(First row)} without adaptivity and \textbf{(Second row)} with self-adaptive angle.
\label{fig:Plane50Error}}
\end{center}
\end{figure}

\indent It is also interesting to see how the errors are distributed over the domain $\Omega$, i.e., where they originate from and how far they spread. Figure~\ref{fig:Plane50Error} shows the error fields $|\psi(t)-\psi^{\textup{ref}}(t)|$ at different times for the ABC$_W^{1/2}$ conditions with and without adaptivity. The first snapshots were taken just as the first wave hits the absorbing boundary. In the subsequent snapshots, we can see how the erroneous reflections travel together with the wave across the absorbing boundary. It is evident that there are fewer reflections present when using the self-adaptive absorbing boundary conditions that take the local angle information into account.

\begin{figure}[h!]
\begin{center}
\input{images/Plane50/Cost.tex}\hspace*{2mm} \input{images/Plane50/Angle_Quality.tex}
\caption{\textbf{Inclined plane boundary:} Relative $L^2(\Omega)$ error of the potential $\psi(t)$ over the simulation time with $\theta=50^{\circ}$ . \textbf{(Left)} Nonlinear vs. linear conditions with and without adaptivity, \textbf{(Right)} Performance of adaptive conditions with numerically computed vs. the exact incidence angle $\theta$.
\label{fig:Plane50Cost}}
\end{center}
\end{figure}
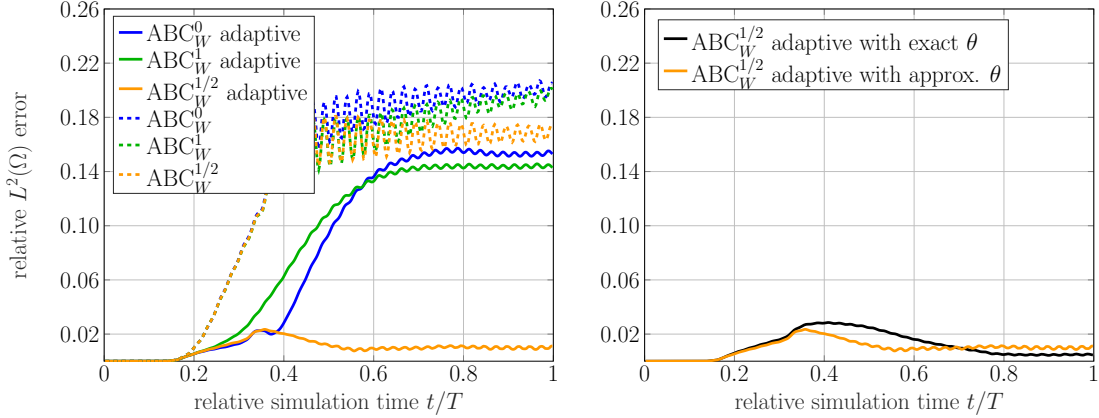

\subsection{Higher source frequency}
So far our experiments have been conducted with the same source excitation. After varying the domain geometry via the angle $\theta$, we investigate the influence of the source frequency on the quality of absorbing conditions. We now test with the excitation frequency $$f =250 \, \textup{kHz}.$$ 

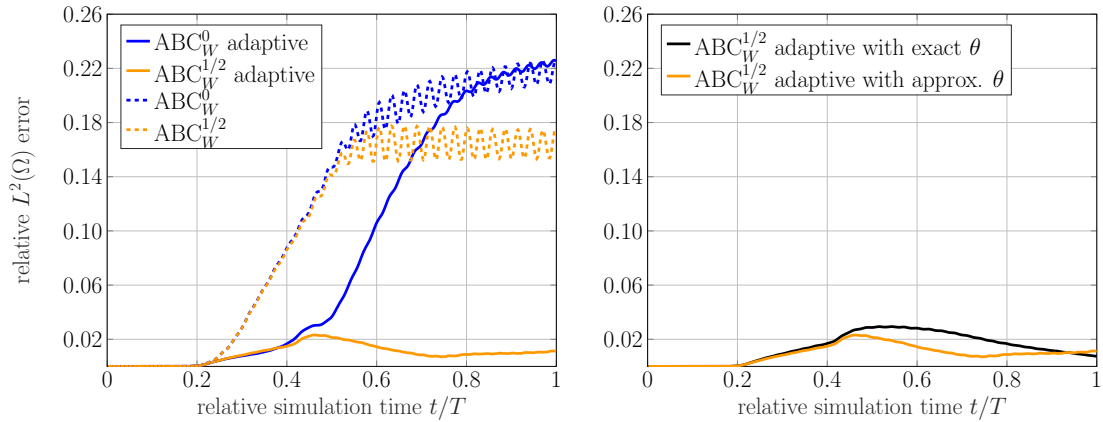
\begin{figure}[h!]
\begin{center}
\input{images/HighFq/Cost.tex}\hspace*{2mm} \input{images/HighFq/Angle_Quality.tex}
\caption{\textbf{Inclined plane boundary:} Relative $L^2(\Omega)$ error of the potential $\psi(t)$ for $t \in [0,T]$ with $\theta=50^{\circ}$ and higher source frequency. \textbf{(Left)} Nonlinear vs. linear conditions with and without adaptivity, \textbf{(Right)} Performance of adaptive conditions with numerically computed vs. the exact incidence angle $\theta$.
\label{fig:HighFqCost}}
\end{center}
\end{figure}

All the remaining parameters being the same as before. The domain again has the upper boundary tilted with the angle of $50\,^{\circ}$. Figure~\ref{fig:HighFqCost} shows the error plots for the higher frequency. We observe that conditions ABC$^{0}_W$ asymptotically show the same poor results as conditions ABC$^0_{W}$ with self-adaptivity. Only our combination of nonlinear conditions with the self-adaptive technique \enlargethispage*{3cm} (\ref{ABC_Westervelt_angle_1/2}) allows for the relative error to stay around $1 \, \%$. By comparing Figure~\ref{fig:Plane50Cost} and Figure~\ref{fig:HighFqCost}, we also notice that the difference in the error between the nonlinear and the corresponding linear conditions increases with the frequency, emphasizing the need to employ nonlinear conditions for high-frequency nonlinear sound waves. 

\FloatBarrier

\FloatBarrier
\newpage

\subsection{Plate with a hole} \label{SubSec:PlateHole}
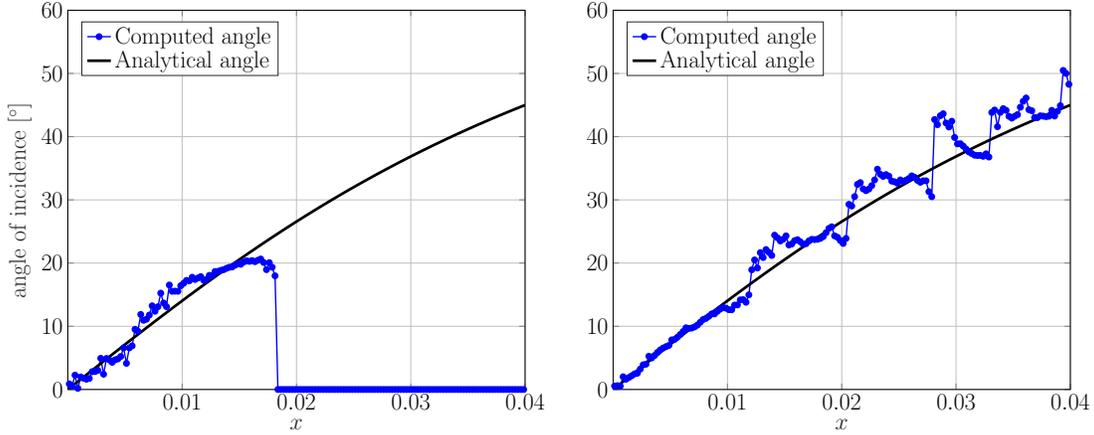
\begin{figure}[h!]
\begin{center}
\input{images/Octant/Angle_Distribution0.tex}\hspace*{2mm} \input{images/Octant/Angle_Distribution.tex}
\caption{Computed versus exact angle of incidence $\theta$. \textbf{(Left)} The angle is only computed on the parts of the absorbing boundary that the wave has reached. \textbf{(Right)} Angle computation towards the end of the simulation.
\label{fig:Octant_AngleDistribution}}
\end{center}
\end{figure}

We previously tested the new absorbing boundary conditions in a domain where the angle of incidence was constant over the absorbing boundary. To show that both our approach of computing the angle of incidence of the wave as well as the nonlinear boundary conditions work in more realistic settings, we now consider the so-called ``\textit{plate with a hole}'' domain. It consists of a square with a circular hole in the center. In our case, the excitation of the wave takes place at the boundary of the hole. By using symmetry, we reduce the simulation of the whole domain to half of one of its quarters; see~Figure \ref{fig:Domains}.\\
\indent An analytical expression for the angle of incidence is also available here which allows us to judge the quality of our angle approximation. If the origin is at the center of the circular hole and the square has sides of length $a$, the angle of incidence is given by
\begin{align} \label{angle_formula}
\theta(x, y, t) = \arccos\left(\frac{a/2}{\sqrt{x^2+(a/2)^2}}\right).
\end{align}
We set $a=0.08$ in the experiments. The domain is resolved with 13820 degrees of freedom in space. For time discretization, we choose 8330 time steps with final time $T = 8.0325\cdot 10^{-5} \,$s.\\
\indent Figure~\ref{fig:Octant_AngleDistribution} shows the angle of incidence computed by our approach at two different points in time and further illustrates our criterion on computing the incidence angle based on the amplitude of the wave on the boundary. Note that formula \eqref{angle_formula} assumes a reflection-free potential field so an exact match of our angle to the analytical angle distribution cannot be expected. In fact, our method tries to compensate also for the waves that originate as spurious reflections from the absorbing boundary.

\FloatBarrier

\begin{figure}[h!]
\begin{center}
\includegraphics[scale=0.09, trim=0cm 0cm 0cm 0cm, clip]{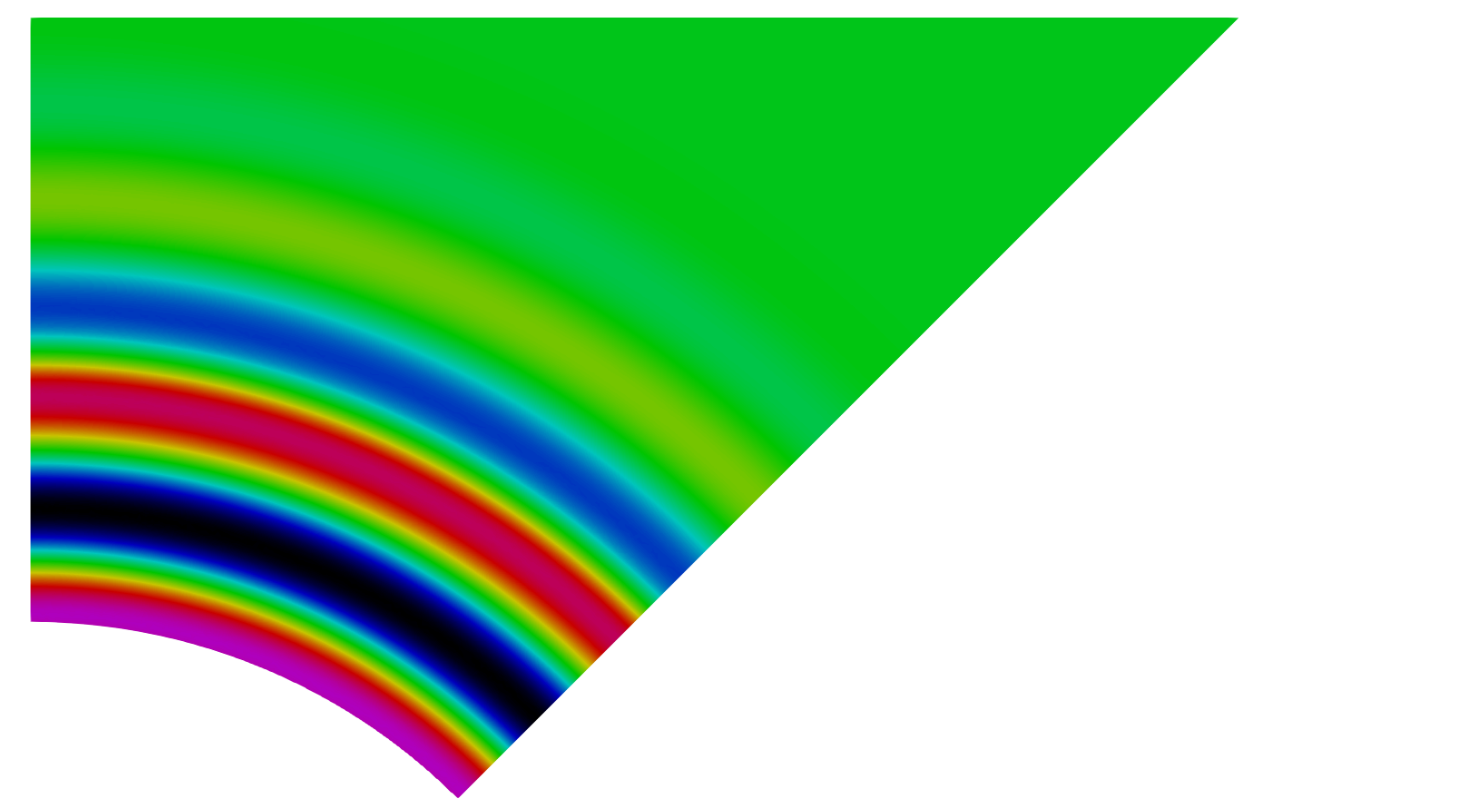}
\includegraphics[scale=0.09, trim=0cm 0cm 8cm 0cm, clip]{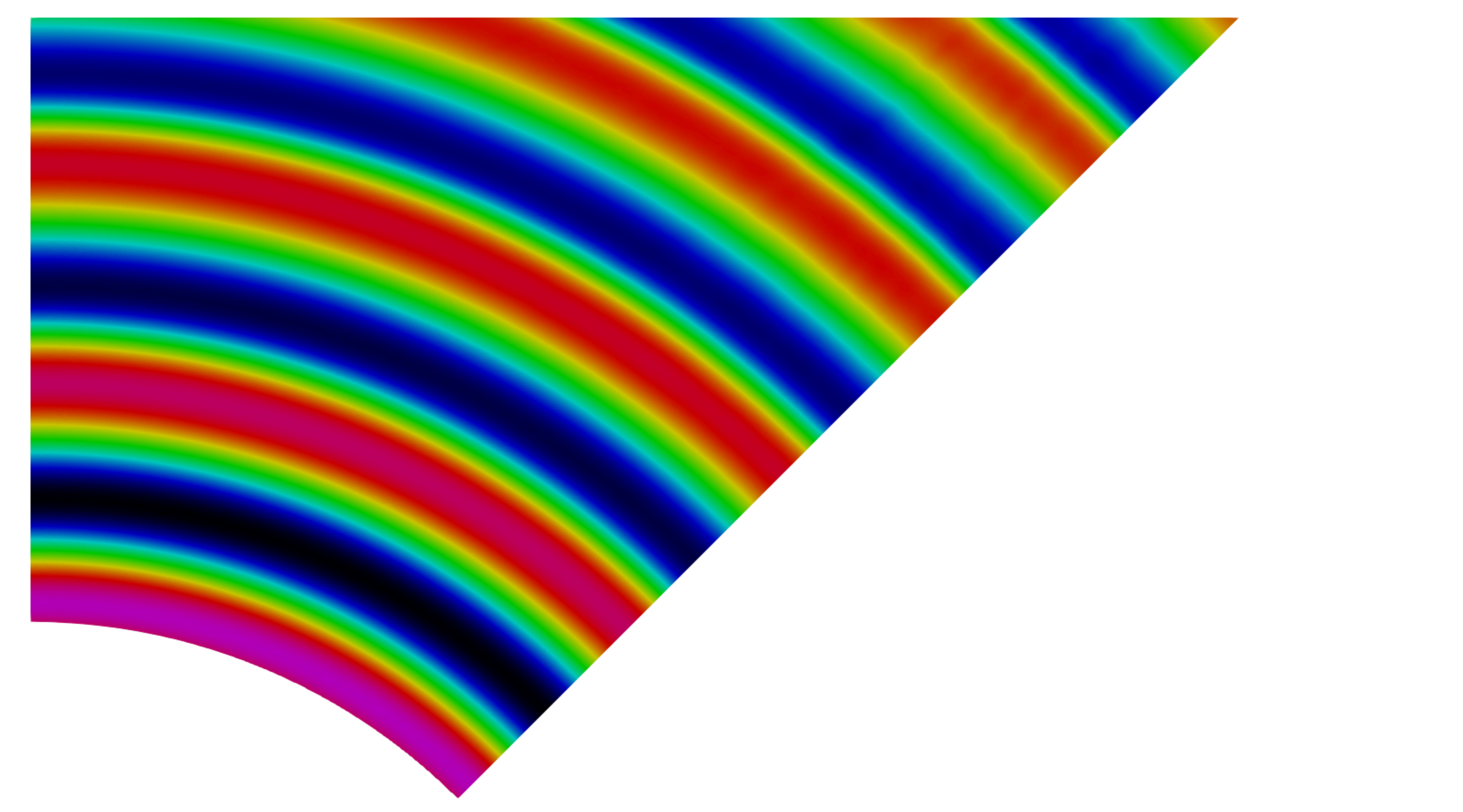}
\caption{Time snapshots of the pressure field $u =\rho\psi_t$.
\label{fig:OctantPressure}}
\end{center}
\end{figure}

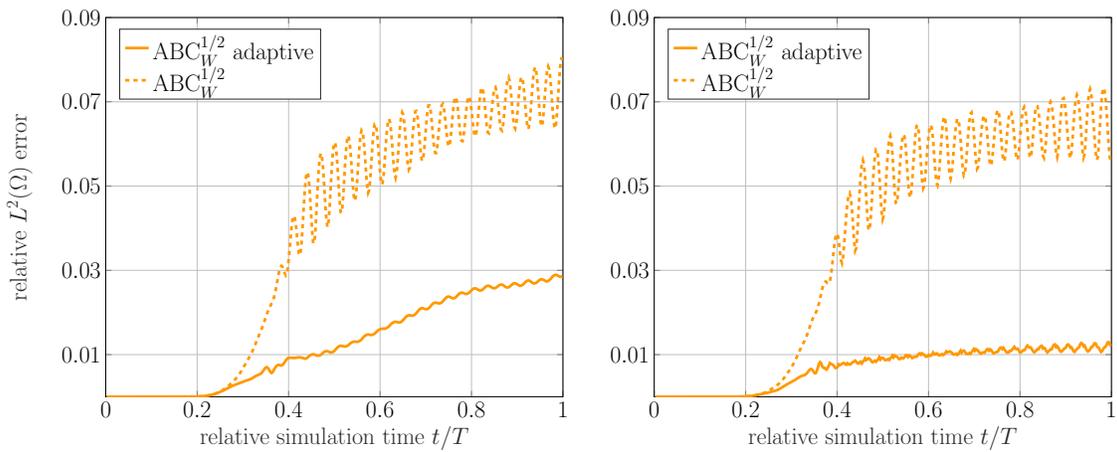
\begin{figure}[h!]
\begin{center}
\input{images/Octant/Cost.tex}\hspace*{2mm} \input{images/Octant/Cost_V.tex}
\caption{\textbf{Plate with a hole:} Relative $L^2(\Omega)$ error of \textbf{(Left)} the potential $\psi(t)$ and \textbf{(Right)} the pressure $u(t)=\rho\psi_t(t)$ over the simulation time.
\label{fig:OctantCost}}
\end{center}
\end{figure}

\indent As in our last experiment, we proceed by showing the wave field at different time steps as well as the error comparison of different absorbing boundary conditions; see Figure~\ref{fig:OctantPressure} and Figure~\ref{fig:OctantCost}. This time we also compute them in terms of the acoustic pressure $u=\rho\psi_t$ due to its practical importance. We also introduce here the relative errors in the $L^2(0,T; L^2(\Omega))$ norm, i.e.  
 
$$e_{\psi}=\frac{\|\psi-\psi^{\textup{ref}}\|_{L^2(0,T; L^2(\Omega))}}{\|\psi^{\textup{ref}}\|_{L^2(0,T; L^2(\Omega))}},\qquad e_{u}=\frac{\|u-u^{\textup{ref}}\|_{L^2(0,T; L^2(\Omega))}}{\|u^{\textup{ref}}\|_{L^2(0,T; L^2(\Omega))}}.$$
 
In the present experiment, those errors amount to $e_{\psi} = 1.82\,$\% and $e_u = 0.93\,$\% for the self-adaptive conditions and $e_{\psi} = 5.41\,$\% and $e_u = 5.17$\,\% for the conditions without adaptivity, giving an overall improvement of $66.36\,$\% in the potential and $82.01\,$\% in the pressure.

\newpage

\begin{figure}[h!]
\begin{center}
\includegraphics[scale=0.09, trim=0cm 0cm 8cm 0cm, clip]{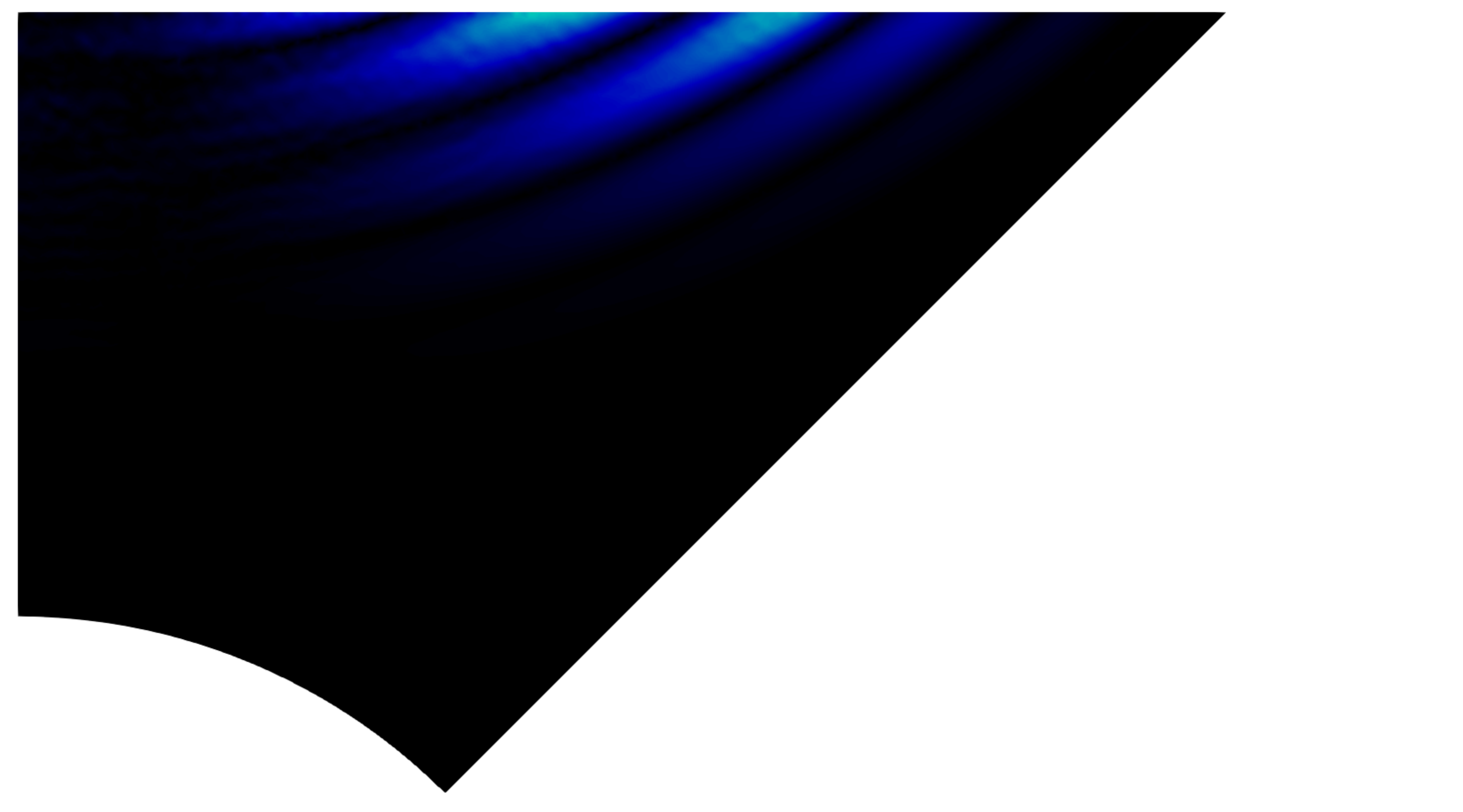}\includegraphics[scale=0.09, trim=0cm 0cm 8cm 0cm, clip]{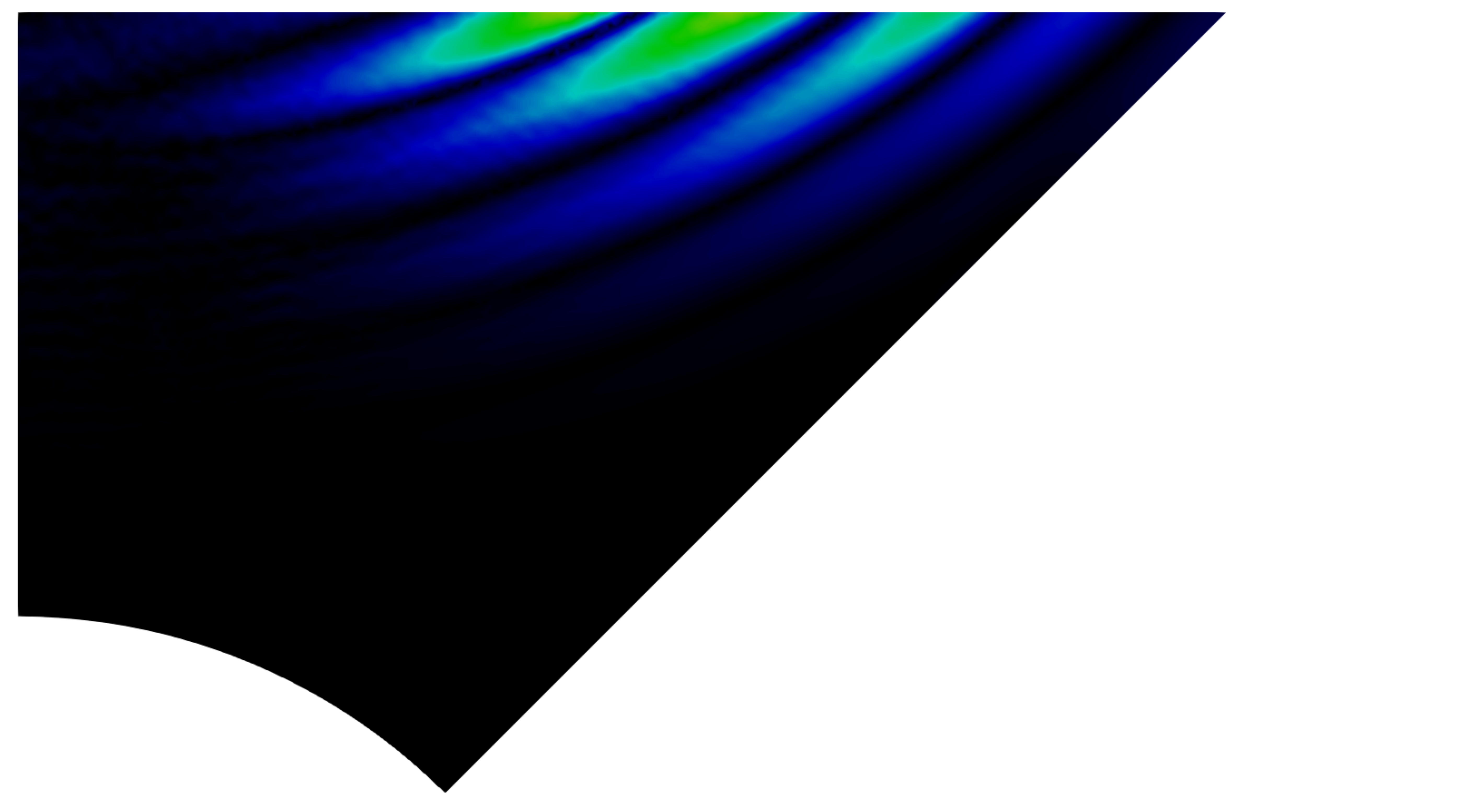}\includegraphics[scale=0.09, trim=0cm 0cm 8cm 0cm, clip]{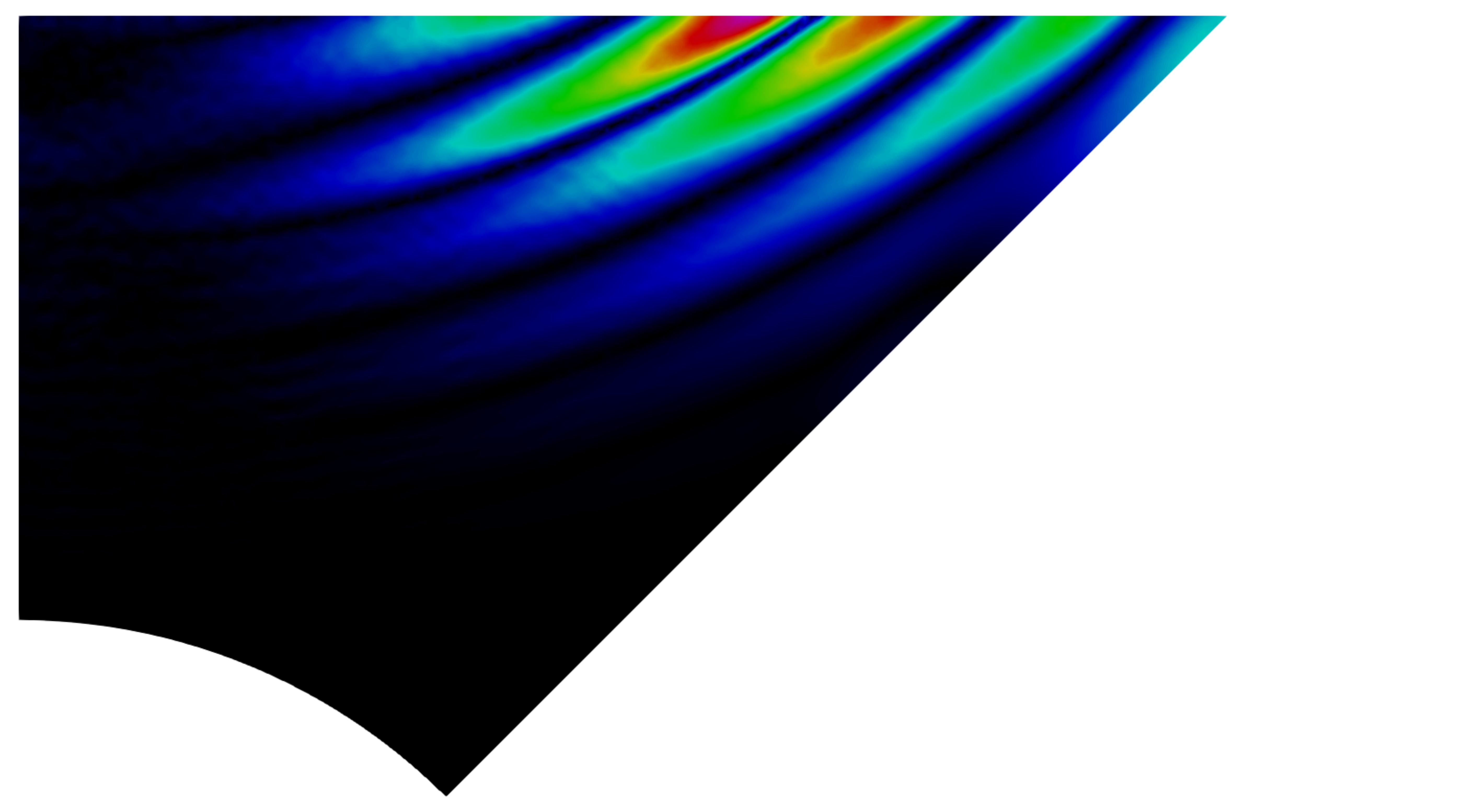}\\
\includegraphics[scale=0.09, trim=0cm 0cm 8cm 0cm, clip]{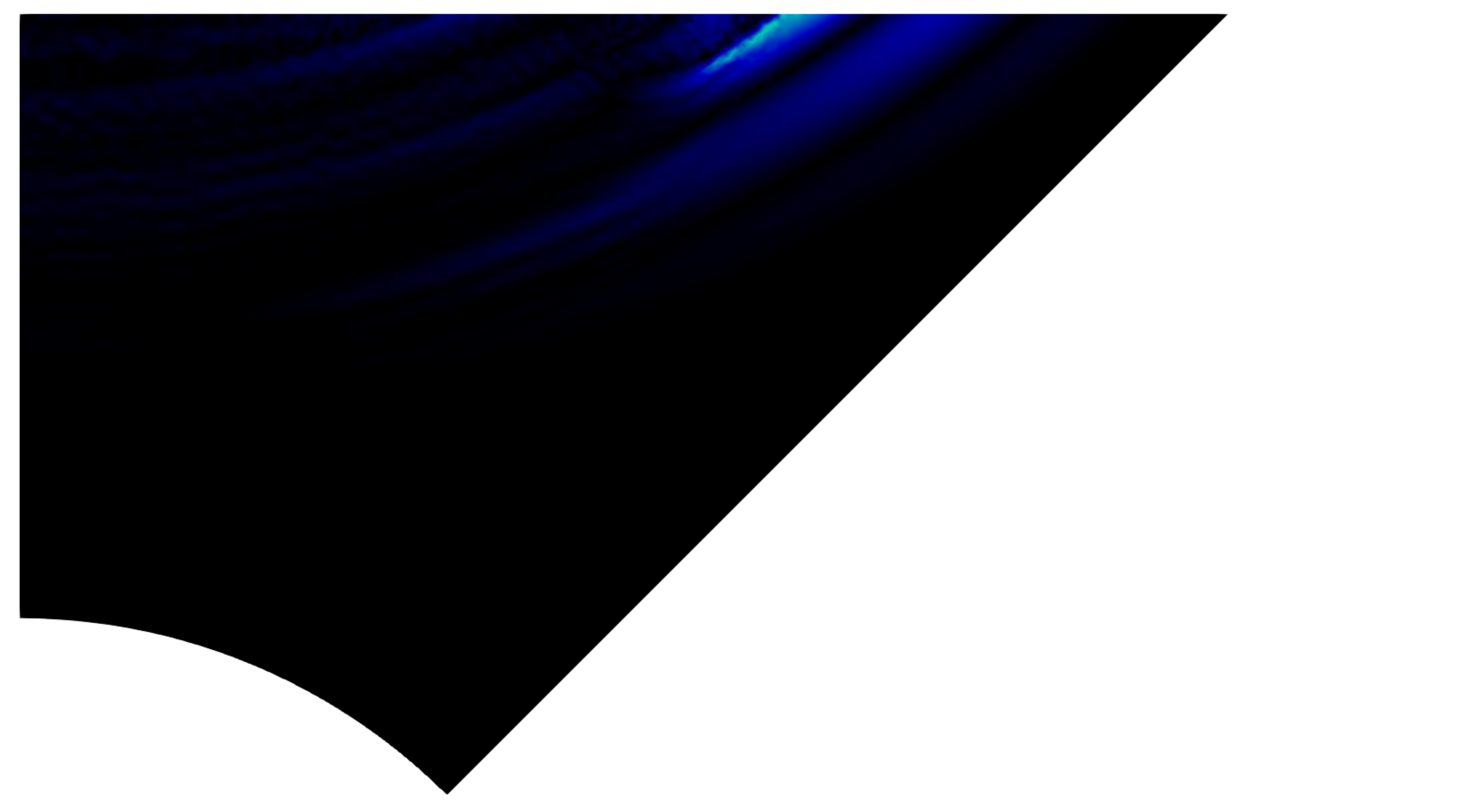}\includegraphics[scale=0.09, trim=0cm 0cm 8cm 0cm, clip]{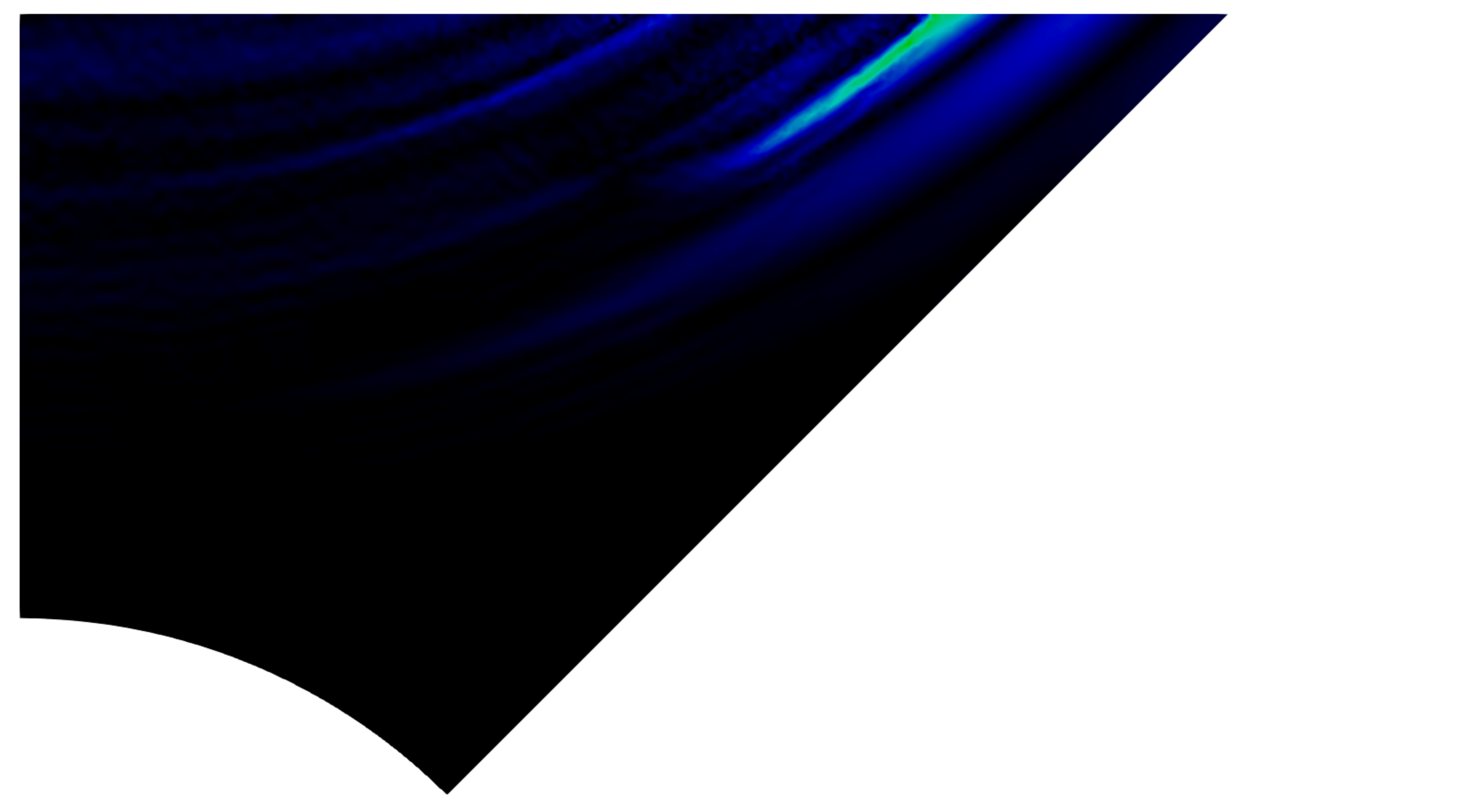}\includegraphics[scale=0.09, trim=0cm 0cm 8cm 0cm, clip]{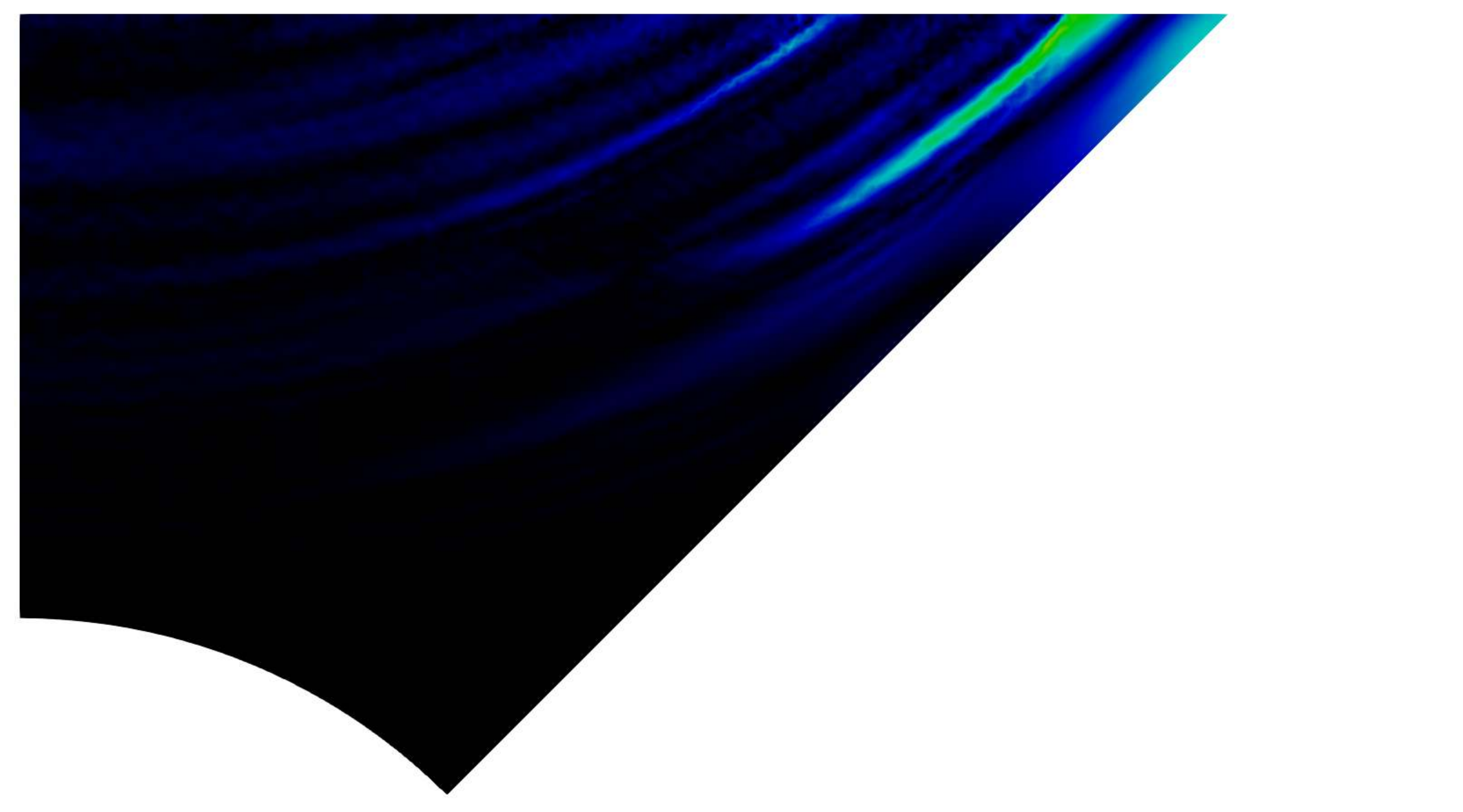}\\
\caption{Pressure difference $|u(t)-u^{\textup{ref}}(t)|$  plotted over time \textbf{(Horizontal)} for ABC$_W^{1/2}$ without \textbf{(First row)} and with \textbf{(Second row)} angle consideration.
\label{fig:OctantErrors}}
\end{center}
\end{figure}

\indent Plots of the difference field $|u(t)-u^{\textup{ref}}(t)|$ at different times are given in Figure~\ref{fig:OctantErrors}. We observe that with the absorbing conditions that are not adaptive, the amplitude of the error increases over the width of the domain as the angle grows.

\clearpage

\subsection{High-intensity focused ultrasound (HIFU)}\label{subsec:Transducer2D}
\begin{figure}[h!]
\begin{center}
\includegraphics[scale=0.1, trim=0cm 0cm 10cm 0cm, clip]{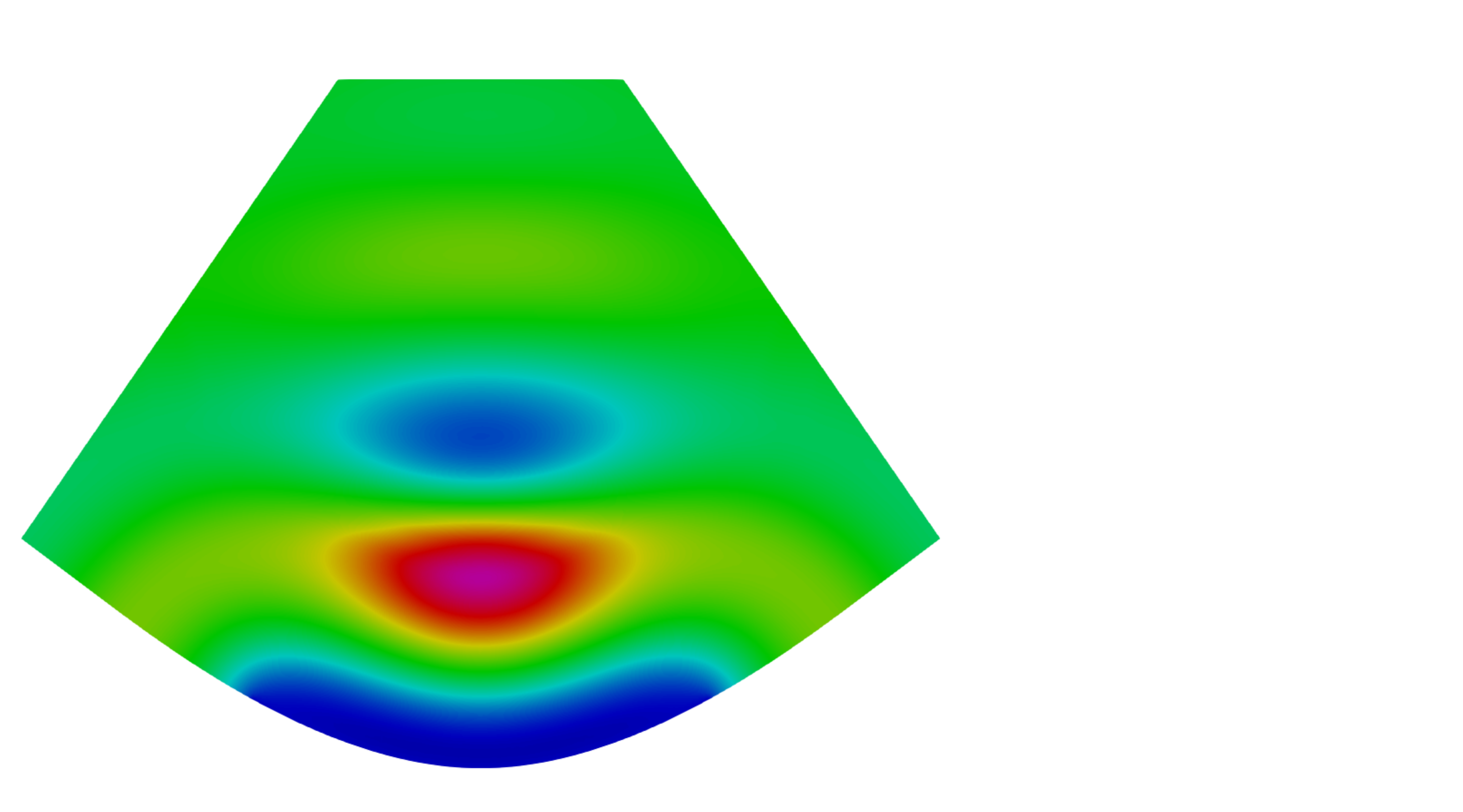}\includegraphics[scale=0.1, trim=0cm 0cm 10cm 0cm, clip]{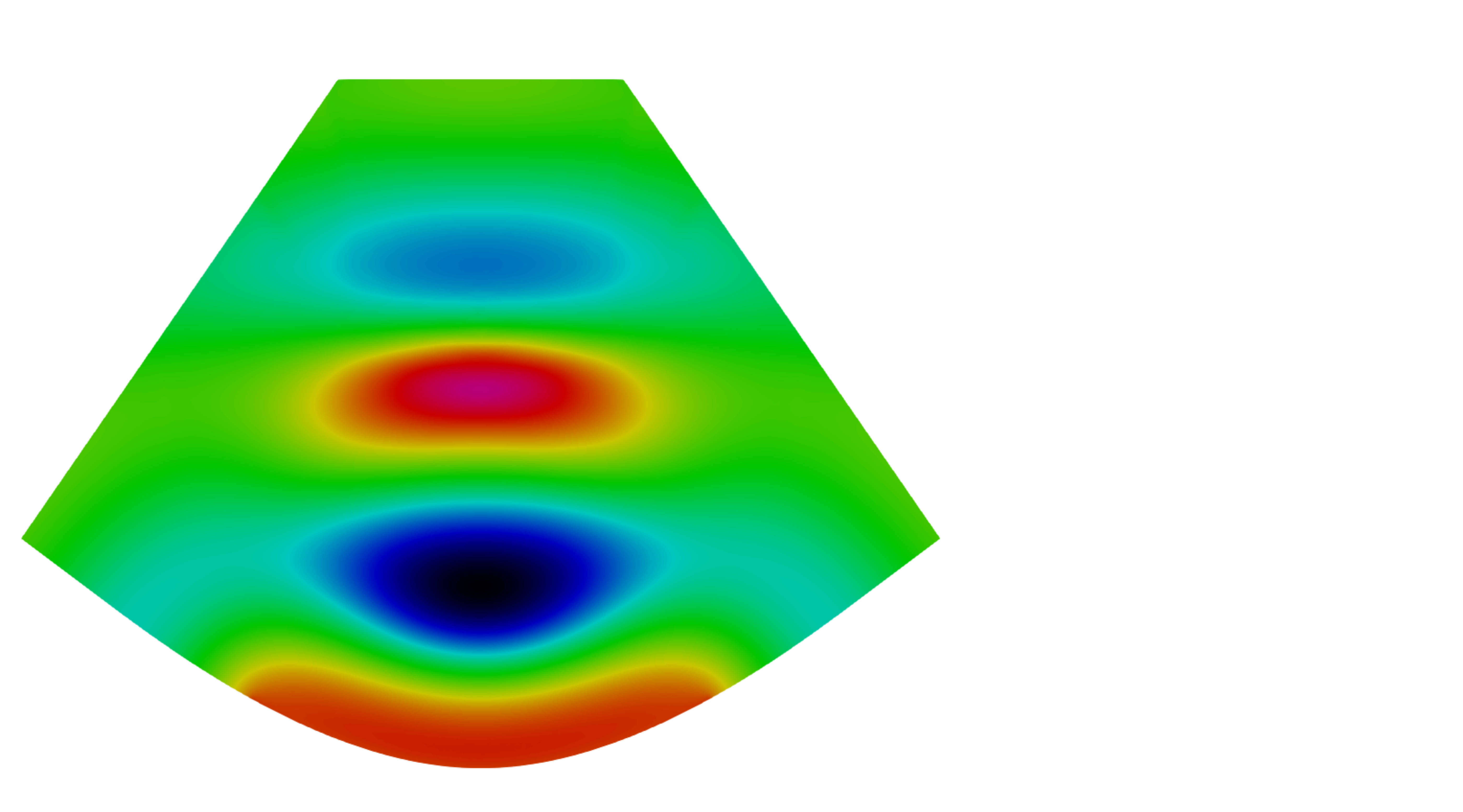}\includegraphics[scale=0.1, trim=0cm 0cm 10cm 0cm, clip]{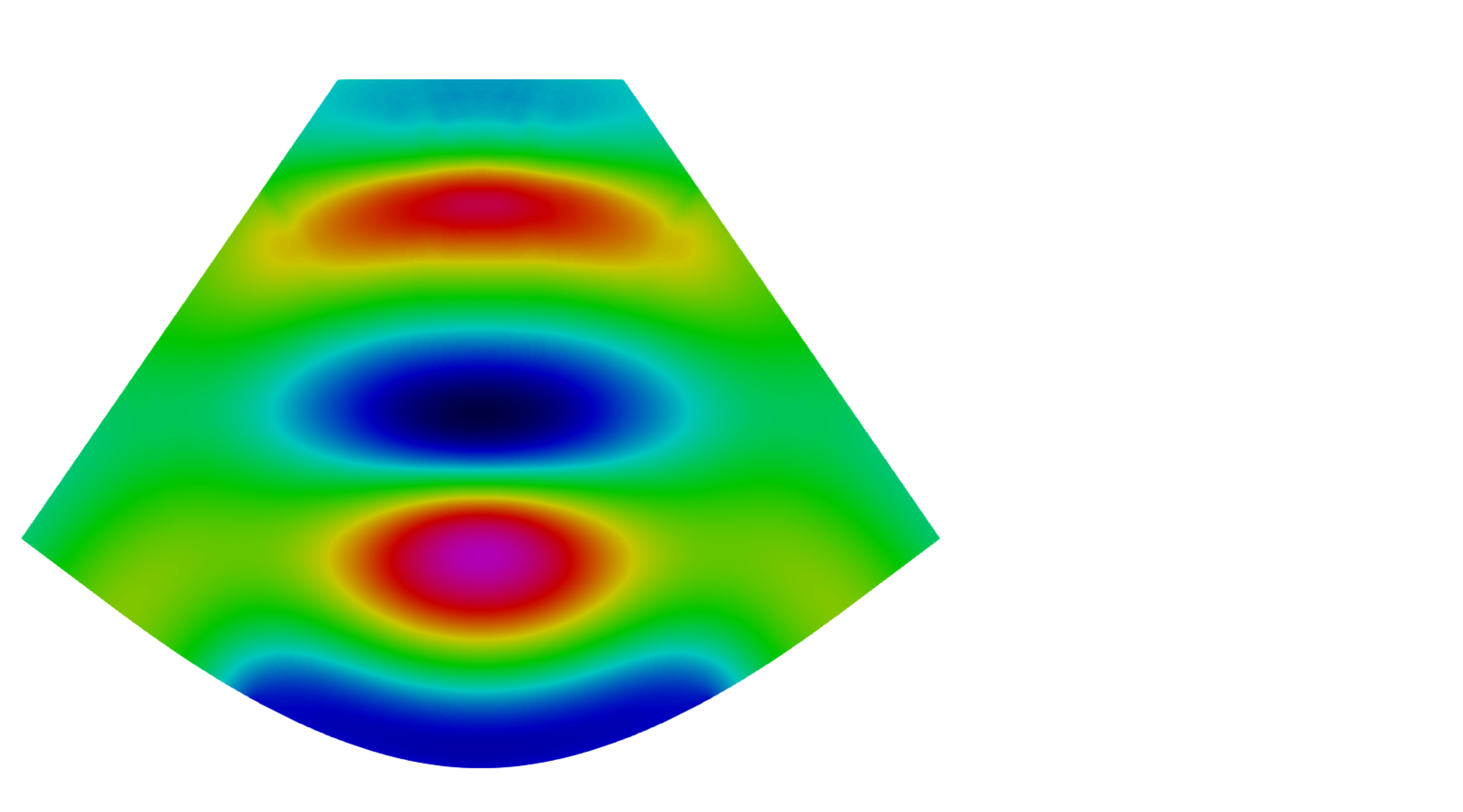}
\caption{Snapshots of the pressure field $u=\rho\psi_t$ of a propagating, self-focusing wave.
\label{fig:Transducer_wave}}
\end{center}
\end{figure}

We next simulate the potential field generated by a piezoelectric transducer. Such devices are made of small plates of a piezoelectric material that are aligned in an arc-shaped array pointing towards a common focal point~\cite{dreyer2000investigations, watanabe1994wave}. When set into motion, those vibrating plates induce acoustic sound waves that propagate towards the focal point. As it travels, the wave is focused more and more the closer it comes to the focal point. This technique of generating high-pressure amplitudes at specific locations is used in medicine to treat kidney stones and certain types of cancer; see~\cite{illing2005safety, kennedy2004high, kennedy2005high, wu2003preliminary, yoshizawa2009high}. The pressure levels in the non-focal region are sufficiently low so that damage to the surrounding tissue is avoided.\\
\indent For this experiment, the medium of propagation is again chosen to be water with the same physical parameters as before. The source frequency is again given by $f=210\,$kHz and the source amplitude by $\mathcal{A}=0.02\,\textup{m}^2/\textup{s}^2$ which increases at the focal point due to focusing. As depicted in Figure \ref{fig:Domains}, the absorbing boundary here consists out of three line segments at the left, right and top. Time-snapshots of the transducer simulation can be seen in Figure~\ref{fig:Transducer_wave}.\\
\indent The computational domain was resolved with 13313 degrees of freedom, while for the time stepping 9800 steps and a final time of $T=4.725\cdot 10^{-5}$\,s were used. \vspace*{-0.25cm}

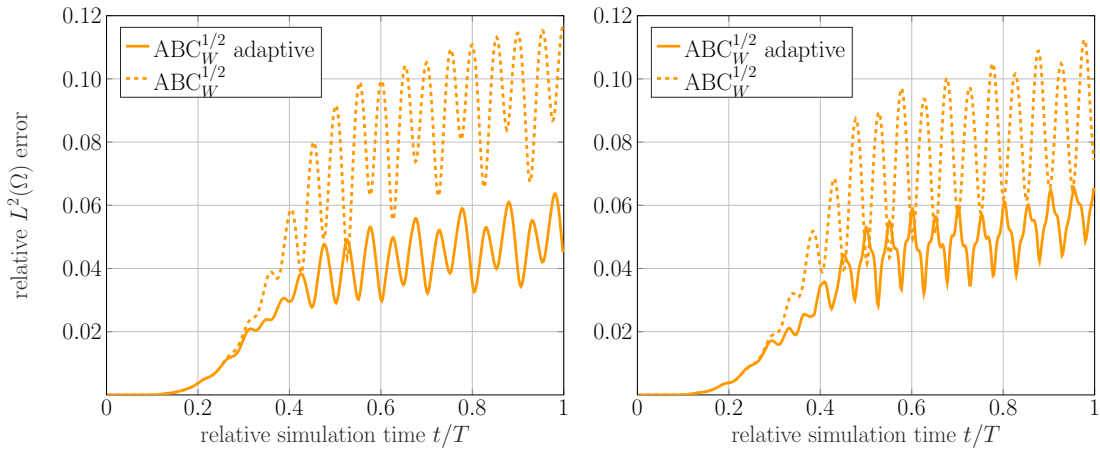
\begin{figure}[h!]
\begin{center}
\input{images/Transducer/Cost.tex}\hspace*{2mm}\input{images/Transducer/Cost_v.tex}
\caption{\textbf{HIFU transducer in 2D:} Relative $L^2(\Omega)$ error of \textbf{(Left)} the potential $\psi(t)$ and \textbf{(Right)} the pressure $u(t)=\rho\psi_t(t)$ over the simulation time.
\label{fig:TransducerCost}}
\end{center}
\end{figure}

A comparison of the relative $L^2(\Omega)$ errors at every time step is displayed in Figure~\ref{fig:TransducerCost}. Due to the relevance of measuring the acoustic pressure in HIFU applications, we again also plot the relative error that the absorbing boundary conditions produce in the pressure $u=\rho\psi_t$. We observe that the qualitative behavior of the errors and especially also the improvement made by the new adaptive conditions remains the same. For the adaptive conditions, the relative errors in the $L^2(0,T; L^2(\Omega))$ norm are $e_{\psi}=4.12\,$\% and $e_{u}=4.46\,$\%, whereas $e_{\psi}=7.94\,$\% and $e_{u}=7.27\,$\% if the adaptivity is not considered, resulting in an improvement of $51.89\,$\% in $\psi$ and even $61.35\,$\% in $u$. The increase in computational time when using the adaptive absorbing conditions amounts to $1.5 \, \%$. \vspace*{-0.2cm}

Since even with the absorbing boundary conditions an impinging wave is not completely absorbed, there are always some spurious reflections that travel back into the interior of the domain. They then interact with the ``main" wave that still travels towards the boundary and other spurious reflections originating from different parts of the boundary. Those effects together can accumulate over time, leading to an increase of the deviation from the reference solution and therefore the error, as can be observed in Figure~\ref{fig:TransducerCost}.  These effects appear to be more pronounced in more complicated geometries as well as with wave focusing.	
	

	\subsection{Multi-source wave superposition}\label{SubSec:MultiSource}
\begin{figure}[h!]
		\begin{center}
			\includegraphics[scale=0.18, trim=20cm 3cm 20cm 0cm, clip]{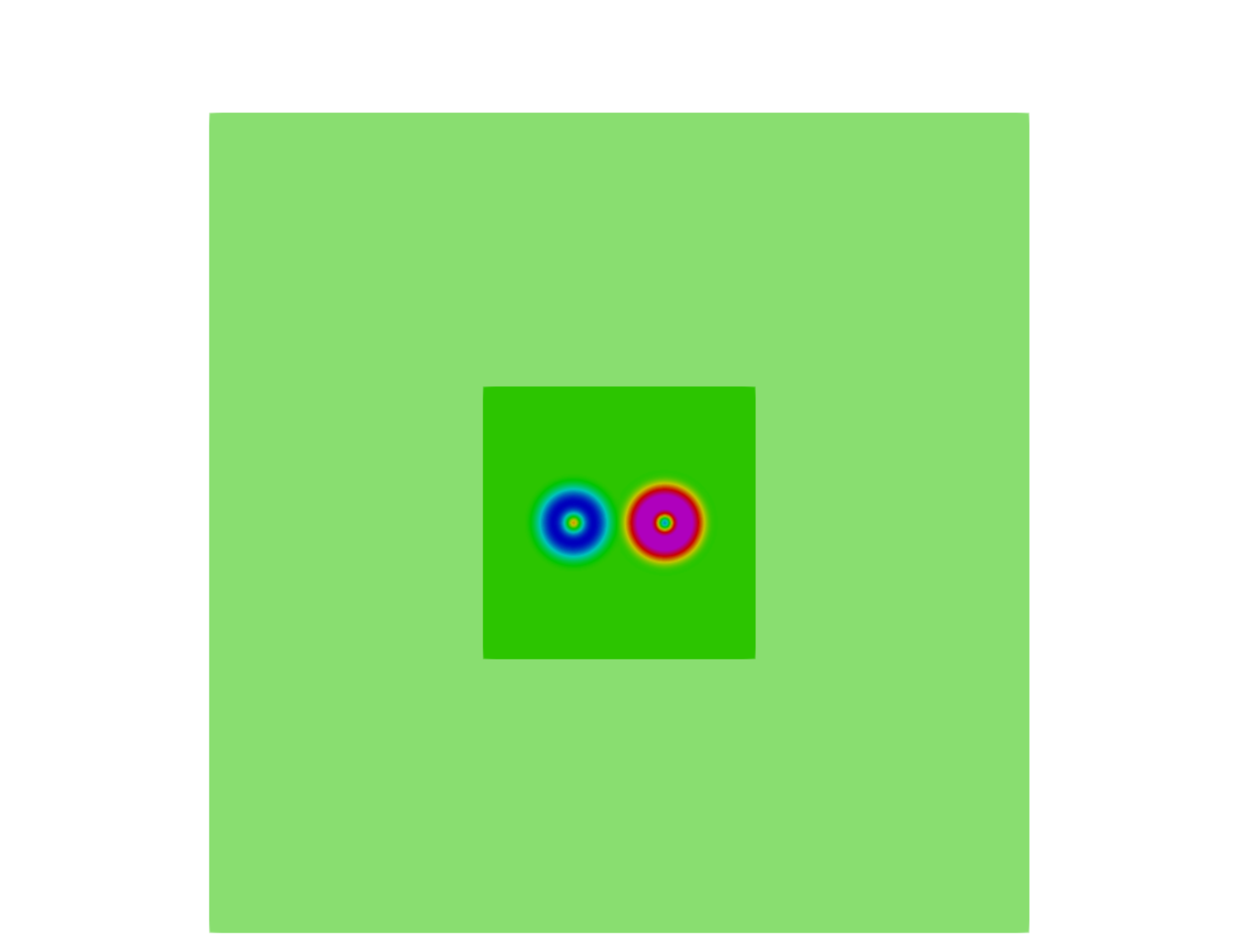}\includegraphics[scale=0.18, trim=20cm 3cm 20cm 0cm, clip]{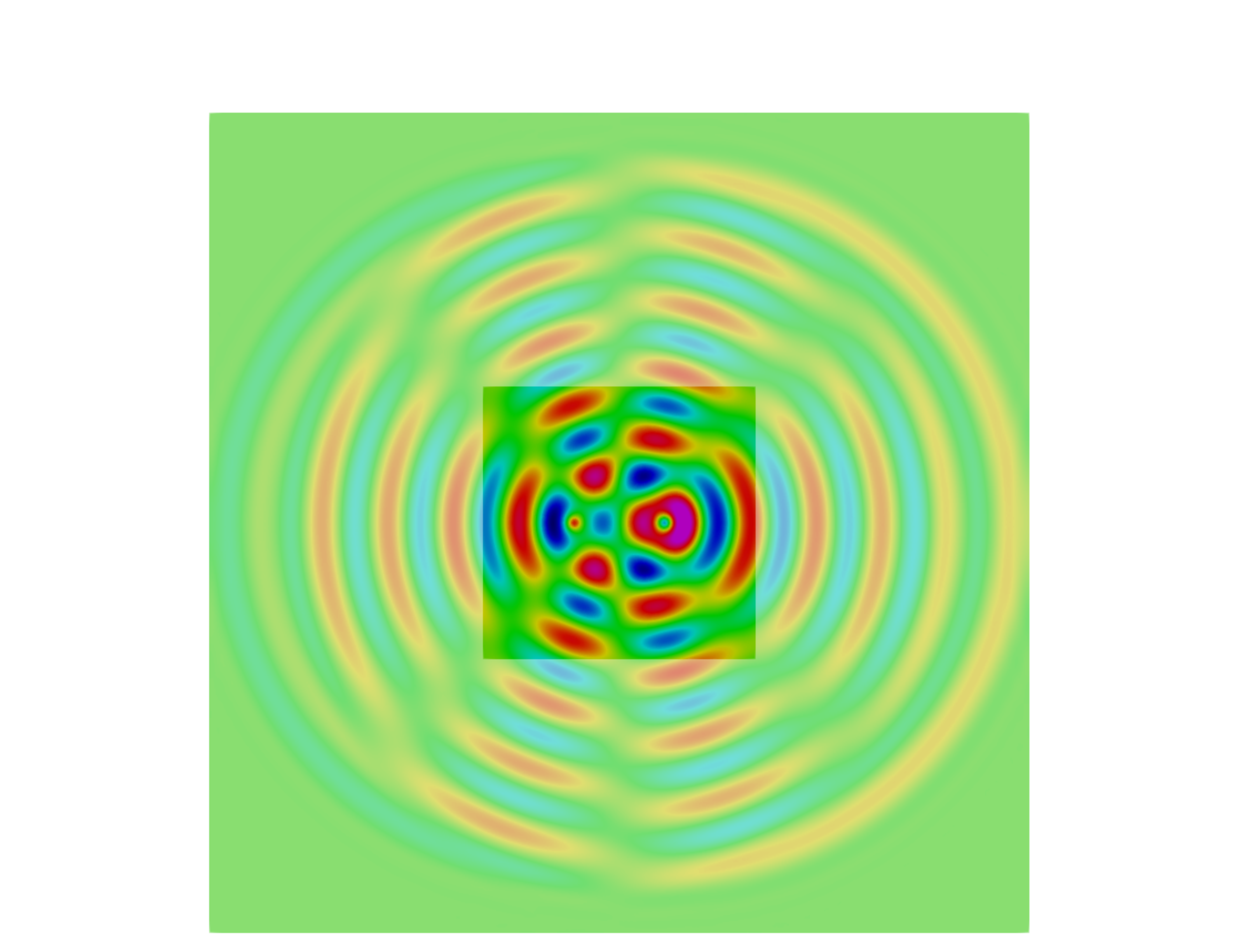}
			\caption{ Wave propagation induced by two different source terms.
				\label{fig:MultiSourceWave}}
		\end{center}
	\end{figure}

\begin{figure}[h!]
		\begin{center}
			\input{images/MultiSource/Cost.tex}
			\caption{ \textbf{Multi-source wave superposition:} Relative $L^2(\Omega)$ error of the potential $\psi(t)$
				\label{fig:MultiSourceCost}}
		\end{center}
	\end{figure}
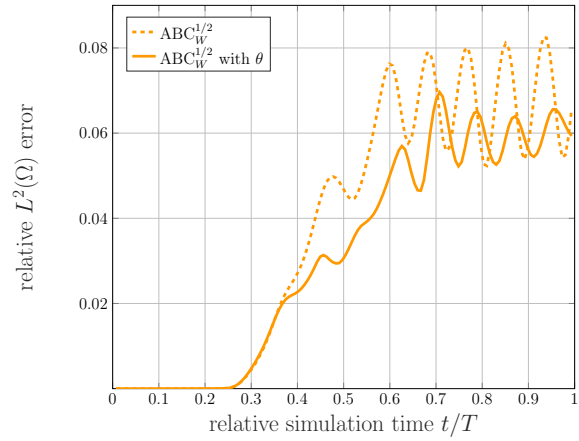

	Our next experiment is intended to illustrate that the adaptive method also works in scenarios with more than one wave source present and when superposition of waves occurs. In such cases, we can expect a less distinct wave propagation direction. In contrast to the other examples, here we use a source term on the right-hand side of the Westervelt equation. The computational domain is given by the square $\Omega = (0,0.03)^2$ and we choose the source term as follows
	\begin{align*}
	    f(x,y,t) &= \mathfrak{A}\sin(\omega t)\, \left[\exp\left({-\left(\frac{x-x_{\textup{mp1}}}{\sigma_x}\right)^2-\left(\frac{y-y_{\textup{mp1}}}{\sigma_y}\right)^2}\right)\right.\\
	    &~~~ \left.- \frac{2}{3}\exp\left({-\left(\frac{x-x_{\textup{mp2}}}{\sigma_x}\right)^2-\left(\frac{y-y_{\textup{mp2}}}{\sigma_y}\right)^2}\right)\right],
	\end{align*}
	with $\sigma_x = \sigma_y = 0.0005$, $x_{\textup{mp1}}=0.02,x_{\textup{mp2}}=0.01$, $y_{\textup{mp1}}=y_{\textup{mp2}}=0.015$, and $\mathfrak{A} = 10^{11}\, \textup{m}^2/\textup{s}^4$. Figure~\ref{fig:MultiSourceWave} depicts wave propagation at two different time induced by two different source terms within the computational domain. The transparent region depicts the reference domain. In Figure~\ref{fig:MultiSourceCost}, we can see the relative $L^2$ errors. We observe that also in this setting the adaptive approach results in a smaller error and a better overall behavior.
	
\subsection{3D Transducer} \label{SubSec:FocusedSound3D}
\begin{figure}[h!]
\begin{center}
\includegraphics[scale=0.1, trim=0cm 0cm 10cm 0cm, clip]{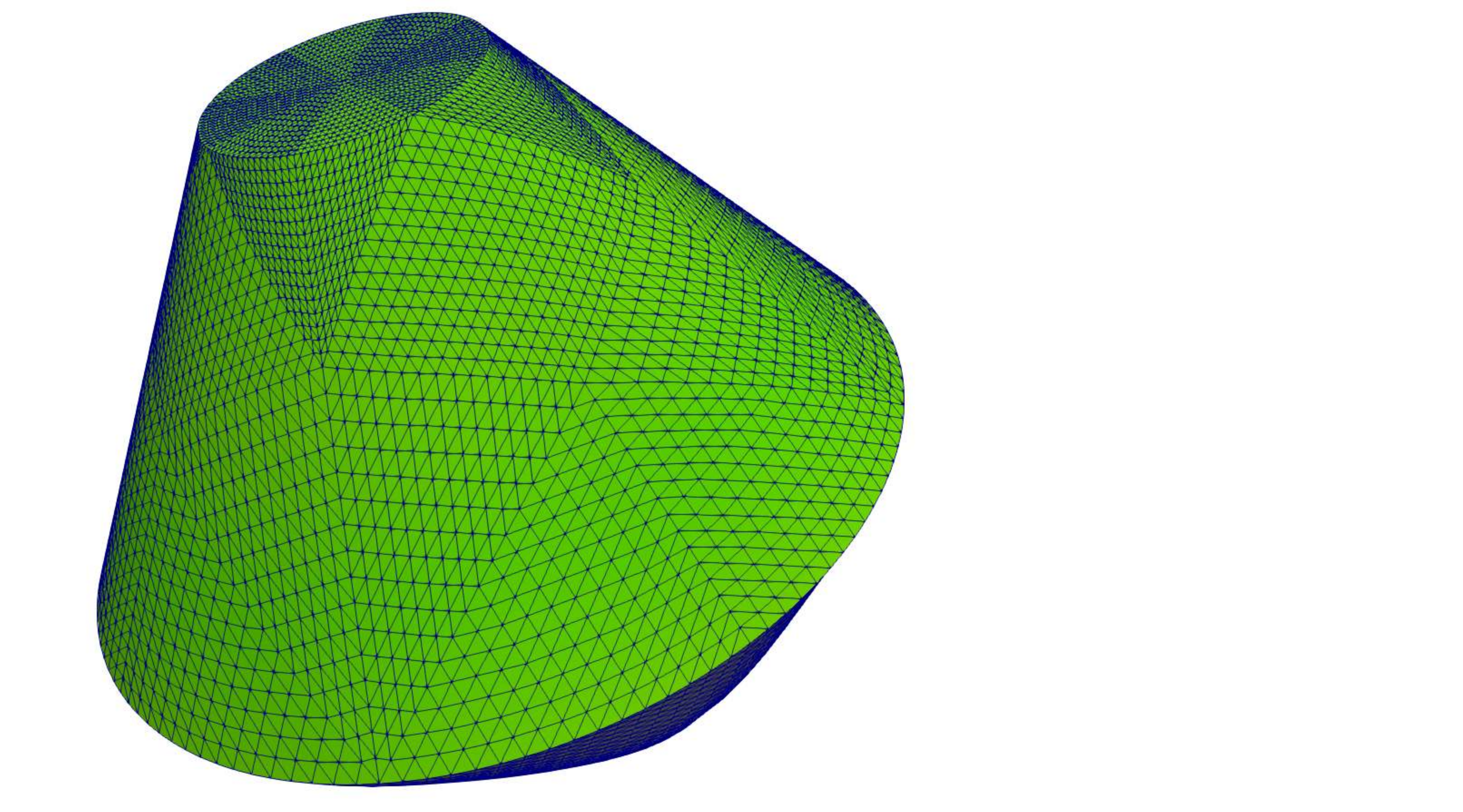}\includegraphics[scale=0.1, trim=0cm 0cm 10cm 0cm, clip]{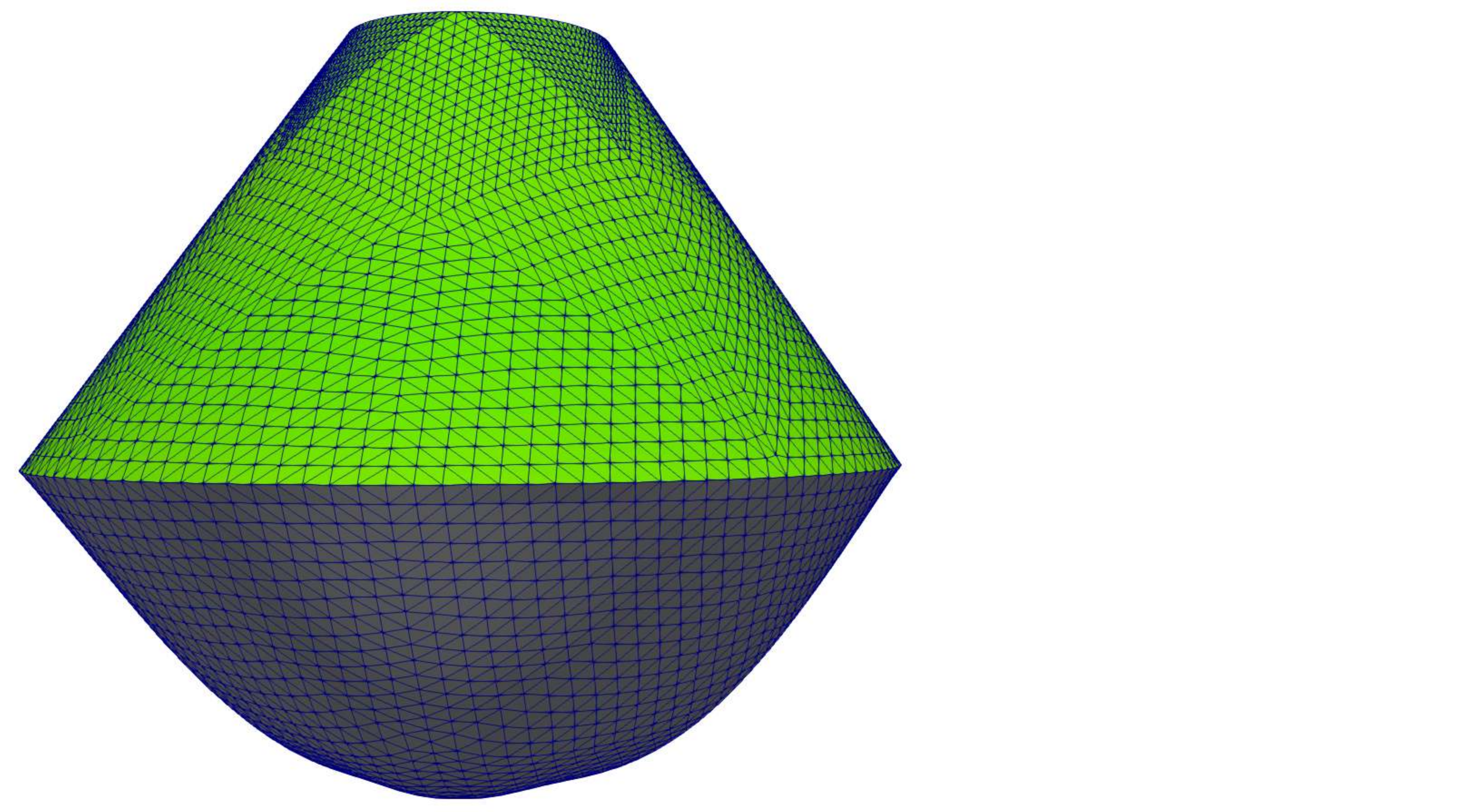}\includegraphics[scale=0.12, trim=0cm 2cm 16cm 0cm, clip]{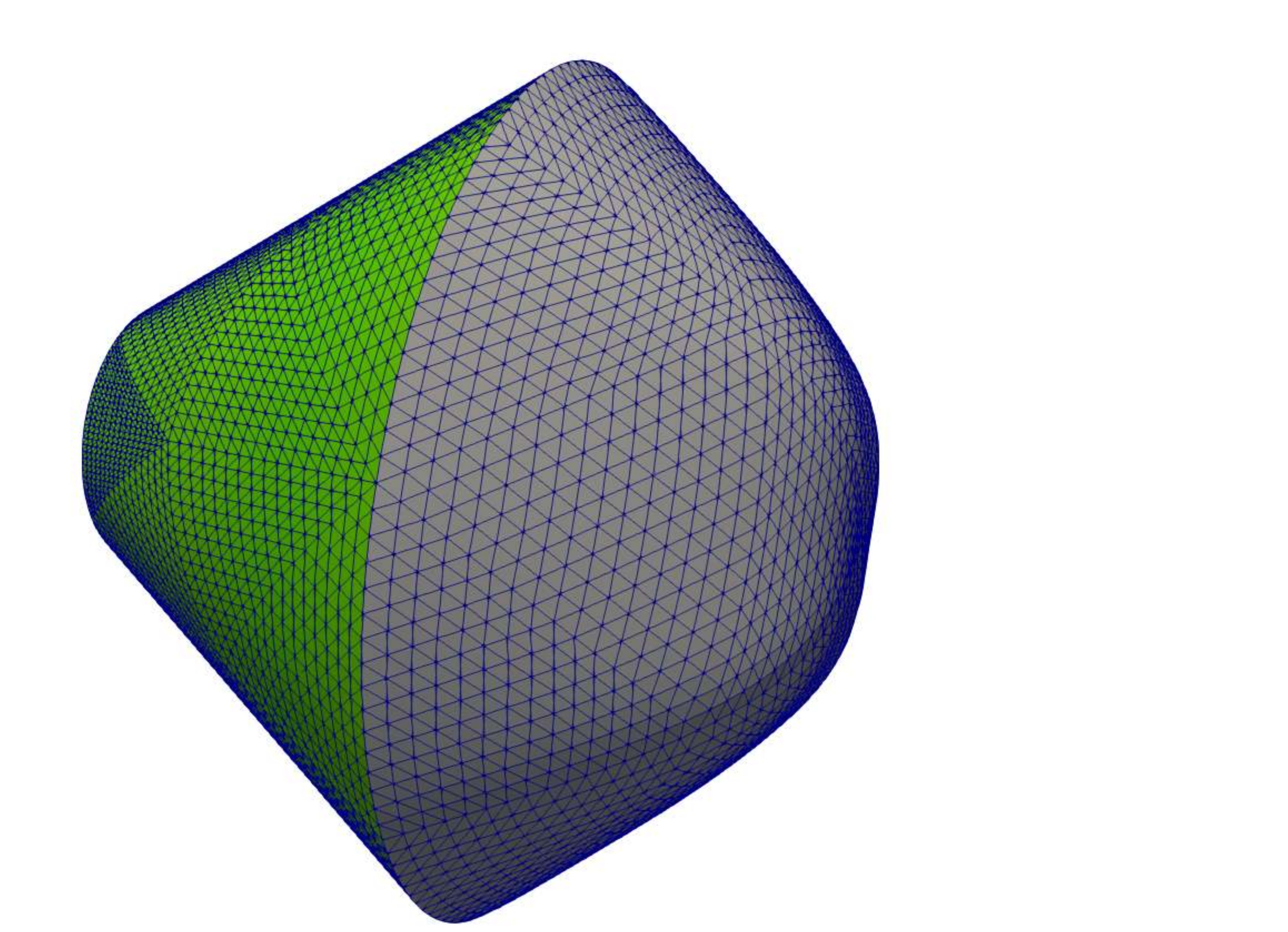}
\caption{Three-dimensional transducer geometry with mesh visible on the surface. Absorbing conditions are employed on the green surface, Dirichlet conditions on the grey.
\label{fig:Transducer3D_mesh}}
\end{center}
\end{figure}

\begin{figure}[h!]
\begin{center}
\includegraphics[scale=0.13, trim=9cm 0cm 14cm 0cm, clip]{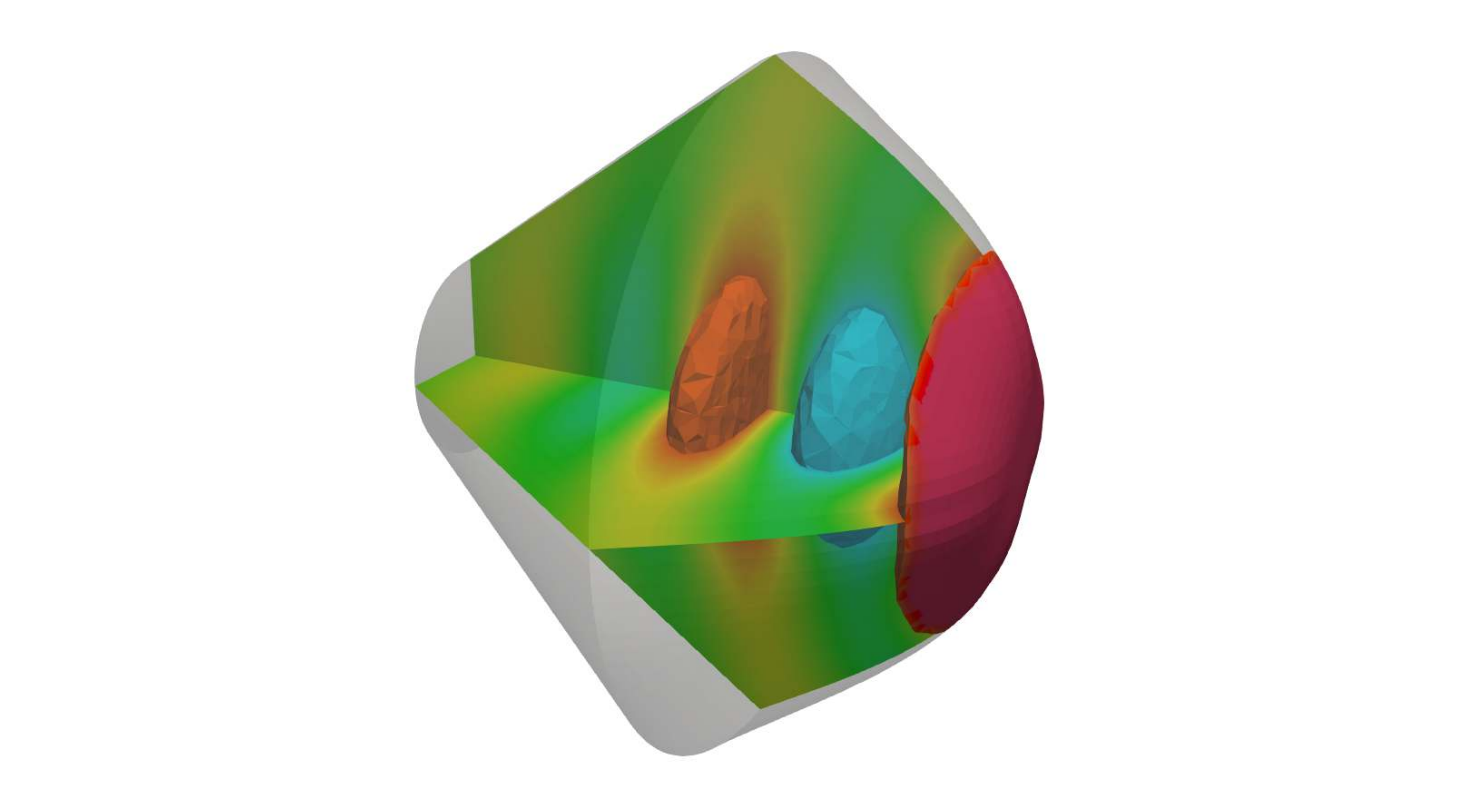}\includegraphics[scale=0.13, trim=7cm 0cm 14cm 0cm, clip]{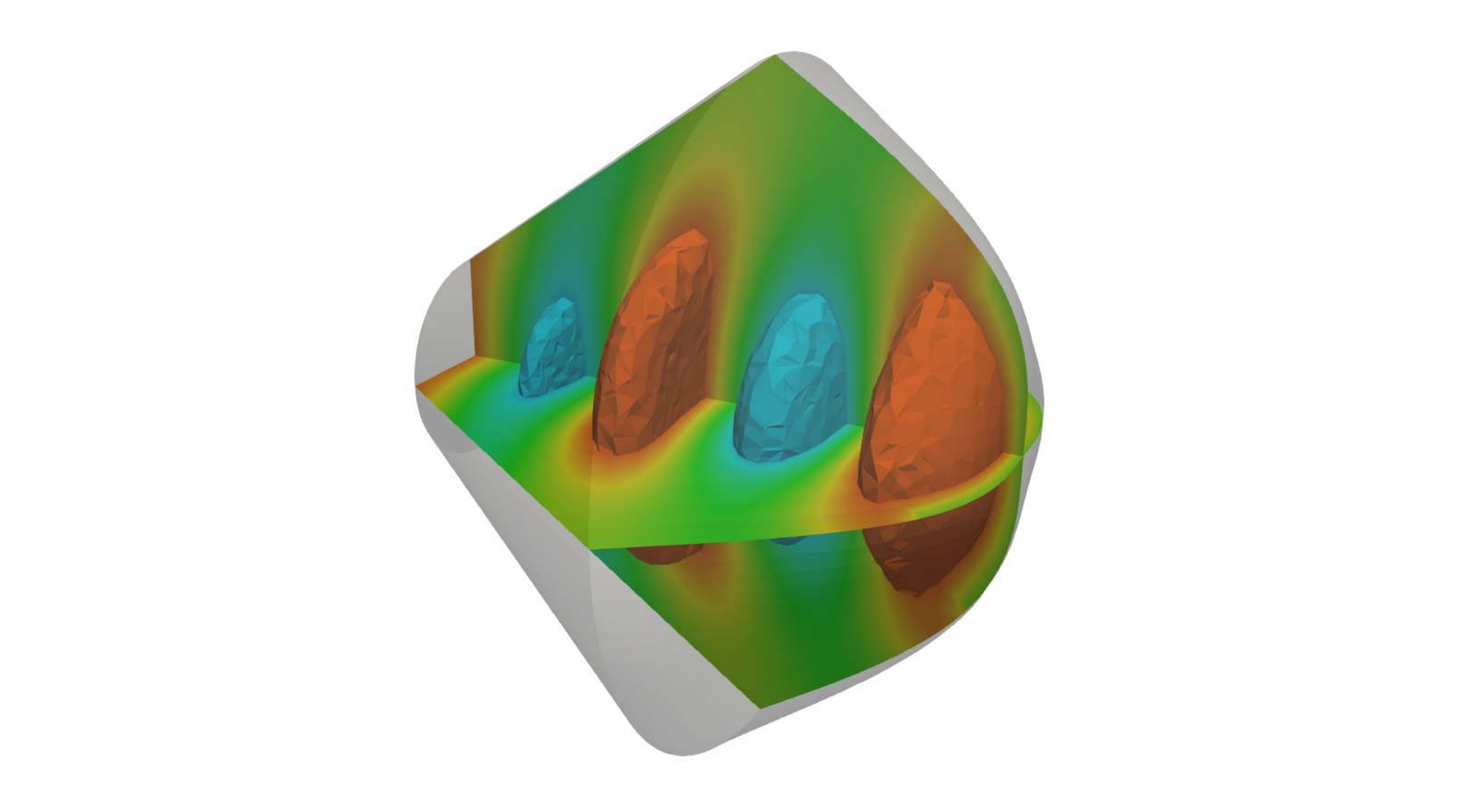}\includegraphics[scale=0.13, trim=7cm 0cm 14cm 0cm, clip]{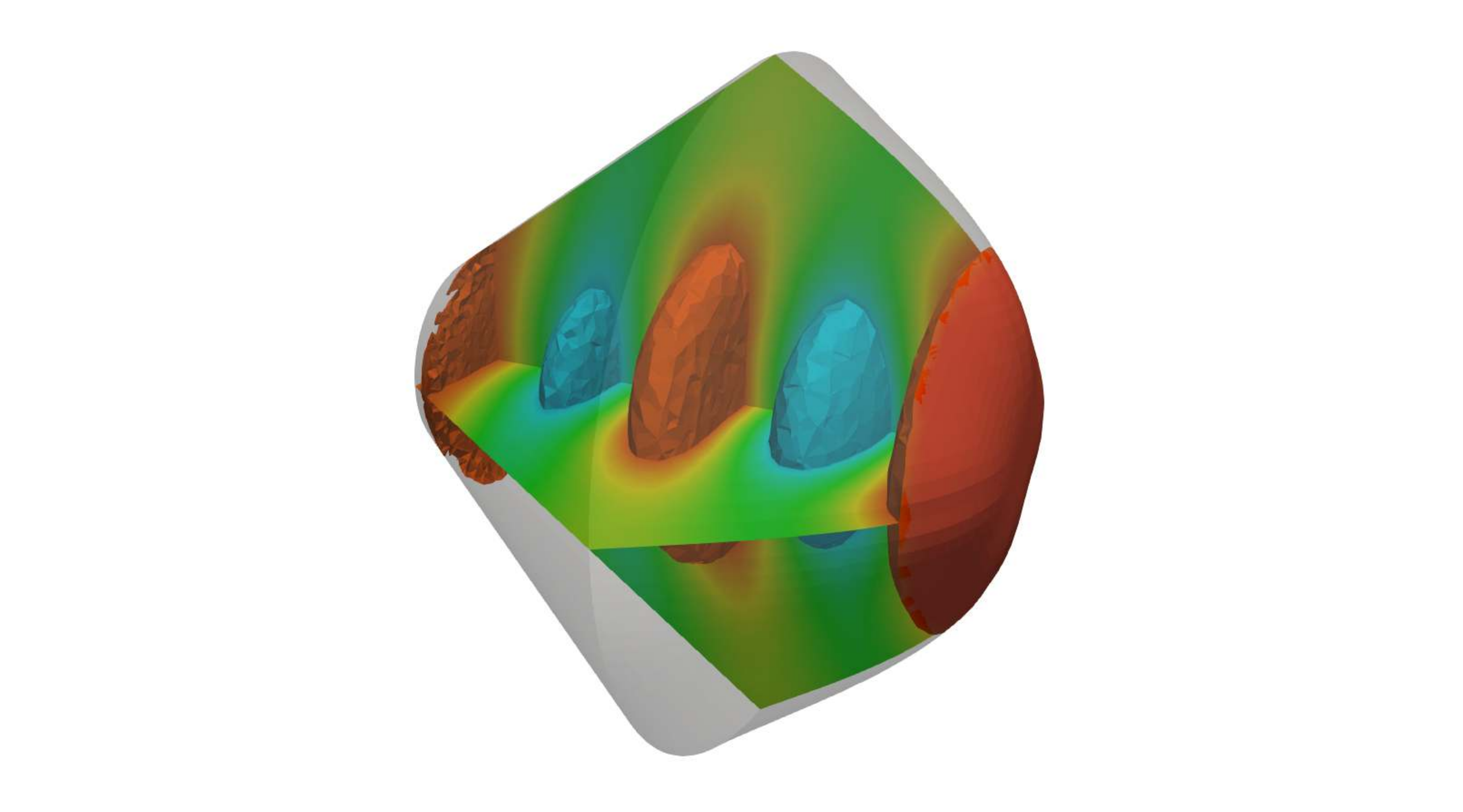}
\caption{Three-dimensional propagation of the pressure wave $u=\rho\psi_t$. The iso-volumes show the regions of highest (in absolute value) pressure amplitudes in 3D, while the two planes show slices through the three-dimensional pressure field.
\label{fig:Transducer3D_wave}}
\end{center}
\end{figure}

To also show the capability of our method in three dimensions, we perform another experiment in a transducer setting, this time in 3D. Due to the high computational costs, especially for the reference solution on the larger domain $\Omega^{\textup{ref}}$, we choose a source amplitude of $\mathfrak{A}=0.002\,\textup{m}^2/\textup{s}^2$, while keeping physical parameters and the excitation frequency the same. Figure~\ref{fig:Transducer3D_mesh} shows the computational domain together with the grid, while Figure~\ref{fig:Transducer3D_wave} depicts the three-dimensional wave propagation in the given domain.

\begin{figure}[h!]
\begin{center}
\input{images/Transducer3D/Cost.tex}\input{images/Transducer3D/Cost_abs.tex}
\caption{\textbf{HIFU transducer 3D:} \textbf{(Left)} Relative $L^2(\Omega)$ error of the potential $\psi(t)$, \textbf{(Right)} Absolute $L^2(\Omega)$ error of the potential $\psi(t)$.
\label{fig:Transducer3DCost}}
\end{center}
\end{figure}
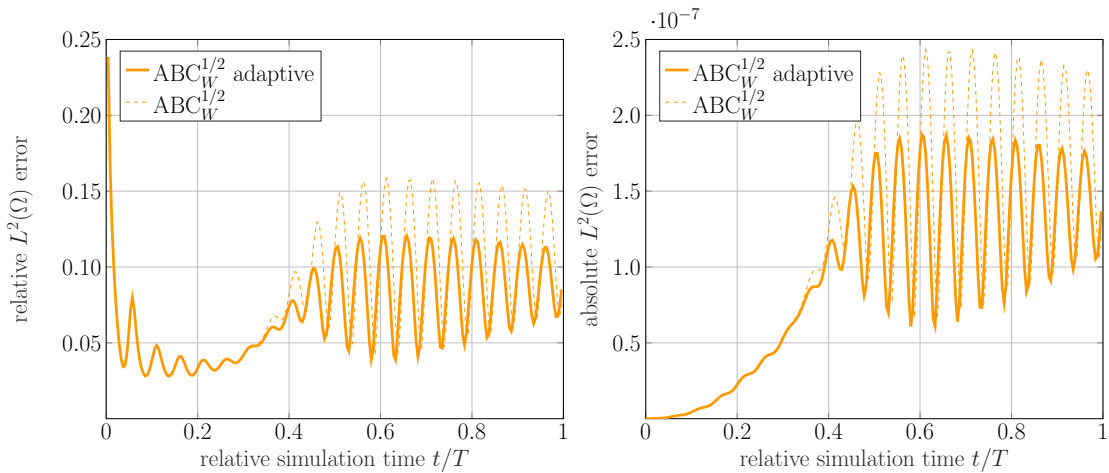

In Figure~\ref{fig:Transducer3DCost}, we compare the adaptive conditions to the nonlinear conditions with the fixed angle $\theta=0 \, ^{\circ}$. 

Note that in 3D, in addition to the discretization error, there is also an interpolation error caused by a mismatch between the grid for the simulation with absorbing conditions and the mesh of the reference solution. Combined with the $L^2$-norm of the reference solution being small at the beginning of the simulation, this results in the initial peak in the relative error for both adaptive and non-adaptive conditions. In 2D, the meshing software was able to avoid the interpolation error. In the long term behavior as well as in the absolute errors, we observe that the adaptive angle information improves the quality of the conditions. The qualitative behavior of the errors of the new adaptive conditions is similar in 3D. The relative errors in the $L^2(0,T; L^2(\Omega))$ norm are $e_{\psi}=8.28\,$\% and $e_{u}=8.3\,$\% if the adaptivity is considered, whereas $e_{\psi}=10.34\,$\% and $e_{u}=10.39\,$\% if the adaptivity is not considered, resulting in an improvement of $19.92\,$\% in $\psi$ and $20.12\,$\% in $u$.

\FloatBarrier

\subsection{3D Acoustic horn}

\begin{figure}[h!]
\begin{center}
\includegraphics[scale=0.16, trim=0cm 2cm 18cm 2cm, clip]{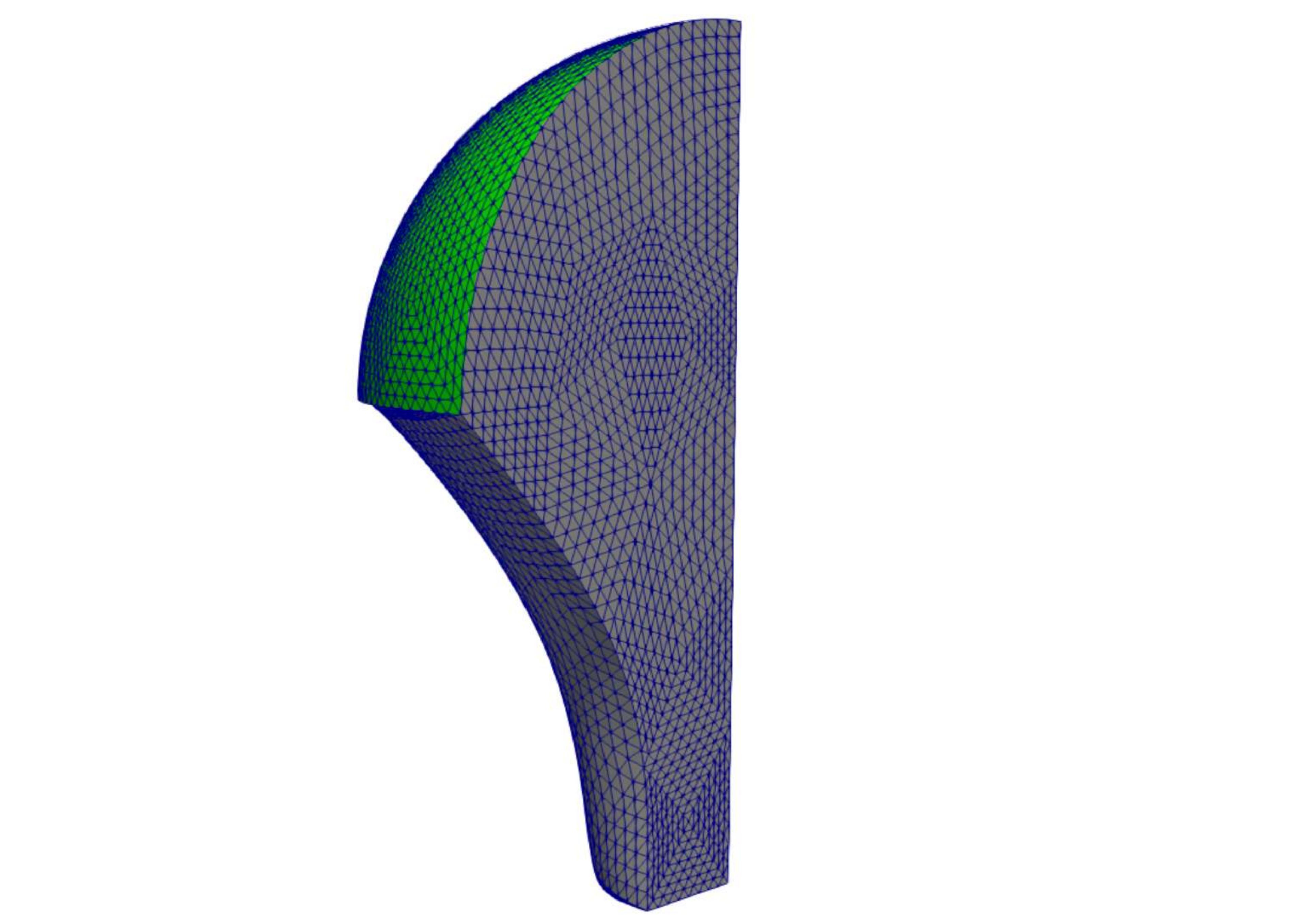}\includegraphics[scale=0.16, trim=15cm 2cm 18cm 3cm, clip]{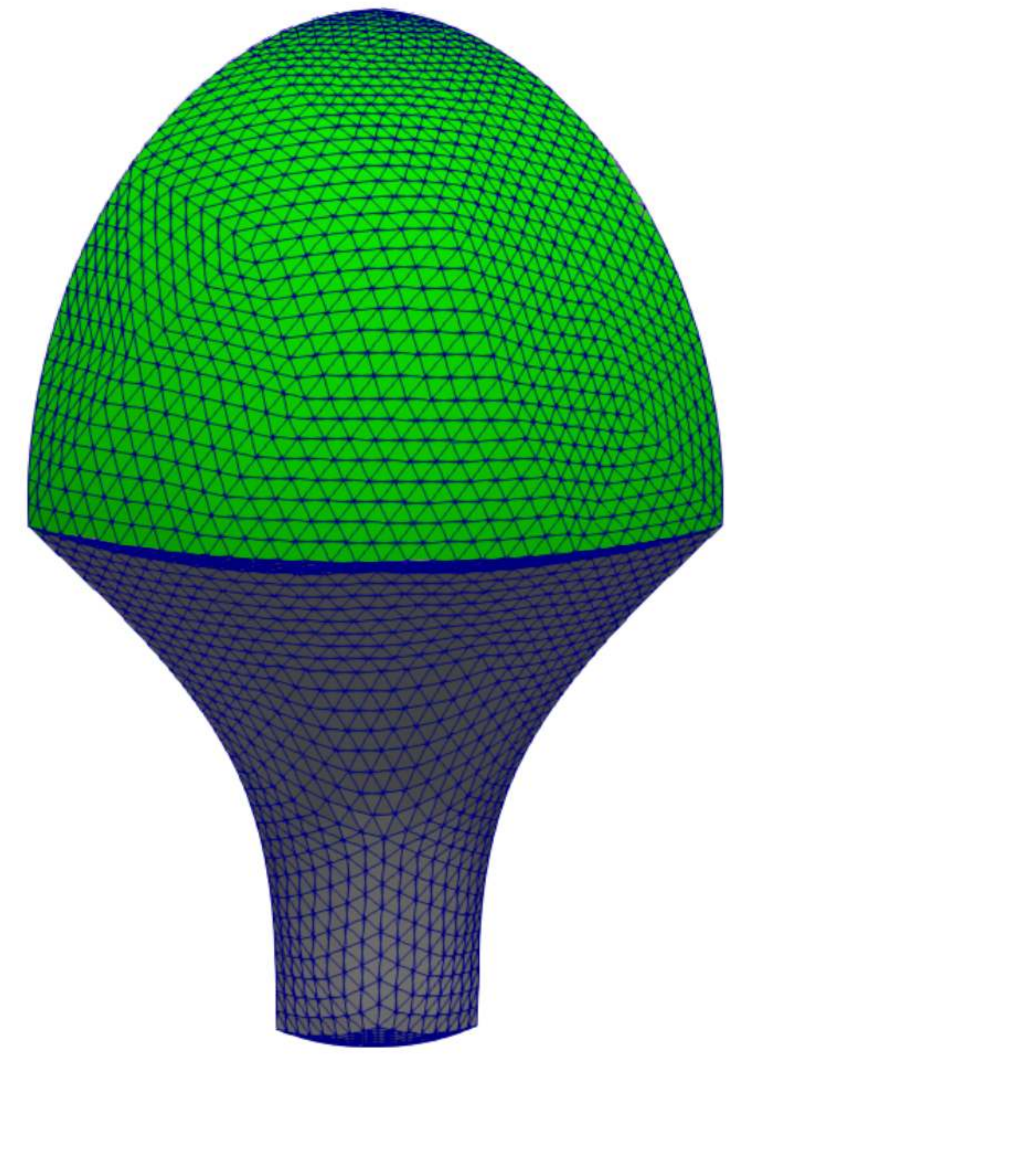}\includegraphics[scale=0.16, trim=11cm 2cm 18cm 2cm, clip]{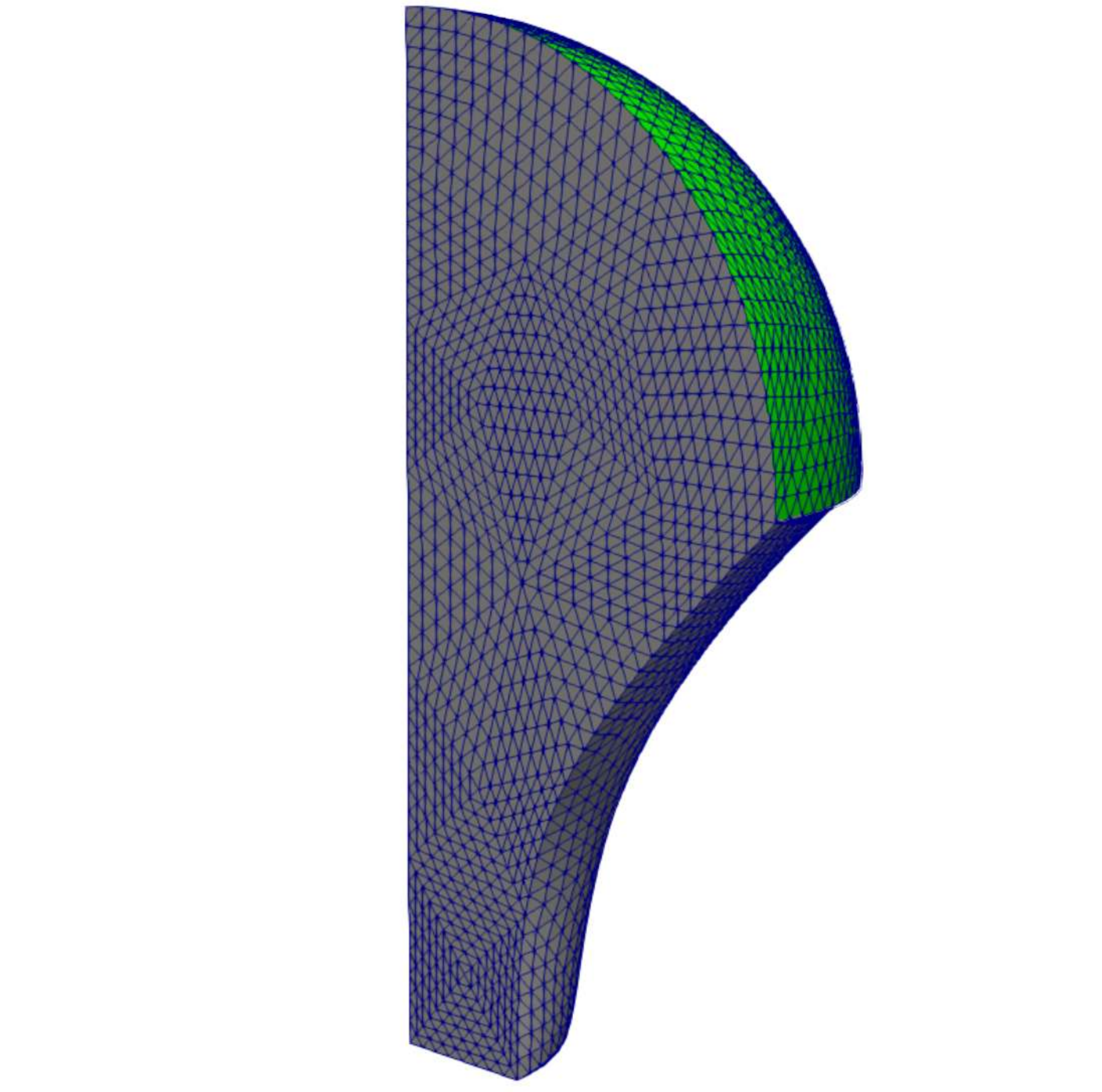}
\caption{Three-dimensional acoustic horn geometry with mesh visible on the surface. Absorbing conditions are employed on the green surface, homogeneous Neumann conditions on the grey and a wave excitation via inhomogeneous Dirichlet conditions at the bottom surface.
\label{fig:Horn3D_mesh}}
\end{center}
\end{figure}

\begin{figure}[h!]
\begin{center}
\includegraphics[scale=0.18, trim=18cm 3cm 18cm 0cm, clip]{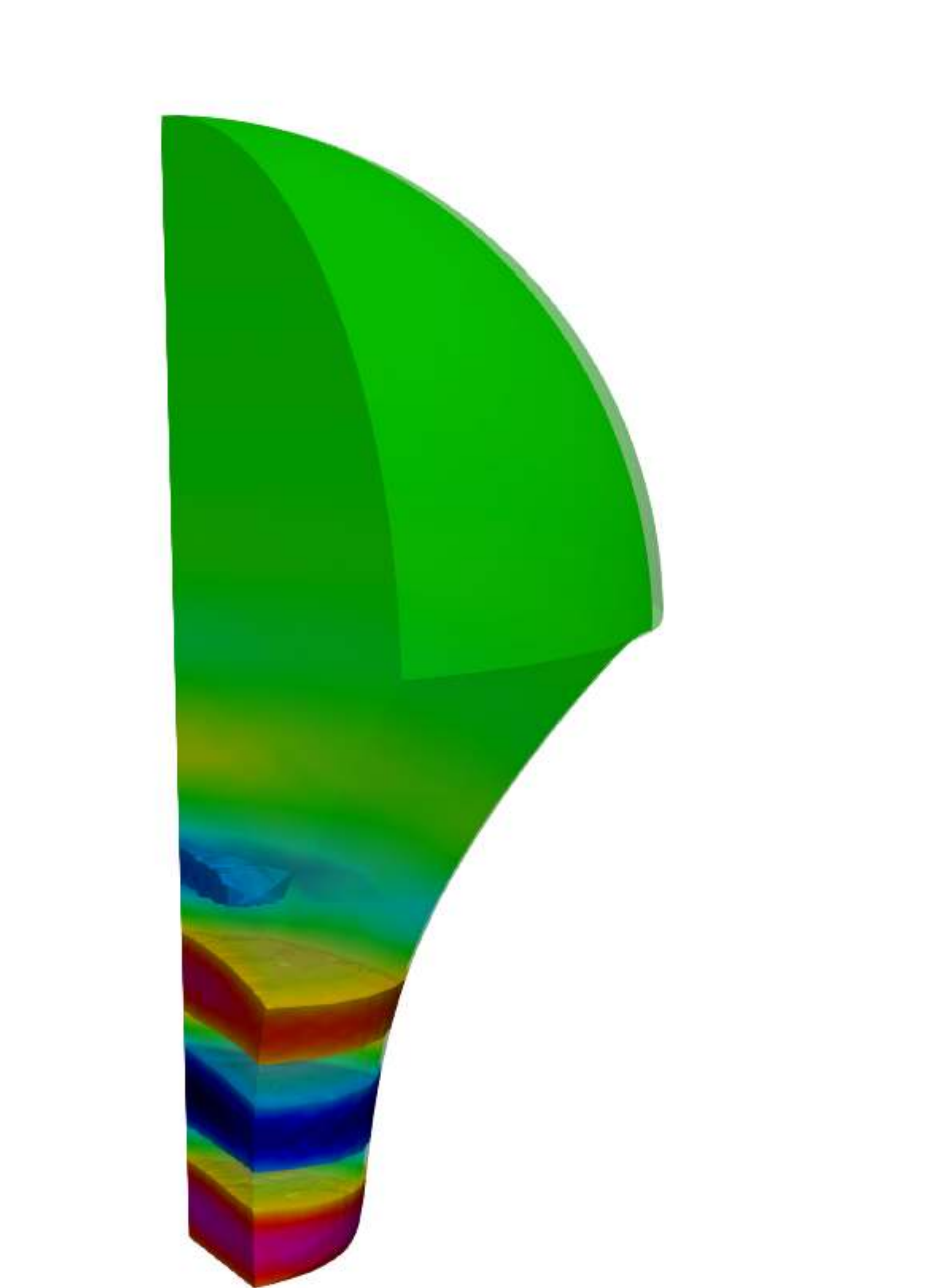}\includegraphics[scale=0.18, trim=18cm 3cm 18cm 0cm, clip]{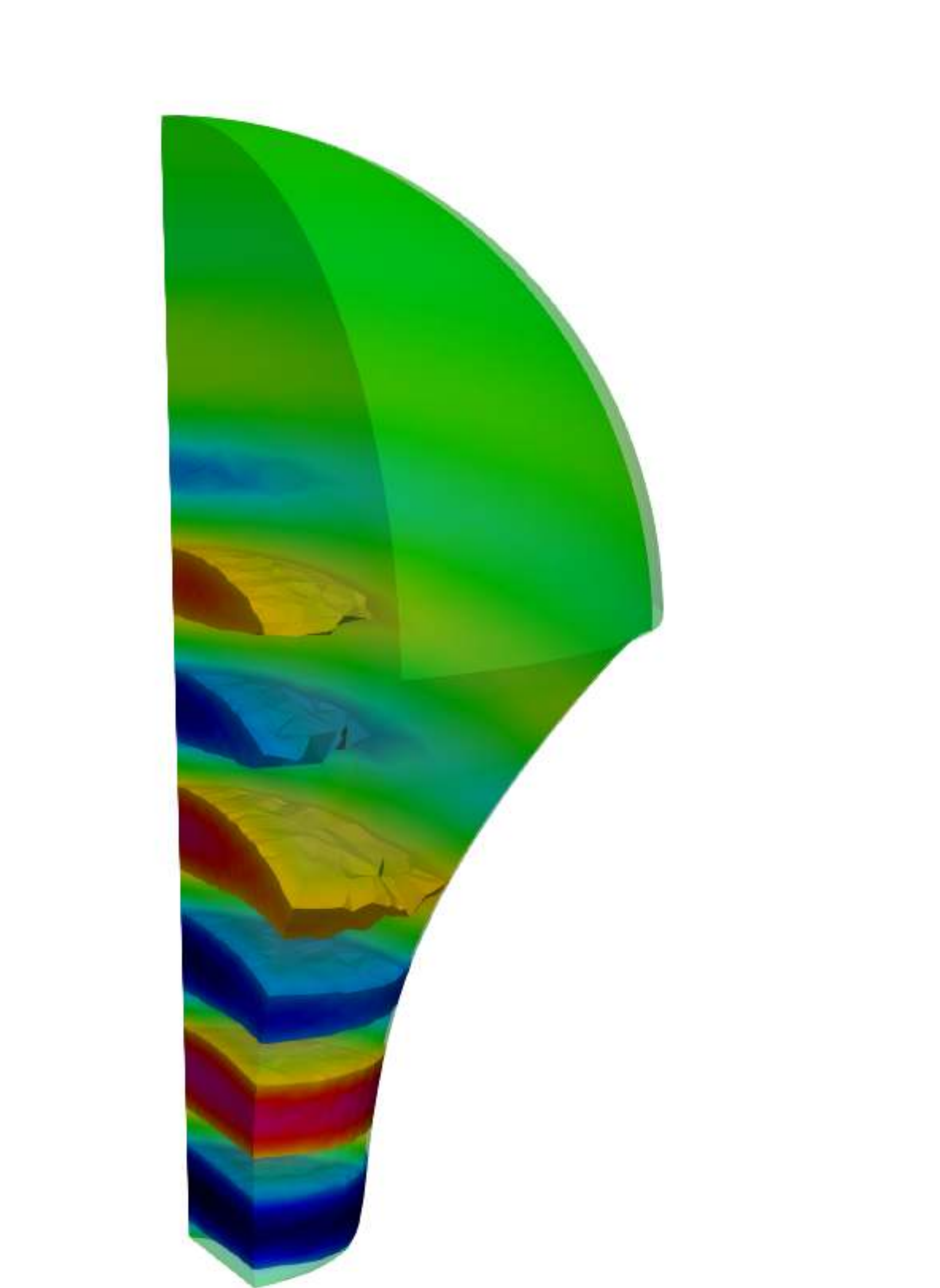}\includegraphics[scale=0.18, trim=18cm 3cm 23cm 0cm, clip]{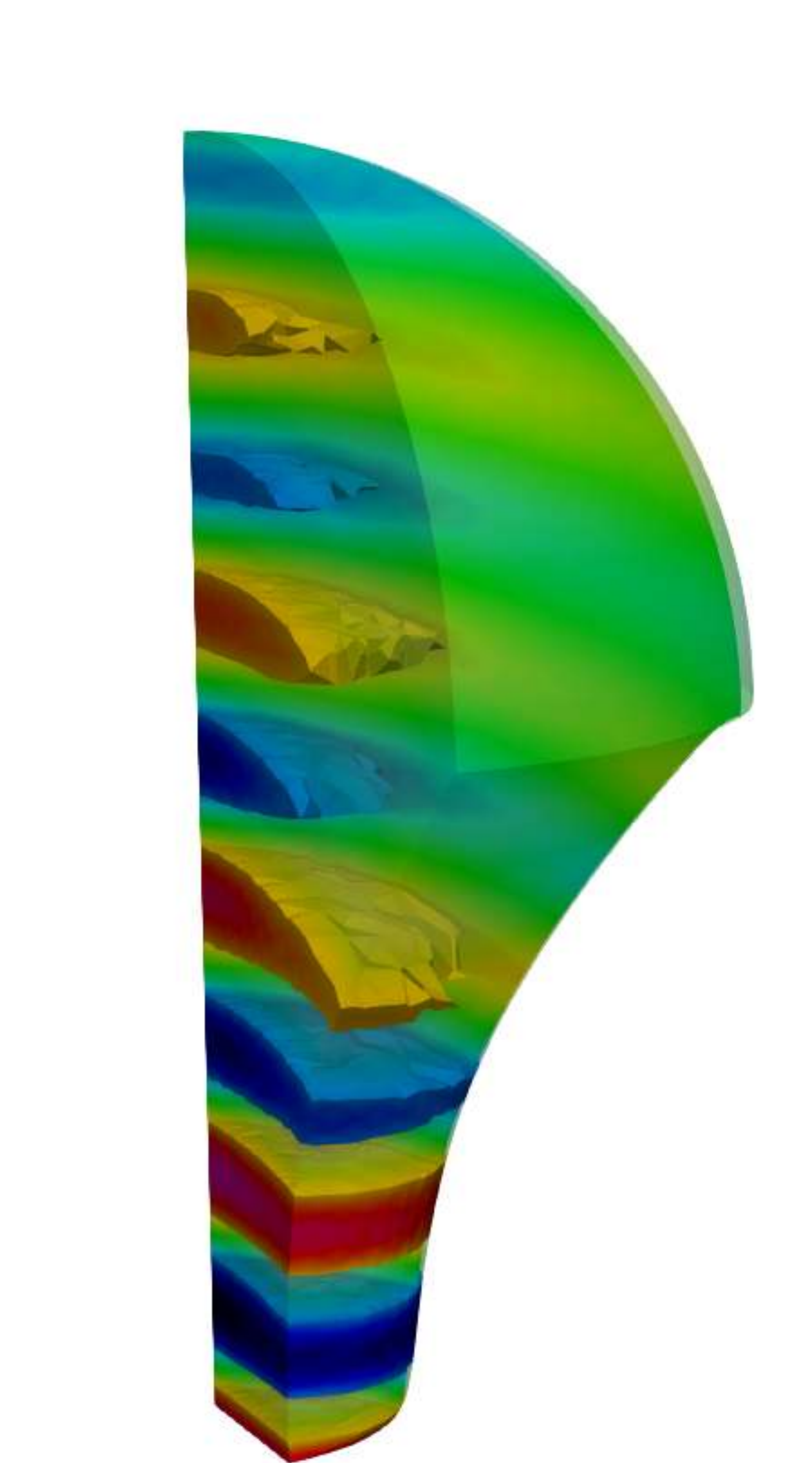}
\caption{Three-dimensional propagation of the pressure wave $u=\rho\psi_t$. The iso-volumes show the regions of highest (in absolute value) pressure amplitudes in 3D.
\label{fig:Horn3D_wave}}
\end{center}
\end{figure}

\begin{figure}[h!]
\begin{center}
\input{images/Horn3D/Cost.tex}\input{images/Horn3D/Cost_abs.tex}
\caption{\textbf{Acoustic horn 3D:} \textbf{(Left)} Relative $L^2(\Omega)$ error of the potential $\psi(t)$, \textbf{(Right)} Absolute $L^2(\Omega)$ error of the potential $\psi(t)$ over the simulation time.
\label{fig:Horn3DCost}}
\end{center}
\end{figure}
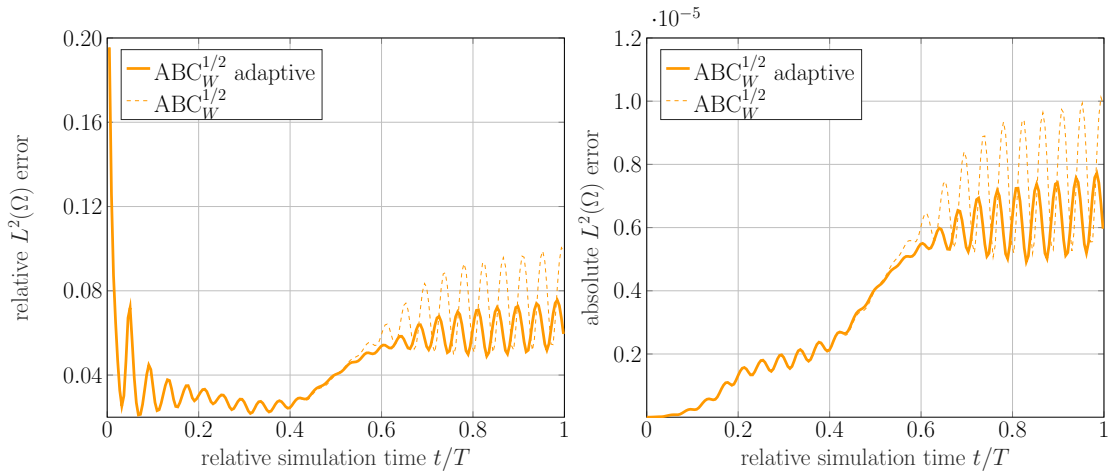

Nonlinear sound propagation has been widely reported to occur in wind instruments; see~\cite{berjamin2017time, fletcher1990nonlinear, myers2012effects, rendon2013nonlinear}. Motivated by this, for our final experiment, we consider a numerical simulation of an acoustic horn.\\
\indent The excitation takes place at the bottom of the domain and the waves then travel through a waveguide with an increasing diameter; see Figure~\ref{fig:Horn3D_wave}. At the end of the waveguide, where the wave starts to propagate into the ambient space, we employ again the absorbing conditions to truncate the domain at a spherical boundary. In this experiment we use the physical parameters of air instead of water:
\begin{align*}
c = 331\,\frac{\textup{m}}{\textup{s}},\ b = 0.00005\,\frac{\textup{kg}}{\textup{m\,s}}, \ B/A = 1.2, \ \varrho=1.29\,\frac{\textup{kg}}{\textup{m}^3};  \end{align*} 
see~\cite{gudra2002parametric}. The excitation has a frequency of $f = 6.5\,\textup{kHz}$ and an amplitude of $\mathfrak{A}=0.01\,\textup{m}^2/\textup{s}^2$. To keep the computational cost reasonable, we again use symmetry to reduce the simulation to a quarter of the actual three-dimensional acoustic horn. The two planes of symmetry are equipped with homogeneous Neumann conditions; see~Figure \ref{fig:Horn3D_mesh}. A comparison of adaptive and non-adaptive conditions is given in~Figure \ref{fig:Horn3DCost}. The relative errors in the $L^2(0,T; L^2(\Omega))$ norm are $e_{\psi}=5.23\,$\% and $e_{u}=5.33\,$\% if the adaptivity is considered, whereas $e_{\psi}=6.11\,$\% and $e_{u}=6.07\,$\% if the adaptivity is not considered, resulting in an improvement of $14.4\,$\% in $\psi$ and $12.19\,$\% in $u$.

~\\
We observe less of a gain in our recent examples compared to the introductory ones in 2D, which has a natural explanation. In our simplest 2D setting in Section~\ref{SubSec:InclinedBoundary}, the angle of incidence is constant over the whole absorbing boundary and large, so the boundary conditions can significantly profit from taking the angle information into account adaptively. In the later, more advanced examples (e.g., in Section~\ref{SubSec:PlateHole}) a large portion of the wave leaves the domain with quite small incidence angles. The same also holds for the application-oriented example of the focusing transducer in Sections~~\ref{subsec:Transducer2D} and~\ref{SubSec:FocusedSound3D}. There, most of the angles of incidence are smaller than in the introductory examples which results smaller improvement compared to the standard non-adaptive conditions.

\section{Conclusion}
We have developed a self-adaptive absorbing technique for sound propagation in the presence of nonlinearities. Within our approach, the angle of incidence of the wave is computed locally by employing the information given by the gradient of the wave field. The absorbing conditions are then updated in real time with the angle values. \\
\indent The method offers three fundamental advantages. It is sufficiently accurate over a range of angles of incidence, and it is easy to implement. Moreover, by only relying on the gradient of the wave field which is readily available in finite element simulations, we can keep the additional computational efforts low.
\section*{Acknowledgements}
 We thank Dr. Igor Shevchenko for helpful comments. The funds provided by the Deutsche Forschungsgemeinschaft under the grant number WO 671/11-1 are gratefully acknowledged.
 
 
\bibliography{references}{}
\bibliographystyle{siam} 
\end{document}

%% file: images/Plane20/Angle_Quality.tex
%
%
\definecolor{mycolor1}{rgb}{1.00000,0.60000,0.00000}%
\pgfplotsset{scaled y ticks=false}
\begin{tikzpicture}[scale = 0.52, font=\huge]

\begin{axis}[%
width=4.521in,
height=3.548in,
at={(0.758in,0.499in)},
scale only axis,
xmin=0,
xmax=1,
xlabel style={font=\huge\color{white!15!black}},
xlabel={relative simulation time $t/T$}, xtick={0, 0.2, 0.4, 0.6, 0.8, 1},
y tick label style={
        /pgf/number format/.cd,
            fixed,
            fixed zerofill,
            precision=2,
        /tikz/.cd
            },
ymin=0,
ymax=0.14, ytick={0.02, 0.04, 0.06, 0.08, 0.1, 0.12, 0.14},
ylabel style={at={(-0.12,0.5)}, font=\huge\color{white!15!black}},
ylabel={},
axis background/.style={fill=white},
xmajorgrids,
ymajorgrids,
legend style={at={(0.03,0.97)}, anchor=north west, legend cell align=left, align=left, draw=white!15!black}
]
\addplot [color=black, line width=2.0pt]
  table[row sep=crcr]{%
0.00214285714285711	1.11022302462516e-16\\
0.204183673469388	5.12904183074259e-05\\
0.224591836734694	0.000276855824008937\\
0.236836734693878	0.000612273677256914\\
0.265408163265306	0.0014650637711292\\
0.281734693877551	0.00181064613458115\\
0.345	0.00348059266504597\\
0.37765306122449	0.00416244084085904\\
0.385816326530612	0.00436591419522026\\
0.402142857142857	0.00456409292354631\\
0.41030612244898	0.00475710334226798\\
0.428673469387755	0.00499287753796551\\
0.434795918367347	0.00512403935306993\\
0.440918367346939	0.00514990339746291\\
0.447040816326531	0.00516881930185287\\
0.453163265306122	0.0052923410435155\\
0.459285714285714	0.00542242261220505\\
0.463367346938776	0.00542330039063332\\
0.473571428571429	0.00532663789487431\\
0.47969387755102	0.00543462955961405\\
0.485816326530612	0.00553504111043079\\
0.489897959183674	0.00550087408973554\\
0.498061224489796	0.00534873446175654\\
0.502142857142857	0.00540639726082792\\
0.51030612244898	0.00568080775632784\\
0.514387755102041	0.00566345328675832\\
0.524591836734694	0.00537546201988137\\
0.528673469387755	0.00547919039878209\\
0.534795918367347	0.00579537490427406\\
0.536836734693878	0.00583865698278185\\
0.540918367346939	0.00576569187136\\
0.549081632653061	0.0053652318488121\\
0.551122448979592	0.00534927172975075\\
0.553163265306122	0.00540327976284283\\
0.557244897959184	0.00565340283363458\\
0.561326530612245	0.00590485694516518\\
0.563367346938776	0.0059567750971542\\
0.565408163265306	0.00592850727663463\\
0.569489795918367	0.00567665445646925\\
0.573571428571429	0.00537319603131137\\
0.575612244897959	0.00528697099634257\\
0.57765306122449	0.00528163419925709\\
0.57969387755102	0.00535559035107358\\
0.583775510204082	0.00568523581502989\\
0.587857142857143	0.00598828231698512\\
0.589897959183673	0.00603791736794923\\
0.591938775510204	0.00600433153804125\\
0.593979591836735	0.0059016306655939\\
0.598061224489796	0.00554884900707298\\
0.600102040816327	0.0053683846146243\\
0.602142857142857	0.00525553454578864\\
0.604183673469388	0.0052404282077313\\
0.606224489795918	0.00534058034531959\\
0.614387755102041	0.00612882399764847\\
0.616428571428571	0.00620931097685651\\
0.618469387755102	0.0061911979784699\\
0.620510204081633	0.00605627035261036\\
0.624591836734694	0.00552542489864005\\
0.626632653061225	0.00526805350991066\\
0.628673469387755	0.00512270828554351\\
0.630714285714286	0.00514947897469242\\
0.632755102040816	0.00533420300296117\\
0.640918367346939	0.00637404245781215\\
0.642959183673469	0.0064650004097363\\
0.645	0.00638772714580149\\
0.647040816326531	0.00613818228676521\\
0.651122448979592	0.00539137246780519\\
0.653163265306122	0.00510716077787587\\
0.655204081632653	0.00500599545936109\\
0.657244897959184	0.00511768756162145\\
0.659285714285714	0.00537999650754051\\
0.665408163265306	0.00643590027707497\\
0.667448979591837	0.00666270411825887\\
0.669489795918367	0.00669013673272956\\
0.671530612244898	0.00648472780073717\\
0.673571428571429	0.00608854727623342\\
0.67765306122449	0.00515787806006485\\
0.67969387755102	0.00489886080774848\\
0.681734693877551	0.00492708662246244\\
0.683775510204082	0.00519101178927039\\
0.685816326530612	0.00558265844800099\\
0.689897959183673	0.00649514919053018\\
0.691938775510204	0.00685550108816935\\
0.693979591836735	0.00701330062150451\\
0.696020408163265	0.00689762304094299\\
0.698061224489796	0.00652773671332907\\
0.702142857142857	0.00544083780597859\\
0.704183673469388	0.00501961211088242\\
0.706224489795918	0.00488878425787076\\
0.708265306122449	0.0050656630261966\\
0.71030612244898	0.00543646827994015\\
0.716428571428571	0.00696481631726142\\
0.718469387755102	0.00727550365013852\\
0.720510204081633	0.0072864293393996\\
0.722551020408163	0.0069847987768511\\
0.724591836734694	0.00644831663733403\\
0.726632653061224	0.00581079610569357\\
0.728673469387755	0.005242825299558\\
0.730714285714286	0.00493634575606361\\
0.732755102040816	0.00499017367882104\\
0.734795918367347	0.00530808231498558\\
0.736836734693878	0.00579736294696465\\
0.740918367346939	0.00699598242026689\\
0.742959183673469	0.00745692003390508\\
0.745	0.00761577574547689\\
0.747040816326531	0.00742118078143794\\
0.749081632653061	0.00692844250615199\\
0.753163265306122	0.00557280906555113\\
0.755204081632653	0.00507733286320045\\
0.757244897959184	0.00494150795639403\\
0.759285714285714	0.0051612985766355\\
0.761326530612245	0.00561798905053545\\
0.763367346938776	0.00624501504485275\\
0.765408163265306	0.00693634151813416\\
0.767448979591837	0.00753999192723864\\
0.769489795918367	0.0078594088137427\\
0.771530612244898	0.00779783040729742\\
0.773571428571429	0.00738238167928462\\
0.775612244897959	0.00670072327749549\\
0.77765306122449	0.00592582487341897\\
0.77969387755102	0.0052542165019599\\
0.781734693877551	0.00490264536497009\\
0.783775510204082	0.00499747107786697\\
0.785816326530612	0.00540021029949089\\
0.787857142857143	0.00603001214446319\\
0.791938775510204	0.00752733597362321\\
0.793979591836735	0.00801540165089509\\
0.796020408163265	0.00811916772782673\\
0.798061224489796	0.00781684673233873\\
0.800102040816327	0.0071809986119955\\
0.804183673469388	0.00556194763203499\\
0.806224489795918	0.00500435127179721\\
0.808265306122449	0.00490257642089231\\
0.81030612244898	0.00520202823735783\\
0.81234693877551	0.00578285054046457\\
0.816428571428571	0.00741194259038869\\
0.818469387755102	0.00806942571700109\\
0.820510204081633	0.00835927087997801\\
0.822551020408163	0.0082036965397092\\
0.824591836734694	0.00765852672444922\\
0.826632653061224	0.00683362718507041\\
0.828673469387755	0.0059288397667141\\
0.830714285714286	0.00517074736090528\\
0.832755102040816	0.00483243271041622\\
0.834795918367347	0.00498813206921089\\
0.836836734693878	0.00547411640154982\\
0.838877551020408	0.00624645557208825\\
0.840918367346939	0.00716166528364726\\
0.842959183673469	0.00798344550925933\\
0.845	0.00846510438189685\\
0.847040816326531	0.00848643243896274\\
0.849081632653061	0.00806733669555493\\
0.851122448979592	0.00731462014096063\\
0.855204081632653	0.0055403201920583\\
0.857244897959184	0.00499248763158111\\
0.859285714285714	0.00495000576961357\\
0.861326530612245	0.00531910533676394\\
0.863367346938775	0.00601764562839113\\
0.867448979591837	0.0078488887650725\\
0.869489795918367	0.00849388260086248\\
0.871530612244898	0.00869429064802674\\
0.873571428571429	0.00840272171048895\\
0.875612244897959	0.00772660027398453\\
0.87969387755102	0.00585579820909066\\
0.881734693877551	0.0051564640463847\\
0.883775510204082	0.00494132018183502\\
0.885816326530612	0.00520545607038214\\
0.887857142857143	0.00581358884999761\\
0.889897959183673	0.00670070331511508\\
0.891938775510204	0.00766249015236653\\
0.893979591836735	0.0084243693783288\\
0.896020408163265	0.00878851653354185\\
0.898061224489796	0.00865431429995622\\
0.900102040816326	0.00810198383504157\\
0.902142857142857	0.00726037089315912\\
0.904183673469388	0.00629459849669201\\
0.906224489795918	0.00547393443823541\\
0.908265306122449	0.00502695338947734\\
0.91030612244898	0.00509777813767964\\
0.91234693877551	0.0055769565282372\\
0.914387755102041	0.00639458863435982\\
0.916428571428571	0.00737412051032194\\
0.918469387755102	0.00825205529673301\\
0.920510204081633	0.00877857860489561\\
0.922551020408163	0.00881303428071012\\
0.924591836734694	0.00837206380951083\\
0.926632653061224	0.00759542630897325\\
0.928673469387755	0.00662442405401475\\
0.930714285714286	0.00572999782772454\\
0.932755102040816	0.00516366147416647\\
0.934795918367347	0.00507503143552102\\
0.936836734693878	0.00543938766862218\\
0.938877551020408	0.00614053632740574\\
0.942959183673469	0.00802073439472661\\
0.945	0.00867989510497735\\
0.947040816326531	0.00889751558079954\\
0.949081632653061	0.00862127691496561\\
0.951122448979592	0.00796678602009726\\
0.955204081632653	0.00612887460734668\\
0.957244897959184	0.00544089850794827\\
0.959285714285714	0.00513329135746365\\
0.961326530612245	0.00533786911084289\\
0.963367346938775	0.00589853181908739\\
0.965408163265306	0.00677130653912883\\
0.967448979591837	0.00774613595532381\\
0.969489795918367	0.0085263234490125\\
0.971530612244898	0.00890126039327022\\
0.973571428571429	0.00877771175119157\\
0.975612244897959	0.00822609404339081\\
0.97765306122449	0.00738641917848581\\
0.97969387755102	0.00643225776695844\\
0.981734693877551	0.0056713682450682\\
0.983775510204082	0.00523899561257768\\
0.985816326530612	0.00530802928119156\\
0.987857142857143	0.00574901291364438\\
0.989897959183674	0.00650997922427921\\
0.991938775510204	0.00747343147144353\\
0.993979591836735	0.00834507640737336\\
0.996020408163265	0.00885499488131081\\
0.998061224489796	0.00890995027808894\\
};
\addlegendentry{ABC$_W^{1/2}$ adaptive with exact $\theta$}

\addplot [color=mycolor1, line width=2.0pt]
  table[row sep=crcr]{%
0.00214285714285711	1.11022302462516e-16\\
0.204183673469388	5.12904183074259e-05\\
0.224591836734694	0.000274578951980797\\
0.281734693877551	0.00144146858701533\\
0.293979591836735	0.00161427818245707\\
0.326632653061224	0.00211343663489816\\
0.336836734693878	0.00229899982300585\\
0.361326530612245	0.00250850014198811\\
0.373571428571429	0.00258193705756027\\
0.387857142857143	0.00276268435631088\\
0.406224489795918	0.00282660504083621\\
0.416428571428571	0.00285727916116707\\
0.426632653061225	0.00288912005786102\\
0.438877551020408	0.00295352920102232\\
0.453163265306122	0.00294634353913059\\
0.461326530612245	0.00296716115334628\\
0.469489795918367	0.00284770149558677\\
0.475612244897959	0.00277762354610145\\
0.491938775510204	0.00268241203594755\\
0.498061224489796	0.00257081425872585\\
0.502142857142857	0.00257917127894292\\
0.514387755102041	0.0027382852826402\\
0.526632653061225	0.00259029625146034\\
0.530714285714286	0.0027038871893944\\
0.536836734693878	0.00291613329595908\\
0.540918367346939	0.00292187172620983\\
0.551122448979592	0.00269129160562764\\
0.555204081632653	0.00279163758414669\\
0.563367346938776	0.00317626360644829\\
0.567448979591837	0.00320383859243012\\
0.573571428571429	0.00302199290411687\\
0.57765306122449	0.00294993079847294\\
0.581734693877551	0.00311531114246832\\
0.587857142857143	0.00355575986099343\\
0.589897959183673	0.00364422865944969\\
0.593979591836735	0.00363154582576841\\
0.602142857142857	0.00338562438317502\\
0.606224489795918	0.00348670179709243\\
0.61030612244898	0.00378089521500857\\
0.614387755102041	0.00406983435341435\\
0.618469387755102	0.00415297820565752\\
0.622551020408163	0.00405197045235528\\
0.626632653061225	0.0038994571332599\\
0.630714285714286	0.00389986642827522\\
0.632755102040816	0.00396528597076529\\
0.638877551020408	0.00445345278646436\\
0.642959183673469	0.00461041747286528\\
0.645	0.00461906675693591\\
0.649081632653061	0.0044875135604181\\
0.653163265306122	0.00434098240736036\\
0.657244897959184	0.00431205381589095\\
0.659285714285714	0.00442480431943804\\
0.665408163265306	0.00490735524505992\\
0.669489795918367	0.00507190289575588\\
0.673571428571429	0.00497060164093033\\
0.67765306122449	0.00472979386132322\\
0.681734693877551	0.00458506971355821\\
0.683775510204082	0.00461082406419955\\
0.691938775510204	0.00519624635960314\\
0.693979591836735	0.00529449731043352\\
0.696020408163265	0.00529918975716748\\
0.698061224489796	0.00524016139520289\\
0.706224489795918	0.00471472670915207\\
0.708265306122449	0.00463097041606664\\
0.71234693877551	0.00493533307208238\\
0.716428571428571	0.00526052974571467\\
0.718469387755102	0.00541622925177054\\
0.720510204081633	0.00547963640698468\\
0.722551020408163	0.00544982157473795\\
0.724591836734694	0.00532064725053971\\
0.728673469387755	0.00497330287124043\\
0.732755102040816	0.00466612220856333\\
0.734795918367347	0.00471438362946297\\
0.745	0.00550959113718363\\
0.747040816326531	0.00551920940224671\\
0.749081632653061	0.00543929406966326\\
0.753163265306122	0.00507044360416675\\
0.755204081632653	0.00492875959946504\\
0.757244897959184	0.00471111040372652\\
0.759285714285714	0.00464778333783644\\
0.765408163265306	0.00522192479332539\\
0.769489795918367	0.0056271560968193\\
0.771530612244898	0.00569837301525289\\
0.773571428571429	0.00566955987627316\\
0.77969387755102	0.00510592553778677\\
0.783775510204082	0.00467601981885268\\
0.785816326530612	0.00478728927437977\\
0.791938775510204	0.00536256056352291\\
0.793979591836735	0.00559167027525409\\
0.796020408163265	0.00572462045378619\\
0.798061224489796	0.00575614749896269\\
0.800102040816327	0.00564039835904495\\
0.802142857142857	0.00540920777655385\\
0.806224489795918	0.00503153124108491\\
0.808265306122449	0.00472320203436172\\
0.81030612244898	0.00474475191296364\\
0.816428571428571	0.00535950257387097\\
0.820510204081633	0.00583458737701947\\
0.822551020408163	0.00591605489818436\\
0.824591836734694	0.00586059778375503\\
0.830714285714286	0.00519429686106587\\
0.832755102040816	0.0048247860091164\\
0.834795918367347	0.00464813803559205\\
0.838877551020408	0.00490925303547307\\
0.842959183673469	0.00537381307265583\\
0.845	0.00565281787604521\\
0.847040816326531	0.00579611972241789\\
0.849081632653061	0.00581785144492342\\
0.851122448979592	0.00566758970386705\\
0.857244897959184	0.00490830065220504\\
0.859285714285714	0.00456426007355792\\
0.861326530612245	0.00463490376978271\\
0.863367346938775	0.00479148828070419\\
0.867448979591837	0.00526870214607433\\
0.869489795918367	0.00559405306106986\\
0.871530612244898	0.00583074375951942\\
0.873571428571429	0.00590348711792932\\
0.875612244897959	0.00581632965793744\\
0.881734693877551	0.00507467825848062\\
0.883775510204082	0.004625839847199\\
0.885816326530612	0.00453541173067851\\
0.887857142857143	0.00464265125667718\\
0.889897959183673	0.00482864452171949\\
0.893979591836735	0.00540051999701796\\
0.896020408163265	0.0056949431010741\\
0.898061224489796	0.00582313797574696\\
0.900102040816326	0.00581466557877641\\
0.902142857142857	0.00562443585789385\\
0.906224489795918	0.00514913451669252\\
0.908265306122449	0.00470071616865397\\
0.91030612244898	0.00443767099808978\\
0.91234693877551	0.00455498415847311\\
0.914387755102041	0.00473669455099002\\
0.918469387755102	0.00532788826683284\\
0.920510204081633	0.00567939281023111\\
0.922551020408163	0.00588407106092381\\
0.924591836734694	0.0059248029921527\\
0.926632653061224	0.00578633618391244\\
0.928673469387755	0.00550379282228641\\
0.930714285714286	0.00529526890812371\\
0.932755102040816	0.00495102365270461\\
0.934795918367347	0.0045125632500147\\
0.936836734693878	0.00452258560151886\\
0.938877551020408	0.00462591016270286\\
0.940918367346939	0.00486807010723378\\
0.947040816326531	0.00575327490772737\\
0.949081632653061	0.00586050994528375\\
0.951122448979592	0.00580330776512861\\
0.957244897959184	0.00509120887909542\\
0.959285714285714	0.00457554229036805\\
0.961326530612245	0.00448668687356368\\
0.963367346938775	0.00460688815806798\\
0.965408163265306	0.00481493712065018\\
0.971530612244898	0.00575402355585131\\
0.973571428571429	0.00590816741019673\\
0.975612244897959	0.00590413458547723\\
0.97765306122449	0.00570500084551584\\
0.97969387755102	0.00544001841168162\\
0.981734693877551	0.00526994370158484\\
0.983775510204082	0.00481934538672435\\
0.985816326530612	0.00452149767876187\\
0.989897959183674	0.00467960615264906\\
0.996020408163265	0.00559793006640297\\
0.998061224489796	0.00580179185247498\\
};
\addlegendentry{ABC$_W^{1/2}$ adaptive with approx. $\theta$}

\end{axis}
\end{tikzpicture}%

%% file: images/Plane50/Angle_Quality.tex
%
%
\definecolor{mycolor1}{rgb}{1.00000,0.60000,0.00000}%
\pgfplotsset{scaled y ticks=false}
\begin{tikzpicture}[scale = 0.52, font=\huge]

\begin{axis}[%
width=4.521in,
height=3.548in,
at={(0.758in,0.499in)},
scale only axis,
xmin=0,
xmax=1,
xlabel style={font=\huge\color{white!15!black}},
xlabel={relative simulation time $t/T$}, xtick={0, 0.2, 0.4, 0.6, 0.8, 1},
y tick label style={
        /pgf/number format/.cd,
            fixed,
            fixed zerofill,
            precision=2,
        /tikz/.cd
            },
ymin=0,
ymax=0.26, ytick={0.02, 0.06, 0.1, 0.14, 0.18, 0.22, 0.26},
ylabel style={at={(-0.15,0.5)}, font=\huge\color{white!15!black}},
ylabel={},
axis background/.style={fill=white},
xmajorgrids,
ymajorgrids,
legend style={at={(0.03,0.97)}, anchor=north west, legend cell align=left, align=left, draw=white!15!black}
]
\addplot [color=black, line width=2.0pt]
  table[row sep=crcr]{%
0.00214285714285711	4.44089209850063e-16\\
0.138877551020408	5.68269401037824e-05\\
0.147040816326531	0.000308035076785007\\
0.161326530612245	0.00091690947229206\\
0.165408163265306	0.00132943697110954\\
0.173571428571429	0.00241583436263793\\
0.185816326530612	0.00420947587771026\\
0.200102040816327	0.00597450435251623\\
0.208265306122449	0.00708852055367037\\
0.214387755102041	0.00764575476028917\\
0.220510204081633	0.00817236505759766\\
0.226632653061224	0.00891067081178432\\
0.232755102040816	0.00967176912438672\\
0.238877551020408	0.0101857943299526\\
0.245	0.0106640204898805\\
0.251122448979592	0.0113842634350265\\
0.257244897959184	0.0121540889082635\\
0.261326530612245	0.0125100590916535\\
0.271530612244898	0.0132128580040098\\
0.27765306122449	0.0139527169946672\\
0.283775510204082	0.0146822644304501\\
0.287857142857143	0.0150026910306196\\
0.298061224489796	0.0156441007032917\\
0.302142857142857	0.016087545999394\\
0.314387755102041	0.0177893897179439\\
0.318469387755102	0.0185469935465132\\
0.322551020408163	0.0196467288604257\\
0.332755102040816	0.0228859182684613\\
0.336836734693878	0.0237476969514329\\
0.342959183673469	0.0246128431186081\\
0.347040816326531	0.0252472124925092\\
0.355204081632653	0.0267332626557024\\
0.359285714285714	0.0271692309345257\\
0.363367346938776	0.0272919624931334\\
0.369489795918367	0.0273327777009463\\
0.373571428571429	0.02753195511478\\
0.383775510204082	0.0282440418515941\\
0.387857142857143	0.0282599156935928\\
0.396020408163265	0.0280228354346762\\
0.400102040816327	0.0281271687231973\\
0.408265306122449	0.0285652251190217\\
0.41234693877551	0.028460323387994\\
0.420510204081633	0.0278262180068855\\
0.424591836734694	0.0277947652954477\\
0.432755102040816	0.0281216156529095\\
0.436836734693878	0.027977349903848\\
0.447040816326531	0.0271070906631383\\
0.451122448979592	0.0271217018319855\\
0.457244897959184	0.0273166913724048\\
0.461326530612245	0.0271621575724791\\
0.467448979591837	0.0264980447496784\\
0.471530612244898	0.0261503082242008\\
0.475612244897959	0.0260846523436404\\
0.481734693877551	0.0262043029136882\\
0.485816326530612	0.0260489939047344\\
0.489897959183674	0.0256132625922947\\
0.496020408163265	0.0249104656787082\\
0.500102040816326	0.0247054033548382\\
0.508265306122449	0.0245782777057775\\
0.51234693877551	0.0241951343989126\\
0.522551020408163	0.0228301839492742\\
0.526632653061225	0.0226293646731816\\
0.532755102040816	0.0224442059776542\\
0.536836734693878	0.0220753250387942\\
0.547040816326531	0.0206908430085326\\
0.551122448979592	0.0204226919521199\\
0.557244897959184	0.0201980489634496\\
0.561326530612245	0.0198331964398445\\
0.567448979591837	0.0189291111529726\\
0.571530612244898	0.018380046226838\\
0.575612244897959	0.0180687813937955\\
0.583775510204082	0.0177341561758073\\
0.587857142857143	0.0173330872690626\\
0.596020408163265	0.0163454060982465\\
0.600102040816327	0.0160872211529379\\
0.608265306122449	0.0158366734186317\\
0.61234693877551	0.0154617834573026\\
0.622551020408163	0.0142523326974636\\
0.626632653061225	0.0140809439685515\\
0.632755102040816	0.0139698715851431\\
0.636836734693878	0.0136884376766355\\
0.647040816326531	0.0126262715703338\\
0.651122448979592	0.0124702952844619\\
0.659285714285714	0.0123214974453767\\
0.663367346938776	0.0119645704897146\\
0.671530612244898	0.0111048710363408\\
0.675612244897959	0.0109623051663795\\
0.683775510204082	0.0109561798521719\\
0.687857142857143	0.0106743110312845\\
0.696020408163265	0.00982644521339848\\
0.700102040816327	0.00965510117718449\\
0.708265306122449	0.0096124155200894\\
0.71234693877551	0.00931814785090757\\
0.722551020408163	0.00822876799993888\\
0.726632653061224	0.00815868582654544\\
0.732755102040816	0.00826303206258749\\
0.736836734693878	0.0080816939278191\\
0.740918367346939	0.00768165940885734\\
0.747040816326531	0.00705332174240725\\
0.751122448979592	0.00689969889438591\\
0.759285714285714	0.00681331773882021\\
0.763367346938776	0.00643082745939272\\
0.771530612244898	0.00546621714500917\\
0.773571428571429	0.00535614745752488\\
0.77765306122449	0.00542318958533672\\
0.781734693877551	0.00564048397401162\\
0.785816326530612	0.00562800406873143\\
0.791938775510204	0.00523240095479471\\
0.796020408163265	0.00493318318713132\\
0.800102040816327	0.0048986945846996\\
0.804183673469388	0.00510134666380579\\
0.808265306122449	0.00536317261876662\\
0.81030612244898	0.00538118339229454\\
0.814387755102041	0.00514435752606879\\
0.822551020408163	0.0045060610937705\\
0.824591836734694	0.00450082218957615\\
0.828673469387755	0.00474872699668549\\
0.832755102040816	0.00508162940762247\\
0.836836734693878	0.00506213195575222\\
0.842959183673469	0.00469615115212962\\
0.847040816326531	0.0044498087864564\\
0.851122448979592	0.00449408796109541\\
0.853163265306122	0.00459497299612199\\
0.857244897959184	0.0049482072118423\\
0.859285714285714	0.00502959027833483\\
0.863367346938775	0.00491776413208012\\
0.869489795918367	0.00448302253203103\\
0.871530612244898	0.00436198531248966\\
0.873571428571429	0.00433949945338286\\
0.87765306122449	0.00456432997781209\\
0.881734693877551	0.00502103270679399\\
0.883775510204082	0.00512323343027743\\
0.887857142857143	0.00503615518714229\\
0.898061224489796	0.00430698171590294\\
0.902142857142857	0.00435404843562281\\
0.908265306122449	0.00494097114102388\\
0.91030612244898	0.00503171540786385\\
0.914387755102041	0.00489495851453048\\
0.920510204081633	0.00443580417946909\\
0.922551020408163	0.00434321735443088\\
0.926632653061224	0.00441875257650237\\
0.932755102040816	0.0051143952218593\\
0.934795918367347	0.00521664787936382\\
0.938877551020408	0.00513631465753883\\
0.947040816326531	0.00457648852551185\\
0.953163265306122	0.00467348272838097\\
0.957244897959184	0.00509313836583858\\
0.961326530612245	0.00528695942918855\\
0.965408163265306	0.00508031139129095\\
0.971530612244898	0.00456931220273638\\
0.975612244897959	0.00445927312608818\\
0.97765306122449	0.00454719136796511\\
0.981734693877551	0.00505022237246355\\
0.985816326530612	0.00529404670051026\\
0.989897959183674	0.00517676944391143\\
0.996020408163265	0.0046958602961894\\
0.998061224489796	0.00461877637848818\\
};
\addlegendentry{ABC$_W^{1/2}$ adaptive with exact $\theta$}

\addplot [color=mycolor1, line width=2.0pt]
  table[row sep=crcr]{%
0.00214285714285711	4.44089209850063e-16\\
0.138877551020408	5.68269401037824e-05\\
0.147040816326531	0.000308035076785007\\
0.161326530612245	0.00091690947229206\\
0.167448979591837	0.001568774244041\\
0.17765306122449	0.00292523993406113\\
0.183775510204082	0.00371060268048551\\
0.191938775510204	0.00450764755760258\\
0.198061224489796	0.00516026695426564\\
0.21030612244898	0.00662549211135077\\
0.224591836734694	0.00785601189647278\\
0.232755102040816	0.00873456224402824\\
0.238877551020408	0.00913784415701968\\
0.245	0.00952961190364077\\
0.249081632653061	0.00993823724151155\\
0.259285714285714	0.0111031423899371\\
0.265408163265306	0.0114991048583176\\
0.271530612244898	0.0119276418128968\\
0.27765306122449	0.0125968743633909\\
0.283775510204082	0.0132568455144382\\
0.289897959183673	0.0136368760284457\\
0.298061224489796	0.0140843643071737\\
0.302142857142857	0.014494538727709\\
0.316428571428572	0.0164290100392172\\
0.318469387755102	0.0168006365386051\\
0.320510204081633	0.0172648400786739\\
0.324591836734694	0.0184726872111283\\
0.330714285714286	0.0205361247761945\\
0.332755102040816	0.0211069346639262\\
0.334795918367347	0.0215574290435776\\
0.336836734693878	0.0218823735419835\\
0.338877551020408	0.0220823544839287\\
0.349081632653061	0.0227111244011992\\
0.355204081632653	0.0234078029550743\\
0.357244897959184	0.0234967801162597\\
0.359285714285714	0.023447432519601\\
0.363367346938776	0.0230414154642851\\
0.367448979591837	0.0225219520086147\\
0.371530612244898	0.0222399184441889\\
0.37969387755102	0.0219421708943545\\
0.383775510204082	0.0215105121038378\\
0.389897959183674	0.0207789977971672\\
0.393979591836735	0.0205502781374632\\
0.400102040816327	0.0203185724531872\\
0.406224489795918	0.0197812573146058\\
0.41234693877551	0.0193515114309105\\
0.416428571428571	0.0190936182738817\\
0.420510204081633	0.0188843052251272\\
0.424591836734694	0.0183382052166443\\
0.428673469387755	0.0177505723845259\\
0.434795918367347	0.0170776845401005\\
0.442959183673469	0.0166996515078597\\
0.447040816326531	0.0162412066026204\\
0.455204081632653	0.0148844349196391\\
0.457244897959184	0.0147167657596533\\
0.459285714285714	0.0146623021939208\\
0.463367346938776	0.01468742491405\\
0.467448979591837	0.0146112051563801\\
0.469489795918367	0.0144905509849994\\
0.471530612244898	0.0142531217651299\\
0.47765306122449	0.0130023494056811\\
0.481734693877551	0.0126293911550444\\
0.483775510204082	0.0125200206215431\\
0.485816326530612	0.0125671458248671\\
0.489897959183674	0.0128369197428446\\
0.491938775510204	0.0128263152263931\\
0.496020408163265	0.0125152247786543\\
0.498061224489796	0.0121847516068695\\
0.504183673469388	0.0109725945025672\\
0.506224489795918	0.0107488127112646\\
0.508265306122449	0.0106844681000507\\
0.51234693877551	0.0109753398148253\\
0.516428571428571	0.0112178004064311\\
0.520510204081633	0.0111122470010098\\
0.522551020408163	0.0108072134049507\\
0.526632653061225	0.0096765681946821\\
0.528673469387755	0.00926230103274805\\
0.530714285714286	0.00904298789943259\\
0.532755102040816	0.0089206482229488\\
0.534795918367347	0.00890575546740524\\
0.536836734693878	0.00912449684585093\\
0.540918367346939	0.00970616363208288\\
0.542959183673469	0.00972226925900066\\
0.545	0.00960781799344201\\
0.547040816326531	0.00940198987010543\\
0.551122448979592	0.00845421505612143\\
0.553163265306122	0.00797547134079957\\
0.555204081632653	0.00762500426109425\\
0.557244897959184	0.00758147424574318\\
0.559285714285714	0.00780656093176102\\
0.565408163265306	0.00898306396964521\\
0.567448979591837	0.00925265908752482\\
0.569489795918367	0.00933952427772089\\
0.571530612244898	0.0092960212985771\\
0.573571428571429	0.00895403799010996\\
0.57765306122449	0.0078926937821332\\
0.57969387755102	0.00767375342538756\\
0.581734693877551	0.00774317869350116\\
0.583775510204082	0.00793084942121036\\
0.585816326530612	0.00824587100494378\\
0.591938775510204	0.00960863852920535\\
0.596020408163265	0.00974671485130307\\
0.598061224489796	0.00963915863674203\\
0.604183673469388	0.00872684229214937\\
0.606224489795918	0.00861517838499792\\
0.608265306122449	0.00876282582792887\\
0.61234693877551	0.00944073924119937\\
0.616428571428571	0.0100491243773602\\
0.618469387755102	0.0101993502836213\\
0.620510204081633	0.0101939193182364\\
0.622551020408163	0.0100293456288955\\
0.624591836734694	0.00957706984712281\\
0.626632653061225	0.00902700921845723\\
0.628673469387755	0.00862721542589384\\
0.630714285714286	0.00857540806337287\\
0.634795918367347	0.00905028775025862\\
0.640918367346939	0.0101502888396353\\
0.642959183673469	0.0103304895145349\\
0.647040816326531	0.0100750240774654\\
0.649081632653061	0.00971782206338678\\
0.653163265306122	0.00885414029228437\\
0.655204081632653	0.00858823493195249\\
0.657244897959184	0.00863143159272972\\
0.659285714285714	0.00889650745519543\\
0.663367346938776	0.00963864762098343\\
0.667448979591837	0.0101945418922267\\
0.669489795918367	0.0102825387095655\\
0.671530612244898	0.0102300251207379\\
0.673571428571429	0.00997624452965296\\
0.67765306122449	0.00901372380138632\\
0.67969387755102	0.00874676108084638\\
0.681734693877551	0.00884036313212777\\
0.685816326530612	0.0093218055379295\\
0.691938775510204	0.010389494262652\\
0.693979591836735	0.0104969884029943\\
0.696020408163265	0.0103795604107871\\
0.698061224489796	0.0101691208600608\\
0.704183673469388	0.00897052580488178\\
0.706224489795918	0.00881246994994023\\
0.708265306122449	0.00895881916791597\\
0.71030612244898	0.00926584655192664\\
0.716428571428571	0.0103742121778346\\
0.718469387755102	0.0106812104175127\\
0.720510204081633	0.0107432147171649\\
0.722551020408163	0.0106713730435058\\
0.724591836734694	0.0103933950037673\\
0.728673469387755	0.00955768918062827\\
0.730714285714286	0.00944387367376121\\
0.734795918367347	0.00983343008615267\\
0.736836734693878	0.010101300497792\\
0.740918367346939	0.0108897053255109\\
0.742959183673469	0.0112895487875239\\
0.745	0.0113469751401455\\
0.747040816326531	0.0112586139756417\\
0.749081632653061	0.0110634844976205\\
0.755204081632653	0.0100264117241106\\
0.757244897959184	0.00994344006792258\\
0.759285714285714	0.0100462129067347\\
0.763367346938776	0.0105962694844697\\
0.769489795918367	0.0114657103516759\\
0.771530612244898	0.0114564184630247\\
0.773571428571429	0.0113303785750435\\
0.775612244897959	0.0110084773210495\\
0.77969387755102	0.0101013208479495\\
0.781734693877551	0.0100151952565896\\
0.785816326530612	0.0101833530874016\\
0.787857142857143	0.0103756741701378\\
0.789897959183673	0.0106708116134275\\
0.791938775510204	0.0111177900236145\\
0.793979591836735	0.0114146108012608\\
0.798061224489796	0.0112541990819836\\
0.800102040816327	0.0109997494663078\\
0.802142857142857	0.0106188976457476\\
0.804183673469388	0.0101421419757818\\
0.806224489795918	0.00982638374962885\\
0.808265306122449	0.00972362346556821\\
0.81030612244898	0.00976370529802717\\
0.81234693877551	0.00998011445661151\\
0.816428571428571	0.0106489220586672\\
0.818469387755102	0.0110589998231683\\
0.820510204081633	0.0111869606318326\\
0.822551020408163	0.0111622435520903\\
0.824591836734694	0.0109867026002768\\
0.826632653061224	0.010591093882498\\
0.828673469387755	0.0100170408894465\\
0.830714285714286	0.00960864612973822\\
0.834795918367347	0.00959023447327545\\
0.836836734693878	0.00966346344814828\\
0.838877551020408	0.00990215809790385\\
0.840918367346939	0.0103086636229253\\
0.842959183673469	0.0108787728767088\\
0.845	0.0110868587574557\\
0.849081632653061	0.0108342111237235\\
0.851122448979592	0.0104984672653369\\
0.857244897959184	0.00893874864465838\\
0.859285714285714	0.00876374603557062\\
0.861326530612245	0.00880504683430217\\
0.863367346938775	0.00911063731441919\\
0.865408163265306	0.00952353256190663\\
0.867448979591837	0.0101410204862636\\
0.869489795918367	0.0106596715189066\\
0.871530612244898	0.0108250037665276\\
0.873571428571429	0.0108510796459147\\
0.875612244897959	0.0106608387639057\\
0.87765306122449	0.0101932135665186\\
0.87969387755102	0.00955035138028137\\
0.881734693877551	0.0092046256851942\\
0.885816326530612	0.00909210880238054\\
0.887857142857143	0.00917557743424846\\
0.889897959183673	0.00948876130533161\\
0.893979591836735	0.0107305437586281\\
0.896020408163265	0.0108655183172437\\
0.898061224489796	0.0108352758441306\\
0.900102040816326	0.0106883840023422\\
0.902142857142857	0.0103350974657658\\
0.906224489795918	0.00917615335237387\\
0.908265306122449	0.00885902368060498\\
0.91030612244898	0.00872429050057\\
0.91234693877551	0.00881888868622838\\
0.914387755102041	0.00914951573537193\\
0.916428571428571	0.00963266030482957\\
0.918469387755102	0.0103277632583685\\
0.920510204081633	0.0107195467430437\\
0.922551020408163	0.010833209733783\\
0.924591836734694	0.0108267371845895\\
0.926632653061224	0.010559786494214\\
0.930714285714286	0.00941921245713728\\
0.932755102040816	0.00917834590125521\\
0.936836734693878	0.00907829328515319\\
0.938877551020408	0.00920881223162551\\
0.940918367346939	0.00958155355289225\\
0.942959183673469	0.0102888138818529\\
0.945	0.0107586323767835\\
0.949081632653061	0.0107330809635535\\
0.951122448979592	0.0105171663916795\\
0.953163265306122	0.0100928483714632\\
0.955204081632653	0.00946039088324346\\
0.957244897959184	0.00894090827056371\\
0.959285714285714	0.00870902863617684\\
0.961326530612245	0.00865833786720549\\
0.963367346938775	0.00886258326427214\\
0.965408163265306	0.00925208429705437\\
0.969489795918367	0.010585674260553\\
0.971530612244898	0.0108654567811836\\
0.973571428571429	0.0109765348948022\\
0.975612244897959	0.0109088413408652\\
0.97765306122449	0.0105825508208163\\
0.981734693877551	0.00944748504640713\\
0.983775510204082	0.00931835556590688\\
0.987857142857143	0.00929790304676492\\
0.989897959183674	0.00952828442816578\\
0.991938775510204	0.0100533904659993\\
0.993979591836735	0.0108406807966139\\
0.996020408163265	0.0111851210440393\\
0.998061224489796	0.0112042826394964\\
};
\addlegendentry{ABC$_W^{1/2}$ adaptive with approx. $\theta$}

\end{axis}
\end{tikzpicture}%

%% file: images/HighFq/Cost.tex
%
%
\definecolor{mycolor1}{rgb}{1.00000,0.60000,0.00000}%
\pgfplotsset{scaled y ticks=false}
\begin{tikzpicture}[scale = 0.52, font=\huge]

\begin{axis}[%
width=4.521in,
height=3.548in,
at={(0.758in,0.499in)},
scale only axis,
xmin=0,
xmax=1,
xlabel style={font=\huge\color{white!15!black}},
xlabel={relative simulation time $t/T$}, xtick={0, 0.2, 0.4, 0.6, 0.8, 1},
y tick label style={
        /pgf/number format/.cd,
            fixed,
            fixed zerofill,
            precision=2,
        /tikz/.cd
            },
ymin=0,
ymax=0.26, ytick={0.02, 0.06, 0.1, 0.14, 0.18, 0.22, 0.26},
ylabel style={at={(-0.15,0.5)}, font=\huge\color{white!15!black}},
ylabel={relative $L^{2}(\Omega)$ error}, ytick={0.02,  0.06, 0.1, 0.14, 0.18, 0.22, 0.26},
axis background/.style={fill=white},
xmajorgrids,
ymajorgrids,
legend style={at={(0.03,0.97)}, anchor=north west, legend cell align=left, align=left, draw=white!15!black}
]
\addplot [color=blue, line width=2.0pt]
  table[row sep=crcr]{%
0.00190476190476185	4.44089209850063e-16\\
0.185124716553288	7.87776842238186e-05\\
0.199637188208617	0.000453000394169956\\
0.206893424036281	0.000671496867433286\\
0.21233560090703	0.00104854829340828\\
0.223219954648526	0.0020505798442807\\
0.235918367346939	0.00323742916243464\\
0.2540589569161	0.00470114860269244\\
0.26312925170068	0.00548115266875704\\
0.293968253968254	0.00743694240278314\\
0.306666666666667	0.00802128462404716\\
0.319365079365079	0.00880235364486492\\
0.333877551020408	0.0094681985138263\\
0.362902494331066	0.0117444785655866\\
0.371972789115646	0.0128317473282318\\
0.379229024943311	0.0136648289664357\\
0.388299319727891	0.0146742683685638\\
0.393741496598639	0.0154946839748628\\
0.400997732426304	0.0168661430884721\\
0.413696145124717	0.0194410470397774\\
0.417324263038549	0.0203979167653165\\
0.420952380952381	0.0215293286046895\\
0.431836734693878	0.0252429191972392\\
0.43546485260771	0.0261889935620869\\
0.440907029478458	0.0273099579859456\\
0.453605442176871	0.0294709176826226\\
0.457233560090703	0.0299385504882506\\
0.460861678004535	0.0302238088738942\\
0.464489795918367	0.0303370728316564\\
0.471746031746032	0.0304286629226406\\
0.475374149659864	0.0306771975233229\\
0.479002267573696	0.0312267537983837\\
0.484444444444444	0.0324380575332536\\
0.497142857142857	0.0355648223446855\\
0.498956916099773	0.0361876154111263\\
0.500770975056689	0.0369263727774425\\
0.502585034013605	0.0377900840683443\\
0.506213151927438	0.0398314573777598\\
0.517097505668934	0.0465259948671887\\
0.526167800453515	0.0514765157879766\\
0.529795918367347	0.0538302090689908\\
0.533424036281179	0.0565533370321628\\
0.54249433106576	0.0637216306093213\\
0.546122448979592	0.0662210012086357\\
0.549750566893424	0.0684342501528721\\
0.553378684807256	0.070593001191572\\
0.557006802721088	0.073053191678299\\
0.560634920634921	0.0759416955338534\\
0.566077097505669	0.080737092125279\\
0.569705215419501	0.0837964354224635\\
0.571519274376417	0.0851658230359428\\
0.575147392290249	0.0875646431794084\\
0.58421768707483	0.0930259415323262\\
0.586031746031746	0.0943475741504939\\
0.589659863945578	0.097320546782016\\
0.596916099773243	0.103677227651333\\
0.598730158730159	0.105099268154273\\
0.600544217687075	0.106379026941851\\
0.602358276643991	0.10750965005809\\
0.605986394557823	0.109407356028313\\
0.609614512471655	0.111151365109934\\
0.611428571428571	0.112138378312852\\
0.613242630385488	0.113263654928765\\
0.615056689342404	0.114540287680367\\
0.618684807256236	0.1174793736328\\
0.624126984126984	0.122249220776552\\
0.6259410430839	0.123684449086114\\
0.627755102040816	0.124955217518481\\
0.629569160997732	0.126048527291154\\
0.631383219954649	0.126987973366198\\
0.635011337868481	0.128554328761861\\
0.638639455782313	0.130138917269199\\
0.640453514739229	0.131092749236634\\
0.642267573696145	0.132203525967263\\
0.644081632653061	0.133476453728476\\
0.647709750566893	0.136402532667218\\
0.651337868480726	0.139497858154432\\
0.653151927437642	0.140933497896591\\
0.654965986394558	0.142224385211882\\
0.656780045351474	0.143324075948735\\
0.65859410430839	0.144223774710728\\
0.660408163265306	0.144944977745664\\
0.664036281179138	0.146055212938599\\
0.665850340136054	0.146605786181047\\
0.66766439909297	0.147272864188086\\
0.669478458049887	0.148109330885496\\
0.671292517006803	0.14913693308382\\
0.673106575963719	0.15034654804743\\
0.676734693877551	0.153158812501011\\
0.680362811791383	0.156067627076715\\
0.682176870748299	0.157347313915846\\
0.683990929705215	0.158427610645445\\
0.685804988662131	0.159287032511777\\
0.687619047619048	0.159948165587855\\
0.69124716553288	0.160863166362741\\
0.694875283446712	0.1617611422216\\
0.696689342403628	0.162401019015192\\
0.698503401360544	0.163233473740216\\
0.70031746031746	0.164271021778611\\
0.702131519274376	0.165487649951767\\
0.705759637188209	0.168267745958664\\
0.709387755102041	0.17097081321704\\
0.711201814058957	0.172099188363605\\
0.713015873015873	0.173005140732014\\
0.714829931972789	0.173672542130833\\
0.716643990929705	0.174124036954586\\
0.720272108843537	0.174622681722925\\
0.722086167800454	0.174854973152092\\
0.72390022675737	0.175213538981842\\
0.725714285714286	0.175764578314235\\
0.727528344671202	0.176539629197955\\
0.729342403628118	0.177532464381023\\
0.731156462585034	0.178705765823957\\
0.736598639455782	0.18260435073892\\
0.738412698412698	0.183720616023416\\
0.740226757369615	0.184613441824218\\
0.742040816326531	0.185251808274318\\
0.743854875283447	0.185653082817524\\
0.745668934240363	0.185861184730992\\
0.751111111111111	0.186191712134919\\
0.752925170068027	0.18651622564957\\
0.754739229024943	0.187054113545425\\
0.756553287981859	0.187824777606657\\
0.758367346938775	0.188804971836423\\
0.761995464852608	0.191186397821741\\
0.76562358276644	0.193525572631895\\
0.767437641723356	0.194454910519806\\
0.769251700680272	0.195137381909822\\
0.771065759637188	0.195549723739176\\
0.772879818594104	0.195716448372026\\
0.77469387755102	0.195698467149047\\
0.778321995464853	0.195482323410412\\
0.780136054421769	0.195511955536584\\
0.781950113378685	0.195751369952785\\
0.783764172335601	0.196239560534524\\
0.785578231292517	0.196971887557823\\
0.787392290249433	0.19791205638282\\
0.794648526077097	0.202141185606097\\
0.796462585034014	0.202822004234425\\
0.79827664399093	0.203222563279338\\
0.800090702947846	0.203357016842873\\
0.801904761904762	0.203271835680419\\
0.80734693877551	0.202635270756213\\
0.809160997732426	0.202639650821977\\
0.810975056689342	0.202874633288824\\
0.812789115646259	0.203363895287166\\
0.814603174603175	0.204086810954719\\
0.816417233560091	0.204999148581802\\
0.821859410430839	0.207974525719275\\
0.823673469387755	0.20870747269722\\
0.825487528344671	0.209180363961573\\
0.827301587301587	0.209365864725317\\
0.829115646258503	0.209285497540959\\
0.83092970521542	0.209004519375766\\
0.834557823129252	0.208227539027179\\
0.836371882086168	0.207973627137729\\
0.838185941043084	0.207942451884126\\
0.84	0.208175477273636\\
0.841814058956916	0.208670582987479\\
0.843628117913832	0.20939887154549\\
0.847256235827664	0.211262673100953\\
0.849070294784581	0.212191823874747\\
0.850884353741497	0.212961377265133\\
0.852698412698413	0.213486230583991\\
0.854512471655329	0.213723222945491\\
0.856326530612245	0.21367947403985\\
0.858140589569161	0.213400111564006\\
0.863582766439909	0.212105978596468\\
0.865396825396825	0.211888138695615\\
0.867210884353741	0.211912105481851\\
0.869024943310658	0.212205014599041\\
0.870839002267574	0.212750391680403\\
0.87265306122449	0.213506179387181\\
0.878095238095238	0.216121264072206\\
0.879909297052154	0.21674959174642\\
0.88172335600907	0.217118658913013\\
0.883537414965986	0.217194701578415\\
0.885351473922902	0.216993820234914\\
0.887165532879819	0.216580787063374\\
0.890793650793651	0.215504442142448\\
0.892607709750567	0.215090940900877\\
0.894421768707483	0.214905224023411\\
0.896235827664399	0.214992975380964\\
0.898049886621315	0.215355133143611\\
0.899863945578231	0.215969031317937\\
0.903492063492063	0.217663412994135\\
0.90530612244898	0.218531998238032\\
0.907120181405896	0.21924814394763\\
0.908934240362812	0.219724008966622\\
0.910748299319728	0.219911651622476\\
0.912562358276644	0.219809224423999\\
0.91437641723356	0.219460610859183\\
0.919818594104308	0.217883716913031\\
0.921632653061224	0.217561176384443\\
0.923446712018141	0.217485422462667\\
0.925260770975057	0.217687832428344\\
0.927074829931973	0.218153912527888\\
0.928888888888889	0.218846996346148\\
0.934331065759637	0.221365520386914\\
0.936145124716553	0.221977695687423\\
0.937959183673469	0.222332407419562\\
0.939773242630386	0.222390674171202\\
0.941587301587302	0.222162267934177\\
0.943401360544218	0.221705749796462\\
0.948843537414966	0.219989424907138\\
0.950657596371882	0.219716118862469\\
0.952471655328798	0.219721329015178\\
0.954285714285714	0.220008161913241\\
0.95609977324263	0.220557976929962\\
0.957913832199547	0.221310679764897\\
0.961541950113379	0.22301568580306\\
0.963356009070295	0.223717234380783\\
0.965170068027211	0.224182037795289\\
0.966984126984127	0.224354979859671\\
0.968798185941043	0.22422622004309\\
0.970612244897959	0.223839912272664\\
0.977868480725624	0.221654722900251\\
0.97968253968254	0.221493459471968\\
0.981496598639456	0.221613406453378\\
0.983310657596372	0.222006296218369\\
0.985124716553288	0.222639760837399\\
0.990566893424036	0.22505111001629\\
0.992380952380952	0.225639500016044\\
0.994195011337868	0.225972028707427\\
0.996009070294785	0.226004702494391\\
0.997823129251701	0.225739177310179\\
0.999637188208617	0.225227065264749\\
};
\addlegendentry{ABC$_W^0$ adaptive}

\addplot [color=mycolor1, line width=2.0pt]
  table[row sep=crcr]{%
0.00190476190476185	4.44089209850063e-16\\
0.185124716553288	7.87497743746091e-05\\
0.199637188208617	0.000452806469949052\\
0.206893424036281	0.000671685805582634\\
0.21233560090703	0.00105023099571266\\
0.223219954648526	0.00205683081791497\\
0.235918367346939	0.00325689151007535\\
0.252244897959184	0.00459607414534413\\
0.26312925170068	0.0056086622727789\\
0.28671201814059	0.00740901105390634\\
0.293968253968254	0.00789968693741838\\
0.304852607709751	0.00854206333893581\\
0.322993197278912	0.00995074675573016\\
0.332063492063492	0.0104247300861359\\
0.339319727891156	0.0110141271694447\\
0.346575963718821	0.0115789321337036\\
0.353832199546485	0.0119108533178272\\
0.36108843537415	0.0122906635855513\\
0.37015873015873	0.013071123368007\\
0.375600907029478	0.0134808234961211\\
0.382857142857143	0.0138111904676127\\
0.391927437641723	0.0142566839639685\\
0.399183673469388	0.0148456422355383\\
0.413696145124717	0.0161700226080294\\
0.417324263038549	0.0167123753006598\\
0.420952380952381	0.0174470844863349\\
0.433650793650794	0.02038035378734\\
0.439092970521542	0.0211953124292996\\
0.446349206349206	0.0219946506516526\\
0.455419501133787	0.0229255695475064\\
0.459047619047619	0.0231455915790678\\
0.462675736961451	0.0231766864789268\\
0.4681179138322	0.0229703244309565\\
0.473560090702948	0.0227422428710317\\
0.47718820861678	0.0227411940092282\\
0.484444444444444	0.0228149946019855\\
0.488072562358277	0.0225817952192507\\
0.500770975056689	0.021277319934207\\
0.506213151927438	0.0212117243002117\\
0.511655328798186	0.0211418706756615\\
0.515283446712018	0.0208704076804916\\
0.526167800453515	0.0198108457046864\\
0.529795918367347	0.0196121727459546\\
0.537052154195011	0.0195284097186231\\
0.54249433106576	0.0191651465699723\\
0.553378684807256	0.0181465268039762\\
0.558820861678004	0.0178859088356064\\
0.564263038548753	0.0176852945109385\\
0.567891156462585	0.017466782892779\\
0.573333333333333	0.0168594289343217\\
0.586031746031746	0.0155290883247046\\
0.589659863945578	0.0153722970249763\\
0.596916099773243	0.0147589549899472\\
0.604172335600907	0.0140707445192881\\
0.611428571428571	0.0134352157133868\\
0.615056689342404	0.0133041233085223\\
0.618684807256236	0.0131810431708944\\
0.6259410430839	0.0127514101421193\\
0.633197278911565	0.012151122739542\\
0.642267573696145	0.011561283614014\\
0.645895691609977	0.0115523982354112\\
0.64952380952381	0.0113582351370177\\
0.65859410430839	0.0105889946028558\\
0.669478458049887	0.00998588082418184\\
0.673106575963719	0.0100102867944711\\
0.678548752834467	0.0097624171099342\\
0.682176870748299	0.00945600191089446\\
0.687619047619048	0.00888355871457769\\
0.693061224489796	0.0086360740836996\\
0.698503401360544	0.00842706505436219\\
0.702131519274376	0.00855185730863861\\
0.705759637188209	0.00838836160686052\\
0.714829931972789	0.00765337625156093\\
0.725714285714286	0.00739620511222738\\
0.729342403628118	0.00761980720281019\\
0.73297052154195	0.00756590474106855\\
0.736598639455782	0.00750192538009753\\
0.743854875283447	0.00700408206918657\\
0.749297052154195	0.00716328225031382\\
0.752925170068027	0.00721199751759172\\
0.754739229024943	0.00733859994513619\\
0.758367346938775	0.00775520150566156\\
0.761995464852608	0.00784709962571806\\
0.772879818594104	0.00767557931790575\\
0.781950113378685	0.00796989834235173\\
0.785578231292517	0.00838408485324615\\
0.791020408163265	0.0086857638325939\\
0.794648526077097	0.00879856010139457\\
0.801904761904762	0.00875013372805489\\
0.805532879818594	0.00880916507469109\\
0.809160997732426	0.00881477327538482\\
0.812789115646259	0.00908470648608228\\
0.816417233560091	0.00928821939536462\\
0.820045351473923	0.0092592770173292\\
0.827301587301587	0.00905363223202704\\
0.84	0.00927249863047719\\
0.843628117913832	0.00942290283920344\\
0.850884353741497	0.00939792482312818\\
0.856326530612245	0.00912662350224136\\
0.861768707482993	0.00931088530991298\\
0.867210884353741	0.00940779668029801\\
0.870839002267574	0.0096820462809194\\
0.874467120181406	0.00976057844999245\\
0.88172335600907	0.00957616021676522\\
0.885351473922902	0.00951386741204874\\
0.888979591836735	0.0095491078773271\\
0.892607709750567	0.0094483455637836\\
0.894421768707483	0.00947487810685166\\
0.898049886621315	0.00970430605154493\\
0.90530612244898	0.00984483300717875\\
0.908934240362812	0.0096774016149217\\
0.912562358276644	0.00948754425839915\\
0.918004535147392	0.00963017284554024\\
0.921632653061224	0.00962577912997842\\
0.923446712018141	0.00970423940828669\\
0.927074829931973	0.0100181871186057\\
0.930702947845805	0.0102189236609443\\
0.934331065759637	0.0102473315537271\\
0.941587301587302	0.0101573587280829\\
0.945215419501134	0.0101691375295897\\
0.948843537414966	0.0100344243788549\\
0.950657596371882	0.0100561353023771\\
0.954285714285714	0.010327292410544\\
0.963356009070295	0.0106640326858271\\
0.968798185941043	0.0104348244295858\\
0.974240362811791	0.0105451133039857\\
0.977868480725624	0.0105369507373079\\
0.97968253968254	0.0106257969084275\\
0.985124716553288	0.0111301045289377\\
0.98875283446712	0.0113551708539015\\
0.994195011337868	0.0113539062413155\\
0.999637188208617	0.0112777311416224\\
};
\addlegendentry{ABC$_W^{1/2}$ adaptive}

\addplot [color=blue, dashed, line width=2.0pt]
  table[row sep=crcr]{%
0.00190476190476185	4.44089209850063e-16\\
0.185124716553288	7.87776842238186e-05\\
0.199637188208617	0.000453000394169956\\
0.206893424036281	0.000671496867433286\\
0.210521541950113	0.000917036425774675\\
0.215963718820862	0.00153225669818835\\
0.237732426303855	0.00458523029897662\\
0.241360544217687	0.00538709639234525\\
0.244988662131519	0.00639507632365466\\
0.250430839002268	0.00815182080978172\\
0.257687074829932	0.0104765655848597\\
0.26312925170068	0.0119529883034962\\
0.266757369614512	0.0129439614761967\\
0.270385487528345	0.0141615692006891\\
0.274013605442177	0.0157088665657845\\
0.279455782312925	0.0183974930390344\\
0.284897959183673	0.0210677483910157\\
0.297596371882086	0.027060264701447\\
0.301224489795918	0.0291086072945755\\
0.304852607709751	0.0314277299419591\\
0.315736961451247	0.0389122695019203\\
0.319365079365079	0.0409876953028491\\
0.322993197278912	0.0427847484148904\\
0.32843537414966	0.045448949110815\\
0.332063492063492	0.0474885792401404\\
0.335691609977324	0.0497749592386152\\
0.342947845804989	0.0545111018063593\\
0.346575963718821	0.0565698573458684\\
0.350204081632653	0.0583336455388372\\
0.355646258503401	0.0608393508981702\\
0.359274376417234	0.0627894687123975\\
0.362902494331066	0.0651174561239819\\
0.371972789115646	0.0715686723315941\\
0.373786848072562	0.0726791814722032\\
0.375600907029478	0.0736725349276065\\
0.379229024943311	0.0753587347496518\\
0.384671201814059	0.0776954311839704\\
0.388299319727891	0.0795742716517535\\
0.391927437641723	0.0818520058199147\\
0.400997732426304	0.0880890534889132\\
0.40281179138322	0.0891440150749381\\
0.406439909297052	0.0909295946424227\\
0.415510204081633	0.0948040949969029\\
0.417324263038549	0.0957802084013548\\
0.419138321995465	0.0968871800557881\\
0.420952380952381	0.098131686969965\\
0.424580498866213	0.100983800667577\\
0.430022675736962	0.105545272174844\\
0.431836734693878	0.106896892936257\\
0.433650793650794	0.108066162474323\\
0.43546485260771	0.109023435860961\\
0.437278911564626	0.109768813287007\\
0.439092970521542	0.110332446235954\\
0.44453514739229	0.111619492742586\\
0.446349206349206	0.112218086783179\\
0.448163265306122	0.113052261514575\\
0.449977324263039	0.114182264040138\\
0.451791383219955	0.115629894125145\\
0.453605442176871	0.117368223694098\\
0.460861678004535	0.125137408160936\\
0.462675736961451	0.126614322385228\\
0.464489795918367	0.127718419342738\\
0.466303854875283	0.128446965274407\\
0.4681179138322	0.12885182628734\\
0.469931972789116	0.129030338058836\\
0.473560090702948	0.129238421247887\\
0.475374149659864	0.129561265522688\\
0.47718820861678	0.130215543267536\\
0.479002267573696	0.131303282083627\\
0.480816326530612	0.132868469305764\\
0.482630385487528	0.134877160487508\\
0.486258503401361	0.139669920614161\\
0.488072562358277	0.142033004678212\\
0.489886621315193	0.144068131771923\\
0.491700680272109	0.145597494713682\\
0.493514739229025	0.146560847911213\\
0.495328798185941	0.147027896023576\\
0.497142857142857	0.147145796490444\\
0.500770975056689	0.147023459317289\\
0.502585034013605	0.147115384230692\\
0.504399092970522	0.147509081208428\\
0.506213151927438	0.148319237673775\\
0.508027210884354	0.149619115055246\\
0.50984126984127	0.151408831390482\\
0.511655328798186	0.153590518940943\\
0.517097505668934	0.160858437031893\\
0.51891156462585	0.162959367394096\\
0.520725623582766	0.164603545234335\\
0.522539682539683	0.165630650796227\\
0.524353741496599	0.165948623571385\\
0.526167800453515	0.165559927098162\\
0.527981859410431	0.164568405320725\\
0.533424036281179	0.160404103314315\\
0.535238095238095	0.1597257834316\\
0.537052154195011	0.159977974365272\\
0.538866213151927	0.161368653599981\\
0.540680272108844	0.163879062959815\\
0.54249433106576	0.16722559802697\\
0.544308390022676	0.170898071945368\\
0.546122448979592	0.174278801143717\\
0.547936507936508	0.176830996378671\\
0.549750566893424	0.178241600979285\\
0.55156462585034	0.178420354867206\\
0.553378684807256	0.177442765668847\\
0.555192743764172	0.175514735853443\\
0.560634920634921	0.167948340798239\\
0.562448979591837	0.166569294102375\\
0.564263038548753	0.166558084481349\\
0.566077097505669	0.168071745496898\\
0.567891156462585	0.170942237719246\\
0.569705215419501	0.174755480227589\\
0.573333333333333	0.182961538568944\\
0.575147392290249	0.186182701105229\\
0.576961451247166	0.188161342589281\\
0.578775510204082	0.188624389265876\\
0.580589569160998	0.187533637096762\\
0.582403628117914	0.18508264059559\\
0.58421768707483	0.181662226508454\\
0.587845804988662	0.1741538588907\\
0.589659863945578	0.171336336210139\\
0.591473922902494	0.169958871218748\\
0.59328798185941	0.170451188048678\\
0.595102040816327	0.172908982497662\\
0.596916099773243	0.177001459085442\\
0.600544217687075	0.18709957857908\\
0.602358276643991	0.19132444936743\\
0.604172335600907	0.194035652792753\\
0.605986394557823	0.194914058960103\\
0.607800453514739	0.193967719866179\\
0.609614512471655	0.191450092291673\\
0.611428571428571	0.187801753596038\\
0.615056689342404	0.179617991090747\\
0.61687074829932	0.176548248230309\\
0.618684807256236	0.175057407885824\\
0.620498866213152	0.175534098589952\\
0.622312925170068	0.177961482402911\\
0.624126984126984	0.181913397611143\\
0.627755102040816	0.191577346175141\\
0.629569160997732	0.195812733758636\\
0.631383219954649	0.19883153549662\\
0.633197278911565	0.200250950589261\\
0.635011337868481	0.199935286078416\\
0.636825396825397	0.198002551593818\\
0.638639455782313	0.19478923225918\\
0.644081632653061	0.182911033610762\\
0.645895691609977	0.18027257221187\\
0.647709750566893	0.179218826356613\\
0.64952380952381	0.180063083552186\\
0.651337868480726	0.182788179136405\\
0.653151927437642	0.186987614084692\\
0.656780045351474	0.196798373735723\\
0.65859410430839	0.200721396891478\\
0.660408163265306	0.203104871065182\\
0.662222222222222	0.203684164648568\\
0.664036281179138	0.202513010736156\\
0.665850340136054	0.199877644194207\\
0.66766439909297	0.196229238042777\\
0.669478458049887	0.192151370937909\\
0.671292517006803	0.188332253818514\\
0.673106575963719	0.18547248581477\\
0.674920634920635	0.184155410645428\\
0.676734693877551	0.18471257245259\\
0.678548752834467	0.187096506912838\\
0.680362811791383	0.190881925464164\\
0.683990929705215	0.199950531044329\\
0.685804988662131	0.203846476899924\\
0.687619047619048	0.206563131874989\\
0.689433106575964	0.207755738686843\\
0.69124716553288	0.207313116278809\\
0.693061224489796	0.205368421015333\\
0.694875283446712	0.20226115535381\\
0.698503401360544	0.194591631819983\\
0.70031746031746	0.191198342522262\\
0.702131519274376	0.188842152059711\\
0.703945578231293	0.187964137429268\\
0.705759637188209	0.188822964735953\\
0.707573696145125	0.191390286027841\\
0.709387755102041	0.195295248392315\\
0.713015873015873	0.204349720291004\\
0.714829931972789	0.207930466671174\\
0.716643990929705	0.210048059569336\\
0.718458049886621	0.210454960431753\\
0.720272108843537	0.20921576100686\\
0.722086167800454	0.206621899100793\\
0.72390022675737	0.203118192690058\\
0.725714285714286	0.199264071652181\\
0.727528344671202	0.195703734727042\\
0.729342403628118	0.193082280608647\\
0.731156462585034	0.191927648534643\\
0.73297052154195	0.19252932814297\\
0.734784580498866	0.19482507009032\\
0.736598639455782	0.198404006493884\\
0.740226757369615	0.206842641017026\\
0.742040816326531	0.210395232440069\\
0.743854875283447	0.212803238798834\\
0.745668934240363	0.213754707168268\\
0.747482993197279	0.213162528493394\\
0.749297052154195	0.211171880572793\\
0.751111111111111	0.208120050588904\\
0.754739229024943	0.200791275163369\\
0.756553287981859	0.197615986542497\\
0.758367346938775	0.195446949332518\\
0.760181405895692	0.194674017297037\\
0.761995464852608	0.19551505941031\\
0.763809523809524	0.197932698200135\\
0.76562358276644	0.20158461646658\\
0.769251700680272	0.210011060003524\\
0.771065759637188	0.213310075759905\\
0.772879818594104	0.215206489983077\\
0.77469387755102	0.215466052630281\\
0.776507936507936	0.21416314158604\\
0.778321995464853	0.211592445358677\\
0.780136054421769	0.208189738518916\\
0.781950113378685	0.204489907714039\\
0.783764172335601	0.201101178926647\\
0.785578231292517	0.198621126993926\\
0.787392290249433	0.197531922252355\\
0.789206349206349	0.19809405618453\\
0.791020408163265	0.200243276343014\\
0.792834467120181	0.203587395625128\\
0.796462585034014	0.211413954045036\\
0.79827664399093	0.214657704680685\\
0.800090702947846	0.21680212013137\\
0.801904761904762	0.217561340099594\\
0.803718820861678	0.216863644324793\\
0.805532879818594	0.214856741164854\\
0.80734693877551	0.211871378472512\\
0.810975056689342	0.204857990337015\\
0.812789115646259	0.201870895310858\\
0.814603174603175	0.199854935673207\\
0.816417233560091	0.199158044690053\\
0.818231292517007	0.199967547693223\\
0.820045351473923	0.202242335680204\\
0.821859410430839	0.205666784612281\\
0.825487528344671	0.213548306790878\\
0.827301587301587	0.216616415768931\\
0.829115646258503	0.218346336532194\\
0.83092970521542	0.218516987378313\\
0.832743764172336	0.217207932887797\\
0.834557823129252	0.214712170766366\\
0.836371882086168	0.211460404824768\\
0.838185941043084	0.207969554728042\\
0.84	0.204810147583857\\
0.841814058956916	0.202534470706374\\
0.843628117913832	0.201586179820152\\
0.845442176870748	0.20220052452812\\
0.847256235827664	0.204305512683055\\
0.849070294784581	0.207518662194727\\
0.852698412698413	0.214913355718962\\
0.854512471655329	0.217914282723236\\
0.856326530612245	0.219838730781943\\
0.858140589569161	0.220424455011997\\
0.859954648526077	0.219617146759192\\
0.861768707482993	0.217575068205565\\
0.863582766439909	0.214633831111747\\
0.867210884353741	0.207939965490297\\
0.869024943310658	0.205192597098761\\
0.870839002267574	0.203423703442209\\
0.87265306122449	0.202935163811105\\
0.874467120181406	0.203870584158485\\
0.876281179138322	0.206159481279082\\
0.878095238095238	0.209486393835626\\
0.88172335600907	0.216945520484332\\
0.883537414965986	0.219762094164469\\
0.885351473922902	0.221257795348999\\
0.887165532879819	0.221235083848478\\
0.888979591836735	0.219798506888899\\
0.890793650793651	0.217262263278092\\
0.896235827664399	0.207839275831646\\
0.898049886621315	0.205856102981248\\
0.899863945578231	0.205197152430304\\
0.901678004535147	0.206043891235799\\
0.903492063492063	0.208285828307481\\
0.90530612244898	0.211523735128476\\
0.908934240362812	0.218637010543262\\
0.910748299319728	0.221373047093918\\
0.912562358276644	0.222996570940905\\
0.91437641723356	0.223286912691297\\
0.916190476190476	0.222227512158863\\
0.918004535147392	0.220004382024094\\
0.919818594104308	0.216975398783887\\
0.923446712018141	0.21043526428239\\
0.925260770975057	0.207908947847204\\
0.927074829931973	0.206409222588142\\
0.928888888888889	0.206181199635591\\
0.930702947845805	0.207309391350377\\
0.932517006802721	0.209677829888568\\
0.934331065759637	0.212955902996668\\
0.936145124716553	0.216613780536978\\
0.937959183673469	0.220015776625478\\
0.939773242630386	0.222555612206899\\
0.941587301587302	0.223778090876583\\
0.943401360544218	0.223509863904772\\
0.945215419501134	0.221881779630412\\
0.94702947845805	0.219238127042793\\
0.950657596371882	0.212840773867577\\
0.952471655328798	0.21013992726573\\
0.954285714285714	0.208430598642059\\
0.95609977324263	0.208056620913384\\
0.957913832199547	0.209142925040958\\
0.959727891156463	0.211539971879687\\
0.961541950113379	0.214826424511075\\
0.963356009070295	0.218416727402474\\
0.965170068027211	0.22171014452893\\
0.966984126984127	0.224203255265359\\
0.968798185941043	0.225542907777614\\
0.970612244897959	0.225551707486687\\
0.972426303854875	0.224251292507474\\
0.974240362811791	0.221852869501768\\
0.976054421768707	0.218734303633353\\
0.977868480725624	0.215383836009994\\
0.97968253968254	0.212307793182423\\
0.981496598639456	0.209954111382751\\
0.983310657596372	0.208658827383454\\
0.985124716553288	0.208619803716745\\
0.986938775510204	0.209874928656867\\
0.98875283446712	0.212275888604814\\
0.990566893424036	0.215486230407503\\
0.992380952380952	0.218996907785589\\
0.994195011337868	0.22220317894853\\
0.996009070294785	0.224535252196379\\
0.997823129251701	0.225563958932727\\
0.999637188208617	0.225124965564847\\
};
\addlegendentry{ABC$_W^0$}

\addplot [color=mycolor1, dashed, line width=2.0pt]
  table[row sep=crcr]{%
0.00190476190476185	4.44089209850063e-16\\
0.185124716553288	7.87497743746091e-05\\
0.199637188208617	0.000452806469949052\\
0.206893424036281	0.000671685805582634\\
0.210521541950113	0.000918143319634912\\
0.215963718820862	0.00153535429954688\\
0.237732426303855	0.00458488204035779\\
0.241360544217687	0.0053765455366771\\
0.244988662131519	0.00637259406821811\\
0.250430839002268	0.00811094564330628\\
0.257687074829932	0.0104123209501711\\
0.268571428571429	0.0134661293136978\\
0.272199546485261	0.01488419169865\\
0.275827664399093	0.0166134251834711\\
0.284897959183673	0.0212080340023207\\
0.290340136054422	0.0237236436663736\\
0.29578231292517	0.0262128608971757\\
0.299410430839002	0.0280577431033316\\
0.303038548752835	0.0301524707754195\\
0.306666666666667	0.0324892651499872\\
0.313922902494331	0.0373594951880835\\
0.317551020408163	0.0395131284292085\\
0.321179138321996	0.0413658816444638\\
0.32843537414966	0.0448844228200198\\
0.332063492063492	0.0469559331022367\\
0.335691609977324	0.0492839210395731\\
0.342947845804989	0.0540831420086233\\
0.346575963718821	0.0561447287739215\\
0.350204081632653	0.0578857206907144\\
0.355646258503401	0.0603092800642656\\
0.359274376417234	0.0621754333932188\\
0.362902494331066	0.064398588660342\\
0.371972789115646	0.070534449429132\\
0.373786848072562	0.0715844446568922\\
0.377414965986395	0.0733590100204509\\
0.386485260770975	0.0772692714974933\\
0.390113378684807	0.0793387699246392\\
0.393741496598639	0.081770590335464\\
0.399183673469388	0.0855610690357549\\
0.40281179138322	0.0877376125806554\\
0.406439909297052	0.0894687800851015\\
0.415510204081633	0.0931091870688069\\
0.417324263038549	0.0940176781466917\\
0.419138321995465	0.0950500968493768\\
0.420952380952381	0.0962140660425539\\
0.424580498866213	0.0988946424023071\\
0.430022675736962	0.103210699527616\\
0.431836734693878	0.104492267885744\\
0.433650793650794	0.105598812795184\\
0.43546485260771	0.106500191023954\\
0.437278911564626	0.107195348260445\\
0.439092970521542	0.107712671858642\\
0.44453514739229	0.108864852590239\\
0.446349206349206	0.109415846875537\\
0.448163265306122	0.11020112779749\\
0.449977324263039	0.111282157961194\\
0.451791383219955	0.112682377472073\\
0.453605442176871	0.11437612616895\\
0.460861678004535	0.12198235405761\\
0.462675736961451	0.123409914551019\\
0.464489795918367	0.124454728325795\\
0.466303854875283	0.125110967862431\\
0.4681179138322	0.125427496156411\\
0.469931972789116	0.125499189439339\\
0.473560090702948	0.125436381213674\\
0.475374149659864	0.125600059817071\\
0.47718820861678	0.126086761477801\\
0.479002267573696	0.127006463723735\\
0.480816326530612	0.128411025585409\\
0.482630385487528	0.130272103462384\\
0.486258503401361	0.13481194346239\\
0.488072562358277	0.137063071209052\\
0.489886621315193	0.138987227400974\\
0.491700680272109	0.140399302301244\\
0.493514739229025	0.141234186082841\\
0.495328798185941	0.141558454893092\\
0.497142857142857	0.141516720625902\\
0.500770975056689	0.141021592378245\\
0.502585034013605	0.140905105839879\\
0.504399092970522	0.141083480762229\\
0.506213151927438	0.141679574646292\\
0.508027210884354	0.142775458757936\\
0.50984126984127	0.144377013812061\\
0.511655328798186	0.14638729942949\\
0.517097505668934	0.153229203265553\\
0.51891156462585	0.155204058918087\\
0.520725623582766	0.156715547263456\\
0.522539682539683	0.157591276409173\\
0.524353741496599	0.157727223812902\\
0.526167800453515	0.157115652143533\\
0.527981859410431	0.155854163180165\\
0.529795918367347	0.154146987873791\\
0.531609977324263	0.152291121867369\\
0.533424036281179	0.150645843816408\\
0.535238095238095	0.149605675601354\\
0.537052154195011	0.149542381506239\\
0.538866213151927	0.150692966107114\\
0.540680272108844	0.153048976804867\\
0.54249433106576	0.156312786442695\\
0.544308390022676	0.159940450373821\\
0.546122448979592	0.163272202181267\\
0.547936507936508	0.16573705245369\\
0.549750566893424	0.167001589165251\\
0.55156462585034	0.166962464175169\\
0.553378684807256	0.165685420865327\\
0.555192743764172	0.16337152703429\\
0.557006802721088	0.160358329215458\\
0.558820861678004	0.157127458267809\\
0.560634920634921	0.154264117851074\\
0.562448979591837	0.15237415414576\\
0.564263038548753	0.151945755739614\\
0.566077097505669	0.153171979953344\\
0.567891156462585	0.155889176027052\\
0.569705215419501	0.159662202703122\\
0.573333333333333	0.167968568814245\\
0.575147392290249	0.171237869022056\\
0.576961451247166	0.173207459477052\\
0.578775510204082	0.173571438805749\\
0.580589569160998	0.172266700990024\\
0.582403628117914	0.169471397664625\\
0.58421768707483	0.165573988838054\\
0.587845804988662	0.156802138238135\\
0.589659863945578	0.153312347988774\\
0.591473922902494	0.15135936346333\\
0.59328798185941	0.151468854769708\\
0.595102040816327	0.153787111501439\\
0.596916099773243	0.157961547713462\\
0.600544217687075	0.168534930001771\\
0.602358276643991	0.172948824855308\\
0.604172335600907	0.175715306661018\\
0.605986394557823	0.176479090125989\\
0.607800453514739	0.175229698318098\\
0.609614512471655	0.172212853942829\\
0.611428571428571	0.167876500323542\\
0.615056689342404	0.157918393768788\\
0.61687074829932	0.153941453819736\\
0.618684807256236	0.151697944362835\\
0.620498866213152	0.151687823748432\\
0.622312925170068	0.153932908364278\\
0.624126984126984	0.157962365577325\\
0.627755102040816	0.168183024883062\\
0.629569160997732	0.172706662450854\\
0.631383219954649	0.175921115491994\\
0.633197278911565	0.177395472547574\\
0.635011337868481	0.1769619016376\\
0.636825396825397	0.174720152813627\\
0.638639455782313	0.171003660898152\\
0.640453514739229	0.166327916855489\\
0.642267573696145	0.161341390887776\\
0.644081632653061	0.156763156231883\\
0.645895691609977	0.153315645979518\\
0.647709750566893	0.151646620413922\\
0.64952380952381	0.152180741310807\\
0.651337868480726	0.154937976412118\\
0.653151927437642	0.159452740558131\\
0.656780045351474	0.170219160187935\\
0.65859410430839	0.174511945040652\\
0.660408163265306	0.177065097874866\\
0.662222222222222	0.177572056051517\\
0.664036281179138	0.176070497691827\\
0.665850340136054	0.172845641764458\\
0.66766439909297	0.16836657816738\\
0.671292517006803	0.158332484903024\\
0.673106575963719	0.154406149321153\\
0.674920634920635	0.152236128030773\\
0.676734693877551	0.152287508948284\\
0.678548752834467	0.154552151580825\\
0.680362811791383	0.158540245080713\\
0.683990929705215	0.168489721869582\\
0.685804988662131	0.172824771611741\\
0.687619047619048	0.175859191567308\\
0.689433106575964	0.177190842762534\\
0.69124716553288	0.176671803804498\\
0.693061224489796	0.174417514586383\\
0.694875283446712	0.170769806841869\\
0.70031746031746	0.157143697035402\\
0.702131519274376	0.153940017269233\\
0.703945578231293	0.152448965807443\\
0.705759637188209	0.153048333541923\\
0.707573696145125	0.155747843639368\\
0.709387755102041	0.160109494075319\\
0.713015873015873	0.170447009666798\\
0.714829931972789	0.174538861268485\\
0.716643990929705	0.176928538977479\\
0.718458049886621	0.177318744342183\\
0.720272108843537	0.175758563963369\\
0.722086167800454	0.172544195990772\\
0.72390022675737	0.16814834661095\\
0.727528344671202	0.158440657567517\\
0.729342403628118	0.154694794765118\\
0.731156462585034	0.152667809211868\\
0.73297052154195	0.152794292876088\\
0.734784580498866	0.155051033267817\\
0.736598639455782	0.158950744014614\\
0.740226757369615	0.168548594962631\\
0.742040816326531	0.172667990711626\\
0.743854875283447	0.175500838732534\\
0.745668934240363	0.176672891614514\\
0.747482993197279	0.176058330783885\\
0.749297052154195	0.173784857340956\\
0.751111111111111	0.170197226431704\\
0.756553287981859	0.157110201254444\\
0.758367346938775	0.154109477107854\\
0.760181405895692	0.15276169668325\\
0.761995464852608	0.153405214082317\\
0.763809523809524	0.156042423428591\\
0.76562358276644	0.160261980912744\\
0.769251700680272	0.170226101552634\\
0.771065759637188	0.174154387285316\\
0.772879818594104	0.17641603682793\\
0.77469387755102	0.176721275151862\\
0.776507936507936	0.17513122101203\\
0.778321995464853	0.17195226866596\\
0.780136054421769	0.167656621494737\\
0.783764172335601	0.158271003162946\\
0.785578231292517	0.154683874203657\\
0.787392290249433	0.152761436487808\\
0.789206349206349	0.152911331812826\\
0.791020408163265	0.155107241621255\\
0.792834467120181	0.158875922610949\\
0.796462585034014	0.168087367957897\\
0.79827664399093	0.172003387174307\\
0.800090702947846	0.174665968915207\\
0.801904761904762	0.175724527291322\\
0.803718820861678	0.175068104640756\\
0.805532879818594	0.172827686853095\\
0.80734693877551	0.169346236349872\\
0.812789115646259	0.156879353007167\\
0.814603174603175	0.154079894870152\\
0.816417233560091	0.152857606803465\\
0.818231292517007	0.153517221563594\\
0.820045351473923	0.156057756247151\\
0.821859410430839	0.160099167562438\\
0.825487528344671	0.169635795270395\\
0.827301587301587	0.173391920919911\\
0.829115646258503	0.175534350284026\\
0.83092970521542	0.175781348450271\\
0.832743764172336	0.174201239675649\\
0.834557823129252	0.17110467704972\\
0.836371882086168	0.166964942278194\\
0.84	0.158043753673226\\
0.841814058956916	0.154689070054657\\
0.843628117913832	0.152941732340251\\
0.845442176870748	0.153180575272718\\
0.847256235827664	0.1553708048999\\
0.849070294784581	0.159048059554556\\
0.852698412698413	0.167900826559084\\
0.854512471655329	0.171603794902853\\
0.856326530612245	0.174070796399659\\
0.858140589569161	0.174975150475697\\
0.859954648526077	0.17422455182332\\
0.861768707482993	0.171962907348542\\
0.863582766439909	0.16854196958866\\
0.867210884353741	0.160335687528609\\
0.869024943310658	0.156752444333777\\
0.870839002267574	0.154268659589846\\
0.87265306122449	0.153326599473184\\
0.874467120181406	0.154174319541324\\
0.876281179138322	0.156771322380121\\
0.878095238095238	0.160737815818797\\
0.88172335600907	0.169857117952969\\
0.883537414965986	0.173347188425046\\
0.885351473922902	0.175226773039832\\
0.887165532879819	0.175237447879347\\
0.888979591836735	0.173477776319731\\
0.890793650793651	0.170288887845994\\
0.892607709750567	0.166171198857586\\
0.894421768707483	0.161729781108579\\
0.896235827664399	0.157651094585598\\
0.898049886621315	0.154626973936411\\
0.899863945578231	0.153230363652016\\
0.901678004535147	0.153771407739496\\
0.903492063492063	0.156161418025053\\
0.90530612244898	0.159909930315045\\
0.908934240362812	0.168513964261325\\
0.910748299319728	0.171939806804027\\
0.912562358276644	0.174082198476975\\
0.91437641723356	0.174661987719541\\
0.916190476190476	0.173631005826114\\
0.918004535147392	0.171167332481994\\
0.919818594104308	0.167652022131989\\
0.923446712018141	0.159664337806533\\
0.925260770975057	0.156385998758902\\
0.927074829931973	0.154284324528605\\
0.928888888888889	0.153726370972104\\
0.930702947845805	0.154878863198903\\
0.932517006802721	0.157637014824393\\
0.934331065759637	0.161599800449031\\
0.936145124716553	0.166097224830011\\
0.937959183673469	0.170322635461871\\
0.939773242630386	0.173506531935909\\
0.941587301587302	0.175068836666146\\
0.943401360544218	0.174777502329446\\
0.945215419501134	0.172762984574192\\
0.94702947845805	0.169405994455628\\
0.952471655328798	0.157034301669718\\
0.954285714285714	0.154331120687692\\
0.95609977324263	0.153292754874881\\
0.957913832199547	0.154153535208749\\
0.959727891156463	0.156767820551887\\
0.961541950113379	0.160611374610142\\
0.963356009070295	0.164942266380589\\
0.965170068027211	0.169000446743012\\
0.966984126984127	0.172157963927387\\
0.968798185941043	0.173979954922929\\
0.970612244897959	0.174237591072599\\
0.972426303854875	0.172930112732498\\
0.974240362811791	0.170267609811029\\
0.976054421768707	0.166657916714861\\
0.97968253968254	0.158865481394478\\
0.981496598639456	0.155853405770166\\
0.983310657596372	0.154073079418047\\
0.985124716553288	0.153826331172102\\
0.986938775510204	0.155214455393648\\
0.98875283446712	0.158085548418688\\
0.990566893424036	0.162029071222117\\
0.992380952380952	0.166402500742937\\
0.994195011337868	0.170437433822122\\
0.996009070294785	0.173406683124445\\
0.997823129251701	0.174760373689215\\
0.999637188208617	0.174276658993459\\
};
\addlegendentry{ABC$_W^{1/2}$}

\end{axis}
\end{tikzpicture}%

%% file: images/HighFq/Angle_Quality.tex
%
%
\definecolor{mycolor1}{rgb}{1.00000,0.60000,0.00000}%
\pgfplotsset{scaled y ticks=false}
\begin{tikzpicture}[scale = 0.52, font=\huge]

\begin{axis}[%
width=4.521in,
height=3.548in,
at={(0.758in,0.499in)},
scale only axis,
xmin=0,
xmax=1,
xlabel style={font=\huge\color{white!15!black}},
xlabel={relative simulation time $t/T$}, xtick={0, 0.2, 0.4, 0.6, 0.8, 1},
y tick label style={
        /pgf/number format/.cd,
            fixed,
            fixed zerofill,
            precision=2,
        /tikz/.cd
            },
ymin=0,
ymax=0.26, ytick={0.02, 0.06, 0.1, 0.14, 0.18, 0.22, 0.26},
ylabel style={font=\huge\color{white!15!black}},
ylabel={}, ytick={0.02,  0.06, 0.1, 0.14, 0.18, 0.22, 0.26},
axis background/.style={fill=white},
xmajorgrids,
ymajorgrids,
legend style={at={(0.03,0.97)}, anchor=north west, legend cell align=left, align=left, draw=white!15!black}
]
\addplot [color=black, line width=2.0pt]
  table[row sep=crcr]{%
0.00190476190476185	4.44089209850063e-16\\
0.185124716553288	7.87497743746091e-05\\
0.199637188208617	0.000452806469949052\\
0.206893424036281	0.000671685805582634\\
0.21233560090703	0.00105498550331373\\
0.22140589569161	0.00195607731512226\\
0.241360544217687	0.00407957225768241\\
0.250430839002268	0.00495953059527376\\
0.26312925170068	0.00626642386235365\\
0.281269841269841	0.00772705319363454\\
0.290340136054422	0.00860791040280429\\
0.299410430839002	0.00919393744800645\\
0.306666666666667	0.00976743565657789\\
0.319365079365079	0.0109408949366643\\
0.335691609977324	0.0120011814626951\\
0.346575963718821	0.013019997821288\\
0.353832199546485	0.0134149560444845\\
0.36108843537415	0.0138695656371036\\
0.37015873015873	0.0147615461523386\\
0.375600907029478	0.0152068498374732\\
0.382857142857143	0.0155661068448437\\
0.390113378684807	0.0159642376045779\\
0.397369614512472	0.0166155409389598\\
0.413696145124717	0.01832310218299\\
0.417324263038549	0.0189250395353805\\
0.420952380952381	0.0197106059183139\\
0.431836734693878	0.0223853981629303\\
0.43546485260771	0.0230497461865009\\
0.449977324263039	0.0253693565948213\\
0.455419501133787	0.0263209608178689\\
0.459047619047619	0.0267852867972842\\
0.462675736961451	0.027057667808234\\
0.473560090702948	0.0276080463016606\\
0.486258503401361	0.0287201898567274\\
0.489886621315193	0.0287618351553998\\
0.500770975056689	0.028698499009668\\
0.506213151927438	0.0289616623709132\\
0.513469387755102	0.0293290729464726\\
0.517097505668934	0.0293447328524242\\
0.524353741496599	0.029078556601509\\
0.529795918367347	0.0289814169194035\\
0.535238095238095	0.0291431850541567\\
0.54249433106576	0.0293809362172165\\
0.546122448979592	0.029320507453276\\
0.558820861678004	0.0288251731719174\\
0.564263038548753	0.0289607932589595\\
0.569705215419501	0.0290852182309093\\
0.573333333333333	0.0289981847174621\\
0.578775510204082	0.0286026515957272\\
0.58421768707483	0.028241071426886\\
0.587845804988662	0.0281668276536107\\
0.600544217687075	0.0282527697572081\\
0.605986394557823	0.0279025353698881\\
0.611428571428571	0.0275322352846555\\
0.615056689342404	0.0274252544645054\\
0.620498866213152	0.0274788360278236\\
0.6259410430839	0.0275295263544959\\
0.629569160997732	0.0274018449158135\\
0.635011337868481	0.0269558729559227\\
0.640453514739229	0.0265286426680303\\
0.645895691609977	0.0263606907766635\\
0.654965986394558	0.0262752534276054\\
0.65859410430839	0.0260487064932966\\
0.671292517006803	0.0249833892570446\\
0.676734693877551	0.0248512351407034\\
0.682176870748299	0.0247321066688521\\
0.685804988662131	0.0245116820445083\\
0.693061224489796	0.0237456676733148\\
0.698503401360544	0.023271743824825\\
0.703945578231293	0.0230349671456855\\
0.711201814058957	0.0227786745743239\\
0.716643990929705	0.022303737052544\\
0.725714285714286	0.0213990665149092\\
0.731156462585034	0.0211143259305326\\
0.740226757369615	0.0207520060570139\\
0.745668934240363	0.020237064283822\\
0.752925170068027	0.0194904677049113\\
0.758367346938775	0.0191614666067103\\
0.769251700680272	0.0186821124141922\\
0.776507936507936	0.0180071174172132\\
0.781950113378685	0.0175555358115627\\
0.787392290249433	0.0173131258279358\\
0.796462585034014	0.0170421465690215\\
0.801904761904762	0.0166157791192142\\
0.809160997732426	0.0159847244681489\\
0.814603174603175	0.0157151719119398\\
0.827301587301587	0.0152369632741143\\
0.84	0.0143289899506248\\
0.847256235827664	0.0141156209627998\\
0.852698412698413	0.013950556798487\\
0.858140589569161	0.0135899971118492\\
0.865396825396825	0.0130655461602721\\
0.87265306122449	0.0128117870287782\\
0.88172335600907	0.0126073226542005\\
0.888979591836735	0.0121682462266415\\
0.896235827664399	0.0117939288686014\\
0.91437641723356	0.0111857503023679\\
0.921632653061224	0.0106997936307964\\
0.928888888888889	0.0104376658240635\\
0.939773242630386	0.0101357600301972\\
0.954285714285714	0.009355178080378\\
0.970612244897959	0.00875288125583085\\
0.977868480725624	0.00834979099769462\\
0.997823129251701	0.00753997530729567\\
0.999637188208617	0.00740145314105567\\
};
\addlegendentry{ABC$_W^{1/2}$ adaptive with exact $\theta$}

\addplot [color=mycolor1, line width=2.0pt]
  table[row sep=crcr]{%
0.00190476190476185	4.44089209850063e-16\\
0.185124716553288	7.87497743746091e-05\\
0.199637188208617	0.000452806469949052\\
0.206893424036281	0.000671685805582634\\
0.21233560090703	0.00105023099571266\\
0.223219954648526	0.00205683081791497\\
0.235918367346939	0.00325689151007535\\
0.252244897959184	0.00459607414534413\\
0.26312925170068	0.0056086622727789\\
0.28671201814059	0.00740901105390634\\
0.293968253968254	0.00789968693741838\\
0.304852607709751	0.00854206333893581\\
0.322993197278912	0.00995074675573016\\
0.332063492063492	0.0104247300861359\\
0.339319727891156	0.0110141271694447\\
0.346575963718821	0.0115789321337036\\
0.353832199546485	0.0119108533178272\\
0.36108843537415	0.0122906635855513\\
0.37015873015873	0.013071123368007\\
0.375600907029478	0.0134808234961211\\
0.382857142857143	0.0138111904676127\\
0.391927437641723	0.0142566839639685\\
0.399183673469388	0.0148456422355383\\
0.413696145124717	0.0161700226080294\\
0.417324263038549	0.0167123753006598\\
0.420952380952381	0.0174470844863349\\
0.433650793650794	0.02038035378734\\
0.439092970521542	0.0211953124292996\\
0.446349206349206	0.0219946506516526\\
0.455419501133787	0.0229255695475064\\
0.459047619047619	0.0231455915790678\\
0.462675736961451	0.0231766864789268\\
0.4681179138322	0.0229703244309565\\
0.473560090702948	0.0227422428710317\\
0.47718820861678	0.0227411940092282\\
0.484444444444444	0.0228149946019855\\
0.488072562358277	0.0225817952192507\\
0.500770975056689	0.021277319934207\\
0.506213151927438	0.0212117243002117\\
0.511655328798186	0.0211418706756615\\
0.515283446712018	0.0208704076804916\\
0.526167800453515	0.0198108457046864\\
0.529795918367347	0.0196121727459546\\
0.537052154195011	0.0195284097186231\\
0.54249433106576	0.0191651465699723\\
0.553378684807256	0.0181465268039762\\
0.558820861678004	0.0178859088356064\\
0.564263038548753	0.0176852945109385\\
0.567891156462585	0.017466782892779\\
0.573333333333333	0.0168594289343217\\
0.586031746031746	0.0155290883247046\\
0.589659863945578	0.0153722970249763\\
0.596916099773243	0.0147589549899472\\
0.604172335600907	0.0140707445192881\\
0.611428571428571	0.0134352157133868\\
0.615056689342404	0.0133041233085223\\
0.618684807256236	0.0131810431708944\\
0.6259410430839	0.0127514101421193\\
0.633197278911565	0.012151122739542\\
0.642267573696145	0.011561283614014\\
0.645895691609977	0.0115523982354112\\
0.64952380952381	0.0113582351370177\\
0.65859410430839	0.0105889946028558\\
0.669478458049887	0.00998588082418184\\
0.673106575963719	0.0100102867944711\\
0.678548752834467	0.0097624171099342\\
0.682176870748299	0.00945600191089446\\
0.687619047619048	0.00888355871457769\\
0.693061224489796	0.0086360740836996\\
0.698503401360544	0.00842706505436219\\
0.702131519274376	0.00855185730863861\\
0.705759637188209	0.00838836160686052\\
0.714829931972789	0.00765337625156093\\
0.725714285714286	0.00739620511222738\\
0.729342403628118	0.00761980720281019\\
0.73297052154195	0.00756590474106855\\
0.736598639455782	0.00750192538009753\\
0.743854875283447	0.00700408206918657\\
0.749297052154195	0.00716328225031382\\
0.752925170068027	0.00721199751759172\\
0.754739229024943	0.00733859994513619\\
0.758367346938775	0.00775520150566156\\
0.761995464852608	0.00784709962571806\\
0.772879818594104	0.00767557931790575\\
0.781950113378685	0.00796989834235173\\
0.785578231292517	0.00838408485324615\\
0.791020408163265	0.0086857638325939\\
0.794648526077097	0.00879856010139457\\
0.801904761904762	0.00875013372805489\\
0.805532879818594	0.00880916507469109\\
0.809160997732426	0.00881477327538482\\
0.812789115646259	0.00908470648608228\\
0.816417233560091	0.00928821939536462\\
0.820045351473923	0.0092592770173292\\
0.827301587301587	0.00905363223202704\\
0.84	0.00927249863047719\\
0.843628117913832	0.00942290283920344\\
0.850884353741497	0.00939792482312818\\
0.856326530612245	0.00912662350224136\\
0.861768707482993	0.00931088530991298\\
0.867210884353741	0.00940779668029801\\
0.870839002267574	0.0096820462809194\\
0.874467120181406	0.00976057844999245\\
0.88172335600907	0.00957616021676522\\
0.885351473922902	0.00951386741204874\\
0.888979591836735	0.0095491078773271\\
0.892607709750567	0.0094483455637836\\
0.894421768707483	0.00947487810685166\\
0.898049886621315	0.00970430605154493\\
0.90530612244898	0.00984483300717875\\
0.908934240362812	0.0096774016149217\\
0.912562358276644	0.00948754425839915\\
0.918004535147392	0.00963017284554024\\
0.921632653061224	0.00962577912997842\\
0.923446712018141	0.00970423940828669\\
0.927074829931973	0.0100181871186057\\
0.930702947845805	0.0102189236609443\\
0.934331065759637	0.0102473315537271\\
0.941587301587302	0.0101573587280829\\
0.945215419501134	0.0101691375295897\\
0.948843537414966	0.0100344243788549\\
0.950657596371882	0.0100561353023771\\
0.954285714285714	0.010327292410544\\
0.963356009070295	0.0106640326858271\\
0.968798185941043	0.0104348244295858\\
0.974240362811791	0.0105451133039857\\
0.977868480725624	0.0105369507373079\\
0.97968253968254	0.0106257969084275\\
0.985124716553288	0.0111301045289377\\
0.98875283446712	0.0113551708539015\\
0.994195011337868	0.0113539062413155\\
0.999637188208617	0.0112777311416224\\
};
\addlegendentry{ABC$_W^{1/2}$ adaptive with approx. $\theta$}

\end{axis}
\end{tikzpicture}%

%% file: images/Octant/Angle_Distribution0.tex
%
%
\pgfplotsset{scaled x ticks=false}
\begin{tikzpicture}[scale = 0.52, font=\huge]

\begin{axis}[%
width=4.602in,
height=3.82in,
at={(0.772in,0.516in)},
scale only axis,
xmin=0,
xmax=0.04, xtick={0.01, 0.02, 0.03, 0.04},
x tick label style={
        /pgf/number format/.cd,
            fixed,
            fixed zerofill,
            precision=2,
        /tikz/.cd
            },
xlabel style={font=\huge\color{white!15!black}},
xlabel={$x$},
ymin=0,
ymax=60,
ylabel style={font=\huge\color{white!15!black}},
ylabel={angle of incidence $[^{\circ}]$},
axis background/.style={fill=white},
xmajorgrids,
ymajorgrids,
legend style={at={(0.03,0.97)}, anchor=north west, legend cell align=left, align=left, draw=white!15!black}
]
\addplot [color=blue, line width=1.0pt, mark=*, mark options={solid, fill=blue, blue}]
  table[row sep=crcr]{%
0.000125000000000597	0.88909851930827\\
0.000375000000001791	0.525020091963317\\
0.000625000000002984	2.2671500532426\\
0.000875000000004178	0.185555158536101\\
0.00112500000000537	1.98070924331135\\
0.00137500000000657	1.77943738283395\\
0.00162500000000776	1.62664669343648\\
0.00187500000000895	1.74758489431558\\
0.00212500000001015	2.8127032594462\\
0.00237500000001134	2.84180451246091\\
0.00262500000000898	3.00227168747427\\
0.00287500000001017	4.93109168605846\\
0.00312500000001137	2.42031428735227\\
0.00337500000001256	4.90976957424554\\
0.0036250000000102	4.67884551985417\\
0.0038750000000114	4.30005760609274\\
0.00412500000000904	4.68669046589444\\
0.00437500000001023	4.8360632109704\\
0.00462500000000787	5.21506098926052\\
0.00487500000000907	6.57558086106654\\
0.00512500000001026	4.13615367371889\\
0.00537500000001145	6.57916053249682\\
0.00562500000001265	6.90480098683609\\
0.00587500000001384	9.51573278480944\\
0.00612500000001504	9.2176883789694\\
0.00637500000001623	11.8720362276249\\
0.00662500000001742	10.946575163601\\
0.00687500000001862	11.1198279695977\\
0.00712500000001981	11.7729804993258\\
0.007375000000021	13.2386652687831\\
0.00762500000001864	12.3732021703988\\
0.00787500000001984	13.087557887916\\
0.00812500000002103	15.2463315001214\\
0.00837500000002223	13.6780431549871\\
0.00862500000002342	13.1053712302532\\
0.00887500000002461	16.5266357850168\\
0.00912500000002581	15.5337676974775\\
0.009375000000027	15.5553731529181\\
0.00962500000002819	15.536234356522\\
0.00987500000002584	16.4134278913607\\
0.010125000000027	16.8001832996072\\
0.0103750000000282	17.2543397981913\\
0.0106250000000259	17.1931790726516\\
0.0108750000000271	17.755312393919\\
0.0111250000000283	17.4082710826144\\
0.0113750000000294	17.6353975351621\\
0.0116250000000306	17.8547475068482\\
0.0118750000000318	17.2887773051841\\
0.012125000000033	17.4164966176771\\
0.0123750000000342	18.0586039956156\\
0.0126250000000354	18.036850568279\\
0.0128750000000366	18.6659333897875\\
0.0131250000000378	18.6825760085489\\
0.0133750000000354	18.8435811866765\\
0.0136250000000366	18.9542047341877\\
0.0138750000000378	19.1724654555216\\
0.0141250000000355	19.3350125724355\\
0.0143750000000367	19.3918553916776\\
0.0146250000000379	19.6653668278314\\
0.0148750000000391	19.918079559741\\
0.0151250000000402	19.8153739094438\\
0.0153750000000414	20.1614353850956\\
0.0156250000000426	20.3514292558777\\
0.0158750000000403	20.2710416179522\\
0.0161250000000415	20.3870827938718\\
0.0163750000000427	20.2217094468971\\
0.0166250000000439	20.4473750632892\\
0.016875000000045	20.6375986220429\\
0.0171250000000462	20.0932851436921\\
0.0173750000000474	18.9566646499738\\
0.0176250000000486	20.0644624364835\\
0.0178750000000463	19.3479590191027\\
0.0181250000000475	17.983766393279\\
0.0183750000000487	0\\
0.0186250000000499	0\\
0.0188750000000475	0\\
0.0191250000000487	0\\
0.0193750000000499	0\\
0.0196250000000511	0\\
0.0198750000000523	0\\
0.0201250000000535	0\\
0.0203750000000511	0\\
0.0206250000000523	0\\
0.0208750000000499	0\\
0.0211250000000476	0\\
0.0213750000000488	0\\
0.02162500000005	0\\
0.0218750000000476	0\\
0.0221250000000488	0\\
0.0223750000000464	0\\
0.0226250000000476	0\\
0.0228750000000488	0\\
0.0231250000000465	0\\
0.0233750000000441	0\\
0.0236250000000453	0\\
0.0238750000000429	0\\
0.0241250000000406	0\\
0.0243750000000418	0\\
0.0246250000000394	0\\
0.0248750000000406	0\\
0.0251250000000418	0\\
0.0253750000000394	0\\
0.0256250000000406	0\\
0.0258750000000383	0\\
0.0261250000000395	0\\
0.0263750000000371	0\\
0.0266250000000383	0\\
0.026875000000036	0\\
0.0271250000000336	0\\
0.0273750000000348	0\\
0.0276250000000324	0\\
0.0278750000000301	0\\
0.0281250000000313	0\\
0.0283750000000289	0\\
0.0286250000000265	0\\
0.0288750000000277	0\\
0.0291250000000254	0\\
0.029375000000023	0\\
0.0296250000000242	0\\
0.0298750000000219	0\\
0.0301250000000195	0\\
0.0303750000000207	0\\
0.0306250000000219	0\\
0.0308750000000195	0\\
0.0311250000000207	0\\
0.0313750000000219	0\\
0.0316250000000196	0\\
0.0318750000000207	0\\
0.0321250000000219	0\\
0.0323750000000196	0\\
0.0326250000000208	0\\
0.032875000000022	0\\
0.0331250000000196	0\\
0.0333750000000208	0\\
0.0336250000000184	0\\
0.0338750000000196	0\\
0.0341250000000173	0\\
0.0343750000000149	0\\
0.0346250000000161	0\\
0.0348750000000138	0\\
0.0351250000000114	0\\
0.0353750000000126	0\\
0.0356250000000102	0\\
0.0358750000000079	0\\
0.0361250000000091	0\\
0.0363750000000067	0\\
0.0366250000000043	0\\
0.0368750000000055	0\\
0.0371250000000032	0\\
0.0373750000000008	0\\
0.037625000000002	0\\
0.0378749999999997	0\\
0.0381250000000009	0\\
0.0383749999999985	0\\
0.0386249999999961	0\\
0.0388749999999973	0\\
0.039124999999995	0\\
0.0393749999999962	0\\
0.0396249999999974	0\\
0.0398749999999986	0\\
};
\addlegendentry{Computed angle}

\addplot [color=black, line width=2.0pt]
  table[row sep=crcr]{%
0	0\\
0.00377358490565882	5.38931175997342\\
0.00578616352201067	8.23097510119033\\
0.00754716981131764	10.6849124000027\\
0.00905660377358686	12.7575321608767\\
0.010566037735849	14.7967622450573\\
0.0118238993710662	16.4675467930027\\
0.0130817610062905	18.1100234999005\\
0.0143396226415078	19.722277764447\\
0.0155974842767321	21.3026557774829\\
0.0168553459119494	22.8497620716492\\
0.0178616352201288	24.0627147257253\\
0.0188679245283012	25.2531633945739\\
0.0198742138364807	26.4207296195553\\
0.0208805031446531	27.5651253052007\\
0.0218867924528325	28.6861475737394\\
0.0228930817610049	29.783673287977\\
0.0238993710691844	30.8576533579673\\
0.0249056603773568	31.9081069356531\\
0.0259119496855362	32.9351155894746\\
0.0269182389937086	33.9388175384082\\
0.0279245283018881	34.9194020124577\\
0.0289308176100604	35.8771037946215\\
0.0299371069182399	36.8121979880955\\
0.0309433962264123	37.7249950421107\\
0.0319496855345918	38.6158360604875\\
0.0329559748427641	39.485088408755\\
0.0339622641509436	40.333141628561\\
0.034968553459116	41.1604036620465\\
0.0359748427672955	41.967297383823\\
0.0369811320754749	42.7542574341049\\
0.0379874213836473	43.5217273433111\\
0.0389937106918268	44.2701569359759\\
0.0399999999999991	45\\
};
\addlegendentry{Analytical angle}

\end{axis}
\end{tikzpicture}%

%% file: images/Octant/Angle_Distribution.tex
%
%
\pgfplotsset{scaled x ticks=false}
\begin{tikzpicture}[scale = 0.52, font=\huge]

\begin{axis}[%
width=4.602in,
height=3.82in,
at={(0.772in,0.516in)},
scale only axis,
xmin=0,
xmax=0.04, xtick={0.01, 0.02, 0.03, 0.04},
x tick label style={
        /pgf/number format/.cd,
            fixed,
            fixed zerofill,
            precision=2,
        /tikz/.cd
            },
xlabel style={font=\huge\color{white!15!black}},
xlabel={$x$},
ymin=0,
ymax=60,
ylabel style={font=\huge\color{white!15!black}},
ylabel={},
axis background/.style={fill=white},
xmajorgrids,
ymajorgrids,
legend style={at={(0.03,0.97)}, anchor=north west, legend cell align=left, align=left, draw=white!15!black}
]
\addplot [color=blue, line width=1.0pt, mark=*, mark options={solid, fill=blue, blue}]
  table[row sep=crcr]{%
0.000124999999997044	0.547390950284239\\
0.000374999999998238	0.608060672506838\\
0.000624999999999432	0.527037574273862\\
0.000875000000000625	2.00995738047405\\
0.00112500000000182	1.62141559547865\\
0.00137500000000301	1.89210977808088\\
0.00162500000000421	2.17462841064386\\
0.0018750000000054	2.46654620853815\\
0.00212500000000659	2.57726833720031\\
0.00237500000000779	3.22949825476024\\
0.00262500000000898	3.90159837693051\\
0.00287500000001017	3.99738249765474\\
0.00312500000001137	5.25878774632045\\
0.00337500000001256	4.98684510531317\\
0.00362500000001376	5.37423517893751\\
0.00387500000000784	5.85365367955893\\
0.00412500000000904	6.26400547145516\\
0.00437500000001023	6.55661630515752\\
0.00462500000001143	6.77285300187128\\
0.00487500000001262	6.96527303684161\\
0.00512500000000671	7.84463871373714\\
0.0053750000000079	7.91860366087815\\
0.00562500000000909	8.27051163499354\\
0.00587500000001029	8.70180983107591\\
0.00612500000001148	9.17897342908671\\
0.00637500000001268	9.76173645220243\\
0.00662500000001387	9.66322203904468\\
0.00687500000001506	9.72288270201477\\
0.00712500000001626	9.89863992895323\\
0.00737500000001745	10.1861519319168\\
0.00762500000001864	10.6324962110614\\
0.00787500000001984	11.067540589457\\
0.00812500000002103	11.2546808170297\\
0.00837500000002223	11.5766958776609\\
0.00862500000002342	11.9590548992221\\
0.00887500000002461	11.9626990474211\\
0.00912500000002581	12.3804601945633\\
0.009375000000027	12.7297946748917\\
0.00962500000002819	12.9702759776872\\
0.00987500000002939	12.8641597415759\\
0.0101250000000235	12.6279060168691\\
0.0103750000000247	12.6163432302943\\
0.0106250000000259	13.3465136390436\\
0.0108750000000271	13.3771467579403\\
0.0111250000000283	14.1753007175487\\
0.0113750000000294	14.2054682301187\\
0.0116250000000306	13.7997514865411\\
0.0118750000000318	14.9752093233218\\
0.012125000000033	18.9288897474458\\
0.0123750000000342	20.4999312765984\\
0.0126250000000354	19.2225301129393\\
0.0128750000000366	21.6557997812151\\
0.0131250000000378	20.8756268517485\\
0.013375000000039	22.1414556684711\\
0.0136250000000402	21.7028647255823\\
0.0138750000000343	21.2033675046657\\
0.0141250000000355	24.41596825123\\
0.0143750000000367	24.0084612204426\\
0.0146250000000379	23.499223418459\\
0.0148750000000391	23.7694823768145\\
0.0151250000000402	24.3022525025947\\
0.0153750000000414	22.8677454616249\\
0.0156250000000426	23.0205763468612\\
0.0158750000000438	23.5467731773921\\
0.016125000000045	23.6982074644841\\
0.0163750000000462	23.3705503239274\\
0.0166250000000474	22.9696338557334\\
0.0168750000000486	23.0649941624361\\
0.0171250000000498	23.5041241460465\\
0.017375000000051	23.7544903400721\\
0.0176250000000451	23.718448405125\\
0.0178750000000463	23.7628193561402\\
0.0181250000000475	23.9033738094237\\
0.0183750000000487	24.175422534135\\
0.0186250000000499	24.8022297702788\\
0.018875000000051	25.4829398952802\\
0.0191250000000522	25.7372738434461\\
0.0193750000000534	24.3045879332573\\
0.0196250000000546	24.0830903880126\\
0.0198750000000558	23.5693796937997\\
0.0201250000000499	23.109839625797\\
0.0203750000000511	23.8941793056132\\
0.0206250000000523	29.2939451112026\\
0.0208750000000464	29.0283381271014\\
0.0211250000000476	30.5162749980666\\
0.0213750000000488	32.4488904421263\\
0.02162500000005	32.7310291439001\\
0.0218750000000512	31.7336339273849\\
0.0221250000000452	31.4551870161921\\
0.0223750000000464	31.7162122181906\\
0.0226250000000476	32.2557526403578\\
0.0228750000000488	33.1524649738072\\
0.0231250000000429	34.8520587598634\\
0.0233750000000441	34.0441038596464\\
0.0236250000000453	33.7140195868425\\
0.0238750000000394	34.0011026194599\\
0.0241250000000406	33.7494873574854\\
0.0243750000000418	32.9685981642712\\
0.024625000000043	32.8573030524069\\
0.0248750000000442	32.6370317052436\\
0.0251250000000383	33.1545550064437\\
0.0253750000000394	32.9359482063002\\
0.0256250000000406	33.1137617240714\\
0.0258750000000418	33.3776114213929\\
0.0261250000000359	33.7952211631504\\
0.0263750000000371	33.5689613528095\\
0.0266250000000383	33.0913209824581\\
0.0268750000000324	32.7865499849364\\
0.0271250000000336	33.0126244876507\\
0.0273750000000348	33.0200579507266\\
0.0276250000000289	31.2987959473801\\
0.0278750000000301	30.5025123293601\\
0.0281250000000313	42.7124615461812\\
0.0283750000000254	41.8967087374019\\
0.0286250000000265	43.3140213311494\\
0.0288750000000277	43.6268330126377\\
0.0291250000000218	42.1645786416245\\
0.029375000000023	41.5336882404106\\
0.0296250000000242	42.4446079848606\\
0.0298750000000254	39.8681055075782\\
0.0301250000000195	38.8651779414853\\
0.0303750000000207	38.882483448772\\
0.0306250000000219	38.5141314162219\\
0.0308750000000231	38.0494934738915\\
0.0311250000000243	37.5875034514907\\
0.0313750000000184	37.3283027179913\\
0.0316250000000196	37.0706428734403\\
0.0318750000000207	37.0206149700009\\
0.0321250000000219	37.0380147065354\\
0.0323750000000231	36.8928039157747\\
0.0326250000000172	37.2860310993098\\
0.0328750000000184	36.7820667461707\\
0.0331250000000196	43.8281398842317\\
0.0333750000000208	44.2169518028267\\
0.033625000000022	41.5751406504729\\
0.0338750000000161	43.8469166649566\\
0.0341250000000173	44.431692252561\\
0.0343750000000185	44.1695166461101\\
0.0346250000000126	43.2253455610704\\
0.0348750000000138	42.9407130544744\\
0.0351250000000149	43.2227921763993\\
0.035375000000009	43.4803517742301\\
0.0356250000000102	44.6653710432963\\
0.0358750000000114	45.5957490455819\\
0.0361250000000055	46.1139975897379\\
0.0363750000000067	44.2665034411966\\
0.0366250000000079	44.101251211711\\
0.036875000000002	43.0340460376295\\
0.0371250000000032	43.0189862096751\\
0.0373750000000044	43.2976741622823\\
0.0376249999999985	43.2470526300871\\
0.0378749999999997	43.1683004652695\\
0.0381250000000009	43.2351189854969\\
0.0383749999999949	44.1622373037109\\
0.0386249999999961	43.2701602303911\\
0.0388749999999973	44.027193131617\\
0.0391249999999985	44.9005700591166\\
0.0393749999999997	50.4792600085758\\
0.0396250000000009	50.0044088365229\\
0.0398750000000021	48.2753645949294\\
};
\addlegendentry{Computed angle}

\addplot [color=black, line width=2.0pt]
  table[row sep=crcr]{%
0	0\\
0.00377358490565882	5.38931175997342\\
0.00578616352201067	8.23097510119033\\
0.00754716981131764	10.6849124000027\\
0.00905660377358686	12.7575321608767\\
0.010566037735849	14.7967622450573\\
0.0118238993710662	16.4675467930027\\
0.0130817610062905	18.1100234999005\\
0.0143396226415078	19.722277764447\\
0.0155974842767321	21.3026557774829\\
0.0168553459119494	22.8497620716492\\
0.0178616352201288	24.0627147257253\\
0.0188679245283012	25.2531633945739\\
0.0198742138364807	26.4207296195553\\
0.0208805031446531	27.5651253052007\\
0.0218867924528325	28.6861475737394\\
0.0228930817610049	29.783673287977\\
0.0238993710691844	30.8576533579673\\
0.0249056603773568	31.9081069356531\\
0.0259119496855362	32.9351155894746\\
0.0269182389937086	33.9388175384082\\
0.0279245283018881	34.9194020124577\\
0.0289308176100604	35.8771037946215\\
0.0299371069182399	36.8121979880955\\
0.0309433962264123	37.7249950421107\\
0.0319496855345918	38.6158360604875\\
0.0329559748427641	39.485088408755\\
0.0339622641509436	40.333141628561\\
0.034968553459116	41.1604036620465\\
0.0359748427672955	41.967297383823\\
0.0369811320754749	42.7542574341049\\
0.0379874213836473	43.5217273433111\\
0.0389937106918268	44.2701569359759\\
0.0399999999999991	45\\
};
\addlegendentry{Analytical angle}

\end{axis}
\end{tikzpicture}%

%% file: images/Octant/Cost.tex
%
%
\definecolor{mycolor1}{rgb}{1.00000,0.60000,0.00000}%
\pgfplotsset{scaled y ticks=false}
\begin{tikzpicture}[scale = 0.52, font=\huge]

\begin{axis}[%
width=4.602in,
height=3.82in,
at={(0.772in,0.516in)},
scale only axis,
xmin=0,
xmax=1,
xlabel style={font=\huge\color{white!15!black}},
xlabel={relative simulation time $t/T$}, xtick={0, 0.2, 0.4, 0.6, 0.8, 1},
y tick label style={
        /pgf/number format/.cd,
            fixed,
            fixed zerofill,
            precision=2,
        /tikz/.cd
            },
ymin=0,
ymax=0.09, ytick={0.01, 0.03, 0.05, 0.07, 0.09},
ylabel style={at={(-0.15,0.5)}, font=\huge\color{white!15!black}},
ylabel={relative $L^{2}(\Omega)$ error},
axis background/.style={fill=white},
xmajorgrids,
ymajorgrids,
legend style={at={(0.03,0.97)}, anchor=north west, legend cell align=left, align=left, draw=white!15!black}
]
\addplot [color=mycolor1, line width=2.0pt]
  table[row sep=crcr]{%
0.00252100840336134	5.77315972805081e-15\\
0.194597839135654	3.14823947963339e-05\\
0.206602641056423	8.94620307470007e-05\\
0.218607442977191	0.000212552694798029\\
0.228211284513806	0.000372295126119915\\
0.235414165666266	0.000554837995334623\\
0.242617046818728	0.000807008368315598\\
0.252220888355342	0.00122637854369245\\
0.259423769507803	0.00158553163607444\\
0.269027611044418	0.00214574994849415\\
0.276230492196879	0.00254249869665835\\
0.28343337334934	0.00285433440309946\\
0.295438175270108	0.00334826668370702\\
0.305042016806723	0.00374564551212209\\
0.317046818727491	0.004186529427887\\
0.321848739495798	0.00447237626812425\\
0.329051620648259	0.00497342758281083\\
0.333853541416567	0.00520303004653144\\
0.33625450180072	0.00529643113509115\\
0.338655462184874	0.00543307142997806\\
0.341056422569028	0.00566222830073837\\
0.343457382953181	0.00599467040553114\\
0.348259303721489	0.00674624875891272\\
0.350660264105642	0.00697606728403932\\
0.353061224489796	0.00698723334642248\\
0.35546218487395	0.00674576635775104\\
0.360264105642257	0.00584659005496613\\
0.362665066026411	0.00563453227268307\\
0.365066026410564	0.00584064781414662\\
0.369867947178872	0.00677061834061143\\
0.372268907563025	0.0071314323456152\\
0.374669867947179	0.00740292414034371\\
0.377070828331333	0.00754090896065185\\
0.379471788715486	0.00752785816920398\\
0.38187274909964	0.007443598673052\\
0.384273709483794	0.00741429136016836\\
0.386674669867947	0.00755981069268841\\
0.389075630252101	0.00785016730191712\\
0.396278511404562	0.00896815745820945\\
0.398679471788716	0.00917601782064836\\
0.401080432172869	0.00925682101793701\\
0.40828331332533	0.00923739217194908\\
0.413085234093638	0.00919751578570949\\
0.417887154861945	0.0092970674815791\\
0.420288115246098	0.00933673900682919\\
0.427490996398559	0.00931113848071363\\
0.43469387755102	0.00905621243795696\\
0.437094837935174	0.00905210780014254\\
0.439495798319328	0.00913638477574807\\
0.441896758703481	0.00933919869794753\\
0.444297719087635	0.00959390730989051\\
0.446698679471789	0.0098059435451584\\
0.449099639855942	0.00994304243777433\\
0.451500600240096	0.0100395245709846\\
0.45390156062425	0.0100927568132879\\
0.456302521008403	0.0100314393185058\\
0.461104441776711	0.00978279679626559\\
0.465906362545018	0.00975662761096419\\
0.468307322929172	0.00980647474851315\\
0.470708283313325	0.00995811924982792\\
0.477911164465786	0.010670280397946\\
0.48031212484994	0.0107758893731695\\
0.482713085234094	0.0108108843438482\\
0.485114045618247	0.0107988133616702\\
0.487515006002401	0.0107055622723574\\
0.489915966386555	0.0105774972957458\\
0.492316926770708	0.0105056980418479\\
0.494717887154862	0.0105525123001766\\
0.497118847539016	0.0106454851288307\\
0.499519807923169	0.0108232317223603\\
0.501920768307323	0.011094942024219\\
0.504321728691477	0.0114065331614813\\
0.50672268907563	0.0116823430546238\\
0.509123649459784	0.0118909101608744\\
0.511524609843938	0.0120331777194148\\
0.513925570228091	0.012105197062952\\
0.516326530612245	0.0120757609483102\\
0.521128451380552	0.0119241509422209\\
0.52593037214886	0.0119582347554767\\
0.528331332533013	0.0120689081379278\\
0.530732292917167	0.012254106511224\\
0.535534213685474	0.0127246468237032\\
0.537935174069628	0.0128954214289567\\
0.540336134453781	0.0129955779195136\\
0.542737094837935	0.0130624893172692\\
0.545138055222089	0.0130867590796768\\
0.549939975990396	0.0130234343859394\\
0.55234093637455	0.0130764554986893\\
0.557142857142857	0.0133730969959593\\
0.559543817527011	0.0136127376132308\\
0.564345738295318	0.0142085098313864\\
0.566746698679472	0.0144470403036366\\
0.569147659063625	0.0146199619245765\\
0.571548619447779	0.0147331158395537\\
0.573949579831933	0.0147695892887868\\
0.581152460984394	0.0146563864246607\\
0.583553421368547	0.0147131601334939\\
0.585954381752701	0.0148148771550035\\
0.588355342136855	0.0149957718116143\\
0.595558223289316	0.0157422298746263\\
0.597959183673469	0.0159162324344536\\
0.602761104441777	0.0161144337027267\\
0.60516206482593	0.0161448915868612\\
0.609963985594238	0.0161324285359581\\
0.612364945978391	0.0162140539292045\\
0.614765906362545	0.0163364519032693\\
0.617166866746699	0.0164959588966832\\
0.619567827130852	0.0167254371407765\\
0.62436974789916	0.0172678835455662\\
0.626770708283313	0.0174890964713755\\
0.629171668667467	0.0176572699713298\\
0.631572629051621	0.0177646800349962\\
0.633973589435774	0.0177782269960242\\
0.641176470588235	0.0176171392318841\\
0.643577430972389	0.017647925035809\\
0.645978391356543	0.0177513162475829\\
0.648379351740696	0.0179508201549037\\
0.65078031212485	0.0182313513750961\\
0.655582232893157	0.0188466467672622\\
0.657983193277311	0.0190758465859034\\
0.660384153661465	0.0192363519560595\\
0.662785114045618	0.019326400649649\\
0.665186074429772	0.0193245427626175\\
0.669987995198079	0.0192374166297219\\
0.672388955582233	0.0192846189193632\\
0.674789915966387	0.0193855546139088\\
0.67719087635054	0.0195768820465907\\
0.679591836734694	0.0198704705176653\\
0.684393757503001	0.0205349015037684\\
0.686794717887155	0.0207915002963173\\
0.689195678271309	0.0209608972445607\\
0.691596638655462	0.0210304163647336\\
0.693997599039616	0.0209870820564634\\
0.698799519807923	0.0207907777233525\\
0.701200480192077	0.0207736222099486\\
0.703601440576231	0.0208214523368281\\
0.706002400960384	0.0209722895963842\\
0.708403361344538	0.0212209637036114\\
0.713205282112845	0.0218453450834883\\
0.715606242496999	0.0220938745084704\\
0.718007202881152	0.0222532883983554\\
0.720408163265306	0.0223305555066581\\
0.72280912364946	0.0223260301785465\\
0.727611044417767	0.0221671448529761\\
0.730012004801921	0.0221537600678245\\
0.732412965186074	0.0222291342832232\\
0.734813925570228	0.0223763030373926\\
0.737214885954382	0.0226223902367525\\
0.742016806722689	0.0232598645466202\\
0.744417767106843	0.0235195822101256\\
0.746818727490996	0.0236934275846393\\
0.74921968787515	0.0237670880494139\\
0.751620648259304	0.0237309529707749\\
0.758823529411765	0.0233691908061948\\
0.761224489795918	0.0233703093029937\\
0.763625450180072	0.023453442859412\\
0.766026410564226	0.023647707678477\\
0.773229291716687	0.0244674418676395\\
0.77563025210084	0.0246315146693066\\
0.778031212484994	0.0246963721897082\\
0.780432172869148	0.0246799619216111\\
0.787635054021609	0.0243903989900028\\
0.790036014405762	0.0244101035006615\\
0.792436974789916	0.02450857809429\\
0.79483793517407	0.0246821702049528\\
0.802040816326531	0.0254703387154894\\
0.804441776710684	0.0256423243380908\\
0.806842737094838	0.0257218344480082\\
0.809243697478992	0.0257032174477304\\
0.811644657863145	0.0255834895832822\\
0.816446578631453	0.0252655304376866\\
0.818847539015606	0.0252129797090878\\
0.82124849939976	0.0252441807153246\\
0.823649459783914	0.0253671122424023\\
0.826050420168067	0.0255870605999799\\
0.830852340936375	0.0261048032294755\\
0.833253301320528	0.0262763043836446\\
0.835654261704682	0.0263447602124233\\
0.838055222088836	0.0263137139564661\\
0.840456182472989	0.0262051664315301\\
0.845258103241296	0.025863513461365\\
0.84765906362545	0.0257616008191636\\
0.850060024009604	0.0257550698712504\\
0.852460984393758	0.0258192089162795\\
0.854861944777911	0.0259676051919591\\
0.859663865546219	0.0263949450826104\\
0.862064825930372	0.0265495491435934\\
0.864465786314526	0.0266235580562902\\
0.866866746698679	0.0266037402298488\\
0.869267707082833	0.0264905900037971\\
0.874069627851141	0.0260862253698894\\
0.876470588235294	0.0259546282455797\\
0.878871548619448	0.0259301840683115\\
0.881272509003601	0.0259866837306938\\
0.883673469387755	0.0261446191241241\\
0.890876350540216	0.0268227717566746\\
0.89327731092437	0.0269226978613007\\
0.895678271308523	0.0269168559975509\\
0.898079231692677	0.0268231047295131\\
0.905282112845138	0.0263538001518727\\
0.907683073229292	0.0263422771642128\\
0.910084033613445	0.0264395878704944\\
0.912484993997599	0.0266150743490789\\
0.91968787515006	0.02738902922347\\
0.922088835534214	0.0275409531982092\\
0.924489795918367	0.0275889055726199\\
0.926890756302521	0.0275280426973924\\
0.929291716686675	0.0273705929241383\\
0.931692677070828	0.0271537964593881\\
0.934093637454982	0.0269730609775732\\
0.936494597839136	0.026893611947004\\
0.938895558223289	0.0269179234572848\\
0.941296518607443	0.0270215784300143\\
0.943697478991597	0.0272201192837527\\
0.948499399759904	0.0276884825070347\\
0.950900360144058	0.0278346068657837\\
0.953301320528211	0.0278780389690692\\
0.955702280912365	0.0278240059210492\\
0.958103241296519	0.0277069297170711\\
0.962905162064826	0.027413038832046\\
0.96530612244898	0.0273756473536054\\
0.967707082833133	0.0274646112648242\\
0.970108043217287	0.0276476667596927\\
0.972509003601441	0.0279017753068095\\
0.977310924369748	0.0285110582565131\\
0.979711884753902	0.0287456156118444\\
0.982112845138055	0.0288833167776438\\
0.984513805522209	0.0289153132823908\\
0.986914765906363	0.0288507770483376\\
0.99171668667467	0.0285523484526701\\
0.994117647058824	0.0284680910303587\\
0.996518607442977	0.0284925026034782\\
0.998919567827131	0.0285998431459978\\
};
\addlegendentry{ABC$_W^{1/2}$ adaptive}

\addplot [color=mycolor1, dashed, line width=2.0pt]
  table[row sep=crcr]{%
0.00252100840336134	5.77315972805081e-15\\
0.194597839135654	3.14823947963339e-05\\
0.206602641056423	8.94620307470007e-05\\
0.218607442977191	0.000212125693292986\\
0.225810324129652	0.000342029472075178\\
0.233013205282113	0.000528904210459369\\
0.240216086434574	0.000781068244333327\\
0.247418967587035	0.00110550111146768\\
0.254621848739496	0.0015040441026023\\
0.259423769507803	0.00182396586007727\\
0.26422569027611	0.00219980756048521\\
0.269027611044418	0.00263339661258732\\
0.276230492196879	0.00336134808436073\\
0.28343337334934	0.00416079751857013\\
0.288235294117647	0.00477407034936006\\
0.293037214885954	0.00549854958233897\\
0.297839135654262	0.00633891681731769\\
0.317046818727491	0.00995304184630763\\
0.319447779111645	0.0104725938273745\\
0.321848739495798	0.0110388224163716\\
0.324249699879952	0.0116511394559741\\
0.33625450180072	0.0148888844797103\\
0.341056422569028	0.0161187546426007\\
0.343457382953181	0.0167821500975518\\
0.345858343337335	0.0175005930824282\\
0.353061224489796	0.0198074583024636\\
0.35546218487395	0.0204929394906229\\
0.357863145258103	0.0210836395382255\\
0.365066026410564	0.0226309631488208\\
0.367466986794718	0.0233722985880009\\
0.369867947178872	0.0243738175520533\\
0.372268907563025	0.0256362928884473\\
0.377070828331333	0.028526147310451\\
0.379471788715486	0.029799267165376\\
0.38187274909964	0.0306956285815031\\
0.384273709483794	0.0310618552467207\\
0.386674669867947	0.0308423027932845\\
0.389075630252101	0.0301435894684511\\
0.391476590636255	0.0292891748414481\\
0.393877551020408	0.0288032668204833\\
0.396278511404562	0.0292411871009414\\
0.398679471788716	0.0308833735245597\\
0.401080432172869	0.0335446438738335\\
0.403481392557023	0.036685662958283\\
0.405882352941176	0.0396586169651763\\
0.40828331332533	0.0418825029407001\\
0.410684273709484	0.0429504374897572\\
0.413085234093638	0.0426873270679383\\
0.415486194477791	0.0411647501528559\\
0.417887154861945	0.038719119275943\\
0.420288115246098	0.0359961398777774\\
0.422689075630252	0.0339471436276587\\
0.425090036014406	0.0335861236697369\\
0.427490996398559	0.0354156924710635\\
0.429891956782713	0.0390502785430753\\
0.432292917166867	0.0435397911832288\\
0.43469387755102	0.0478871846774734\\
0.437094837935174	0.0512836396122561\\
0.439495798319328	0.0531643456633425\\
0.441896758703481	0.0532258848364477\\
0.444297719087635	0.0514306286533448\\
0.446698679471789	0.0480243152688772\\
0.451500600240096	0.0392477174558716\\
0.45390156062425	0.0364435570770937\\
0.456302521008403	0.036535189103293\\
0.458703481392557	0.0396388366926206\\
0.461104441776711	0.0445466231641067\\
0.463505402160864	0.0497371184338955\\
0.465906362545018	0.0540638849259006\\
0.468307322929172	0.0568489253205823\\
0.470708283313325	0.0577708978915044\\
0.473109243697479	0.056768037435022\\
0.475510204081633	0.0540187160845336\\
0.477911164465786	0.0499883560409897\\
0.48031212484994	0.0455117619416339\\
0.482713085234094	0.0418186623143448\\
0.485114045618247	0.0402642587789042\\
0.487515006002401	0.0415800109895972\\
0.489915966386555	0.0452887501366551\\
0.492316926770708	0.0501179325440043\\
0.494717887154862	0.0547798373107672\\
0.497118847539016	0.0583394709513559\\
0.499519807923169	0.0602502126672594\\
0.501920768307323	0.0602840773864691\\
0.504321728691477	0.0584687278365447\\
0.50672268907563	0.0550915276859949\\
0.509123649459784	0.0507764116431116\\
0.511524609843938	0.0465664181334328\\
0.513925570228091	0.0438260936728067\\
0.516326530612245	0.0436926973088551\\
0.518727490996399	0.0462712364983681\\
0.521128451380552	0.0505488850715852\\
0.523529411764706	0.0551518815501817\\
0.52593037214886	0.0589752329497963\\
0.528331332533013	0.0613415018332403\\
0.530732292917167	0.061941383831145\\
0.533133253301321	0.0607566747882062\\
0.535534213685474	0.0580337668031501\\
0.537935174069628	0.0543066709724743\\
0.540336134453781	0.0504250489997004\\
0.542737094837935	0.0474855266131085\\
0.545138055222089	0.0465092476175218\\
0.547539015606243	0.0478967900481762\\
0.549939975990396	0.0511561036816469\\
0.55234093637455	0.0552448767435896\\
0.554741896758703	0.0590971254002626\\
0.557142857142857	0.0619138323584446\\
0.559543817527011	0.0632201146370577\\
0.561944777911164	0.0628251707242685\\
0.564345738295318	0.0607949776306158\\
0.566746698679472	0.0574787786104819\\
0.569147659063625	0.0535781981624365\\
0.571548619447779	0.0501568996256344\\
0.573949579831933	0.0484101401735132\\
0.576350540216086	0.0490821604206619\\
0.57875150060024	0.0519471039856463\\
0.583553421368547	0.0600028527669909\\
0.585954381752701	0.0630956943946657\\
0.588355342136855	0.0647230531172168\\
0.590756302521008	0.0646712501403151\\
0.593157262905162	0.0630009999561626\\
0.595558223289316	0.0600365129141933\\
0.597959183673469	0.0563835109802329\\
0.600360144057623	0.0529241123772483\\
0.602761104441777	0.0506754956507802\\
0.60516206482593	0.0504203488889471\\
0.607563025210084	0.0522737425042199\\
0.609963985594238	0.0556155950069216\\
0.612364945978391	0.0594599161751189\\
0.614765906362545	0.0628681224009094\\
0.617166866746699	0.0651594052130929\\
0.619567827130852	0.0659475649528379\\
0.621968787515006	0.0651116236472187\\
0.62436974789916	0.0627880339015499\\
0.626770708283313	0.0594138626950038\\
0.629171668667467	0.0557771987870188\\
0.631572629051621	0.0529550198664593\\
0.633973589435774	0.0519801509765565\\
0.636374549819928	0.0532837442578455\\
0.638775510204082	0.0564022321943125\\
0.641176470588235	0.0602957875193676\\
0.643577430972389	0.0638974940004137\\
0.645978391356543	0.0664325904755696\\
0.648379351740696	0.067472814457499\\
0.65078031212485	0.0668944291510883\\
0.653181272509004	0.064844723834286\\
0.655582232893157	0.0617422291292438\\
0.657983193277311	0.058280721715649\\
0.660384153661465	0.0553675899533426\\
0.662785114045618	0.0539010337842946\\
0.665186074429772	0.0543873543418418\\
0.667587034813926	0.05666076344132\\
0.669987995198079	0.0599899061711356\\
0.672388955582233	0.063445323516321\\
0.674789915966387	0.0662166474602192\\
0.67719087635054	0.0677597180168096\\
0.679591836734694	0.0678187528836734\\
0.681992797118848	0.0664068540402758\\
0.684393757503001	0.0638088336788174\\
0.686794717887155	0.0606117313682679\\
0.689195678271309	0.0576893060457622\\
0.691596638655462	0.0560153607649914\\
0.693997599039616	0.0562498769586207\\
0.696398559423769	0.0583508049694031\\
0.698799519807923	0.0615996306635566\\
0.701200480192077	0.0650070808656856\\
0.703601440576231	0.0677120416564045\\
0.706002400960384	0.0691639371441699\\
0.708403361344538	0.0691383907280948\\
0.710804321728691	0.0677030142459908\\
0.713205282112845	0.0651950620647519\\
0.715606242496999	0.062201258990776\\
0.718007202881152	0.0594968239642927\\
0.720408163265306	0.0578723978039103\\
0.72280912364946	0.0578448025769214\\
0.725210084033613	0.0594112853342615\\
0.727611044417767	0.0620667733769281\\
0.730012004801921	0.0650498302382458\\
0.732412965186074	0.0676105261369427\\
0.734813925570228	0.0691895443955365\\
0.737214885954382	0.0694928321932438\\
0.739615846338535	0.0684949683625978\\
0.742016806722689	0.0664217463836584\\
0.744417767106843	0.0637498306004081\\
0.746818727490996	0.0611875697407125\\
0.74921968787515	0.0595484133289756\\
0.751620648259304	0.0594575804979441\\
0.754021608643457	0.0610238133205779\\
0.756422569027611	0.0637634003388138\\
0.758823529411765	0.0668493311046316\\
0.761224489795918	0.0694526749114679\\
0.763625450180072	0.0709730880651378\\
0.766026410564226	0.0711204362604599\\
0.768427370948379	0.0699122228458589\\
0.770828331332533	0.0676482616056502\\
0.773229291716687	0.0648805317964665\\
0.77563025210084	0.0623401839030129\\
0.778031212484994	0.0607696001640808\\
0.780432172869148	0.0606569272769344\\
0.782833133253301	0.0620170911201573\\
0.785234093637455	0.0644049550757524\\
0.787635054021609	0.0671336683614706\\
0.790036014405762	0.0695104351180648\\
0.792436974789916	0.0709941347651776\\
0.79483793517407	0.0712824119940771\\
0.797238895558223	0.0703425252785044\\
0.799639855942377	0.0683967677595783\\
0.802040816326531	0.0659012378425647\\
0.804441776710684	0.0635115728438251\\
0.806842737094838	0.0619674805096941\\
0.809243697478992	0.0618466799206385\\
0.811644657863145	0.0632837116352033\\
0.814045618247299	0.0658798323098559\\
0.816446578631453	0.0688906973286806\\
0.818847539015606	0.0715155644772404\\
0.82124849939976	0.0731215029563204\\
0.823649459783914	0.0733629505993457\\
0.826050420168067	0.0722106333126845\\
0.828451380552221	0.0699369727496404\\
0.830852340936375	0.0670896652152668\\
0.833253301320528	0.0644209622648058\\
0.835654261704682	0.06271866387983\\
0.838055222088836	0.0625253487746444\\
0.840456182472989	0.0638842343952355\\
0.842857142857143	0.0663322228260304\\
0.845258103241296	0.069142868665658\\
0.84765906362545	0.0715964415297792\\
0.850060024009604	0.0731359322751181\\
0.852460984393758	0.073436485151163\\
0.854861944777911	0.0724427969186379\\
0.857262905162065	0.0703753801840068\\
0.859663865546219	0.0677090313604536\\
0.862064825930372	0.0651253704356871\\
0.864465786314526	0.0633969976689917\\
0.866866746698679	0.0631468740728266\\
0.869267707082833	0.0645558947222016\\
0.871668667466987	0.0672492398889053\\
0.874069627851141	0.0704732375028487\\
0.876470588235294	0.0733920578330738\\
0.878871548619448	0.0753138192691193\\
0.881272509003601	0.0758157183225264\\
0.883673469387755	0.0747954708575422\\
0.886074429771909	0.072474042487825\\
0.890876350540216	0.0663078247710907\\
0.89327731092437	0.0641581803132794\\
0.895678271308523	0.0636396434047174\\
0.898079231692677	0.0649195248302271\\
0.900480192076831	0.0675387197225723\\
0.902881152460984	0.0706761246448309\\
0.905282112845138	0.073497572620626\\
0.907683073229292	0.075356396594346\\
0.910084033613445	0.0758615935059722\\
0.912484993997599	0.0749032364067934\\
0.914885954381753	0.0726672209908081\\
0.917286914765906	0.0696306516887905\\
0.91968787515006	0.0665272415516766\\
0.922088835534214	0.064238773929357\\
0.924489795918367	0.0635447724471189\\
0.926890756302521	0.0647706414144833\\
0.929291716686675	0.0675949862354825\\
0.931692677070828	0.0712110053517271\\
0.934093637454982	0.0746767883554442\\
0.936494597839136	0.0771864778477569\\
0.938895558223289	0.0782028423341677\\
0.941296518607443	0.0775151943480235\\
0.943697478991597	0.0752582498616386\\
0.94609843937575	0.0719102710938078\\
0.948499399759904	0.0682742092967176\\
0.950900360144058	0.0653773663447339\\
0.953301320528211	0.0641778336347665\\
0.955702280912365	0.0651124925808323\\
0.958103241296519	0.0678343589159854\\
0.960504201680672	0.0714309703751639\\
0.962905162064826	0.0748808070935073\\
0.96530612244898	0.0773657998664722\\
0.967707082833133	0.0783711906707842\\
0.970108043217287	0.0776947816662701\\
0.972509003601441	0.0754506226630807\\
0.974909963985594	0.0720769344206864\\
0.977310924369748	0.0683274599787774\\
0.979711884753902	0.0651981662751117\\
0.982112845138055	0.0636931416657697\\
0.984513805522209	0.0644029182056168\\
0.986914765906363	0.0671614583306547\\
0.989315726290516	0.0711327846390447\\
0.99171668667467	0.0752305443413673\\
0.994117647058824	0.0784875975963828\\
0.996518607442977	0.080224346105977\\
0.998919567827131	0.0801021462063296\\
};
\addlegendentry{ABC$_W^{1/2}$}

\end{axis}
\end{tikzpicture}%

%% file: images/Octant/Cost_V.tex
%
%
\definecolor{mycolor1}{rgb}{1.00000,0.60000,0.00000}%
\pgfplotsset{scaled y ticks=false}
\begin{tikzpicture}[scale = 0.52, font=\huge]

\begin{axis}[%
width=4.602in,
height=3.82in,
at={(0.772in,0.516in)},
scale only axis,
xmin=0,
xmax=1,
xlabel style={font=\huge\color{white!15!black}},
xlabel={relative simulation time $t/T$}, xtick={0, 0.2, 0.4, 0.6, 0.8, 1},
y tick label style={
        /pgf/number format/.cd,
            fixed,
            fixed zerofill,
            precision=2,
        /tikz/.cd
            },
ymin=0,
ymax=0.09, ytick={0.01, 0.03, 0.05, 0.07, 0.09},
ylabel style={font=\huge\color{white!15!black}},
ylabel={},
axis background/.style={fill=white},
xmajorgrids,
ymajorgrids,
legend style={at={(0.03,0.97)}, anchor=north west, legend cell align=left, align=left, draw=white!15!black}
]
\addplot [color=mycolor1, line width=2.0pt]
  table[row sep=crcr]{%
0.00252100840336134	6.21724893790088e-15\\
0.187394957983193	2.78745822634985e-05\\
0.204201680672269	0.000109222392993824\\
0.213805522208884	0.00020491832533287\\
0.223409363745498	0.000384768649810407\\
0.247418967587035	0.000937641767185249\\
0.266626650660264	0.00147868027982734\\
0.271428571428571	0.00169698871579371\\
0.276230492196879	0.00198703600774341\\
0.285834333733493	0.00265202043901147\\
0.293037214885954	0.00309949842339252\\
0.307442977190876	0.00395410451049372\\
0.319447779111645	0.00479523005631766\\
0.324249699879952	0.00500323123443758\\
0.329051620648259	0.00516991644436948\\
0.331452581032413	0.00532619315573646\\
0.333853541416567	0.0055610119978371\\
0.338655462184874	0.00625102231696728\\
0.341056422569028	0.00651817722328052\\
0.343457382953181	0.00658292834789764\\
0.345858343337335	0.0063808166148509\\
0.350660264105642	0.00560821375328779\\
0.353061224489796	0.00558032069289138\\
0.35546218487395	0.0060820863905483\\
0.360264105642257	0.00779429419232613\\
0.362665066026411	0.008216434811953\\
0.365066026410564	0.00796265503760785\\
0.367466986794718	0.00721135726729472\\
0.369867947178872	0.00688335952328412\\
0.372268907563025	0.00687270960751485\\
0.374669867947179	0.0065761885790202\\
0.377070828331333	0.00644490220257177\\
0.379471788715486	0.00687485142166588\\
0.38187274909964	0.00737322698554455\\
0.384273709483794	0.00799738062509048\\
0.386674669867947	0.00781445505736811\\
0.389075630252101	0.00768902787243053\\
0.391476590636255	0.00750853376753546\\
0.393877551020408	0.00720996315405753\\
0.396278511404562	0.00681844745401117\\
0.398679471788716	0.00684708424115121\\
0.401080432172869	0.0070715273304196\\
0.403481392557023	0.00770178897142937\\
0.405882352941176	0.0074887902232198\\
0.413085234093638	0.00798018255755173\\
0.415486194477791	0.00798862467466321\\
0.417887154861945	0.00771084634342467\\
0.420288115246098	0.00760117873816579\\
0.422689075630252	0.00795181147652679\\
0.425090036014406	0.00770259721328925\\
0.427490996398559	0.00763045793606931\\
0.429891956782713	0.00773801057217416\\
0.432292917166867	0.00806317774042675\\
0.437094837935174	0.00834804277675716\\
0.439495798319328	0.00857252212595283\\
0.441896758703481	0.00831427966004405\\
0.444297719087635	0.00794692894197802\\
0.446698679471789	0.00787222283057865\\
0.449099639855942	0.00816642350281804\\
0.451500600240096	0.00832213605618837\\
0.45390156062425	0.00789232613330881\\
0.456302521008403	0.00823675117105482\\
0.458703481392557	0.00861736982370431\\
0.461104441776711	0.00909701069230706\\
0.463505402160864	0.00870930965348549\\
0.468307322929172	0.00911849617538985\\
0.470708283313325	0.0090052894685585\\
0.473109243697479	0.00878182816140738\\
0.475510204081633	0.00831799289998936\\
0.477911164465786	0.00818721270183742\\
0.482713085234094	0.00856603732109606\\
0.485114045618247	0.00831270097042436\\
0.487515006002401	0.00856882337369613\\
0.489915966386555	0.00895075041473647\\
0.492316926770708	0.00940811804170294\\
0.494717887154862	0.00885648292336882\\
0.497118847539016	0.00921450554513126\\
0.499519807923169	0.0092548767410362\\
0.501920768307323	0.00906233597572992\\
0.504321728691477	0.00871947785463623\\
0.50672268907563	0.00867999220015825\\
0.509123649459784	0.00877626582097057\\
0.511524609843938	0.0089352736357009\\
0.513925570228091	0.00873599493540356\\
0.518727490996399	0.00939641554294723\\
0.521128451380552	0.00979831956997768\\
0.523529411764706	0.00933622692687885\\
0.52593037214886	0.00979109076705009\\
0.528331332533013	0.00963227981784265\\
0.530732292917167	0.00952942098429921\\
0.535534213685474	0.0089405272600348\\
0.537935174069628	0.00898143060619361\\
0.540336134453781	0.00929267994313721\\
0.545138055222089	0.0090994852729005\\
0.547539015606243	0.00939111861688124\\
0.549939975990396	0.00984358229767612\\
0.55234093637455	0.00983080685324633\\
0.554741896758703	0.00944167875213353\\
0.557142857142857	0.00981720604770508\\
0.559543817527011	0.00973872684668142\\
0.561944777911164	0.00943609396396006\\
0.564345738295318	0.00919272382321312\\
0.566746698679472	0.00931871224929526\\
0.569147659063625	0.00949034830996653\\
0.573949579831933	0.00949674783119647\\
0.576350540216086	0.00982563278225279\\
0.57875150060024	0.010188014347657\\
0.581152460984394	0.0102974057399992\\
0.583553421368547	0.00998622252305392\\
0.585954381752701	0.0102777851305406\\
0.588355342136855	0.0100450030967862\\
0.590756302521008	0.00994858305207269\\
0.593157262905162	0.00959940301984208\\
0.595558223289316	0.00947804616374615\\
0.597959183673469	0.00959887168343165\\
0.600360144057623	0.00992784891686549\\
0.602761104441777	0.009706578217763\\
0.60516206482593	0.00972247813226834\\
0.607563025210084	0.0100088293356373\\
0.609963985594238	0.0104801414811665\\
0.612364945978391	0.0100324816348438\\
0.614765906362545	0.0100302285917866\\
0.617166866746699	0.0102953797107527\\
0.619567827130852	0.0102038572047437\\
0.621968787515006	0.00987097017652094\\
0.62436974789916	0.00981955252246758\\
0.629171668667467	0.0101299546849052\\
0.631572629051621	0.00989786398790715\\
0.633973589435774	0.0100273560343859\\
0.636374549819928	0.0103399679806595\\
0.638775510204082	0.0107670645678385\\
0.641176470588235	0.0105378915721657\\
0.645978391356543	0.0107580610969874\\
0.648379351740696	0.0106466327641116\\
0.65078031212485	0.0104619667634216\\
0.653181272509004	0.0100628213592631\\
0.655582232893157	0.00999601706036202\\
0.660384153661465	0.010178578734517\\
0.662785114045618	0.00991958663155557\\
0.665186074429772	0.0101256944212486\\
0.669987995198079	0.010947248273853\\
0.672388955582233	0.0104640387847413\\
0.674789915966387	0.0108452547429111\\
0.67719087635054	0.010960851844104\\
0.679591836734694	0.0106630016917604\\
0.681992797118848	0.0102196891925493\\
0.684393757503001	0.0101629185287551\\
0.686794717887155	0.0101860271968701\\
0.689195678271309	0.0103326651513344\\
0.691596638655462	0.0101690393055693\\
0.696398559423769	0.0109191842011129\\
0.698799519807923	0.0113175426511309\\
0.701200480192077	0.0108990378706971\\
0.703601440576231	0.01128633353641\\
0.706002400960384	0.0111297375284866\\
0.708403361344538	0.0109402914222101\\
0.710804321728691	0.0105700554812828\\
0.713205282112845	0.0102376273984218\\
0.715606242496999	0.0102597290237842\\
0.718007202881152	0.010500756717705\\
0.720408163265306	0.010502376133168\\
0.72280912364946	0.0104683648454166\\
0.727611044417767	0.0112557905136207\\
0.730012004801921	0.0112909092901362\\
0.732412965186074	0.0109396590817509\\
0.734813925570228	0.0112601636472143\\
0.737214885954382	0.0110716505323939\\
0.739615846338535	0.0106114425270828\\
0.742016806722689	0.0102896357256853\\
0.746818727490996	0.010541424556917\\
0.751620648259304	0.0107048169062157\\
0.754021608643457	0.0112110889329897\\
0.756422569027611	0.0115833869255112\\
0.758823529411765	0.0116884506719394\\
0.761224489795918	0.0113124848750047\\
0.763625450180072	0.0116213397617806\\
0.770828331332533	0.0105104581376919\\
0.773229291716687	0.0104524098366149\\
0.77563025210084	0.0105704289927065\\
0.778031212484994	0.0109550625501003\\
0.780432172869148	0.0109071464991451\\
0.782833133253301	0.0111440763437223\\
0.785234093637455	0.01153496488143\\
0.787635054021609	0.0118418817756188\\
0.790036014405762	0.0113835667590436\\
0.792436974789916	0.0113106806216471\\
0.79483793517407	0.0114036361974005\\
0.797238895558223	0.0110204560662712\\
0.799639855942377	0.0105418533899287\\
0.802040816326531	0.0105187794075519\\
0.806842737094838	0.011024489348876\\
0.809243697478992	0.0110059224331389\\
0.811644657863145	0.011329998356985\\
0.814045618247299	0.0118697937353038\\
0.816446578631453	0.012145696699803\\
0.818847539015606	0.0117829685746127\\
0.82124849939976	0.0117118378015731\\
0.823649459783914	0.0116995699791299\\
0.828451380552221	0.0108323125854999\\
0.830852340936375	0.0105517781320091\\
0.833253301320528	0.0107197524462378\\
0.835654261704682	0.0109561954849718\\
0.838055222088836	0.0113023673671383\\
0.840456182472989	0.0112732519969994\\
0.842857142857143	0.0117376595913737\\
0.845258103241296	0.011996671367744\\
0.84765906362545	0.0120784694140881\\
0.850060024009604	0.0114525234902813\\
0.852460984393758	0.0116078323545723\\
0.854861944777911	0.0114341451455332\\
0.857262905162065	0.0108617208606014\\
0.859663865546219	0.0105078763284131\\
0.864465786314526	0.0109808045863244\\
0.866866746698679	0.0112764995232852\\
0.869267707082833	0.0112516880768385\\
0.871668667466987	0.0118552946511046\\
0.874069627851141	0.0123344137921025\\
0.876470588235294	0.0124186199648195\\
0.878871548619448	0.011789444030829\\
0.881272509003601	0.0120171672191343\\
0.883673469387755	0.0116704504430161\\
0.886074429771909	0.0111296025295485\\
0.888475390156062	0.0107368799306035\\
0.890876350540216	0.0106965038510732\\
0.89327731092437	0.0109416374051431\\
0.895678271308523	0.0113147811159187\\
0.900480192076831	0.0116915889540354\\
0.902881152460984	0.0122755702635161\\
0.905282112845138	0.012413174230127\\
0.907683073229292	0.0121849198360964\\
0.910084033613445	0.0116591208378619\\
0.912484993997599	0.0117447530965096\\
0.914885954381753	0.0113099123753191\\
0.917286914765906	0.0107443737687573\\
0.91968787515006	0.0106508233725324\\
0.922088835534214	0.0110519672541147\\
0.926890756302521	0.0115754595517681\\
0.929291716686675	0.0117374358784943\\
0.931692677070828	0.0124745615124868\\
0.934093637454982	0.0127552619885322\\
0.936494597839136	0.0125655004829892\\
0.938895558223289	0.011982558828495\\
0.941296518607443	0.0121231332606895\\
0.943697478991597	0.0115317944236915\\
0.94609843937575	0.0110268780213585\\
0.948499399759904	0.0107501427558899\\
0.950900360144058	0.0109197419722561\\
0.953301320528211	0.0112243611638528\\
0.955702280912365	0.0117417638948396\\
0.958103241296519	0.0118042743694705\\
0.962905162064826	0.0126688573068027\\
0.96530612244898	0.0126750959194755\\
0.967707082833133	0.0120586540390746\\
0.970108043217287	0.0117505124684093\\
0.972509003601441	0.0116102900709928\\
0.974909963985594	0.0110728683391547\\
0.977310924369748	0.0106445397365285\\
0.979711884753902	0.0108165979096235\\
0.982112845138055	0.0112953239015529\\
0.984513805522209	0.0117019377201046\\
0.986914765906363	0.0118628413144108\\
0.989315726290516	0.0122874332659741\\
0.99171668667467	0.0129026508227554\\
0.994117647058824	0.0129482046473058\\
0.996518607442977	0.0124079096331468\\
0.998919567827131	0.0120455587638693\\
};
\addlegendentry{ABC$_W^{1/2}$ adaptive}

\addplot [color=mycolor1, dashed, line width=2.0pt]
  table[row sep=crcr]{%
0.00252100840336134	6.21724893790088e-15\\
0.187394957983193	2.78745822634985e-05\\
0.204201680672269	0.000109222392993824\\
0.213805522208884	0.000202373101602649\\
0.221008403361345	0.000329710353446977\\
0.230612244897959	0.000567524003315389\\
0.23781512605042	0.000809920605769121\\
0.245018007202881	0.00113615449245741\\
0.252220888355342	0.00155679568805278\\
0.261824729891957	0.00221370711917157\\
0.269027611044418	0.00276179339338656\\
0.273829531812725	0.00319866968240257\\
0.276230492196879	0.00345001789554511\\
0.281032412965186	0.00404274819716477\\
0.285834333733493	0.00473743538973526\\
0.295438175270108	0.00623712239180407\\
0.300240096038415	0.00702575954921736\\
0.305042016806723	0.00791870397714689\\
0.307442977190876	0.00841550537282221\\
0.30984393757503	0.00894923745151344\\
0.314645858343337	0.0100977603801293\\
0.319447779111645	0.0112557606417154\\
0.324249699879952	0.0123208278355775\\
0.329051620648259	0.0133446327847271\\
0.331452581032413	0.0139101018754261\\
0.333853541416567	0.014550488734104\\
0.33625450180072	0.015273509804359\\
0.343457382953181	0.0176440452584923\\
0.345858343337335	0.0183285146489767\\
0.348259303721489	0.0188992495205698\\
0.353061224489796	0.0198783887767127\\
0.35546218487395	0.0204684539187453\\
0.357863145258103	0.0212299561481512\\
0.360264105642257	0.0221858306328423\\
0.362665066026411	0.0232990511153298\\
0.365066026410564	0.0244813514247446\\
0.367466986794718	0.0256112750701413\\
0.369867947178872	0.0265468889389738\\
0.372268907563025	0.0271566349667629\\
0.374669867947179	0.0273680503261657\\
0.377070828331333	0.0272213930925138\\
0.379471788715486	0.0269258804196704\\
0.38187274909964	0.0268607970887172\\
0.384273709483794	0.027442864417922\\
0.386674669867947	0.0288874094397973\\
0.389075630252101	0.0310621249545412\\
0.391476590636255	0.0335707001038132\\
0.393877551020408	0.0359469257108427\\
0.396278511404562	0.037774305672119\\
0.398679471788716	0.0387132913173682\\
0.401080432172869	0.0385078023769275\\
0.403481392557023	0.0370808955579616\\
0.405882352941176	0.0347244782557444\\
0.40828331332533	0.0324660470108071\\
0.410684273709484	0.0317369061425372\\
0.413085234093638	0.033183606998492\\
0.415486194477791	0.0362398198878295\\
0.417887154861945	0.0398981216852591\\
0.420288115246098	0.0433590393428406\\
0.422689075630252	0.0461197913324689\\
0.425090036014406	0.0478559731621793\\
0.427490996398559	0.0483184336211565\\
0.429891956782713	0.0473207850968022\\
0.432292917166867	0.0447987614168109\\
0.43469387755102	0.0410288426248758\\
0.437094837935174	0.0370456140302212\\
0.439495798319328	0.0346676673562952\\
0.441896758703481	0.0353330450413845\\
0.444297719087635	0.0387827979793031\\
0.449099639855942	0.0484524840621301\\
0.451500600240096	0.0524789714250034\\
0.45390156062425	0.0551521740442985\\
0.456302521008403	0.056131324562577\\
0.458703481392557	0.0551614057793266\\
0.461104441776711	0.0522026448036421\\
0.463505402160864	0.0477934120590935\\
0.465906362545018	0.0431991030080123\\
0.468307322929172	0.0400259245673186\\
0.470708283313325	0.0395170211523957\\
0.473109243697479	0.0417817491588551\\
0.475510204081633	0.0457944127003383\\
0.477911164465786	0.0502729743183224\\
0.48031212484994	0.0542627989529924\\
0.482713085234094	0.0571854501579726\\
0.485114045618247	0.0586763261153056\\
0.487515006002401	0.0584817542422148\\
0.489915966386555	0.0564638199038044\\
0.492316926770708	0.0527618612905819\\
0.494717887154862	0.0481435740943874\\
0.497118847539016	0.0440777041656321\\
0.499519807923169	0.0421240117976833\\
0.501920768307323	0.0429768328189051\\
0.504321728691477	0.0461276967872103\\
0.50672268907563	0.0504632347472911\\
0.509123649459784	0.0549147817598453\\
0.511524609843938	0.0586681045206614\\
0.513925570228091	0.0611384668083818\\
0.516326530612245	0.0618948702124509\\
0.518727490996399	0.0606529288652869\\
0.521128451380552	0.0574643334659541\\
0.523529411764706	0.0530281884697322\\
0.52593037214886	0.0486926591896307\\
0.528331332533013	0.0459993648415897\\
0.530732292917167	0.0459042056710202\\
0.533133253301321	0.0481946440042753\\
0.535534213685474	0.051813425369945\\
0.537935174069628	0.0556545267347858\\
0.540336134453781	0.0589549982853943\\
0.542737094837935	0.0612488618601589\\
0.545138055222089	0.0622302913829709\\
0.547539015606243	0.0616716008867202\\
0.549939975990396	0.0594766178681531\\
0.55234093637455	0.0558934853777824\\
0.554741896758703	0.0517906393316544\\
0.557142857142857	0.0485189818107513\\
0.559543817527011	0.047255326731762\\
0.561944777911164	0.0483754998566723\\
0.564345738295318	0.051373784532022\\
0.569147659063625	0.0592436847235741\\
0.571548619447779	0.0624200175566498\\
0.573949579831933	0.064278517495417\\
0.576350540216086	0.0644198544168669\\
0.57875150060024	0.0626437147921562\\
0.581152460984394	0.0591783996140306\\
0.583553421368547	0.0549273073828538\\
0.585954381752701	0.0512671836019486\\
0.588355342136855	0.0494870488447761\\
0.590756302521008	0.0501109380362486\\
0.593157262905162	0.0526401346995274\\
0.597959183673469	0.0594713629495029\\
0.600360144057623	0.0622441928507586\\
0.602761104441777	0.0639994939626055\\
0.60516206482593	0.0644638799605863\\
0.607563025210084	0.0634452588343576\\
0.609963985594238	0.0609107304863107\\
0.612364945978391	0.0572465950102956\\
0.614765906362545	0.0534410501608602\\
0.617166866746699	0.0507822445885413\\
0.619567827130852	0.0501976038567417\\
0.621968787515006	0.051799883120188\\
0.62436974789916	0.0549792738617129\\
0.626770708283313	0.0588229994549201\\
0.629171668667467	0.0624516107722435\\
0.631572629051621	0.0651696448014562\\
0.633973589435774	0.0664787395948429\\
0.636374549819928	0.0660364793567925\\
0.638775510204082	0.0637435215541511\\
0.641176470588235	0.0600138286639476\\
0.643577430972389	0.0559092869045192\\
0.645978391356543	0.0528126673709961\\
0.648379351740696	0.05180108436383\\
0.65078031212485	0.0530691418092001\\
0.653181272509004	0.0559078395841294\\
0.655582232893157	0.0592923516173137\\
0.657983193277311	0.0624027384382553\\
0.660384153661465	0.0647189695464618\\
0.662785114045618	0.0659190914038611\\
0.665186074429772	0.0657881674968727\\
0.667587034813926	0.0642024900456426\\
0.669987995198079	0.0612534241828336\\
0.672388955582233	0.0575421658360842\\
0.674789915966387	0.0542280798165596\\
0.67719087635054	0.0525219030625897\\
0.679591836734694	0.0529927631501461\\
0.681992797118848	0.0553429442296476\\
0.684393757503001	0.0587518777847538\\
0.686794717887155	0.0623195100799471\\
0.689195678271309	0.065294214837304\\
0.691596638655462	0.0671266198972046\\
0.693997599039616	0.0674555181721037\\
0.696398559423769	0.0660803826356022\\
0.698799519807923	0.0631217625783005\\
0.701200480192077	0.0592887129820936\\
0.703601440576231	0.0558504842064705\\
0.706002400960384	0.0540884035427041\\
0.708403361344538	0.0545827412127711\\
0.710804321728691	0.0569007322901141\\
0.713205282112845	0.0600363589582855\\
0.715606242496999	0.0630426593697863\\
0.718007202881152	0.0653250651855668\\
0.720408163265306	0.0665759683106303\\
0.72280912364946	0.0666550928148996\\
0.725210084033613	0.0655151924580906\\
0.727611044417767	0.0632008088837844\\
0.732412965186074	0.056893441487282\\
0.734813925570228	0.0549413394376281\\
0.737214885954382	0.0549321690681128\\
0.739615846338535	0.0567467442052945\\
0.742016806722689	0.0596670499167934\\
0.744417767106843	0.0628483364894359\\
0.746818727490996	0.0656021239301494\\
0.74921968787515	0.0674313889228582\\
0.751620648259304	0.0680079360514213\\
0.754021608643457	0.0671403590851593\\
0.756422569027611	0.0648238481290611\\
0.761224489795918	0.0580621050288034\\
0.763625450180072	0.0559427423134891\\
0.766026410564226	0.0559717947179962\\
0.768427370948379	0.0580019653554691\\
0.773229291716687	0.0641066266745349\\
0.77563025210084	0.0663824213811044\\
0.778031212484994	0.0675566026145662\\
0.780432172869148	0.0675578089620286\\
0.782833133253301	0.066437215713319\\
0.785234093637455	0.0643196220599893\\
0.787635054021609	0.0614658449208079\\
0.790036014405762	0.0584962504125565\\
0.792436974789916	0.0564525692849763\\
0.79483793517407	0.0562331934044601\\
0.797238895558223	0.0578911019468619\\
0.799639855942377	0.0607194129371699\\
0.802040816326531	0.0638203790584894\\
0.804441776710684	0.0664717702812336\\
0.806842737094838	0.06820496656803\\
0.809243697478992	0.0687669711035525\\
0.811644657863145	0.068039036927088\\
0.814045618247299	0.0660243630170821\\
0.816446578631453	0.0629565150416784\\
0.818847539015606	0.0595685268672238\\
0.82124849939976	0.0571218202995526\\
0.823649459783914	0.0567383669052401\\
0.826050420168067	0.0585811105680348\\
0.828451380552221	0.0617843439568908\\
0.830852340936375	0.0651294423388084\\
0.833253301320528	0.0676823923268498\\
0.835654261704682	0.0689823683595355\\
0.838055222088836	0.0689393562327845\\
0.840456182472989	0.067679938217404\\
0.842857142857143	0.0654309399931883\\
0.845258103241296	0.062522955130133\\
0.84765906362545	0.0594842630228779\\
0.850060024009604	0.0572026334114698\\
0.852460984393758	0.0566797872873347\\
0.854861944777911	0.058245999495166\\
0.857262905162065	0.0612566155958705\\
0.859663865546219	0.0646489286901826\\
0.862064825930372	0.067532831966441\\
0.864465786314526	0.0693849044020928\\
0.866866746698679	0.0699751858127773\\
0.869267707082833	0.0692605148930479\\
0.871668667466987	0.0673071883260465\\
0.874069627851141	0.0643142468725554\\
0.876470588235294	0.0608131161242588\\
0.878871548619448	0.0579008324956044\\
0.881272509003601	0.0569018991474716\\
0.883673469387755	0.0584033675123716\\
0.886074429771909	0.0617706345457735\\
0.888475390156062	0.0656777501135815\\
0.890876350540216	0.0688883737506346\\
0.89327731092437	0.070680950088954\\
0.895678271308523	0.0708498865240478\\
0.898079231692677	0.0695385315976531\\
0.900480192076831	0.0670773625171102\\
0.902881152460984	0.063879411758423\\
0.905282112845138	0.0604829012267052\\
0.907683073229292	0.0577056094835036\\
0.910084033613445	0.0566123356056546\\
0.912484993997599	0.0578645751540279\\
0.914885954381753	0.0610463059518364\\
0.917286914765906	0.0649750646093045\\
0.91968787515006	0.0684724105565553\\
0.922088835534214	0.0707909882364415\\
0.924489795918367	0.07162987418042\\
0.926890756302521	0.0709938375522899\\
0.929291716686675	0.0690228854129659\\
0.931692677070828	0.0659472408804875\\
0.934093637454982	0.0621983452071113\\
0.936494597839136	0.0586752963852173\\
0.938895558223289	0.0567758363487356\\
0.941296518607443	0.0575878964502672\\
0.943697478991597	0.0609106626732413\\
0.94609843937575	0.0654149113997667\\
0.948499399759904	0.069554224732244\\
0.950900360144058	0.0722356318461032\\
0.953301320528211	0.0729978268940659\\
0.955702280912365	0.0718962894116415\\
0.958103241296519	0.0692998788658963\\
0.960504201680672	0.0657098024180585\\
0.962905162064826	0.0617380165910225\\
0.96530612244898	0.0582089419810707\\
0.967707082833133	0.0562256910645271\\
0.970108043217287	0.0567882176154193\\
0.972509003601441	0.0598969403973093\\
0.974909963985594	0.0643867521064625\\
0.977310924369748	0.0687548989764646\\
0.979711884753902	0.0718923316818701\\
0.982112845138055	0.0732881266126753\\
0.984513805522209	0.0729030316491391\\
0.986914765906363	0.0709495065143743\\
0.989315726290516	0.0677389058332905\\
0.99171668667467	0.0636946638807885\\
0.994117647058824	0.0595551820846795\\
0.996518607442977	0.0566134025774888\\
0.998919567827131	0.0563429771231591\\
};
\addlegendentry{ABC$_W^{1/2}$}

\end{axis}
\end{tikzpicture}%

%% file: images/Transducer/Cost.tex
%
%
\definecolor{mycolor1}{rgb}{1.00000,0.60000,0.00000}%
\begin{tikzpicture}[scale = 0.52, font=\huge]

\begin{axis}[%
width=4.602in,
height=3.82in,
at={(0.772in,0.516in)},
scale only axis,
xmin=0,
xmax=1,
xlabel style={font=\huge\color{white!15!black}},
y tick label style={
        /pgf/number format/.cd,
            fixed,
            fixed zerofill,
            precision=2,
        /tikz/.cd
            },
xlabel={relative simulation time $t/T$}, xtick={0, 0.2, 0.4, 0.6, 0.8, 1},
ymin=0,
ymax=0.12,  ytick={0.02, 0.04, 0.06, 0.08, 0.1, 0.12},
ylabel style={at={(-0.15,0.5)}, font=\huge\color{white!15!black}},
ylabel={relative $L^{2}(\Omega)$ error},
axis background/.style={fill=white},
xmajorgrids,
ymajorgrids,
legend style={at={(0.03,0.97)}, anchor=north west, legend cell align=left, align=left, draw=white!15!black}
]
\addplot [color=mycolor1, line width=2.0pt]
  table[row sep=crcr]{%
0.00214285714285711	0\\
0.104183673469388	3.09745474589285e-05\\
0.114387755102041	0.000153286872072278\\
0.138877551020408	0.000587671064360085\\
0.151122448979592	0.00088222627119805\\
0.159285714285714	0.0011484594517831\\
0.169489795918367	0.00157636056317689\\
0.17969387755102	0.00209146802563265\\
0.185816326530612	0.00247253526162949\\
0.191938775510204	0.00293786977886024\\
0.198061224489796	0.00350468742929078\\
0.21030612244898	0.00474331140186113\\
0.214387755102041	0.00505879948908894\\
0.222551020408163	0.00554898051120267\\
0.226632653061224	0.00581504138218214\\
0.230714285714286	0.0061481446767041\\
0.234795918367347	0.00657022365660287\\
0.238877551020408	0.00709329064812736\\
0.242959183673469	0.00771851854295613\\
0.247040816326531	0.00843382783264468\\
0.259285714285714	0.0106825259568899\\
0.263367346938776	0.0112462531467162\\
0.267448979591837	0.0116432389154641\\
0.271530612244898	0.0119093284964716\\
0.27765306122449	0.0122716779695091\\
0.27969387755102	0.0124419850716678\\
0.281734693877551	0.0126612266454176\\
0.283775510204082	0.0129450261256957\\
0.285816326530612	0.0133075056714719\\
0.287857142857143	0.013757037204854\\
0.289897959183673	0.0142916902379895\\
0.291938775510204	0.014904257450978\\
0.293979591836735	0.0155829901489701\\
0.298061224489796	0.0170625403091139\\
0.302142857142857	0.0185393425864206\\
0.304183673469388	0.0192067280148449\\
0.306224489795918	0.0197886333033074\\
0.308265306122449	0.0202620414975077\\
0.31030612244898	0.020609218566266\\
0.31234693877551	0.0208232620868859\\
0.314387755102041	0.0209131560577145\\
0.316428571428572	0.0208991366816839\\
0.320510204081633	0.0206825416572947\\
0.324591836734694	0.0204577468523903\\
0.326632653061224	0.0204303851890306\\
0.328673469387755	0.020495968454774\\
0.330714285714286	0.0206687289867489\\
0.332755102040816	0.020948414422946\\
0.334795918367347	0.021320442867665\\
0.338877551020408	0.0222235825699312\\
0.342959183673469	0.0230934314833966\\
0.345	0.0234366621645737\\
0.347040816326531	0.0236885072585583\\
0.349081632653061	0.0238402861102954\\
0.351122448979592	0.0238936102674645\\
0.355204081632653	0.0237991897319966\\
0.357244897959184	0.0237388187963704\\
0.359285714285714	0.0237390012625586\\
0.361326530612245	0.0238499099436658\\
0.363367346938776	0.0241087168792654\\
0.365408163265306	0.0245300542381056\\
0.367448979591837	0.0250930850835884\\
0.369489795918367	0.0257587053274225\\
0.37765306122449	0.0286913140914944\\
0.37969387755102	0.0293198966648706\\
0.381734693877551	0.0298519369588252\\
0.383775510204082	0.0302626926595039\\
0.385816326530612	0.0305214106571394\\
0.387857142857143	0.0306072969874592\\
0.389897959183674	0.0305325875276902\\
0.391938775510204	0.0303420435778996\\
0.398061224489796	0.0296025673405658\\
0.400102040816327	0.0294927909354085\\
0.402142857142857	0.0295476990330161\\
0.404183673469388	0.0297909925155228\\
0.406224489795918	0.0302312845119607\\
0.408265306122449	0.0308824083517288\\
0.41030612244898	0.0317504617769446\\
0.41234693877551	0.0328117143003368\\
0.418469387755102	0.0363684104254408\\
0.420510204081633	0.0373192715810817\\
0.422551020408163	0.037995752073425\\
0.424591836734694	0.0383403205542532\\
0.426632653061225	0.0383340420251922\\
0.428673469387755	0.038001976091434\\
0.430714285714286	0.0373911990527773\\
0.432755102040816	0.0365485942157787\\
0.434795918367347	0.0355105544162018\\
0.436836734693878	0.034307535918914\\
0.438877551020408	0.0329747512142082\\
0.442959183673469	0.0301960352127604\\
0.445	0.0289867696112949\\
0.447040816326531	0.0281169192461128\\
0.449081632653061	0.0277423641476209\\
0.451122448979592	0.0279569921336235\\
0.453163265306122	0.0287832954746328\\
0.455204081632653	0.0301528585583871\\
0.457244897959184	0.0319194099756666\\
0.459285714285714	0.0339245720405096\\
0.463367346938776	0.0382324652486085\\
0.465408163265306	0.04040907884982\\
0.467448979591837	0.0425093911727036\\
0.469489795918367	0.04442300234798\\
0.471530612244898	0.0460112893773988\\
0.473571428571429	0.0471282182553658\\
0.475612244897959	0.0476650144759574\\
0.47765306122449	0.047594089812988\\
0.47969387755102	0.0469635770605813\\
0.481734693877551	0.0458639178173923\\
0.483775510204082	0.0444056734429628\\
0.485816326530612	0.0426873524903442\\
0.487857142857143	0.0407817983255453\\
0.489897959183674	0.0387344222721623\\
0.496020408163265	0.032251584974303\\
0.498061224489796	0.0304748209299629\\
0.500102040816326	0.0293339096881677\\
0.502142857142857	0.0290529694324005\\
0.504183673469388	0.0296584004772998\\
0.506224489795918	0.0309918180930365\\
0.508265306122449	0.0328189040654209\\
0.51030612244898	0.0349226818135088\\
0.514387755102041	0.0394455990068724\\
0.518469387755102	0.0440043246569564\\
0.520510204081633	0.0460521828548393\\
0.522551020408163	0.0477050003069023\\
0.524591836734694	0.0487509509527311\\
0.526632653061225	0.0490503807352908\\
0.528673469387755	0.0485934363510037\\
0.530714285714286	0.0474847251080908\\
0.532755102040816	0.0458900242522753\\
0.534795918367347	0.0439846553812483\\
0.536836734693878	0.0419072751841126\\
0.540918367346939	0.0375300205107676\\
0.542959183673469	0.0353096394605085\\
0.545	0.033200906293167\\
0.547040816326531	0.0314468413320844\\
0.549081632653061	0.030359356024663\\
0.551122448979592	0.0301994883919421\\
0.553163265306122	0.0310514501957461\\
0.555204081632653	0.0327800855505852\\
0.557244897959184	0.035085654183232\\
0.561326530612245	0.0401918851067496\\
0.563367346938776	0.0426285733301454\\
0.565408163265306	0.0449156686464035\\
0.567448979591837	0.047060041032061\\
0.569489795918367	0.0490511984666587\\
0.571530612244898	0.0508150205093574\\
0.573571428571429	0.0522022825004249\\
0.575612244897959	0.0530142361068476\\
0.57765306122449	0.053091454165016\\
0.57969387755102	0.0523809139698943\\
0.581734693877551	0.0509557819310554\\
0.583775510204082	0.0489932412810776\\
0.585816326530612	0.0467040525668133\\
0.587857142857143	0.0442531656072429\\
0.589897959183673	0.041719165517707\\
0.596020408163265	0.0338371102285909\\
0.598061224489796	0.0315968938846173\\
0.600102040816327	0.0301167144602349\\
0.602142857142857	0.0297348348741414\\
0.604183673469388	0.0305373825437566\\
0.606224489795918	0.0322820439019951\\
0.608265306122449	0.0345375184924486\\
0.61030612244898	0.0369255927250295\\
0.61234693877551	0.0392338761159198\\
0.614387755102041	0.0413927016860728\\
0.616428571428571	0.0434140566071343\\
0.618469387755102	0.0453149942144837\\
0.620510204081633	0.0470533346597838\\
0.622551020408163	0.048506210820179\\
0.624591836734694	0.0494958844524375\\
0.626632653061225	0.0498381993646079\\
0.628673469387755	0.0494423355886422\\
0.630714285714286	0.0483637023354672\\
0.632755102040816	0.0467812110191289\\
0.634795918367347	0.0449301797193338\\
0.636836734693878	0.0430192114706323\\
0.638877551020408	0.0411740108333416\\
0.640918367346939	0.0394423268879827\\
0.642959183673469	0.0378364299948936\\
0.645	0.0364156194321027\\
0.647040816326531	0.0353323941110114\\
0.649081632653061	0.0348082178135294\\
0.651122448979592	0.0350527475869226\\
0.653163265306122	0.0361554691782995\\
0.655204081632653	0.0380177643780606\\
0.657244897959184	0.0403586559275222\\
0.659285714285714	0.0428434389913722\\
0.661326530612245	0.0452300174266766\\
0.663367346938776	0.0474103921778626\\
0.665408163265306	0.0493683085941788\\
0.667448979591837	0.0511234326079428\\
0.669489795918367	0.0527024817022728\\
0.671530612244898	0.0540932568937126\\
0.673571428571429	0.0552039429214553\\
0.675612244897959	0.0558450171142246\\
0.67765306122449	0.0558030721501626\\
0.67969387755102	0.0549796459972743\\
0.681734693877551	0.0534493249515539\\
0.683775510204082	0.0514061086336053\\
0.685816326530612	0.0490596349194214\\
0.689897959183673	0.0440373990804442\\
0.693979591836735	0.0389388317202303\\
0.696020408163265	0.0364716270984408\\
0.698061224489796	0.0343028164845302\\
0.700102040816327	0.0327702472715211\\
0.702142857142857	0.0322007784114853\\
0.704183673469388	0.0327273319683562\\
0.706224489795918	0.0341816846466156\\
0.708265306122449	0.036201268262292\\
0.71234693877551	0.0406478136625732\\
0.714387755102041	0.0427672245406386\\
0.716428571428571	0.0447899183956965\\
0.718469387755102	0.046729951906438\\
0.720510204081633	0.0485626153616573\\
0.722551020408163	0.0501807675259736\\
0.724591836734694	0.0514193134794932\\
0.726632653061224	0.0521025388024379\\
0.728673469387755	0.0521179204106686\\
0.730714285714286	0.0514774696832224\\
0.732755102040816	0.0503187870481918\\
0.734795918367347	0.0488452427290166\\
0.742959183673469	0.0424777283176017\\
0.745	0.0409944849987544\\
0.747040816326531	0.0397043987941943\\
0.749081632653061	0.0387951323885033\\
0.751122448979592	0.0384701355883683\\
0.753163265306122	0.0388853924776091\\
0.755204081632653	0.0400660536251836\\
0.757244897959184	0.0418609890409665\\
0.759285714285714	0.0439922771170742\\
0.763367346938776	0.0483645660763766\\
0.765408163265306	0.0504290289214629\\
0.767448979591837	0.0523972177318338\\
0.769489795918367	0.0542581619285409\\
0.771530612244898	0.0559584443346736\\
0.773571428571429	0.0574111277187519\\
0.775612244897959	0.0584837264528014\\
0.77765306122449	0.0589859750757457\\
0.77969387755102	0.0587809164236317\\
0.781734693877551	0.0578667854756501\\
0.783775510204082	0.0563447241966395\\
0.785816326530612	0.0543432228860595\\
0.787857142857143	0.0519894211280666\\
0.789897959183673	0.049404283649303\\
0.791938775510204	0.0466638491557942\\
0.793979591836735	0.0437876794423782\\
0.798061224489796	0.0378720633212464\\
0.800102040816327	0.035278396635446\\
0.802142857142857	0.0334205100002383\\
0.804183673469388	0.0326042031836284\\
0.806224489795918	0.0328724631254584\\
0.808265306122449	0.0339979647761136\\
0.81030612244898	0.0356691272958569\\
0.81234693877551	0.0376314296802946\\
0.814387755102041	0.0397576742523051\\
0.818469387755102	0.044292088108536\\
0.820510204081633	0.0465906258596539\\
0.822551020408163	0.0487698540906903\\
0.824591836734694	0.0506672396909448\\
0.826632653061224	0.0520950040844277\\
0.828673469387755	0.052888722076033\\
0.830714285714286	0.0529910701506398\\
0.832755102040816	0.052481288287232\\
0.834795918367347	0.0515088151818999\\
0.836836734693878	0.0502341122262366\\
0.838877551020408	0.0487883123371766\\
0.840918367346939	0.0472337898908841\\
0.847040816326531	0.0423497630674664\\
0.849081632653061	0.0410103263210492\\
0.851122448979592	0.0401367466507601\\
0.853163265306122	0.0399002532502265\\
0.855204081632653	0.0403930542311709\\
0.857244897959184	0.0415785862514282\\
0.859285714285714	0.0432843270953922\\
0.861326530612245	0.0452805981884774\\
0.863367346938775	0.0474065664731477\\
0.869489795918367	0.0539652375668722\\
0.871530612244898	0.0560225853274624\\
0.873571428571429	0.0578869132333166\\
0.875612244897959	0.0594546649686852\\
0.87765306122449	0.0605535855859863\\
0.87969387755102	0.0609997196828665\\
0.881734693877551	0.0607145169265023\\
0.883775510204082	0.0597510936139095\\
0.885816326530612	0.0582068538250865\\
0.887857142857143	0.0561644101204014\\
0.889897959183673	0.0537143568106627\\
0.891938775510204	0.0509638575755689\\
0.893979591836735	0.0479905998495621\\
0.898061224489796	0.0417155938556942\\
0.900102040816326	0.0387759909880847\\
0.902142857142857	0.0363674705731408\\
0.904183673469388	0.0348096532661702\\
0.906224489795918	0.0342918523000479\\
0.908265306122449	0.0347741861503155\\
0.91030612244898	0.0360193015111224\\
0.91234693877551	0.03775545334837\\
0.914387755102041	0.0397981826223633\\
0.916428571428571	0.0420426414638119\\
0.918469387755102	0.0444275048646282\\
0.922551020408163	0.0493370832663871\\
0.924591836734694	0.0515524547840304\\
0.926632653061224	0.0533769602771171\\
0.928673469387755	0.0546333363470258\\
0.930714285714286	0.0551814499095258\\
0.932755102040816	0.0550497213562204\\
0.934795918367347	0.0543594904492091\\
0.936836734693878	0.0532532441691606\\
0.938877551020408	0.0518933658969918\\
0.940918367346939	0.0503806191112801\\
0.942959183673469	0.0487505453363297\\
0.947040816326531	0.0453250713489399\\
0.949081632653061	0.0437957760686128\\
0.951122448979592	0.0426480720008608\\
0.953163265306122	0.0420370061714566\\
0.955204081632653	0.0421042501027805\\
0.957244897959184	0.0429186229775927\\
0.959285714285714	0.0443601670113216\\
0.961326530612245	0.0462475206789589\\
0.963367346938775	0.0484007600860774\\
0.971530612244898	0.0574626998690502\\
0.973571428571429	0.0595192166678419\\
0.975612244897959	0.0613216126920081\\
0.97765306122449	0.0627382643182456\\
0.97969387755102	0.0635689806826346\\
0.981734693877551	0.0636760724257757\\
0.983775510204082	0.0630583895631002\\
0.985816326530612	0.0617604078685879\\
0.987857142857143	0.0598506636650438\\
0.989897959183674	0.0574332026484655\\
0.991938775510204	0.0546375050269192\\
0.993979591836735	0.0515773958381807\\
0.998061224489796	0.0450607748758793\\
};
\addlegendentry{ABC$_W^{1/2}$ adaptive}

\addplot [color=mycolor1, dashed, line width=2.0pt]
  table[row sep=crcr]{%
0.00214285714285711	0\\
0.104183673469388	3.09745474589285e-05\\
0.114387755102041	0.000153286872072278\\
0.138877551020408	0.000587671064360085\\
0.151122448979592	0.00088222627119805\\
0.159285714285714	0.0011484594517831\\
0.169489795918367	0.00157636056317689\\
0.17969387755102	0.00209146802563265\\
0.185816326530612	0.00247253526162949\\
0.191938775510204	0.00293786977886024\\
0.198061224489796	0.00350468742929078\\
0.21030612244898	0.00474331140186113\\
0.214387755102041	0.0050587797194126\\
0.222551020408163	0.00554601412633038\\
0.226632653061224	0.0058096171197175\\
0.230714285714286	0.00614208307378983\\
0.234795918367347	0.00656842017866488\\
0.238877551020408	0.00710260262814899\\
0.242959183673469	0.00774913101384234\\
0.247040816326531	0.00850040690525511\\
0.253163265306122	0.00976191641280588\\
0.259285714285714	0.011021784795583\\
0.263367346938776	0.0117545260163957\\
0.267448979591837	0.0123609787107651\\
0.271530612244898	0.0128633244283374\\
0.27765306122449	0.0135881172538257\\
0.281734693877551	0.0142002789090763\\
0.283775510204082	0.0145808738419653\\
0.285816326530612	0.0150227326803378\\
0.287857142857143	0.0155307310343108\\
0.289897959183673	0.0161059959033748\\
0.291938775510204	0.0167455356228151\\
0.293979591836735	0.0174420047202666\\
0.298061224489796	0.0189551050501916\\
0.302142857142857	0.0205076580360035\\
0.304183673469388	0.0212445112976372\\
0.306224489795918	0.021925430380541\\
0.308265306122449	0.0225310285123802\\
0.31030612244898	0.0230469301105034\\
0.31234693877551	0.0234659417790283\\
0.314387755102041	0.0237898030290863\\
0.316428571428572	0.0240302205581459\\
0.320510204081633	0.0243564599718707\\
0.322551020408163	0.0245098765859674\\
0.324591836734694	0.0247094230317072\\
0.326632653061224	0.0249944639940577\\
0.328673469387755	0.0253992372787543\\
0.330714285714286	0.0259490414349743\\
0.332755102040816	0.0266574630110129\\
0.334795918367347	0.0275249383881944\\
0.336836734693878	0.0285391044397588\\
0.338877551020408	0.0296761323926025\\
0.342959183673469	0.0321791323175593\\
0.347040816326531	0.0346980424176758\\
0.349081632653061	0.0358455616492186\\
0.351122448979592	0.0368569057997387\\
0.353163265306123	0.0376913738634942\\
0.355204081632653	0.0383161542589999\\
0.357244897959184	0.038709852065268\\
0.359285714285714	0.0388659259115306\\
0.361326530612245	0.0387959425520096\\
0.363367346938776	0.038532190440462\\
0.365408163265306	0.0381285399137925\\
0.369489795918367	0.0372205469304017\\
0.371530612244898	0.0369128017672415\\
0.373571428571429	0.0368410551395665\\
0.375612244897959	0.0370957738400127\\
0.37765306122449	0.0377412944781051\\
0.37969387755102	0.0388038026531589\\
0.381734693877551	0.0402692722030824\\
0.383775510204082	0.0420854984810962\\
0.385816326530612	0.0441708105700379\\
0.389897959183674	0.0487447004442552\\
0.391938775510204	0.0510206103336667\\
0.393979591836735	0.0531531730169347\\
0.396020408163265	0.055048502241276\\
0.398061224489796	0.0566207380675271\\
0.400102040816327	0.0577930279723716\\
0.402142857142857	0.0584982701417336\\
0.404183673469388	0.0586805456884906\\
0.406224489795918	0.0582990321005215\\
0.408265306122449	0.0573323061619163\\
0.41030612244898	0.0557837259629139\\
0.41234693877551	0.0536916863611422\\
0.414387755102041	0.051138594247977\\
0.416428571428571	0.0482619566597678\\
0.418469387755102	0.0452679339379445\\
0.420510204081633	0.04244041524022\\
0.422551020408163	0.0401359482772775\\
0.424591836734694	0.038744394429424\\
0.426632653061225	0.0385984774856324\\
0.428673469387755	0.0398545179553635\\
0.430714285714286	0.0424331685481072\\
0.432755102040816	0.0460696482855274\\
0.434795918367347	0.0504213593658527\\
0.440918367346939	0.0646250225885154\\
0.442959183673469	0.0689185265698993\\
0.445	0.0726805494838213\\
0.447040816326531	0.0757723783902081\\
0.449081632653061	0.078082344971579\\
0.451122448979592	0.0795246738757068\\
0.453163265306122	0.0800414912686231\\
0.455204081632653	0.0796061827215059\\
0.457244897959184	0.0782252424816634\\
0.459285714285714	0.0759453997139782\\
0.461326530612245	0.072858307407932\\
0.463367346938776	0.0691083064580372\\
0.465408163265306	0.064901330666883\\
0.467448979591837	0.0605161612095633\\
0.469489795918367	0.0563162729514818\\
0.471530612244898	0.0527467092935023\\
0.473571428571429	0.0502944098062571\\
0.475612244897959	0.0493879919150841\\
0.47765306122449	0.0502540809586608\\
0.47969387755102	0.0528131891844237\\
0.481734693877551	0.0567199719459791\\
0.483775510204082	0.0615058743696378\\
0.487857142857143	0.0719271205737186\\
0.489897959183674	0.0768712862549313\\
0.491938775510204	0.0813139790622864\\
0.493979591836735	0.0850920785609368\\
0.496020408163265	0.0880867718903966\\
0.498061224489796	0.0902106224708131\\
0.500102040816326	0.0913992028762827\\
0.502142857142857	0.0916041587637192\\
0.504183673469388	0.0907945089240776\\
0.506224489795918	0.0889562806110912\\
0.508265306122449	0.0860978760825817\\
0.51030612244898	0.0822624657619024\\
0.51234693877551	0.0775352407040923\\
0.514387755102041	0.0720563441187109\\
0.516428571428571	0.0660406661516113\\
0.518469387755102	0.0597980372253739\\
0.520510204081633	0.0537624292884443\\
0.522551020408163	0.0485244251716983\\
0.524591836734694	0.0448115450157893\\
0.526632653061225	0.0433235807497735\\
0.528673469387755	0.0444104058234818\\
0.530714285714286	0.0478551817694481\\
0.532755102040816	0.0530303519545571\\
0.534795918367347	0.0592281633977968\\
0.538877551020408	0.0724613403099711\\
0.540918367346939	0.0787232228625241\\
0.542959183673469	0.0843977194493788\\
0.545	0.089301513425068\\
0.547040816326531	0.0932955226512195\\
0.549081632653061	0.0962761876136307\\
0.551122448979592	0.0981727606953733\\
0.553163265306122	0.098945164147917\\
0.555204081632653	0.0985869974962156\\
0.557244897959184	0.0971287755155859\\
0.559285714285714	0.0946393914885179\\
0.561326530612245	0.0912313119094317\\
0.563367346938776	0.0870711353015557\\
0.565408163265306	0.0823866571994215\\
0.567448979591837	0.0774704214201188\\
0.569489795918367	0.072682245424753\\
0.571530612244898	0.0684431136824185\\
0.573571428571429	0.0652024603624223\\
0.575612244897959	0.0633695821660013\\
0.57765306122449	0.0632086030391086\\
0.57969387755102	0.0647364310017109\\
0.581734693877551	0.0677061600919316\\
0.583775510204082	0.0717099050844314\\
0.585816326530612	0.0763028976901592\\
0.587857142857143	0.0810727834138488\\
0.589897959183673	0.0856769527933754\\
0.591938775510204	0.0898538979800843\\
0.593979591836735	0.0934172943388031\\
0.596020408163265	0.0962371771365262\\
0.598061224489796	0.0982208909164823\\
0.600102040816327	0.0992983132922033\\
0.602142857142857	0.099417427770251\\
0.604183673469388	0.0985488709632372\\
0.606224489795918	0.0966929233762635\\
0.608265306122449	0.0938776995358559\\
0.61030612244898	0.0901610205364414\\
0.61234693877551	0.0856434699059883\\
0.614387755102041	0.0804880088278096\\
0.618469387755102	0.0693033426713339\\
0.620510204081633	0.0640191474336115\\
0.622551020408163	0.0595903153512523\\
0.624591836734694	0.0565601253442579\\
0.626632653061225	0.0553952693898438\\
0.628673469387755	0.0563070319049993\\
0.630714285714286	0.0591393481591486\\
0.632755102040816	0.0634587607033129\\
0.634795918367347	0.068741340975534\\
0.640918367346939	0.0859312842314642\\
0.642959183673469	0.0910416310115634\\
0.645	0.0954774927457646\\
0.647040816326531	0.0990953551706397\\
0.649081632653061	0.101789841826972\\
0.651122448979592	0.103490497848631\\
0.653163265306122	0.10416170884773\\
0.655204081632653	0.103799761678778\\
0.657244897959184	0.102434208957614\\
0.659285714285714	0.100135381659957\\
0.661326530612245	0.0970228389227552\\
0.663367346938776	0.0932746614511248\\
0.665408163265306	0.0891279788680702\\
0.667448979591837	0.0848713696708946\\
0.669489795918367	0.0808420516799734\\
0.671530612244898	0.0774218530216507\\
0.673571428571429	0.0749967158321611\\
0.675612244897959	0.0738682057531262\\
0.67765306122449	0.0741761351449983\\
0.67969387755102	0.0758767805690983\\
0.681734693877551	0.078753917403267\\
0.683775510204082	0.0824632036108459\\
0.687857142857143	0.0908079239555996\\
0.689897959183673	0.094774393618731\\
0.691938775510204	0.0982879278083905\\
0.693979591836735	0.101189447777534\\
0.696020408163265	0.103359377567662\\
0.698061224489796	0.104707419645486\\
0.700102040816327	0.105176458694812\\
0.702142857142857	0.104737477440885\\
0.704183673469388	0.103372315385621\\
0.706224489795918	0.101083743360095\\
0.708265306122449	0.0979183796652051\\
0.71030612244898	0.0939659303252208\\
0.71234693877551	0.0893523842437685\\
0.714387755102041	0.0842560185852003\\
0.718469387755102	0.0736919856428255\\
0.720510204081633	0.0689942312419827\\
0.722551020408163	0.0653467597238152\\
0.724591836734694	0.0632230916196894\\
0.726632653061224	0.0629179170049325\\
0.728673469387755	0.0644604294747914\\
0.730714285714286	0.0676158195491483\\
0.732755102040816	0.0719860688131841\\
0.734795918367347	0.0771238877185013\\
0.738877551020408	0.0880796845926622\\
0.740918367346939	0.0932697460827112\\
0.742959183673469	0.0979625456278298\\
0.745	0.101988069640586\\
0.747040816326531	0.105217019295693\\
0.749081632653061	0.107557011937342\\
0.751122448979592	0.108950625516962\\
0.753163265306122	0.109370890672456\\
0.755204081632653	0.108818835788916\\
0.757244897959184	0.107324861527835\\
0.759285714285714	0.104961387834497\\
0.761326530612245	0.101855744220093\\
0.763367346938776	0.098196736264154\\
0.767448979591837	0.09023586725507\\
0.769489795918367	0.0865459664992391\\
0.771530612244898	0.0835093681946752\\
0.773571428571429	0.0814791056649313\\
0.775612244897959	0.0807323191581298\\
0.77765306122449	0.0813725600875783\\
0.77969387755102	0.0833213004032295\\
0.781734693877551	0.0863528054701128\\
0.783775510204082	0.0901172738554413\\
0.787857142857143	0.0982346408022891\\
0.789897959183673	0.101956547318007\\
0.791938775510204	0.105195781117788\\
0.793979591836735	0.107812696624434\\
0.796020408163265	0.109679696852177\\
0.798061224489796	0.110690840035143\\
0.800102040816327	0.110792602936564\\
0.802142857142857	0.109988375350697\\
0.804183673469388	0.108299181328785\\
0.806224489795918	0.105729278237544\\
0.808265306122449	0.102292897360402\\
0.81030612244898	0.0980547591102219\\
0.81234693877551	0.0931460801945223\\
0.814387755102041	0.0877694407002073\\
0.816428571428571	0.0821920927576929\\
0.818469387755102	0.0767478173952688\\
0.820510204081633	0.0718703566243954\\
0.822551020408163	0.0680974120720077\\
0.824591836734694	0.0659392230176772\\
0.826632653061224	0.0656862351506351\\
0.828673469387755	0.0673249151375014\\
0.830714285714286	0.0705880627028366\\
0.832755102040816	0.0750635824398329\\
0.834795918367347	0.0803015941752788\\
0.838877551020408	0.0914587924007476\\
0.840918367346939	0.0967334040227368\\
0.842959183673469	0.101496669842966\\
0.845	0.105581025279024\\
0.847040816326531	0.108861623957525\\
0.849081632653061	0.111249754233806\\
0.851122448979592	0.112682318454009\\
0.853163265306122	0.113129531795621\\
0.855204081632653	0.112595393936608\\
0.857244897959184	0.111114508557983\\
0.859285714285714	0.108760085724559\\
0.861326530612245	0.105648038163852\\
0.863367346938775	0.101954137821146\\
0.867448979591837	0.0938371669946108\\
0.869489795918367	0.090037790792777\\
0.871530612244898	0.0868720138767002\\
0.873571428571429	0.084711120049788\\
0.875612244897959	0.0838633787961403\\
0.87765306122449	0.0844512685936631\\
0.87969387755102	0.0863814673841723\\
0.881734693877551	0.0894165415706827\\
0.883775510204082	0.0932252374414555\\
0.887857142857143	0.101597301212761\\
0.889897959183673	0.105448876119623\\
0.891938775510204	0.108766942150889\\
0.893979591836735	0.11141892705758\\
0.896020408163265	0.113304995255764\\
0.898061224489796	0.11435139692808\\
0.900102040816326	0.114505528292377\\
0.902142857142857	0.113754510912534\\
0.904183673469388	0.112114887737301\\
0.906224489795918	0.109603689085013\\
0.908265306122449	0.106244411946587\\
0.91030612244898	0.102086281318\\
0.91234693877551	0.0972213980458415\\
0.914387755102041	0.091805990784528\\
0.918469387755102	0.0803421186274289\\
0.920510204081633	0.075075247509775\\
0.922551020408163	0.0708266277340085\\
0.924591836734694	0.0681354924790342\\
0.926632653061224	0.0673413385111554\\
0.928673469387755	0.068470326029731\\
0.930714285714286	0.0712873785378541\\
0.932755102040816	0.0754414771505875\\
0.934795918367347	0.0805180658132842\\
0.940918367346939	0.0972488262377965\\
0.942959183673469	0.102241210422162\\
0.945	0.106566393099335\\
0.947040816326531	0.11010341213743\\
0.949081632653061	0.112782097529975\\
0.951122448979592	0.114534026475958\\
0.953163265306122	0.115302739024344\\
0.955204081632653	0.115086933027369\\
0.957244897959184	0.11392101778346\\
0.959285714285714	0.111849010852262\\
0.961326530612245	0.108977229487004\\
0.963367346938775	0.105477778646345\\
0.965408163265306	0.101577346472053\\
0.967448979591837	0.097553649451923\\
0.969489795918367	0.0937181031615339\\
0.971530612244898	0.0904201119836308\\
0.973571428571429	0.0880273390323011\\
0.975612244897959	0.0868618625465827\\
0.97765306122449	0.0871067493524504\\
0.97969387755102	0.0887318168717329\\
0.981734693877551	0.0915358673794542\\
0.983775510204082	0.0952055180503717\\
0.987857142857143	0.10353825812817\\
0.989897959183674	0.107445566442856\\
0.991938775510204	0.110841488069718\\
0.993979591836735	0.113602115594535\\
0.996020408163265	0.115643267359182\\
0.998061224489796	0.116889596141874\\
};
\addlegendentry{ABC$_W^{1/2}$}

\end{axis}
\end{tikzpicture}%

%% file: images/Transducer/Cost_v.tex
%
%
\definecolor{mycolor1}{rgb}{1.00000,0.60000,0.00000}%
\begin{tikzpicture}[scale = 0.52, font=\huge]

\begin{axis}[%
width=4.602in,
height=3.82in,
at={(0.772in,0.516in)},
scale only axis,
xmin=0,
xmax=1,
xlabel style={font=\huge\color{white!15!black}},
y tick label style={
        /pgf/number format/.cd,
            fixed,
            fixed zerofill,
            precision=2,
        /tikz/.cd
            },
xlabel={relative simulation time $t/T$},xtick={0, 0.2, 0.4, 0.6, 0.8, 1},
ymin=0,
ymax=0.12,  ytick={0.02, 0.04, 0.06, 0.08, 0.1, 0.12},
ylabel style={font=\huge\color{white!15!black}},
ylabel={},
axis background/.style={fill=white},
xmajorgrids,
ymajorgrids,
legend style={at={(0.03,0.97)}, anchor=north west, legend cell align=left, align=left, draw=white!15!black}
]
\addplot [color=mycolor1, line width=2.0pt]
  table[row sep=crcr]{%
0.00214285714285711	0\\
0.0980612244897959	3.40778019684018e-05\\
0.108265306122449	0.000139090589205471\\
0.118469387755102	0.000317282202808133\\
0.132755102040816	0.000666066434471535\\
0.142959183673469	0.000878333007428234\\
0.151122448979592	0.00106634514608983\\
0.157244897959184	0.00128897272323225\\
0.163367346938776	0.00158795363598863\\
0.169489795918367	0.00196629232223733\\
0.175612244897959	0.00242782027059485\\
0.187857142857143	0.0034047058152572\\
0.191938775510204	0.00360992075861577\\
0.198061224489796	0.00377518160738\\
0.204183673469388	0.00393594349419246\\
0.208265306122449	0.00415743468281726\\
0.21234693877551	0.00448722774071719\\
0.216428571428571	0.00492920988297751\\
0.220510204081633	0.00549698037969071\\
0.224591836734694	0.00616153047728119\\
0.234795918367347	0.00787956112795973\\
0.238877551020408	0.00843629417983427\\
0.242959183673469	0.00887372305568179\\
0.247040816326531	0.00920021691267958\\
0.259285714285714	0.0100752566505223\\
0.261326530612245	0.0102899965672788\\
0.263367346938776	0.0105480694306312\\
0.265408163265306	0.0108540216647925\\
0.267448979591837	0.0112128241447735\\
0.269489795918367	0.0116301832094613\\
0.271530612244898	0.0121033895543428\\
0.275612244897959	0.0131936732593355\\
0.285816326530612	0.0161385134831343\\
0.287857142857143	0.0165686410980493\\
0.289897959183673	0.0168820931470963\\
0.291938775510204	0.0170782388255694\\
0.293979591836735	0.0171492734484079\\
0.296020408163265	0.0171008399181598\\
0.298061224489796	0.0169556471249505\\
0.300102040816327	0.0167315924162329\\
0.304183673469388	0.0162079481499585\\
0.306224489795918	0.0160123598182137\\
0.308265306122449	0.0159134002657573\\
0.31030612244898	0.0159515750786442\\
0.31234693877551	0.0161655247510369\\
0.314387755102041	0.016544266505533\\
0.316428571428572	0.0170829980345228\\
0.320510204081633	0.0184384910298766\\
0.324591836734694	0.0198128153774267\\
0.326632653061224	0.0203935830735565\\
0.328673469387755	0.0208313952612142\\
0.330714285714286	0.0210656560411018\\
0.332755102040816	0.0210909995748839\\
0.334795918367347	0.0208862779072727\\
0.336836734693878	0.0204977356024841\\
0.340918367346939	0.0195419321735099\\
0.342959183673469	0.0191664228736178\\
0.345	0.0189618514831392\\
0.347040816326531	0.0190541840126819\\
0.349081632653061	0.0193923625631125\\
0.351122448979592	0.0200924792048544\\
0.353163265306123	0.0210972171596385\\
0.357244897959184	0.0232661585857318\\
0.359285714285714	0.02421327906852\\
0.361326530612245	0.0249651280883583\\
0.363367346938776	0.0254801672120267\\
0.365408163265306	0.0255554559649508\\
0.37765306122449	0.0247281532961003\\
0.385816326530612	0.0251843990014\\
0.387857142857143	0.0260462720126872\\
0.389897959183674	0.0276291507145315\\
0.391938775510204	0.0292712539378492\\
0.393979591836735	0.0307435280317935\\
0.396020408163265	0.0321401458602901\\
0.398061224489796	0.0334261110138732\\
0.400102040816327	0.0345354942248522\\
0.402142857142857	0.0349644930962728\\
0.404183673469388	0.0351397472438585\\
0.408265306122449	0.0357866272204769\\
0.41030612244898	0.0357145849192974\\
0.41234693877551	0.0350980263462154\\
0.414387755102041	0.0339806136430598\\
0.416428571428571	0.0323774218492398\\
0.418469387755102	0.0305768647139282\\
0.420510204081633	0.0290316039175383\\
0.422551020408163	0.0277340250105117\\
0.424591836734694	0.0272067747344742\\
0.426632653061225	0.0277264157251043\\
0.428673469387755	0.0287888265144536\\
0.432755102040816	0.0310992104000001\\
0.434795918367347	0.0322186411809824\\
0.436836734693878	0.0335166484360736\\
0.438877551020408	0.0351185077362413\\
0.440918367346939	0.0369881412355727\\
0.442959183673469	0.0390981183919257\\
0.445	0.0411144115887168\\
0.447040816326531	0.0427374319872349\\
0.449081632653061	0.0437151005256845\\
0.451122448979592	0.0443501912869014\\
0.453163265306122	0.0443500448619792\\
0.455204081632653	0.0436696983139095\\
0.457244897959184	0.0427937265699176\\
0.459285714285714	0.0421385956432614\\
0.463367346938776	0.0412855097140618\\
0.465408163265306	0.0403652302605502\\
0.467448979591837	0.0386545343037242\\
0.469489795918367	0.0361946661149033\\
0.471530612244898	0.033188756519742\\
0.473571428571429	0.0306876333477508\\
0.475612244897959	0.0300104716363573\\
0.47765306122449	0.0314507694322862\\
0.47969387755102	0.0339695267064669\\
0.481734693877551	0.0367486512770594\\
0.483775510204082	0.0390521058576729\\
0.485816326530612	0.0407334056917957\\
0.487857142857143	0.0418900894236983\\
0.489897959183674	0.0428595495296347\\
0.491938775510204	0.0441405058928425\\
0.493979591836735	0.0463094087180523\\
0.496020408163265	0.0491779993001331\\
0.498061224489796	0.0516703512830454\\
0.500102040816326	0.0530584363104051\\
0.502142857142857	0.0528643752388823\\
0.504183673469388	0.0515012547958156\\
0.508265306122449	0.0475774706342963\\
0.51030612244898	0.0461532766938225\\
0.51234693877551	0.0454743549639208\\
0.514387755102041	0.0450924527810166\\
0.516428571428571	0.0441284182835213\\
0.518469387755102	0.0419179377211276\\
0.520510204081633	0.0383717277836748\\
0.522551020408163	0.0333970956945419\\
0.524591836734694	0.0292053682249052\\
0.526632653061225	0.0283862383276534\\
0.528673469387755	0.0309925141239484\\
0.530714285714286	0.0350103629588202\\
0.532755102040816	0.0388566302738852\\
0.534795918367347	0.041734714473293\\
0.536836734693878	0.0434829974449071\\
0.538877551020408	0.0444850181801071\\
0.540918367346939	0.0453762578589997\\
0.542959183673469	0.0470157162060458\\
0.547040816326531	0.0521736653531875\\
0.549081632653061	0.0540411662608938\\
0.551122448979592	0.0548846389444417\\
0.553163265306122	0.0544945306503635\\
0.555204081632653	0.0528336278014624\\
0.557244897959184	0.0504737682327085\\
0.559285714285714	0.0482396695768786\\
0.561326530612245	0.0468062505260122\\
0.563367346938776	0.0462721486533119\\
0.565408163265306	0.0462925718624861\\
0.567448979591837	0.0459784128703766\\
0.569489795918367	0.0447304756697703\\
0.571530612244898	0.0419629763389635\\
0.573571428571429	0.0375896172633182\\
0.575612244897959	0.033951487834964\\
0.57765306122449	0.0334133818964104\\
0.57969387755102	0.0361941214599451\\
0.581734693877551	0.0405510232016031\\
0.583775510204082	0.0446966930702363\\
0.585816326530612	0.0474638291578983\\
0.587857142857143	0.0486554841162156\\
0.589897959183673	0.0488880851282429\\
0.591938775510204	0.0492116426920971\\
0.593979591836735	0.0508018038815684\\
0.596020408163265	0.0536557984238132\\
0.598061224489796	0.0566854215833441\\
0.600102040816327	0.0586548493259612\\
0.602142857142857	0.0588789304539947\\
0.604183673469388	0.0569802491620565\\
0.606224489795918	0.0535260501809952\\
0.608265306122449	0.0498417636463272\\
0.61030612244898	0.0471356409121755\\
0.61234693877551	0.045855450576697\\
0.614387755102041	0.0456758930294959\\
0.616428571428571	0.0456826390527069\\
0.618469387755102	0.0448073841386514\\
0.620510204081633	0.0425750449868563\\
0.622551020408163	0.0390609309629049\\
0.624591836734694	0.0347047970804052\\
0.626632653061225	0.0321778110980405\\
0.628673469387755	0.0334757338564082\\
0.630714285714286	0.0376112997584335\\
0.632755102040816	0.0422109246931184\\
0.634795918367347	0.0456535043493733\\
0.636836734693878	0.0475243272715955\\
0.638877551020408	0.0482079435185411\\
0.640918367346939	0.0483788434975194\\
0.642959183673469	0.0491190810027814\\
0.645	0.0508414579418592\\
0.647040816326531	0.0531520885226006\\
0.649081632653061	0.0552439449403824\\
0.651122448979592	0.0564204616865182\\
0.653163265306122	0.0563912843512557\\
0.655204081632653	0.0547124780898193\\
0.657244897959184	0.0517638910609689\\
0.659285714285714	0.0489553653719156\\
0.661326530612245	0.047268839232298\\
0.663367346938776	0.0468115217440926\\
0.665408163265306	0.0471242918784152\\
0.667448979591837	0.0474839908338892\\
0.669489795918367	0.0469309437432045\\
0.671530612244898	0.0450508596752301\\
0.673571428571429	0.0414438378978007\\
0.675612244897959	0.0372253453430219\\
0.67765306122449	0.0355363386884781\\
0.67969387755102	0.0378474110921527\\
0.681734693877551	0.0423247124377091\\
0.683775510204082	0.0464795617465515\\
0.685816326530612	0.0491971298752829\\
0.687857142857143	0.0503927552630615\\
0.689897959183673	0.0504644520083525\\
0.691938775510204	0.050443743743889\\
0.693979591836735	0.0514090522790256\\
0.696020408163265	0.0536752170865765\\
0.698061224489796	0.0564943872159406\\
0.700102040816327	0.0586130652785316\\
0.702142857142857	0.0592911944863348\\
0.704183673469388	0.0579068259547114\\
0.706224489795918	0.054816915526806\\
0.708265306122449	0.0512487657046148\\
0.71030612244898	0.0484881638762983\\
0.71234693877551	0.0472046719821404\\
0.714387755102041	0.0470916527884891\\
0.716428571428571	0.0471771690823533\\
0.718469387755102	0.0466516969979539\\
0.720510204081633	0.0446826442192516\\
0.722551020408163	0.0415914161002938\\
0.724591836734694	0.0377078922624784\\
0.726632653061224	0.0347156904813624\\
0.728673469387755	0.034758599136983\\
0.730714285714286	0.0378565474360131\\
0.732755102040816	0.0419167464352511\\
0.734795918367347	0.0451916078685195\\
0.736836734693878	0.0470500307476289\\
0.738877551020408	0.0478031627169646\\
0.740918367346939	0.0481277783737981\\
0.742959183673469	0.0487389102726352\\
0.745	0.0503377757418644\\
0.747040816326531	0.0526464999222156\\
0.749081632653061	0.0548402942274504\\
0.751122448979592	0.0564634143752311\\
0.753163265306122	0.0571223676486392\\
0.755204081632653	0.0564554839637906\\
0.757244897959184	0.0543986223986883\\
0.759285714285714	0.0519430959297351\\
0.761326530612245	0.0502308724773295\\
0.763367346938776	0.0495305875071722\\
0.765408163265306	0.0496449878302793\\
0.767448979591837	0.0499494532368645\\
0.769489795918367	0.0493268059920177\\
0.771530612244898	0.0475941998346564\\
0.773571428571429	0.0447798570443015\\
0.775612244897959	0.0406631583508482\\
0.77765306122449	0.0376535711186498\\
0.77969387755102	0.037935817150892\\
0.781734693877551	0.0409401342801939\\
0.783775510204082	0.0446488196219563\\
0.785816326530612	0.0476204713549109\\
0.787857142857143	0.0494139045382912\\
0.789897959183673	0.0502430733165142\\
0.791938775510204	0.0506674550607173\\
0.793979591836735	0.0518498163525609\\
0.796020408163265	0.0540384868882201\\
0.798061224489796	0.0570887504367024\\
0.800102040816327	0.0597229392388524\\
0.802142857142857	0.0611768584711264\\
0.804183673469388	0.061009525946369\\
0.806224489795918	0.0590991768420818\\
0.808265306122449	0.056465924813889\\
0.81030612244898	0.0539738242367876\\
0.81234693877551	0.0526250767524072\\
0.814387755102041	0.0521784213131516\\
0.816428571428571	0.0518056428230026\\
0.818469387755102	0.0512752615099282\\
0.820510204081633	0.0494527194937646\\
0.822551020408163	0.0467951572250257\\
0.824591836734694	0.0433815525596526\\
0.826632653061224	0.0398193632334004\\
0.828673469387755	0.0380650756221305\\
0.830714285714286	0.039575033832894\\
0.834795918367347	0.0453981616853022\\
0.836836734693878	0.0477662945784463\\
0.838877551020408	0.0488364993318103\\
0.842959183673469	0.0504331877726503\\
0.845	0.0520066976766687\\
0.847040816326531	0.0545777408504526\\
0.849081632653061	0.0570551863764264\\
0.851122448979592	0.0588697075858632\\
0.853163265306122	0.0600469342912154\\
0.855204081632653	0.0601707480032039\\
0.857244897959184	0.0590716705260879\\
0.859285714285714	0.05709432657366\\
0.861326530612245	0.0556022121893546\\
0.863367346938775	0.0551981788764445\\
0.865408163265306	0.0551888629858268\\
0.867448979591837	0.0552698689539797\\
0.869489795918367	0.054790711300298\\
0.871530612244898	0.0533660940171439\\
0.873571428571429	0.0510181921488941\\
0.87765306122449	0.0436017900507489\\
0.87969387755102	0.042367548715682\\
0.881734693877551	0.0441692764900181\\
0.883775510204082	0.046889981030007\\
0.885816326530612	0.0493691197477543\\
0.887857142857143	0.0515691166397251\\
0.889897959183673	0.05320442745059\\
0.891938775510204	0.0540858924973823\\
0.893979591836735	0.0551751129524408\\
0.896020408163265	0.0572369592335077\\
0.900102040816326	0.0625341883614033\\
0.902142857142857	0.0641993342457029\\
0.904183673469388	0.0647577684365135\\
0.906224489795918	0.0639616884653116\\
0.91030612244898	0.0595145096076386\\
0.91234693877551	0.0580744134271657\\
0.914387755102041	0.0573472189657167\\
0.918469387755102	0.0566059886924576\\
0.920510204081633	0.0550177158584289\\
0.922551020408163	0.0521942620877437\\
0.924591836734694	0.0498349323907985\\
0.926632653061224	0.0470249041901823\\
0.928673469387755	0.0441328691690758\\
0.930714285714286	0.0446389467728762\\
0.932755102040816	0.0467543453922414\\
0.934795918367347	0.0487868094205057\\
0.936836734693878	0.0512995620747083\\
0.938877551020408	0.0530814302124126\\
0.942959183673469	0.0553511637562324\\
0.945	0.0569137468589896\\
0.949081632653061	0.0613561020685645\\
0.951122448979592	0.0628553313088448\\
0.953163265306122	0.0642788127400329\\
0.955204081632653	0.0653411160414041\\
0.957244897959184	0.0644980896490815\\
0.961326530612245	0.0621443883367438\\
0.963367346938775	0.0610583468414633\\
0.965408163265306	0.0605387872890953\\
0.967448979591837	0.0603287235913438\\
0.969489795918367	0.0601761971934184\\
0.971530612244898	0.0589662868588092\\
0.973571428571429	0.0568594954598993\\
0.975612244897959	0.0545413522110446\\
0.97765306122449	0.0509013017852382\\
0.97969387755102	0.0486795603572031\\
0.981734693877551	0.0495466532988313\\
0.983775510204082	0.0514260695182647\\
0.985816326530612	0.0538423928220955\\
0.987857142857143	0.0563303948145483\\
0.989897959183674	0.0583204811286506\\
0.993979591836735	0.0611054338358626\\
0.996020408163265	0.0627724771296516\\
0.998061224489796	0.0655344183443399\\
};
\addlegendentry{ABC$_W^{1/2}$ adaptive}

\addplot [color=mycolor1, dashed, line width=2.0pt]
  table[row sep=crcr]{%
0.00214285714285711	0\\
0.0980612244897959	3.40778019684018e-05\\
0.108265306122449	0.000139090589205471\\
0.118469387755102	0.000317282202808133\\
0.132755102040816	0.000666066434471535\\
0.142959183673469	0.000878333007428234\\
0.151122448979592	0.00106634514608983\\
0.157244897959184	0.00128897272323225\\
0.163367346938776	0.00158795363598863\\
0.169489795918367	0.00196629232223733\\
0.175612244897959	0.00242782027059485\\
0.187857142857143	0.0034047058152572\\
0.191938775510204	0.00360992075861577\\
0.198061224489796	0.00377518160738\\
0.204183673469388	0.00393594349419246\\
0.208265306122449	0.00415743468281726\\
0.21234693877551	0.00448722774071719\\
0.216428571428571	0.00491959630461813\\
0.220510204081633	0.00544479248428964\\
0.226632653061224	0.00636443023659083\\
0.234795918367347	0.00764007363577401\\
0.238877551020408	0.00819594656447631\\
0.242959183673469	0.00866504563876513\\
0.249081632653061	0.00924833124170976\\
0.255204081632653	0.00982467571777768\\
0.259285714285714	0.0102967999638548\\
0.263367346938776	0.0109037478751162\\
0.265408163265306	0.0112734658717214\\
0.267448979591837	0.0116936443920187\\
0.269489795918367	0.0121672579543706\\
0.271530612244898	0.0126940532302079\\
0.275612244897959	0.0138917124469143\\
0.285816326530612	0.0171278366414258\\
0.287857142857143	0.0176662608421407\\
0.289897959183673	0.0181269834584356\\
0.291938775510204	0.0185014883714437\\
0.293979591836735	0.018789203193614\\
0.296020408163265	0.0189979541643227\\
0.300102040816327	0.0192509249914855\\
0.304183673469388	0.0194700330711045\\
0.306224489795918	0.0196502660035243\\
0.308265306122449	0.0199221268618951\\
0.31030612244898	0.0203136613716957\\
0.31234693877551	0.0208440776010103\\
0.314387755102041	0.0215233522538189\\
0.316428571428572	0.0223481742072245\\
0.318469387755102	0.0233046500411649\\
0.320510204081633	0.0243672047922177\\
0.328673469387755	0.0289121513259862\\
0.330714285714286	0.0298973732749956\\
0.332755102040816	0.0307343270733426\\
0.334795918367347	0.0313903214703785\\
0.336836734693878	0.0318431241094835\\
0.338877551020408	0.0320842747846326\\
0.340918367346939	0.0321243922383952\\
0.342959183673469	0.0319919604966835\\
0.345	0.0317357685449107\\
0.349081632653061	0.0311399730831325\\
0.351122448979592	0.0309762434139694\\
0.353163265306123	0.0310250662331328\\
0.355204081632653	0.0313633491177139\\
0.357244897959184	0.0320463089539518\\
0.359285714285714	0.0330903155714766\\
0.361326530612245	0.0344832653833214\\
0.363367346938776	0.0361747862613929\\
0.365408163265306	0.0380988279076943\\
0.369489795918367	0.0423041214674272\\
0.371530612244898	0.0444076133006474\\
0.373571428571429	0.0463908033492748\\
0.375612244897959	0.048171587853444\\
0.37765306122449	0.0496676598279563\\
0.37969387755102	0.0508159687692379\\
0.381734693877551	0.0515590845305502\\
0.383775510204082	0.0518547111007861\\
0.385816326530612	0.0516785476260945\\
0.387857142857143	0.0510307263700154\\
0.389897959183674	0.049917509104382\\
0.391938775510204	0.0483804641578031\\
0.393979591836735	0.0464971140066751\\
0.398061224489796	0.0423473293592607\\
0.400102040816327	0.0405717534719118\\
0.402142857142857	0.0393742841766669\\
0.404183673469388	0.0389820211888789\\
0.406224489795918	0.0395499535006202\\
0.408265306122449	0.0410608849452447\\
0.41030612244898	0.0433982487128893\\
0.41234693877551	0.0463630424208026\\
0.414387755102041	0.0497092347883318\\
0.418469387755102	0.056716315624851\\
0.420510204081633	0.0599913260232413\\
0.422551020408163	0.0629093663633948\\
0.424591836734694	0.0653258758336414\\
0.426632653061225	0.0671282985514332\\
0.428673469387755	0.0682118118073671\\
0.430714285714286	0.0684953766829541\\
0.432755102040816	0.067923433337864\\
0.434795918367347	0.0664738883593969\\
0.436836734693878	0.0641468722214544\\
0.438877551020408	0.0610030235481537\\
0.440918367346939	0.0571920708713938\\
0.442959183673469	0.052915861125836\\
0.445	0.0485154728340007\\
0.447040816326531	0.0444561317666379\\
0.449081632653061	0.0413023661165089\\
0.451122448979592	0.0396569578361573\\
0.453163265306122	0.0399609691427963\\
0.455204081632653	0.0422372459664911\\
0.457244897959184	0.0461375465547412\\
0.459285714285714	0.0511519361441336\\
0.461326530612245	0.0567663661962298\\
0.463367346938776	0.0625636859988101\\
0.465408163265306	0.0682201836176753\\
0.467448979591837	0.0734807002975004\\
0.469489795918367	0.0781310653877962\\
0.471530612244898	0.0819872406744259\\
0.473571428571429	0.084889856619015\\
0.475612244897959	0.0867277740638162\\
0.47765306122449	0.0874256768405243\\
0.47969387755102	0.0869705697563397\\
0.481734693877551	0.0853988200686875\\
0.483775510204082	0.0827692970308507\\
0.485816326530612	0.0791319442935203\\
0.487857142857143	0.0745506727010665\\
0.489897959183674	0.0691651493903853\\
0.491938775510204	0.0632243514946812\\
0.493979591836735	0.0570621263585255\\
0.496020408163265	0.0511813411703613\\
0.498061224489796	0.0462081289818598\\
0.500102040816326	0.0427999602417659\\
0.502142857142857	0.0416102982701699\\
0.504183673469388	0.042822601873996\\
0.506224489795918	0.0461028024842257\\
0.508265306122449	0.0508468092287594\\
0.51030612244898	0.0563857437680831\\
0.514387755102041	0.0679387699088402\\
0.516428571428571	0.0732956273411458\\
0.518469387755102	0.0780970228236274\\
0.520510204081633	0.0821995456050283\\
0.522551020408163	0.0854704936392099\\
0.524591836734694	0.0878094047356205\\
0.526632653061225	0.0891372443204345\\
0.528673469387755	0.0893884512266323\\
0.530714285714286	0.0885204145836874\\
0.532755102040816	0.0865150885508952\\
0.534795918367347	0.0833889139933249\\
0.536836734693878	0.0792194479437398\\
0.538877551020408	0.0741232884698659\\
0.540918367346939	0.0682769157150204\\
0.545	0.0556972274353502\\
0.547040816326531	0.0499975523834283\\
0.549081632653061	0.0456644115274557\\
0.551122448979592	0.0434413627278336\\
0.553163265306122	0.0437837268557852\\
0.555204081632653	0.046639025176571\\
0.557244897959184	0.0514179174908259\\
0.559285714285714	0.0573540137071058\\
0.565408163265306	0.0767353306857399\\
0.567448979591837	0.0825190241217761\\
0.569489795918367	0.0875489456442812\\
0.571530612244898	0.0916387636708957\\
0.573571428571429	0.0946458541603722\\
0.575612244897959	0.0964624510435479\\
0.57765306122449	0.096996774482513\\
0.57969387755102	0.0962186907888302\\
0.581734693877551	0.0942290501459789\\
0.583775510204082	0.0911776506008632\\
0.585816326530612	0.0871451405339676\\
0.587857142857143	0.0822547451758283\\
0.589897959183673	0.0766672827806768\\
0.593979591836735	0.0644679068928407\\
0.596020408163265	0.0587185945627851\\
0.598061224489796	0.0539421057998833\\
0.600102040816327	0.0507317178372867\\
0.602142857142857	0.0495865425194377\\
0.604183673469388	0.0506882042182156\\
0.606224489795918	0.0537134174398001\\
0.608265306122449	0.0580806033800517\\
0.61030612244898	0.0632291866010891\\
0.614387755102041	0.0741041263234663\\
0.616428571428571	0.0791769444926506\\
0.618469387755102	0.0836904029024663\\
0.620510204081633	0.0874846877211831\\
0.622551020408163	0.0904616997261988\\
0.624591836734694	0.092529009987273\\
0.626632653061225	0.0935548544825917\\
0.628673469387755	0.0935001670966992\\
0.630714285714286	0.0923859192807297\\
0.632755102040816	0.0902424802025743\\
0.634795918367347	0.0870943338672696\\
0.636836734693878	0.083011525638417\\
0.638877551020408	0.0781031762766342\\
0.640918367346939	0.072587120433675\\
0.642959183673469	0.0668159526585516\\
0.645	0.0612214603887654\\
0.647040816326531	0.0563330277737971\\
0.649081632653061	0.0527499809974058\\
0.651122448979592	0.0510616193073927\\
0.653163265306122	0.0515945051374177\\
0.655204081632653	0.0542249542829151\\
0.657244897959184	0.0584536188138609\\
0.659285714285714	0.0637362653726827\\
0.661326530612245	0.0696656048184723\\
0.663367346938776	0.0758158865123806\\
0.665408163265306	0.0817525166833032\\
0.667448979591837	0.0871674000275275\\
0.669489795918367	0.0918669892742341\\
0.671530612244898	0.095704281366576\\
0.673571428571429	0.0984713418462916\\
0.675612244897959	0.100014958745488\\
0.67765306122449	0.100233812258437\\
0.67969387755102	0.099138962368424\\
0.681734693877551	0.0969210916555829\\
0.683775510204082	0.0936886137443992\\
0.685816326530612	0.0895340330697647\\
0.687857142857143	0.0846661508766186\\
0.689897959183673	0.0792824954309197\\
0.691938775510204	0.0736308115179396\\
0.693979591836735	0.0681094488824864\\
0.696020408163265	0.0631406336010745\\
0.698061224489796	0.0593207365863938\\
0.700102040816327	0.0571078996459264\\
0.702142857142857	0.0567798343684224\\
0.704183673469388	0.0583345566235941\\
0.706224489795918	0.061473543361327\\
0.708265306122449	0.0657250361318148\\
0.71030612244898	0.0706473960039417\\
0.71234693877551	0.0757947566871253\\
0.714387755102041	0.0808203235985567\\
0.716428571428571	0.0854789417715196\\
0.718469387755102	0.0896028528597479\\
0.720510204081633	0.0929767493792706\\
0.722551020408163	0.0954366151891346\\
0.724591836734694	0.0969130993117612\\
0.726632653061224	0.0974025177113937\\
0.728673469387755	0.0968880334482914\\
0.730714285714286	0.095401236995827\\
0.732755102040816	0.0929113391958243\\
0.734795918367347	0.0894905739227245\\
0.736836734693878	0.0853354848986129\\
0.738877551020408	0.0806035800034183\\
0.740918367346939	0.0753855895089829\\
0.742959183673469	0.0700269307620055\\
0.745	0.0650810225573749\\
0.747040816326531	0.061065220113239\\
0.749081632653061	0.0584124853107504\\
0.751122448979592	0.0574531726530951\\
0.753163265306122	0.0583413484249006\\
0.755204081632653	0.0609495398319689\\
0.757244897959184	0.0649822170927571\\
0.759285714285714	0.070110476387522\\
0.761326530612245	0.0758379109598933\\
0.763367346938776	0.0817367032512319\\
0.765408163265306	0.0874182872079715\\
0.767448979591837	0.0926800101797343\\
0.769489795918367	0.0971430105887069\\
0.771530612244898	0.10074380670208\\
0.773571428571429	0.103399399773664\\
0.775612244897959	0.104792618195993\\
0.77765306122449	0.104814914018641\\
0.77969387755102	0.103559173099756\\
0.781734693877551	0.101106345620945\\
0.783775510204082	0.0975018481366473\\
0.785816326530612	0.0930419280903593\\
0.787857142857143	0.0881797699592966\\
0.789897959183673	0.0831135416638503\\
0.793979591836735	0.0723356028775909\\
0.796020408163265	0.0672771714412523\\
0.798061224489796	0.063650451889367\\
0.800102040816327	0.0621332681561594\\
0.802142857142857	0.0623907046889645\\
0.804183673469388	0.0639987367450683\\
0.806224489795918	0.0667442381591563\\
0.808265306122449	0.0706372246407702\\
0.81030612244898	0.0754045957445417\\
0.814387755102041	0.085595964283909\\
0.816428571428571	0.0902338364905586\\
0.818469387755102	0.0943092447580895\\
0.820510204081633	0.0976035586753214\\
0.822551020408163	0.10000101359816\\
0.824591836734694	0.101328173794948\\
0.826632653061224	0.101602536152708\\
0.828673469387755	0.100968363852201\\
0.830714285714286	0.0994090201856489\\
0.832755102040816	0.0968734620153305\\
0.834795918367347	0.0934290550505098\\
0.836836734693878	0.0891630500457526\\
0.840918367346939	0.0790941107833191\\
0.842959183673469	0.0737998614267301\\
0.845	0.0689047971381017\\
0.847040816326531	0.0650332927544669\\
0.849081632653061	0.0623264001781957\\
0.851122448979592	0.0613841930600788\\
0.853163265306122	0.0622965463667636\\
0.855204081632653	0.064740529945569\\
0.857244897959184	0.0685663550747084\\
0.859285714285714	0.0733184698426883\\
0.861326530612245	0.0789416569220868\\
0.865408163265306	0.090737798549031\\
0.867448979591837	0.0960298594863647\\
0.869489795918367	0.100582357304318\\
0.871530612244898	0.104490435054058\\
0.873571428571429	0.107347604957447\\
0.875612244897959	0.108771828557374\\
0.87765306122449	0.108817969258993\\
0.87969387755102	0.107687240202062\\
0.881734693877551	0.105447015806446\\
0.883775510204082	0.102065014937845\\
0.885816326530612	0.0974614879302035\\
0.887857142857143	0.0921688542955872\\
0.893979591836735	0.0759511610069218\\
0.896020408163265	0.0714776813154127\\
0.898061224489796	0.0678609401676399\\
0.900102040816326	0.0658628754719628\\
0.902142857142857	0.0655334705029451\\
0.904183673469388	0.0666301719762153\\
0.906224489795918	0.069210925962656\\
0.908265306122449	0.072595551160705\\
0.91030612244898	0.0769964352512017\\
0.916428571428571	0.0922694056466341\\
0.918469387755102	0.0966379028706028\\
0.920510204081633	0.100101704465635\\
0.922551020408163	0.102559525876924\\
0.924591836734694	0.104202357598625\\
0.926632653061224	0.104827878023455\\
0.928673469387755	0.104362450807085\\
0.930714285714286	0.103370714225475\\
0.932755102040816	0.101222673895459\\
0.934795918367347	0.0977622401258647\\
0.936836734693878	0.0935224324317606\\
0.938877551020408	0.0886611576902501\\
0.940918367346939	0.083547420321112\\
0.942959183673469	0.078573722321309\\
0.945	0.0741166017577268\\
0.947040816326531	0.0704321508559536\\
0.949081632653061	0.0672243182872306\\
0.951122448979592	0.0653388436211239\\
0.953163265306122	0.0659063894338453\\
0.955204081632653	0.067959090324468\\
0.957244897959184	0.0706945353413478\\
0.959285714285714	0.0754248512034178\\
0.965408163265306	0.092710335275763\\
0.967448979591837	0.0979761817756836\\
0.969489795918367	0.1028492257822\\
0.971530612244898	0.10693576802646\\
0.973571428571429	0.109925807314584\\
0.975612244897959	0.111924089063216\\
0.97765306122449	0.11241532575973\\
0.97969387755102	0.11171124399473\\
0.981734693877551	0.109911831256387\\
0.983775510204082	0.106760437879068\\
0.985816326530612	0.102399977212957\\
0.989897959183674	0.0918745443462987\\
0.991938775510204	0.0869544628681902\\
0.993979591836735	0.08240876096195\\
0.996020408163265	0.0780055444018674\\
0.998061224489796	0.0744861717883306\\
};
\addlegendentry{ABC$_W^{1/2}$}

\end{axis}
\end{tikzpicture}%

%% file: images/MultiSource/Cost.tex
%
%
\definecolor{mycolor1}{rgb}{1.00000,0.60000,0.00000}%
\pgfplotsset{scaled y ticks=false}
\begin{tikzpicture}[scale = 0.52, font=\Large]

\begin{axis}[%
width=4.65in,
height=3.861in,
at={(0.78in,0.521in)},
scale only axis,
xmin=0,
xmax=1,
xlabel style={font=\huge\color{white!15!black}},
y tick label style={
        /pgf/number format/.cd,
            fixed,
            fixed zerofill,
            precision=2,
        /tikz/.cd
            },
xlabel={\huge relative simulation time $t/T$},
ymin=0,
ymax=0.09, ytick={0.02, 0.04, 0.06, 0.08},
ylabel style={at={(-0.15, 0.5)}, font=\color{white!15!black}},
ylabel={\huge relative $L^{2}(\Omega)$ error},
axis background/.style={fill=white},
xmajorgrids,
ymajorgrids,
legend style={at={(0.03,0.97)}, anchor=north west, legend cell align=left, align=left, draw=white!15!black}
]
\addplot [color=mycolor1, dashed, line width=2.0pt]
  table[row sep=crcr]{%
0.00714285714285712	5.6621374255883e-15\\
0.224829931972789	1.8161840355968e-05\\
0.245238095238095	8.9330159971146e-05\\
0.252040816326531	0.000164939186650281\\
0.258843537414966	0.000331963053036177\\
0.265646258503401	0.000628623020592078\\
0.272448979591837	0.00108311075075262\\
0.279251700680272	0.00170549893909233\\
0.286054421768708	0.00248165085048679\\
0.292857142857143	0.00338301755666259\\
0.299659863945578	0.00438680156475513\\
0.306462585034014	0.00548777328250249\\
0.313265306122449	0.00669865759993693\\
0.320068027210884	0.00804502641295535\\
0.32687074829932	0.00955694694119591\\
0.333673469387755	0.0112553481434577\\
0.34047619047619	0.0131372148403592\\
0.347278911564626	0.0151591150799364\\
0.354081632653061	0.0172246046794081\\
0.360884353741497	0.0191975498870248\\
0.367687074829932	0.0209515967740183\\
0.374489795918367	0.0224335220543023\\
0.381292517006803	0.0236906352110078\\
0.388095238095238	0.0248433259141246\\
0.394897959183674	0.0260357409479699\\
0.401700680272109	0.0274001933743162\\
0.408503401360544	0.0290422184865453\\
0.41530612244898	0.0310211316039321\\
0.422108843537415	0.0333388099551643\\
0.42891156462585	0.0359541295901639\\
0.435714285714286	0.0387884523889979\\
0.442517006802721	0.0416974908426305\\
0.449319727891156	0.0444541829424242\\
0.456122448979592	0.0467969473526226\\
0.462925170068027	0.0485189478108862\\
0.469727891156463	0.0495251574668023\\
0.476530612244898	0.0498222343666458\\
0.483333333333333	0.0494787326966444\\
0.490136054421769	0.0486095389502699\\
0.496938775510204	0.0473899903011765\\
0.503741496598639	0.0460683835116478\\
0.510544217687075	0.0449622238123429\\
0.51734693877551	0.0444215686123764\\
0.524149659863946	0.0447329236136189\\
0.530952380952381	0.0459993058025183\\
0.537755102040816	0.0480986261589231\\
0.544557823129252	0.0507783692301945\\
0.551360544217687	0.0538115431127162\\
0.558163265306122	0.0570876038883477\\
0.564965986394558	0.0605962496718858\\
0.571768707482993	0.0643315292208079\\
0.578571428571429	0.0681724586631953\\
0.585374149659864	0.0718044611519281\\
0.592176870748299	0.0747135226645007\\
0.598979591836735	0.0762600781241463\\
0.60578231292517	0.0758534019556527\\
0.612585034013605	0.073237258994385\\
0.619387755102041	0.0687956371268919\\
0.626190476190476	0.0636660119683781\\
0.632993197278912	0.0594654754621087\\
0.639795918367347	0.0576342038047116\\
0.646598639455782	0.0587263544102994\\
0.653401360544218	0.0622384210939692\\
0.660204081632653	0.0670880699230055\\
0.667006802721088	0.0721287723540156\\
0.673809523809524	0.0763350414422402\\
0.680612244897959	0.0788090440896245\\
0.687414965986395	0.0788330025109396\\
0.69421768707483	0.0760723125641073\\
0.701020408163265	0.0708799335614069\\
0.707823129251701	0.0644894795267872\\
0.714625850340136	0.0588360746262853\\
0.721428571428571	0.0558653119869976\\
0.728231292517007	0.0565206044660221\\
0.735034013605442	0.0602702876029507\\
0.741836734693878	0.0657343412358556\\
0.748639455782313	0.0715049794085422\\
0.755442176870748	0.076420849035146\\
0.762244897959184	0.0795197202220554\\
0.769047619047619	0.0800132317017265\\
0.775850340136054	0.0774310046279445\\
0.78265306122449	0.0719240934133568\\
0.789455782312925	0.0645588004536268\\
0.796258503401361	0.057340517443866\\
0.803061224489796	0.0527096014071357\\
0.809863945578231	0.0523396356098789\\
0.816666666666667	0.0560115298762399\\
0.823469387755102	0.062093220218183\\
0.830272108843537	0.0688359418861825\\
0.837074829931973	0.0748967986893013\\
0.843877551020408	0.0792532725302625\\
0.850680272108844	0.0810817768533353\\
0.857482993197279	0.0798172802817828\\
0.864285714285714	0.0754081618213438\\
0.87108843537415	0.0686392128651914\\
0.877891156462585	0.0612605840611012\\
0.88469387755102	0.0556515177149114\\
0.891496598639456	0.0538369976127179\\
0.898299319727891	0.0562426704375899\\
0.905102040816327	0.0616196940373737\\
0.911904761904762	0.0682236407810853\\
0.918707482993197	0.0746046556494772\\
0.925510204081633	0.0796583485214334\\
0.932312925170068	0.0824870189096627\\
0.939115646258503	0.082388219059133\\
0.945918367346939	0.0790478861076477\\
0.952721088435374	0.0728650643260788\\
0.95952380952381	0.0651909790691358\\
0.966326530612245	0.0582241870922467\\
0.97312925170068	0.054344470079937\\
0.979931972789116	0.0548578527046396\\
0.986734693877551	0.0591824912064867\\
0.993537414965986	0.0656213396322739\\
};
\addlegendentry{ABC$_W^{1/2}$}

\addplot [color=mycolor1, line width=2.0pt]
  table[row sep=crcr]{%
0.00714285714285712	5.6621374255883e-15\\
0.224829931972789	1.8161840355968e-05\\
0.245238095238095	8.9330159971146e-05\\
0.252040816326531	0.000165834078648563\\
0.258843537414966	0.000340980679634062\\
0.265646258503401	0.000654669221037185\\
0.272448979591837	0.00113495389701079\\
0.279251700680272	0.00179402514264559\\
0.286054421768708	0.00261860337545305\\
0.292857142857143	0.00358214071047591\\
0.299659863945578	0.00465455414211158\\
0.306462585034014	0.00582147500679797\\
0.313265306122449	0.00708434427687887\\
0.320068027210884	0.00845807181026381\\
0.32687074829932	0.00996829798656407\\
0.333673469387755	0.0116281908095127\\
0.34047619047619	0.0134215834959196\\
0.347278911564626	0.0152856824122579\\
0.354081632653061	0.017105144103345\\
0.360884353741497	0.018725685858191\\
0.367687074829932	0.020020256083717\\
0.374489795918367	0.0209378984431499\\
0.381292517006803	0.0215398693519863\\
0.388095238095238	0.0219539982373891\\
0.394897959183674	0.0223303789370438\\
0.401700680272109	0.0228159650147984\\
0.408503401360544	0.0235012746671366\\
0.41530612244898	0.0243961506595273\\
0.422108843537415	0.02547016429871\\
0.42891156462585	0.0267243643358173\\
0.435714285714286	0.0281887699382535\\
0.442517006802721	0.0297689890284082\\
0.449319727891156	0.0309877569371655\\
0.456122448979592	0.0313592570951985\\
0.462925170068027	0.0309851303634776\\
0.469727891156463	0.0302622013799708\\
0.476530612244898	0.0296948405447537\\
0.483333333333333	0.0294001798243065\\
0.490136054421769	0.0294861258425076\\
0.496938775510204	0.0301493035066734\\
0.503741496598639	0.0313152147611281\\
0.510544217687075	0.0327714980015088\\
0.524149659863946	0.0359295657606518\\
0.530952380952381	0.0372691328257972\\
0.537755102040816	0.0381883257767436\\
0.551360544217687	0.039337927161062\\
0.558163265306122	0.040055658534441\\
0.564965986394558	0.0409841859260384\\
0.571768707482993	0.0422660427153729\\
0.578571428571429	0.0438796989963406\\
0.585374149659864	0.0456805206372795\\
0.592176870748299	0.0476639976201733\\
0.598979591836735	0.0497753830378523\\
0.60578231292517	0.0519444836433552\\
0.612585034013605	0.0540784495541798\\
0.619387755102041	0.0559688680209988\\
0.626190476190476	0.0569618691863137\\
0.632993197278912	0.0561553513939326\\
0.639795918367347	0.0536809248635344\\
0.646598639455782	0.0506923702340167\\
0.653401360544218	0.0481708050177084\\
0.660204081632653	0.0465838552481852\\
0.667006802721088	0.0464806097730619\\
0.673809523809524	0.0489958232878305\\
0.680612244897959	0.0541750380088383\\
0.687414965986395	0.0600578392329725\\
0.69421768707483	0.0649148001725037\\
0.701020408163265	0.0682203493933039\\
0.707823129251701	0.0695908900807306\\
0.714625850340136	0.0685260245319876\\
0.721428571428571	0.0652272145202677\\
0.728231292517007	0.0608535938193301\\
0.735034013605442	0.0566981979735762\\
0.741836734693878	0.0535935190998744\\
0.748639455782313	0.0522208792326511\\
0.755442176870748	0.0532048859312546\\
0.762244897959184	0.0563347284015817\\
0.769047619047619	0.0600813742637299\\
0.775850340136054	0.063064314593253\\
0.78265306122449	0.0647451059644866\\
0.789455782312925	0.0650298088584349\\
0.796258503401361	0.0638177903775845\\
0.803061224489796	0.0611616007685455\\
0.809863945578231	0.0579124845425296\\
0.816666666666667	0.0550919699172414\\
0.823469387755102	0.0532251001753244\\
0.830272108843537	0.0526082311264886\\
0.837074829931973	0.0534625721247488\\
0.843877551020408	0.0557511985655957\\
0.850680272108844	0.0588333249246087\\
0.857482993197279	0.0616607280263778\\
0.864285714285714	0.0634802861836464\\
0.87108843537415	0.0639332305114578\\
0.877891156462585	0.063082500604447\\
0.88469387755102	0.0611759341470689\\
0.891496598639456	0.0587239071300917\\
0.898299319727891	0.056463303785694\\
0.905102040816327	0.054915298913504\\
0.911904761904762	0.0544541058891881\\
0.918707482993197	0.0551983369735796\\
0.925510204081633	0.0569772186071569\\
0.932312925170068	0.0594959571280821\\
0.939115646258503	0.062165953972985\\
0.945918367346939	0.0643887869753528\\
0.952721088435374	0.0655819840396216\\
0.95952380952381	0.0655661202433515\\
0.966326530612245	0.0647178077045292\\
0.97312925170068	0.0633575347927906\\
0.979931972789116	0.0618027314996666\\
0.986734693877551	0.0603539298707445\\
0.993537414965986	0.0593456115774784\\
};
\addlegendentry{ABC$_W^{1/2}$ with $\theta$}

\end{axis}
\end{tikzpicture}%

%% file: images/Transducer3D/Cost.tex
%
%
\definecolor{mycolor1}{rgb}{1.00000,0.60000,0.00000}%
\begin{tikzpicture}[scale = 0.52, font=\huge]

\begin{axis}[%
width=4.602in,
height=3.82in,
at={(0.772in,0.516in)},
scale only axis,
xmin=0,
xmax=1,
xlabel style={font=\huge\color{white!15!black}},
xlabel={relative simulation time $t/T$}, xtick={0, 0.2, 0.4, 0.6, 0.8, 1},
y tick label style={
        /pgf/number format/.cd,
            fixed,
            fixed zerofill,
            precision=2,
        /tikz/.cd
            },
ymin=0,
ymax=0.25,
ylabel style={at={(-0.15,0.5)}, font=\huge\color{white!15!black}}, ytick={0.05, 0.1, 0.15, 0.2, 0.25},
ylabel={relative $L^{2}(\Omega)$ error},
axis background/.style={fill=white},
xmajorgrids,
ymajorgrids,
legend style={at={(0.03,0.97)}, anchor=north west, legend cell align=left, align=left, draw=white!15!black}
]
\addplot [color=mycolor1, line width=2.0pt]
  table[row sep=crcr]{%
0.00428571428571434	0.238478639468435\\
0.0083673469387755	0.168763083192172\\
0.0124489795918368	0.128174631181275\\
0.0165306122448979	0.10042060772048\\
0.0206122448979592	0.0786225385159446\\
0.0246938775510204	0.061540691666722\\
0.0287755102040816	0.0484187566354961\\
0.0328571428571428	0.0391683342443725\\
0.0369387755102041	0.0342785111763484\\
0.0410204081632654	0.0349845399952116\\
0.0451020408163265	0.0422719757480307\\
0.0491836734693878	0.0557844101344758\\
0.053265306122449	0.072010287288562\\
0.0573469387755102	0.0796627288052413\\
0.0614285714285714	0.071729866321614\\
0.0655102040816327	0.0581782862338296\\
0.0695918367346939	0.046611054195958\\
0.0736734693877551	0.0382385000979831\\
0.0777551020408164	0.0326190654295183\\
0.0818367346938775	0.029333484688961\\
0.0859183673469388	0.028160023821665\\
0.09	0.0289870172694253\\
0.0940816326530612	0.0316975939006691\\
0.0981632653061224	0.0360618587561147\\
0.102244897959184	0.041426926264664\\
0.106326530612245	0.0462329129437584\\
0.110408163265306	0.0481945771082549\\
0.114489795918367	0.0462615406323309\\
0.118571428571429	0.0418268830477195\\
0.12265306122449	0.0369308897104553\\
0.126734693877551	0.0327782345710139\\
0.130816326530612	0.0298398971128608\\
0.134897959183674	0.0282781704099571\\
0.138979591836735	0.0281454659582105\\
0.143061224489796	0.0294102126868113\\
0.147142857142857	0.0319054448953359\\
0.15530612244898	0.0385938588344792\\
0.159387755102041	0.0408815093030445\\
0.163469387755102	0.0411232272240996\\
0.167551020408163	0.0392919964826973\\
0.171632653061225	0.0362796508439063\\
0.175714285714286	0.0331331152584142\\
0.179795918367347	0.0305686847237966\\
0.183877551020408	0.0289689946238786\\
0.187959183673469	0.0285065949174078\\
0.192040816326531	0.029206447139932\\
0.196122448979592	0.0309485572278004\\
0.200204081632653	0.0334377533161774\\
0.204285714285714	0.0361095225241088\\
0.208367346938776	0.0380179469753575\\
0.212448979591837	0.0387491438210659\\
0.216530612244898	0.0383104958997088\\
0.220612244897959	0.0369028018119995\\
0.22469387755102	0.0350958052530668\\
0.228775510204082	0.0334330366649588\\
0.232857142857143	0.0322643704772914\\
0.236938775510204	0.0318647095787216\\
0.241020408163265	0.0322931626868802\\
0.245102040816327	0.0334416683111717\\
0.249183673469388	0.0350704321894266\\
0.253265306122449	0.0368296040537043\\
0.25734693877551	0.0383217698267809\\
0.261428571428571	0.0392329392881968\\
0.265510204081633	0.0393828735222219\\
0.269591836734694	0.0388754148300368\\
0.277755102040816	0.0371912412143026\\
0.281836734693878	0.0367448657929713\\
0.285918367346939	0.0369232898224489\\
0.29	0.0377753225602719\\
0.294081632653061	0.0392014547719969\\
0.298163265306122	0.0410092215545178\\
0.302244897959184	0.0429818613854428\\
0.306326530612245	0.0448240845277629\\
0.310408163265306	0.0462944553817399\\
0.314489795918367	0.0472552908748776\\
0.318571428571429	0.0477120803103749\\
0.330816326530612	0.0481413317937629\\
0.334897959183673	0.0487893208639474\\
0.338979591836735	0.0499422083984086\\
0.343061224489796	0.0516163591333958\\
0.347142857142857	0.0537102515718734\\
0.351224489795918	0.0559883856577973\\
0.35530612244898	0.0581029247634566\\
0.359387755102041	0.0596800456744567\\
0.363469387755102	0.0604493510915067\\
0.367551020408163	0.0603659255831976\\
0.375714285714286	0.0589397700397406\\
0.379795918367347	0.058832710543911\\
0.383877551020408	0.0599486307462861\\
0.387959183673469	0.0624893901524743\\
0.392040816326531	0.0661814793052892\\
0.396122448979592	0.0703855154272961\\
0.400204081632653	0.0742901209288528\\
0.404285714285714	0.0770571279355299\\
0.408367346938776	0.0779408413793\\
0.412448979591837	0.0765153825879593\\
0.416530612244898	0.0729863516161056\\
0.420612244897959	0.0684268210740838\\
0.42469387755102	0.0647030969244642\\
0.428775510204082	0.0638495279404555\\
0.432857142857143	0.0669485997015085\\
0.436938775510204	0.0734364808227033\\
0.445102040816327	0.0896829977643928\\
0.449183673469388	0.0959696224465456\\
0.453265306122449	0.0991895336076132\\
0.45734693877551	0.0983786279260596\\
0.461428571428571	0.0930752767521682\\
0.465510204081633	0.0835713033225064\\
0.469591836734694	0.0712840370614631\\
0.473673469387755	0.0592853962286097\\
0.477755102040816	0.0526317570861924\\
0.481836734693878	0.0556845026939536\\
0.485918367346939	0.066903160763958\\
0.49	0.0813400874350112\\
0.494081632653061	0.0952085035206635\\
0.498163265306122	0.10615798356044\\
0.502244897959184	0.112605781849668\\
0.506326530612245	0.113457253301408\\
0.510408163265306	0.108166326794908\\
0.514489795918367	0.096923146263048\\
0.518571428571429	0.0808969479128924\\
0.52265306122449	0.0627126972493743\\
0.526734693877551	0.0477015153312615\\
0.530816326530612	0.0447068334177539\\
0.534897959183673	0.056120722122863\\
0.543061224489796	0.0912472961328528\\
0.547142857142857	0.105685773193329\\
0.551224489795918	0.115294578856706\\
0.55530612244898	0.118903518731648\\
0.559387755102041	0.115829618988139\\
0.563469387755102	0.106022925594049\\
0.567551020408163	0.0902685559867688\\
0.575714285714286	0.05053160993265\\
0.579795918367347	0.0389973558007368\\
0.583877551020408	0.0456826327785252\\
0.587959183673469	0.0638049446006385\\
0.592040816326531	0.0836017318667309\\
0.596122448979592	0.100863216652614\\
0.600204081632653	0.113531882701782\\
0.604285714285714	0.12029895617395\\
0.608367346938775	0.120320858973772\\
0.612448979591837	0.11332433650445\\
0.616530612244898	0.0997716684401868\\
0.620612244897959	0.0810546226214061\\
0.62469387755102	0.0599072312934281\\
0.628775510204082	0.0420693001807951\\
0.632857142857143	0.0389711792410938\\
0.636938775510204	0.053179125594782\\
0.641020408163265	0.0730654846540015\\
0.645102040816326	0.0919589809191468\\
0.649183673469388	0.10708797970361\\
0.653265306122449	0.116891690062618\\
0.65734693877551	0.120338216913523\\
0.661428571428571	0.116891715538702\\
0.665510204081633	0.106666699312994\\
0.669591836734694	0.0906234523550143\\
0.677755102040816	0.051470594664268\\
0.681836734693878	0.0409761268002198\\
0.685918367346939	0.0478407669910083\\
0.69	0.0653439063669109\\
0.694081632653061	0.0844330435565137\\
0.698163265306122	0.101004725891272\\
0.702244897959184	0.113015071333982\\
0.706326530612245	0.119181266329123\\
0.710408163265306	0.118738687110098\\
0.714489795918367	0.111540437321817\\
0.718571428571429	0.0982189232462569\\
0.72265306122449	0.0804240892333113\\
0.726734693877551	0.0613467126894279\\
0.730816326530612	0.0471285014992748\\
0.734897959183674	0.0466271736145791\\
0.738979591836735	0.0597131446884948\\
0.743061224489796	0.0776119305903427\\
0.747142857142857	0.0946953602965647\\
0.751224489795918	0.108248568466455\\
0.75530612244898	0.116693754055558\\
0.759387755102041	0.119044181020373\\
0.763469387755102	0.114860241484797\\
0.767551020408163	0.104418058851917\\
0.771632653061225	0.0889230581448049\\
0.775714285714286	0.0708820889978049\\
0.779795918367347	0.0549555390110992\\
0.783877551020408	0.0487098197551051\\
0.787959183673469	0.0563515873581105\\
0.792040816326531	0.0719534769843717\\
0.796122448979592	0.0887926993999252\\
0.800204081632653	0.103272454826974\\
0.804285714285714	0.113427072208667\\
0.808367346938776	0.118057255426674\\
0.812448979591837	0.116519555552512\\
0.816530612244898	0.108815809842092\\
0.820612244897959	0.095789492177123\\
0.828775510204082	0.0633345456075379\\
0.832857142857143	0.0534985458310496\\
0.836938775510204	0.0555189131049265\\
0.841020408163265	0.067431579224454\\
0.845102040816327	0.082876245143076\\
0.849183673469388	0.0974747188846171\\
0.853265306122449	0.108799331500532\\
0.85734693877551	0.115393703934775\\
0.861428571428571	0.11637233603313\\
0.865510204081633	0.111434629807582\\
0.869591836734694	0.101050848595376\\
0.873673469387755	0.0867470031826003\\
0.877755102040816	0.0714726025768035\\
0.881836734693878	0.0599946652875959\\
0.885918367346939	0.0577025977746951\\
0.89	0.0654385177730612\\
0.894081632653061	0.0785320750317579\\
0.898163265306122	0.092323855376905\\
0.902244897959184	0.103984315440978\\
0.906326530612245	0.111824128716573\\
0.910408163265306	0.11479219749725\\
0.914489795918367	0.112368294659359\\
0.918571428571429	0.104717866762101\\
0.92265306122449	0.0929457233018821\\
0.926734693877551	0.0793872100203764\\
0.930816326530612	0.067808729361601\\
0.934897959183673	0.0627211265069426\\
0.938979591836735	0.0663388656675392\\
0.943061224489796	0.0762722922372601\\
0.947142857142857	0.0884509088093586\\
0.951224489795918	0.0998156962422444\\
0.95530612244898	0.108419348182112\\
0.959387755102041	0.112994827435622\\
0.963469387755102	0.112769861366292\\
0.967551020408163	0.107566383490523\\
0.971632653061224	0.0980249413604587\\
0.975714285714286	0.0859324672148714\\
0.979795918367347	0.0743724542346773\\
0.983877551020408	0.0673111213890609\\
0.987959183673469	0.0676222953172169\\
0.992040816326531	0.0746389322630066\\
0.996122448979592	0.0851647029362714\\
};
\addlegendentry{ABC$_W^{1/2}$ adaptive}

\addplot [color=mycolor1, dashed]
  table[row sep=crcr]{%
0.00428571428571434	0.238478639468435\\
0.0083673469387755	0.168763083192172\\
0.0124489795918368	0.128174631181275\\
0.0165306122448979	0.10042060772048\\
0.0206122448979592	0.0786225385159446\\
0.0246938775510204	0.061540691666722\\
0.0287755102040816	0.0484187566354961\\
0.0328571428571428	0.0391683342443725\\
0.0369387755102041	0.0342785111763484\\
0.0410204081632654	0.0349845399952116\\
0.0451020408163265	0.0422719757480307\\
0.0491836734693878	0.0557844101344758\\
0.053265306122449	0.072010287288562\\
0.0573469387755102	0.0796627288052413\\
0.0614285714285714	0.071729866321614\\
0.0655102040816327	0.0581782862338296\\
0.0695918367346939	0.046611054195958\\
0.0736734693877551	0.0382385000979831\\
0.0777551020408164	0.0326190654295183\\
0.0818367346938775	0.029333484688961\\
0.0859183673469388	0.028160023821665\\
0.09	0.0289870172694253\\
0.0940816326530612	0.0316975939006691\\
0.0981632653061224	0.0360618587561147\\
0.102244897959184	0.041426926264664\\
0.106326530612245	0.0462329129437584\\
0.110408163265306	0.0481945771082549\\
0.114489795918367	0.0462615406323309\\
0.118571428571429	0.0418268830477195\\
0.12265306122449	0.0369308897104553\\
0.126734693877551	0.0327782345710139\\
0.130816326530612	0.0298398971128608\\
0.134897959183674	0.0282781704099571\\
0.138979591836735	0.0281454659582105\\
0.143061224489796	0.0294102126868113\\
0.147142857142857	0.0319054448953359\\
0.15530612244898	0.0385938588344792\\
0.159387755102041	0.0408815093030445\\
0.163469387755102	0.0411232272240996\\
0.167551020408163	0.0392919964826973\\
0.171632653061225	0.0362796508439063\\
0.175714285714286	0.0331331152584142\\
0.179795918367347	0.0305686847237966\\
0.183877551020408	0.0289689946238786\\
0.187959183673469	0.0285065949174078\\
0.192040816326531	0.029206447139932\\
0.196122448979592	0.0309485572278004\\
0.200204081632653	0.0334377533161774\\
0.204285714285714	0.0361095225241088\\
0.208367346938776	0.0380179469753575\\
0.212448979591837	0.0387491438210659\\
0.216530612244898	0.0383104958997088\\
0.220612244897959	0.0369028018119995\\
0.22469387755102	0.0350958052530668\\
0.228775510204082	0.0334330413916376\\
0.232857142857143	0.0322639298321331\\
0.236938775510204	0.0318585727067447\\
0.241020408163265	0.0322689670208413\\
0.245102040816327	0.033379891195631\\
0.249183673469388	0.0349493312349175\\
0.253265306122449	0.0366345269465848\\
0.25734693877551	0.0380524761219521\\
0.261428571428571	0.038904979668646\\
0.265510204081633	0.039026902348719\\
0.269591836734694	0.0385314238018781\\
0.277755102040816	0.0369678270468431\\
0.281836734693878	0.0366039409601033\\
0.285918367346939	0.0368632036583918\\
0.29	0.0377872516808964\\
0.294081632653061	0.0392743904957035\\
0.298163265306122	0.041137599552894\\
0.302244897959184	0.0431746651064966\\
0.306326530612245	0.0451099399254552\\
0.310408163265306	0.0467099276768801\\
0.314489795918367	0.0478261408042123\\
0.318571428571429	0.0484486458613792\\
0.32265306122449	0.0487314324355348\\
0.326734693877551	0.0489085756686759\\
0.330816326530612	0.049262030090569\\
0.334897959183673	0.050035410816099\\
0.338979591836735	0.0513985124250331\\
0.343061224489796	0.0534277137347177\\
0.347142857142857	0.056070019395227\\
0.35530612244898	0.0622434178378134\\
0.359387755102041	0.0649629425962003\\
0.363469387755102	0.0668707481646761\\
0.367551020408163	0.0676980590203254\\
0.371632653061225	0.0674339410883928\\
0.379795918367347	0.0653236213104014\\
0.383877551020408	0.0650578483058427\\
0.387959183673469	0.0664480267325586\\
0.392040816326531	0.0699387845196791\\
0.396122448979592	0.0753595475310095\\
0.404285714285714	0.0886072353623638\\
0.408367346938776	0.0940172231185483\\
0.412448979591837	0.096973450632448\\
0.416530612244898	0.0966295295519614\\
0.420612244897959	0.0928019750069364\\
0.42469387755102	0.0861758478352496\\
0.428775510204082	0.0784258052261207\\
0.432857142857143	0.0722434998576075\\
0.436938775510204	0.0708001699798356\\
0.441020408163265	0.0759645569179296\\
0.445102040816327	0.0867758966469692\\
0.453265306122449	0.113815314323569\\
0.45734693877551	0.124325917436835\\
0.461428571428571	0.129779065806568\\
0.465510204081633	0.128701018717868\\
0.469591836734694	0.120593656535052\\
0.473673469387755	0.106167184762324\\
0.481836734693878	0.0689047612399919\\
0.485918367346939	0.0579562757754948\\
0.49	0.0628982758641757\\
0.494081632653061	0.0808815718055604\\
0.498163265306122	0.103212899919264\\
0.502244897959184	0.124126129723148\\
0.506326530612245	0.140121718882929\\
0.510408163265306	0.148805194457589\\
0.514489795918367	0.148689903286734\\
0.518571428571429	0.139367842314379\\
0.52265306122449	0.121668972499956\\
0.526734693877551	0.0976976108082107\\
0.530816326530612	0.0711891543108173\\
0.534897959183673	0.0499754573917501\\
0.538979591836735	0.0493147633143824\\
0.543061224489796	0.0699239792484254\\
0.547142857142857	0.0966799899592544\\
0.551224489795918	0.12174133834553\\
0.55530612244898	0.14145155629353\\
0.559387755102041	0.153445284886568\\
0.563469387755102	0.156159608759804\\
0.567551020408163	0.148990855268271\\
0.571632653061224	0.132503265674759\\
0.575714285714286	0.108475697975949\\
0.579795918367347	0.07999105053765\\
0.583877551020408	0.0528951703653575\\
0.587959183673469	0.0421741293180908\\
0.592040816326531	0.0593792250908296\\
0.600204081632653	0.114726939794313\\
0.604285714285714	0.13724299458889\\
0.608367346938775	0.152359414229791\\
0.612448979591837	0.158369752833204\\
0.616530612244898	0.154442903107131\\
0.620612244897959	0.140835923613422\\
0.62469387755102	0.11895419413266\\
0.628775510204082	0.0913225817469635\\
0.632857142857143	0.0621674644953132\\
0.636938775510204	0.0417309089337984\\
0.641020408163265	0.0488327397356144\\
0.645102040816326	0.0748753518391408\\
0.649183673469388	0.103358223790212\\
0.653265306122449	0.128318749125937\\
0.65734693877551	0.146830515088763\\
0.661428571428571	0.156905202869075\\
0.665510204081633	0.157318915494104\\
0.669591836734694	0.147836900481231\\
0.673673469387755	0.129378193731259\\
0.677755102040816	0.10404672791753\\
0.681836734693878	0.0753682363039944\\
0.685918367346939	0.050469909342156\\
0.69	0.0457968376393848\\
0.694081632653061	0.0661232013707852\\
0.698163265306122	0.0939916407245488\\
0.702244897959184	0.120280883358871\\
0.706326530612245	0.141179960719191\\
0.710408163265306	0.15435370888423\\
0.714489795918367	0.158293656195616\\
0.718571428571429	0.152421209113735\\
0.72265306122449	0.137285632978486\\
0.726734693877551	0.114622271977197\\
0.730816326530612	0.0874220225719728\\
0.734897959183674	0.0610011921174983\\
0.738979591836735	0.0473162300652056\\
0.743061224489796	0.0587766530204171\\
0.747142857142857	0.0839462885760224\\
0.751224489795918	0.110450674626357\\
0.75530612244898	0.133131275651897\\
0.759387755102041	0.149196918811157\\
0.763469387755102	0.156787698670897\\
0.767551020408163	0.154898832469012\\
0.771632653061225	0.143582230311338\\
0.775714285714286	0.124104023325447\\
0.779795918367347	0.0990225023359929\\
0.783877551020408	0.0727162829439125\\
0.787959183673469	0.0539555887431602\\
0.792040816326531	0.0561092470119623\\
0.796122448979592	0.0766331798367359\\
0.800204081632653	0.10216059795075\\
0.804285714285714	0.125648546213118\\
0.808367346938776	0.143650992994034\\
0.812448979591837	0.153997817414687\\
0.816530612244898	0.155374520081285\\
0.820612244897959	0.147459790234163\\
0.82469387755102	0.131107284474186\\
0.828775510204082	0.108448045366474\\
0.832857142857143	0.0831479028109555\\
0.836938775510204	0.061821519669759\\
0.841020408163265	0.0559300776007421\\
0.845102040816327	0.0699540358118294\\
0.853265306122449	0.116604804471037\\
0.85734693877551	0.136143718357384\\
0.861428571428571	0.149159042767103\\
0.865510204081633	0.153960116667081\\
0.869591836734694	0.149754677876775\\
0.873673469387755	0.136904373837613\\
0.877755102040816	0.117095739978417\\
0.881836734693878	0.0934762983477423\\
0.885918367346939	0.0713498953265901\\
0.89	0.0597722752122404\\
0.894081632653061	0.0665086004006515\\
0.898163265306122	0.0858612832387364\\
0.902244897959184	0.108345302993039\\
0.906326530612245	0.128600702210461\\
0.910408163265306	0.143597852957999\\
0.914489795918367	0.151330369937195\\
0.918571428571429	0.150635095382787\\
0.92265306122449	0.141430422361491\\
0.926734693877551	0.124936531056412\\
0.934897959183673	0.0822785674919003\\
0.938979591836735	0.0672814236215546\\
0.943061224489796	0.0667034634604277\\
0.947142857142857	0.0803629467221837\\
0.95530612244898	0.120529108462298\\
0.959387755102041	0.137061314389698\\
0.963469387755102	0.14752976601615\\
0.967551020408163	0.150272857507118\\
0.971632653061224	0.144600686166152\\
0.975714285714286	0.131142019268703\\
0.979795918367347	0.112006729078805\\
0.983877551020408	0.0908791777493447\\
0.987959183673469	0.0734691435182191\\
0.992040816326531	0.0673512303398152\\
0.996122448979592	0.0758884187046904\\
};
\addlegendentry{ABC$_W^{1/2}$}

\end{axis}
\end{tikzpicture}%

%% file: images/Transducer3D/Cost_abs.tex
%
%
\definecolor{mycolor1}{rgb}{1.00000,0.60000,0.00000}%
\begin{tikzpicture}[scale = 0.52, font=\huge]

\begin{axis}[%
width=4.602in,
height=3.82in,
at={(0.772in,0.516in)},
scale only axis,
xmin=0,
xmax=1,
xlabel style={font=\huge\color{white!15!black}},
xlabel={relative simulation time $t/T$}, xtick={0, 0.2, 0.4, 0.6, 0.8, 1},
y tick label style={
        /pgf/number format/.cd,
            fixed,
            fixed zerofill,
            precision=1,
        /tikz/.cd
            },
ymin=0,
ymax=2.5e-07,
ylabel style={font=\huge\color{white!15!black}},  ytick={0.00000005, 0.0000001, 0.00000015, 0.0000002, 0.00000025},
ylabel={absolute $L^{2}(\Omega)$ error},
axis background/.style={fill=white},
xmajorgrids,
ymajorgrids,
legend style={at={(0.03,0.97)}, anchor=north west, legend cell align=left, align=left, draw=white!15!black}
]
\addplot [color=mycolor1, line width=2.0pt]
  table[row sep=crcr]{%
0.000204081632653064	0\\
0.0124489795918368	9.30855392766716e-11\\
0.0328571428571428	3.99368538239742e-10\\
0.0410204081632654	5.06275354972274e-10\\
0.0491836734693878	8.78358163980408e-10\\
0.0655102040816327	1.93574511975214e-09\\
0.0736734693877551	2.27844820788192e-09\\
0.0818367346938775	2.55820098438164e-09\\
0.0859183673469388	2.77207146126557e-09\\
0.09	3.09003478360381e-09\\
0.0940816326530612	3.52595530550559e-09\\
0.102244897959184	4.67235605761118e-09\\
0.110408163265306	5.89675630546793e-09\\
0.114489795918367	6.42556330365807e-09\\
0.118571428571429	6.85945533707155e-09\\
0.126734693877551	7.46617345726008e-09\\
0.134897959183674	8.08788958028828e-09\\
0.138979591836735	8.61585669476028e-09\\
0.143061224489796	9.3677386958646e-09\\
0.147142857142857	1.03329809153507e-08\\
0.159387755102041	1.3692421640954e-08\\
0.163469387755102	1.46427363567625e-08\\
0.167551020408163	1.53962981253031e-08\\
0.171632653061225	1.59507518304736e-08\\
0.183877551020408	1.72680331145258e-08\\
0.187959183673469	1.8068830764939e-08\\
0.192040816326531	1.92461321324444e-08\\
0.196122448979592	2.07907282501196e-08\\
0.204285714285714	2.44544181571271e-08\\
0.208367346938776	2.60143436792148e-08\\
0.212448979591837	2.72503473119912e-08\\
0.216530612244898	2.82275093299589e-08\\
0.220612244897959	2.88892372246607e-08\\
0.232857142857143	3.02137258545443e-08\\
0.236938775510204	3.10374328460483e-08\\
0.241020408163265	3.22895070681284e-08\\
0.245102040816327	3.3927917253429e-08\\
0.253265306122449	3.76819010439178e-08\\
0.25734693877551	3.93496296657858e-08\\
0.261428571428571	4.06679814268784e-08\\
0.265510204081633	4.156106137021e-08\\
0.269591836734694	4.21129432437439e-08\\
0.277755102040816	4.29896149789499e-08\\
0.281836734693878	4.37988726398331e-08\\
0.285918367346939	4.51360763209507e-08\\
0.29	4.70218175507853e-08\\
0.294081632653061	4.93297926951186e-08\\
0.302244897959184	5.44123371915006e-08\\
0.306326530612245	5.67777309612438e-08\\
0.310408163265306	5.88239928989154e-08\\
0.314489795918367	6.05273466902645e-08\\
0.318571428571429	6.19566546955497e-08\\
0.32265306122449	6.32879149087273e-08\\
0.326734693877551	6.47110780516158e-08\\
0.330816326530612	6.6422445543779e-08\\
0.334897959183673	6.85484385920532e-08\\
0.338979591836735	7.11198773206689e-08\\
0.343061224489796	7.40935909293583e-08\\
0.351224489795918	8.05185147445542e-08\\
0.35530612244898	8.33412192591254e-08\\
0.359387755102041	8.54548763973284e-08\\
0.363469387755102	8.66763730789089e-08\\
0.367551020408163	8.70631804428257e-08\\
0.371632653061225	8.69656190394252e-08\\
0.375714285714286	8.70389573748298e-08\\
0.379795918367347	8.81285857623482e-08\\
0.383877551020408	9.09520603276093e-08\\
0.387959183673469	9.56869024015816e-08\\
0.392040816326531	1.01800185059453e-07\\
0.396122448979592	1.08223633388249e-07\\
0.400204081632653	1.13715974858764e-07\\
0.404285714285714	1.17152904710593e-07\\
0.408367346938776	1.17704969104793e-07\\
0.412448979591837	1.15067690042991e-07\\
0.416530612244898	1.09761406075215e-07\\
0.420612244897959	1.03385265570743e-07\\
0.42469387755102	9.85791139829217e-08\\
0.428775510204082	9.82625998346975e-08\\
0.432857142857143	1.04008859636018e-07\\
0.436938775510204	1.14860075561829e-07\\
0.441020408163265	1.28031197199618e-07\\
0.445102040816327	1.40412190319239e-07\\
0.449183673469388	1.49423248108072e-07\\
0.453265306122449	1.53217109555293e-07\\
0.45734693877551	1.50720887992506e-07\\
0.461428571428571	1.41705362888267e-07\\
0.465510204081633	1.26915413733997e-07\\
0.469591836734694	1.08473774740325e-07\\
0.473673469387755	9.07648406434092e-08\\
0.477755102040816	8.12695774010663e-08\\
0.481836734693878	8.67465904663689e-08\\
0.485918367346939	1.04953145507736e-07\\
0.49	1.28049637782013e-07\\
0.494081632653061	1.49778178082371e-07\\
0.498163265306122	1.66233798615245e-07\\
0.502244897959184	1.75042448979923e-07\\
0.506326530612245	1.74939513541972e-07\\
0.510408163265306	1.65663662787807e-07\\
0.514489795918367	1.47934871019828e-07\\
0.518571428571429	1.23582554634183e-07\\
0.52265306122449	9.62863343501041e-08\\
0.526734693877551	7.3818082646504e-08\\
0.530816326530612	6.9790008883075e-08\\
0.534897959183673	8.82622366304986e-08\\
0.538979591836735	1.16459394572388e-07\\
0.543061224489796	1.44218753161951e-07\\
0.547142857142857	1.66457157835254e-07\\
0.551224489795918	1.80398136540028e-07\\
0.55530612244898	1.84556025240035e-07\\
0.559387755102041	1.78459781285945e-07\\
0.563469387755102	1.62576463735853e-07\\
0.567551020408163	1.38313707043203e-07\\
0.571632653061224	1.08365073137051e-07\\
0.575714285714286	7.82230338280954e-08\\
0.579795918367347	6.08662761303336e-08\\
0.583877551020408	7.18400325894208e-08\\
0.587959183673469	1.00839280925413e-07\\
0.592040816326531	1.32287628185068e-07\\
0.596122448979592	1.59152436562948e-07\\
0.600204081632653	1.78041018594755e-07\\
0.604285714285714	1.87147474850491e-07\\
0.608367346938775	1.85720596013184e-07\\
0.612448979591837	1.73941758108143e-07\\
0.616530612244898	1.52850708867724e-07\\
0.620612244897959	1.24468099826913e-07\\
0.62469387755102	9.25433715215362e-08\\
0.628775510204082	6.55092149415282e-08\\
0.632857142857143	6.11731379995817e-08\\
0.636938775510204	8.39880486269706e-08\\
0.641020408163265	1.15723212346097e-07\\
0.645102040816326	1.45493029846122e-07\\
0.649183673469388	1.6863820495594e-07\\
0.653265306122449	1.82756396793771e-07\\
0.65734693877551	1.866502599146e-07\\
0.661428571428571	1.8008610802589e-07\\
0.665510204081633	1.63721967649799e-07\\
0.669591836734694	1.39140124422177e-07\\
0.673673469387755	1.09257492164616e-07\\
0.677755102040816	7.99285737596023e-08\\
0.681836734693878	6.41440182080899e-08\\
0.685918367346939	7.54091035171101e-08\\
0.69	1.03421377350976e-07\\
0.694081632653061	1.33675404301314e-07\\
0.698163265306122	1.59343981676052e-07\\
0.702244897959184	1.77118858579739e-07\\
0.706326530612245	1.85282035825018e-07\\
0.710408163265306	1.83206165793592e-07\\
0.714489795918367	1.71226443801409e-07\\
0.718571428571429	1.50586446645562e-07\\
0.72265306122449	1.23664152584801e-07\\
0.726734693877551	9.49296031782509e-08\\
0.730816326530612	7.35235428095393e-08\\
0.734897959183674	7.33197819124953e-08\\
0.738979591836735	9.44453648621391e-08\\
0.743061224489796	1.23047792288666e-07\\
0.747142857142857	1.49893673251356e-07\\
0.751224489795918	1.7045871925081e-07\\
0.75530612244898	1.82369481405154e-07\\
0.759387755102041	1.84540715375547e-07\\
0.763469387755102	1.76892622438807e-07\\
0.767551020408163	1.60291097617282e-07\\
0.771632653061225	1.366376599643e-07\\
0.775714285714286	1.09454669772191e-07\\
0.779795918367347	8.55059926196944e-08\\
0.783877551020408	7.64146116383202e-08\\
0.787959183673469	8.90081645055218e-08\\
0.792040816326531	1.14079102053743e-07\\
0.796122448979592	1.40746864829566e-07\\
0.800204081632653	1.63020873911179e-07\\
0.804285714285714	1.77776300569477e-07\\
0.808367346938776	1.8347772268168e-07\\
0.812448979591837	1.79704720881801e-07\\
0.816530612244898	1.67001235951147e-07\\
0.820612244897959	1.46889302055087e-07\\
0.82469387755102	1.22216541553755e-07\\
0.828775510204082	9.81613975659101e-08\\
0.832857142857143	8.36534654879983e-08\\
0.836938775510204	8.75493715302156e-08\\
0.841020408163265	1.0697928387593e-07\\
0.845102040816327	1.31772542655284e-07\\
0.849183673469388	1.54644642869073e-07\\
0.853265306122449	1.71552996142132e-07\\
0.85734693877551	1.80371144020697e-07\\
0.861428571428571	1.80242210268133e-07\\
0.865510204081633	1.71353951916586e-07\\
0.869591836734694	1.54882079383967e-07\\
0.873673469387755	1.3318410951868e-07\\
0.877755102040816	1.10421765397994e-07\\
0.881836734693878	9.35415318537025e-08\\
0.885918367346939	9.08421663448067e-08\\
0.89	1.03824028374611e-07\\
0.894081632653061	1.2509202818034e-07\\
0.898163265306122	1.46937033984607e-07\\
0.902244897959184	1.64580415518145e-07\\
0.906326530612245	1.75396978052689e-07\\
0.910408163265306	1.7819089825899e-07\\
0.914489795918367	1.72846769896928e-07\\
0.918571428571429	1.60221055867105e-07\\
0.92265306122449	1.42208946396494e-07\\
0.926734693877551	1.22130397572917e-07\\
0.930816326530612	1.05320984489232e-07\\
0.934897959183673	9.85196108027608e-08\\
0.938979591836735	1.05261653726707e-07\\
0.943061224489796	1.21820933784988e-07\\
0.947142857142857	1.41477771276044e-07\\
0.951224489795918	1.59004973854238e-07\\
0.95530612244898	1.71194969755817e-07\\
0.959387755102041	1.76366909854231e-07\\
0.963469387755102	1.74002375286086e-07\\
0.967551020408163	1.6459238949551e-07\\
0.971632653061224	1.49581390829212e-07\\
0.975714285714286	1.31631549438183e-07\\
0.979795918367347	1.15006576706378e-07\\
0.983877551020408	1.05401594230337e-07\\
0.987959183673469	1.07204258914528e-07\\
0.992040816326531	1.1940338229266e-07\\
0.996122448979592	1.36723441679365e-07\\
};
\addlegendentry{ABC$_W^{1/2}$ adaptive}

\addplot [color=mycolor1, dashed]
  table[row sep=crcr]{%
0.000204081632653064	0\\
0.0124489795918368	9.30855392766716e-11\\
0.0328571428571428	3.99368538239742e-10\\
0.0410204081632654	5.06275354972274e-10\\
0.0491836734693878	8.78358163980408e-10\\
0.0655102040816327	1.93574511975214e-09\\
0.0736734693877551	2.27844820788192e-09\\
0.0818367346938775	2.55820098438164e-09\\
0.0859183673469388	2.77207146126557e-09\\
0.09	3.09003478360381e-09\\
0.0940816326530612	3.52595530550559e-09\\
0.102244897959184	4.67235605761118e-09\\
0.110408163265306	5.89675630546793e-09\\
0.114489795918367	6.42556330365807e-09\\
0.118571428571429	6.85945533707155e-09\\
0.126734693877551	7.46617345726008e-09\\
0.134897959183674	8.08788958028828e-09\\
0.138979591836735	8.61585669476028e-09\\
0.143061224489796	9.3677386958646e-09\\
0.147142857142857	1.03329809153507e-08\\
0.159387755102041	1.3692421640954e-08\\
0.163469387755102	1.46427363567625e-08\\
0.167551020408163	1.53962981253031e-08\\
0.171632653061225	1.59507518304736e-08\\
0.183877551020408	1.72680331145258e-08\\
0.187959183673469	1.8068830764939e-08\\
0.192040816326531	1.92461321324444e-08\\
0.196122448979592	2.07907282501196e-08\\
0.204285714285714	2.44544181571271e-08\\
0.208367346938776	2.60143436792148e-08\\
0.212448979591837	2.72503473119912e-08\\
0.216530612244898	2.82275093299589e-08\\
0.220612244897959	2.88892372246607e-08\\
0.232857142857143	3.02133132956683e-08\\
0.236938775510204	3.10314552942614e-08\\
0.241020408163265	3.22653141981988e-08\\
0.245102040816327	3.38652419440422e-08\\
0.253265306122449	3.74823094784205e-08\\
0.25734693877551	3.90731130783806e-08\\
0.261428571428571	4.03280258076677e-08\\
0.265510204081633	4.11854020887503e-08\\
0.269591836734694	4.1740304768112e-08\\
0.273673469387755	4.21764932978519e-08\\
0.277755102040816	4.2731368998794e-08\\
0.281836734693878	4.36308941198504e-08\\
0.285918367346939	4.50626252979092e-08\\
0.29	4.7036666561695e-08\\
0.294081632653061	4.9421572501096e-08\\
0.306326530612245	5.71398179882721e-08\\
0.310408163265306	5.93519122737973e-08\\
0.314489795918367	6.12585248083164e-08\\
0.318571428571429	6.29131239326952e-08\\
0.326734693877551	6.60820598241685e-08\\
0.330816326530612	6.79687160998199e-08\\
0.334897959183673	7.02991808099185e-08\\
0.338979591836735	7.31937175313391e-08\\
0.343061224489796	7.66937310325844e-08\\
0.347142857142857	8.07147064740832e-08\\
0.35530612244898	8.92802273222415e-08\\
0.359387755102041	9.30193696513371e-08\\
0.363469387755102	9.5883807116337e-08\\
0.367551020408163	9.76380014616751e-08\\
0.371632653061225	9.82610754984847e-08\\
0.379795918367347	9.78516595528944e-08\\
0.383877551020408	9.87035945199111e-08\\
0.387959183673469	1.01748566128101e-07\\
0.392040816326531	1.07579662400603e-07\\
0.396122448979592	1.15871624850605e-07\\
0.400204081632653	1.25453051724733e-07\\
0.404285714285714	1.34712975663476e-07\\
0.408367346938776	1.41983254819955e-07\\
0.412448979591837	1.45833564224596e-07\\
0.416530612244898	1.45317484712137e-07\\
0.420612244897959	1.40213394117517e-07\\
0.42469387755102	1.31294159433359e-07\\
0.428775510204082	1.2069507415724e-07\\
0.432857142857143	1.1223481999334e-07\\
0.436938775510204	1.1073669081707e-07\\
0.441020408163265	1.19137539833503e-07\\
0.445102040816327	1.35860687255906e-07\\
0.449183673469388	1.56335058787427e-07\\
0.453265306122449	1.75809411029348e-07\\
0.45734693877551	1.904734094893e-07\\
0.461428571428571	1.97586193295685e-07\\
0.465510204081633	1.954515772562e-07\\
0.469591836734694	1.83508814588507e-07\\
0.473673469387755	1.6253998835758e-07\\
0.477755102040816	1.35183630400881e-07\\
0.481836734693878	1.07341411204054e-07\\
0.485918367346939	9.09178785590825e-08\\
0.49	9.9017614774155e-08\\
0.494081632653061	1.2723962694583e-07\\
0.498163265306122	1.61622064065448e-07\\
0.502244897959184	1.92950498356659e-07\\
0.506326530612245	2.16053400059479e-07\\
0.510408163265306	2.27904693450576e-07\\
0.514489795918367	2.269470452676e-07\\
0.518571428571429	2.12905856611378e-07\\
0.52265306122449	1.86805222512909e-07\\
0.526734693877551	1.51187027563893e-07\\
0.530816326530612	1.11130476820875e-07\\
0.534897959183673	7.85974498818121e-08\\
0.538979591836735	7.79042192977641e-08\\
0.543061224489796	1.10516689688644e-07\\
0.547142857142857	1.52272873310544e-07\\
0.551224489795918	1.90485197060219e-07\\
0.55530612244898	2.19553948244133e-07\\
0.559387755102041	2.36414590792933e-07\\
0.563469387755102	2.39456483863698e-07\\
0.567551020408163	2.28290762849248e-07\\
0.571632653061224	2.03733714565502e-07\\
0.575714285714286	1.6792059864823e-07\\
0.579795918367347	1.24848397176081e-07\\
0.583877551020408	8.31823940838916e-08\\
0.587959183673469	6.66532805615461e-08\\
0.592040816326531	9.39590206749941e-08\\
0.596122448979592	1.37808816469409e-07\\
0.600204081632653	1.79915110476436e-07\\
0.604285714285714	2.13507088453468e-07\\
0.608367346938775	2.35173530760413e-07\\
0.612448979591837	2.43082060702449e-07\\
0.616530612244898	2.3660732129116e-07\\
0.620612244897959	2.16268723818303e-07\\
0.62469387755102	1.83757819405272e-07\\
0.628775510204082	1.4220513788743e-07\\
0.632857142857143	9.75843933437304e-08\\
0.636938775510204	6.59073945419308e-08\\
0.641020408163265	7.73426952482836e-08\\
0.645102040816326	1.18464142273567e-07\\
0.649183673469388	1.62764722699826e-07\\
0.653265306122449	2.00622236046222e-07\\
0.65734693877551	2.27741066005649e-07\\
0.661428571428571	2.41731821559732e-07\\
0.665510204081633	2.41467698502973e-07\\
0.669591836734694	2.26983680184922e-07\\
0.673673469387755	1.99493571440001e-07\\
0.677755102040816	1.61573936008352e-07\\
0.681836734693878	1.17981417413482e-07\\
0.685918367346939	7.95532942854393e-08\\
0.69	7.24837599941353e-08\\
0.694081632653061	1.04687043145368e-07\\
0.698163265306122	1.48280213108087e-07\\
0.702244897959184	1.88505944631601e-07\\
0.706326530612245	2.19481730123761e-07\\
0.710408163265306	2.38157856324861e-07\\
0.714489795918367	2.42997611255014e-07\\
0.718571428571429	2.33687842632335e-07\\
0.72265306122449	2.11097341940025e-07\\
0.726734693877551	1.77369679943418e-07\\
0.730816326530612	1.36384069371331e-07\\
0.734897959183674	9.59224795149538e-08\\
0.738979591836735	7.48377704251624e-08\\
0.743061224489796	9.31858973318711e-08\\
0.747142857142857	1.3287892364211e-07\\
0.751224489795918	1.73926369750887e-07\\
0.75530612244898	2.08058108075448e-07\\
0.759387755102041	2.31283090901258e-07\\
0.763469387755102	2.41463772199246e-07\\
0.767551020408163	2.37783618572962e-07\\
0.771632653061225	2.20626015012826e-07\\
0.775714285714286	1.91638890423285e-07\\
0.779795918367347	1.54070317615052e-07\\
0.783877551020408	1.14075283996762e-07\\
0.787959183673469	8.52236493598113e-08\\
0.792040816326531	8.89587659091973e-08\\
0.796122448979592	1.21472597092165e-07\\
0.800204081632653	1.61265750753259e-07\\
0.804285714285714	1.96931237650411e-07\\
0.808367346938776	2.23254021602948e-07\\
0.812448979591837	2.37506353917283e-07\\
0.816530612244898	2.38455578616126e-07\\
0.820612244897959	2.26123609059492e-07\\
0.82469387755102	2.0176633641622e-07\\
0.828775510204082	1.6808223046727e-07\\
0.832857142857143	1.3001493992526e-07\\
0.836938775510204	9.74881331217148e-08\\
0.841020408163265	8.87323078346824e-08\\
0.845102040816327	1.11226336918513e-07\\
0.849183673469388	1.47716565312983e-07\\
0.853265306122449	1.83860537572578e-07\\
0.85734693877551	2.1280535589252e-07\\
0.861428571428571	2.31023595920199e-07\\
0.865510204081633	2.36745744786049e-07\\
0.869591836734694	2.29531134321448e-07\\
0.873673469387755	2.10191550720218e-07\\
0.877755102040816	1.80907338886982e-07\\
0.881836734693878	1.45744894042288e-07\\
0.885918367346939	1.12327335433093e-07\\
0.89	9.48340306283768e-08\\
0.894081632653061	1.05940097805401e-07\\
0.898163265306122	1.3665159714904e-07\\
0.902244897959184	1.7148273667722e-07\\
0.906326530612245	2.01711158442919e-07\\
0.910408163265306	2.22905658819883e-07\\
0.914489795918367	2.32778878617168e-07\\
0.918571428571429	2.30475608153924e-07\\
0.92265306122449	2.16391573881758e-07\\
0.926734693877551	1.92204112070549e-07\\
0.930816326530612	1.61198180470556e-07\\
0.934897959183673	1.29239586366126e-07\\
0.938979591836735	1.06757235718291e-07\\
0.943061224489796	1.0653774740188e-07\\
0.947142857142857	1.2854102626747e-07\\
0.951224489795918	1.59931830667759e-07\\
0.95530612244898	1.90316372661314e-07\\
0.959387755102041	2.13930858916633e-07\\
0.963469387755102	2.27636439475987e-07\\
0.967551020408163	2.29939576712468e-07\\
0.971632653061224	2.20653758709055e-07\\
0.975714285714286	2.00883644496663e-07\\
0.979795918367347	1.73202708064757e-07\\
0.983877551020408	1.4230650180469e-07\\
0.987959183673469	1.16473495070935e-07\\
0.992040816326531	1.07744905553098e-07\\
0.996122448979592	1.21831291943764e-07\\
};
\addlegendentry{ABC$_W^{1/2}$}

\end{axis}
\end{tikzpicture}%

%% file: images/Horn3D/Cost.tex
%
%
\definecolor{mycolor1}{rgb}{1.00000,0.60000,0.00000}%
\begin{tikzpicture}[scale = 0.52, font=\huge]

\begin{axis}[%
width=4.602in,
height=3.82in,
at={(0.772in,0.516in)},
scale only axis,
xmin=0,
xmax=1,
xlabel style={font=\huge\color{white!15!black}},
xlabel={relative simulation time $t/T$}, xtick={0, 0.2, 0.4, 0.6, 0.8, 1},
y tick label style={
        /pgf/number format/.cd,
            fixed,
            fixed zerofill,
            precision=2,
        /tikz/.cd
            },
ymin=0.02,
ymax=0.2,
ylabel style={at={(-0.15,0.5)}, font=\huge\color{white!15!black}},  ytick={0.04, 0.08, 0.12, 0.16, 0.2},
ylabel={relative $L^{2}(\Omega)$ error}, 
axis background/.style={fill=white},
xmajorgrids,
ymajorgrids,
legend style={at={(0.03,0.97)}, anchor=north west, legend cell align=left, align=left, draw=white!15!black}
]
\addplot [color=mycolor1, line width=2.0pt]
  table[row sep=crcr]{%
0.00476190476190474	0.195526817761151\\
0.00929705215419496	0.125782721718886\\
0.0138321995464853	0.0871449782559917\\
0.0183673469387755	0.0650653146692521\\
0.0229024943310657	0.046905034648709\\
0.0274376417233561	0.0327622147153339\\
0.0319727891156463	0.0253369009768525\\
0.0365079365079365	0.0301743786905428\\
0.0410430839002267	0.0466264595907681\\
0.045578231292517	0.0671688859905816\\
0.0501133786848073	0.0719654159978033\\
0.0546485260770975	0.0558905461464337\\
0.0591836734693878	0.038714115985078\\
0.063718820861678	0.0268896959646023\\
0.0682539682539682	0.0210559726973538\\
0.0727891156462586	0.0213368616105069\\
0.0773242630385488	0.0260491433984756\\
0.081859410430839	0.0329712554515228\\
0.0863945578231292	0.0402395433476007\\
0.0909297052154195	0.0447560133642134\\
0.0954648526077098	0.0426223621632168\\
0.10453514739229	0.027020968283881\\
0.109070294784581	0.023110783553804\\
0.113605442176871	0.0235365580583089\\
0.118140589569161	0.0266914627118093\\
0.122675736961451	0.0309315571009762\\
0.127210884353742	0.0350828527497536\\
0.131746031746032	0.0377952730023066\\
0.136281179138322	0.0374062592968304\\
0.140816326530612	0.0334854686330347\\
0.145351473922902	0.028349515555351\\
0.149886621315193	0.025097317162956\\
0.154421768707483	0.0248845077064632\\
0.158956916099773	0.0268927103509734\\
0.168027210884354	0.0328600088908604\\
0.172562358276644	0.034934519836673\\
0.177097505668934	0.0349809256972601\\
0.181632653061225	0.0324582224031906\\
0.186167800453515	0.0291550368970745\\
0.190702947845805	0.0277713684901617\\
0.195238095238095	0.0286276914984496\\
0.204308390022676	0.0318663861336048\\
0.208843537414966	0.0331121380085903\\
0.213378684807256	0.0332932985342009\\
0.217913832199547	0.0321609986548876\\
0.226984126984127	0.0273083013783635\\
0.231519274376417	0.0263339205586449\\
0.236054421768708	0.026989729901972\\
0.245124716553288	0.0297191258405567\\
0.249659863945578	0.030578449177978\\
0.254195011337868	0.0304491673236403\\
0.258730158730159	0.0291775766067498\\
0.263265306122449	0.0270215123760247\\
0.267800453514739	0.0247846611060698\\
0.272335600907029	0.023581948253489\\
0.27687074829932	0.0241039603449043\\
0.28140589569161	0.0257547638423863\\
0.2859410430839	0.0274888649716792\\
0.290476190476191	0.02857921104523\\
0.295011337868481	0.0285101004498025\\
0.299546485260771	0.0272293696292595\\
0.304081632653061	0.0251489391852009\\
0.308616780045352	0.0229648705171024\\
0.313151927437642	0.021790817048831\\
0.317687074829932	0.0223005765213546\\
0.322222222222222	0.0240479531355959\\
0.326757369614513	0.0260182227273259\\
0.331292517006803	0.0273115360604821\\
0.335827664399093	0.0275225978154532\\
0.340362811791383	0.026663104284359\\
0.344897959183673	0.0250554181870897\\
0.349433106575964	0.0233726075081772\\
0.353968253968254	0.0224754064494548\\
0.358503401360544	0.0228989205055885\\
0.363038548752834	0.0243507936803364\\
0.367573696145125	0.0260385412175629\\
0.372108843537415	0.0273328821171365\\
0.376643990929705	0.0278320346246679\\
0.381179138321995	0.0273892669497042\\
0.385714285714286	0.0262760915688017\\
0.390249433106576	0.0250162328987371\\
0.394784580498866	0.0242740283987818\\
0.399319727891156	0.0244365834702162\\
0.403854875283447	0.0255407001618224\\
0.408390022675737	0.0270005382635078\\
0.412925170068027	0.0282403703293306\\
0.417460317460318	0.0289884591999628\\
0.421995464852608	0.0291390050422513\\
0.431065759637188	0.0286354814775666\\
0.435600907029479	0.0286787494423123\\
0.440136054421769	0.0293245884953978\\
0.444671201814059	0.0304919108968691\\
0.449206349206349	0.031804840533472\\
0.453741496598639	0.0329440781239788\\
0.45827664399093	0.0337997440487031\\
0.46281179138322	0.0343264046032019\\
0.46734693877551	0.0347032769636975\\
0.4718820861678	0.0351466962278555\\
0.476417233560091	0.0358158326319267\\
0.480952380952381	0.0367598176912401\\
0.490022675736961	0.0388917988332773\\
0.494557823129252	0.0397027749639518\\
0.499092970521542	0.0402477675052245\\
0.508163265306122	0.0410444784077327\\
0.512698412698413	0.0416820675591132\\
0.517233560090703	0.0425932416665119\\
0.526303854875283	0.0447775505061676\\
0.530839002267574	0.0456191504976341\\
0.535374149659864	0.0460423717327773\\
0.544444444444444	0.0462835081277931\\
0.548979591836735	0.0466885894149626\\
0.553514739229025	0.0475585007784912\\
0.562585034013605	0.0498519082313499\\
0.567120181405896	0.0505142302592061\\
0.571655328798186	0.0505645852405276\\
0.576190476190476	0.0502648935776492\\
0.580725623582766	0.0500585198906824\\
0.585260770975057	0.0503324102243484\\
0.589795918367347	0.0512239069973273\\
0.594331065759637	0.0524867027614575\\
0.598866213151927	0.0536291287311474\\
0.603401360544218	0.0542552034471885\\
0.607936507936508	0.05416506155171\\
0.617006802721088	0.0528642917422317\\
0.621541950113379	0.0528286863154386\\
0.626077097505669	0.0537056624290225\\
0.630612244897959	0.0553500984056332\\
0.635147392290249	0.0571891915661336\\
0.63968253968254	0.0584019659158929\\
0.64421768707483	0.0584032191424715\\
0.64875283446712	0.0570723521683788\\
0.65328798185941	0.054815642161778\\
0.657823129251701	0.0528255257304501\\
0.662358276643991	0.0525100149424957\\
0.666893424036281	0.0544097663735668\\
0.671428571428571	0.0578521130971256\\
0.675963718820862	0.0614639606258134\\
0.680498866213152	0.0638128704884268\\
0.685034013605442	0.0639035512293634\\
0.689569160997732	0.0615543441240308\\
0.694104308390023	0.0575011924824425\\
0.698639455782313	0.0535367104774926\\
0.703174603174603	0.0520187261715513\\
0.707709750566893	0.0540982797660966\\
0.712244897959184	0.0586834523857307\\
0.716780045351474	0.0636564418101853\\
0.721315192743764	0.067124715359649\\
0.725850340136054	0.0678489179712842\\
0.730385487528345	0.0654198175468487\\
0.734920634920635	0.0603791544376299\\
0.739455782312925	0.0544758529539209\\
0.743990929705215	0.0506452372128421\\
0.748526077097506	0.0513760422756193\\
0.753061224489796	0.0562345014947638\\
0.757596371882086	0.0625645838992738\\
0.762131519274376	0.0677530420644066\\
0.766666666666667	0.0700791022132224\\
0.771201814058957	0.068773682154961\\
0.775736961451247	0.064082547517283\\
0.780272108843537	0.0574004133304833\\
0.784807256235828	0.0515074482989218\\
0.789342403628118	0.0498521093134745\\
0.793877551020408	0.0535548444170065\\
0.798412698412698	0.0602356932420676\\
0.802947845804989	0.066713150852708\\
0.807482993197279	0.0707091995690412\\
0.812018140589569	0.0709698075859464\\
0.816553287981859	0.0672998822528742\\
0.82108843537415	0.0606185918935465\\
0.82562358276644	0.0532696991785577\\
0.83015873015873	0.0489808033825188\\
0.83469387755102	0.0506259416083619\\
0.839229024943311	0.0569572880643401\\
0.843764172335601	0.0643591988079445\\
0.848299319727891	0.0699521641821639\\
0.852834467120181	0.0721062770889076\\
0.857369614512472	0.0702979220745705\\
0.861904761904762	0.0650161451274073\\
0.866439909297052	0.0578381797703946\\
0.870975056689342	0.0516584300925744\\
0.875510204081633	0.0500090802217636\\
0.880045351473923	0.0539852956473188\\
0.884580498866213	0.0610276678692755\\
0.889115646258503	0.0678263150978318\\
0.893650793650794	0.0720892419436223\\
0.898185941043084	0.072642299633832\\
0.902721088435374	0.0693566348684067\\
0.907256235827664	0.0630800776556786\\
0.911791383219955	0.05583172997125\\
0.916326530612245	0.0510395463197888\\
0.920861678004535	0.0517993753468208\\
0.925396825396825	0.0575236049098283\\
0.929931972789116	0.0649282068511686\\
0.934467120181406	0.071021488091578\\
0.939002267573696	0.0739603805351658\\
0.943537414965986	0.0729791356620365\\
0.948072562358277	0.0683504386040831\\
0.952607709750567	0.061380756633308\\
0.957142857142857	0.0546616299113554\\
0.961678004535147	0.051835661114885\\
0.966213151927438	0.0548579540109212\\
0.970748299319728	0.0617525870671171\\
0.975283446712018	0.0690401732098603\\
0.979818594104308	0.0741353905746842\\
0.984353741496599	0.0755960165828879\\
0.988888888888889	0.073024801733045\\
0.993424036281179	0.0670495927638652\\
0.997959183673469	0.0595290546533354\\
};
\addlegendentry{ABC$_W^{1/2}$ adaptive}

\addplot [color=mycolor1, dashed]
  table[row sep=crcr]{%
0.00476190476190474	0.195526817761151\\
0.00929705215419496	0.125782721718886\\
0.0138321995464853	0.0871449782559917\\
0.0183673469387755	0.0650653146692521\\
0.0229024943310657	0.046905034648709\\
0.0274376417233561	0.0327622147153339\\
0.0319727891156463	0.0253369009768525\\
0.0365079365079365	0.0301743786905428\\
0.0410430839002267	0.0466264595907681\\
0.045578231292517	0.0671688859905816\\
0.0501133786848073	0.0719654159978033\\
0.0546485260770975	0.0558905461464337\\
0.0591836734693878	0.038714115985078\\
0.063718820861678	0.0268896959646023\\
0.0682539682539682	0.0210559726973538\\
0.0727891156462586	0.0213368616105069\\
0.0773242630385488	0.0260491433984756\\
0.081859410430839	0.0329712554515228\\
0.0863945578231292	0.0402395433476007\\
0.0909297052154195	0.0447560133642134\\
0.0954648526077098	0.0426223621632168\\
0.10453514739229	0.027020968283881\\
0.109070294784581	0.023110783553804\\
0.113605442176871	0.0235365580583089\\
0.118140589569161	0.0266914627118093\\
0.122675736961451	0.0309315571009762\\
0.127210884353742	0.0350828527497536\\
0.131746031746032	0.0377952730023066\\
0.136281179138322	0.0374062592968304\\
0.140816326530612	0.0334854686330347\\
0.145351473922902	0.028349515555351\\
0.149886621315193	0.025097317162956\\
0.154421768707483	0.0248845077064632\\
0.158956916099773	0.0268927103509734\\
0.168027210884354	0.0328600088908604\\
0.172562358276644	0.034934519836673\\
0.177097505668934	0.0349809256972601\\
0.181632653061225	0.0324582224031906\\
0.186167800453515	0.0291550368970745\\
0.190702947845805	0.0277713684901617\\
0.195238095238095	0.0286276914984496\\
0.204308390022676	0.0318663861336048\\
0.208843537414966	0.0331121380085903\\
0.213378684807256	0.0332932985342009\\
0.217913832199547	0.0321609986548876\\
0.226984126984127	0.0273083013783635\\
0.231519274376417	0.0263339205586449\\
0.236054421768708	0.026989729901972\\
0.245124716553288	0.0297191258405567\\
0.249659863945578	0.030578449177978\\
0.254195011337868	0.0304491673236403\\
0.258730158730159	0.0291775766067498\\
0.263265306122449	0.0270215123760247\\
0.267800453514739	0.0247846611060698\\
0.272335600907029	0.023581948253489\\
0.27687074829932	0.0241039603449043\\
0.28140589569161	0.0257547638423863\\
0.2859410430839	0.0274888649716792\\
0.290476190476191	0.02857921104523\\
0.295011337868481	0.0285101004498025\\
0.299546485260771	0.0272293696292595\\
0.304081632653061	0.0251489391852009\\
0.308616780045352	0.0229648705171024\\
0.313151927437642	0.021790817048831\\
0.317687074829932	0.0223005765213546\\
0.322222222222222	0.0240479531355959\\
0.326757369614513	0.0260182227273259\\
0.331292517006803	0.0273115360604821\\
0.335827664399093	0.0275225978154532\\
0.340362811791383	0.026663104284359\\
0.344897959183673	0.0250554181870897\\
0.349433106575964	0.0233726075081772\\
0.353968253968254	0.0224754064494548\\
0.358503401360544	0.0228989205055885\\
0.363038548752834	0.0243507936803364\\
0.367573696145125	0.0260385412175629\\
0.372108843537415	0.0273328821171365\\
0.376643990929705	0.0278320407439292\\
0.381179138321995	0.0273885410132495\\
0.385714285714286	0.0262680142979413\\
0.390249433106576	0.0249847526557083\\
0.394784580498866	0.0241978664066594\\
0.399319727891156	0.024304196004742\\
0.403854875283447	0.0253633572853487\\
0.408390022675737	0.0268086576532571\\
0.412925170068027	0.028039428152676\\
0.417460317460318	0.0287722417210882\\
0.421995464852608	0.0288884128190837\\
0.426530612244898	0.0285844841533198\\
0.431065759637188	0.0281798055422734\\
0.435600907029479	0.0280388938506763\\
0.440136054421769	0.0284927632896199\\
0.444671201814059	0.0295245825429722\\
0.449206349206349	0.0307959673567223\\
0.453741496598639	0.0319825171605118\\
0.45827664399093	0.0329453594013743\\
0.46281179138322	0.0335611729712351\\
0.4718820861678	0.0342865247629138\\
0.476417233560091	0.0347934788043748\\
0.480952380952381	0.0355977528725392\\
0.485487528344671	0.0366716112678941\\
0.490022675736961	0.0378805285404653\\
0.494557823129252	0.0390114652183842\\
0.499092970521542	0.039956854457631\\
0.512698412698413	0.0422269862909813\\
0.517233560090703	0.0432288372634544\\
0.521768707482993	0.0444128638210933\\
0.526303854875283	0.0456652838760759\\
0.530839002267574	0.0468096617824243\\
0.535374149659864	0.0476549598662657\\
0.539909297052154	0.0482159582568162\\
0.544444444444444	0.0486607057736339\\
0.548979591836735	0.049216202452528\\
0.553514739229025	0.0501531228742937\\
0.558049886621315	0.0515301527917738\\
0.562585034013605	0.0531530904688785\\
0.567120181405896	0.0545960033871667\\
0.571655328798186	0.055390733855693\\
0.576190476190476	0.0553343197779397\\
0.585260770975057	0.054032354221786\\
0.589795918367347	0.0542430333461106\\
0.594331065759637	0.0558017417442234\\
0.598866213151927	0.0585425522955357\\
0.603401360544218	0.0616126907820714\\
0.607936507936508	0.0638244181600885\\
0.612471655328798	0.0641761995456758\\
0.617006802721088	0.062326721974142\\
0.621541950113379	0.0589223609016466\\
0.626077097505669	0.0556231420711852\\
0.630612244897959	0.0545672127649951\\
0.635147392290249	0.0570647624504146\\
0.63968253968254	0.0624565655035397\\
0.64421768707483	0.0685914040088574\\
0.64875283446712	0.073103100370193\\
0.65328798185941	0.074156933243831\\
0.657823129251701	0.071032958754816\\
0.662358276643991	0.0644788573490966\\
0.666893424036281	0.0569491754136559\\
0.671428571428571	0.0523619495592963\\
0.675963718820862	0.0540978409601096\\
0.680498866213152	0.0616944873153566\\
0.685034013605442	0.0714539638309373\\
0.689569160997732	0.0796101171863886\\
0.694104308390023	0.0834731468456171\\
0.698639455782313	0.0816394750877281\\
0.703174603174603	0.0742721219657458\\
0.707709750566893	0.0634594589810256\\
0.712244897959184	0.0535916580306598\\
0.716780045351474	0.0507980731166392\\
0.721315192743764	0.0576167123238487\\
0.725850340136054	0.0695945890160299\\
0.730385487528345	0.0810128086451097\\
0.734920634920635	0.0880258803603099\\
0.739455782312925	0.088556922856757\\
0.743990929705215	0.0822265534223572\\
0.748526077097506	0.0705729882360016\\
0.753061224489796	0.0574632586599807\\
0.757596371882086	0.0495924548090126\\
0.762131519274376	0.0530148887913856\\
0.766666666666667	0.0650845199335991\\
0.771201814058957	0.0787144180912829\\
0.775736961451247	0.0889257711756393\\
0.780272108843537	0.0929109323889293\\
0.784807256235828	0.089542743994935\\
0.789342403628118	0.0795037881282186\\
0.793877551020408	0.0656055390636103\\
0.798412698412698	0.0534174468383679\\
0.802947845804989	0.0507377838886321\\
0.807482993197279	0.0597475815434345\\
0.812018140589569	0.0737454468229691\\
0.816553287981859	0.0861848372225706\\
0.82108843537415	0.09327001242491\\
0.82562358276644	0.0931620890607462\\
0.83015873015873	0.0857297543197089\\
0.83469387755102	0.0728166647793138\\
0.839229024943311	0.0587279547361581\\
0.843764172335601	0.0507475246410747\\
0.848299319727891	0.0549784902536106\\
0.852834467120181	0.067859291901411\\
0.857369614512472	0.0818934657800382\\
0.861904761904762	0.0921048350231313\\
0.866439909297052	0.0957630208430964\\
0.870975056689342	0.0918061430836228\\
0.875510204081633	0.080976955523397\\
0.880045351473923	0.0662885986582132\\
0.884580498866213	0.0538266268748199\\
0.889115646258503	0.0518981737668077\\
0.893650793650794	0.0620186478210466\\
0.898185941043084	0.0767342338845525\\
0.902721088435374	0.0894577657405122\\
0.907256235827664	0.0964698632687611\\
0.911791383219955	0.0959330679614212\\
0.916326530612245	0.0877466360628858\\
0.920861678004535	0.0739125144374037\\
0.925396825396825	0.0591784482045472\\
0.929931972789116	0.0516401819009964\\
0.934467120181406	0.057279882892866\\
0.939002267573696	0.0712819415120717\\
0.943537414965986	0.0857530456186991\\
0.948072562358277	0.0958253851464042\\
0.952607709750567	0.0988842729884294\\
0.957142857142857	0.0939810892486704\\
0.961678004535147	0.0820916535802589\\
0.966213151927438	0.0667067232099452\\
0.970748299319728	0.0546721964131539\\
0.975283446712018	0.0544413311130899\\
0.979818594104308	0.0660343458226288\\
0.984353741496599	0.0813483476299742\\
0.988888888888889	0.0940099682436331\\
0.993424036281179	0.100475561704407\\
0.997959183673469	0.0990377806624997\\
};
\addlegendentry{ABC$_W^{1/2}$}

\end{axis}
\end{tikzpicture}%

%% file: images/Horn3D/Cost_abs.tex
%
%
\definecolor{mycolor1}{rgb}{1.00000,0.60000,0.00000}%
\begin{tikzpicture}[scale = 0.52, font=\huge]

\begin{axis}[%
width=4.602in,
height=3.82in,
at={(0.772in,0.516in)},
scale only axis,
xmin=0,
xmax=1,
xlabel style={font=\huge\color{white!15!black}},
xlabel={relative simulation time $t/T$}, xtick={0, 0.2, 0.4, 0.6, 0.8, 1},
y tick label style={
        /pgf/number format/.cd,
            fixed,
            fixed zerofill,
            precision=1,
        /tikz/.cd
            },
ymin=0,
ymax=1.2e-05,
ylabel style={font=\huge\color{white!15!black}}, ytick={0.000002, 0.000004, 0.000006, 0.000008, 0.00001, 0.000012},
ylabel={absolute $L^{2}(\Omega)$ error},
axis background/.style={fill=white},
xmajorgrids,
ymajorgrids,
legend style={at={(0.03,0.97)}, anchor=north west, legend cell align=left, align=left, draw=white!15!black}
]
\addplot [color=mycolor1, line width=2.0pt]
  table[row sep=crcr]{%
0.000226757369614528	0\\
0.00929705215419496	2.77857203911935e-09\\
0.0274376417233561	1.45164456011315e-08\\
0.0319727891156463	1.4862199582133e-08\\
0.0365079365079365	2.07625842074677e-08\\
0.045578231292517	4.98861908448944e-08\\
0.0501133786848073	6.44927802273543e-08\\
0.0546485260770975	7.47726257577241e-08\\
0.0591836734693878	7.85579905571154e-08\\
0.0682539682539682	7.77235572568813e-08\\
0.0727891156462586	9.42000372150531e-08\\
0.0773242630385488	1.28235129071008e-07\\
0.0863945578231292	2.1032581742908e-07\\
0.0909297052154195	2.38932429175165e-07\\
0.0954648526077098	2.50722057337782e-07\\
0.10453514739229	2.38312900746607e-07\\
0.109070294784581	2.48327996277276e-07\\
0.113605442176871	2.94605337569642e-07\\
0.118140589569161	3.70081384870602e-07\\
0.122675736961451	4.53218788498688e-07\\
0.127210884353742	5.24059345585748e-07\\
0.131746031746032	5.68468046524551e-07\\
0.136281179138322	5.80551164586218e-07\\
0.145351473922902	5.49135144800594e-07\\
0.149886621315193	5.64248462064931e-07\\
0.154421768707483	6.38197102942328e-07\\
0.158956916099773	7.6105095458967e-07\\
0.163492063492063	8.98726536191496e-07\\
0.168027210884354	1.01636340332778e-06\\
0.172562358276644	1.08876842463967e-06\\
0.177097505668934	1.1006133716629e-06\\
0.186167800453515	1.01033811883422e-06\\
0.190702947845805	1.04483511309539e-06\\
0.195238095238095	1.16830400809231e-06\\
0.199773242630386	1.3180132405477e-06\\
0.204308390022676	1.45459500844325e-06\\
0.208843537414966	1.5492133235151e-06\\
0.213378684807256	1.57283085266968e-06\\
0.217913832199547	1.52578984324858e-06\\
0.222448979591837	1.42591275908988e-06\\
0.226984126984127	1.3384842617592e-06\\
0.231519274376417	1.33881434161154e-06\\
0.236054421768708	1.43392288343147e-06\\
0.240589569160998	1.57287582758237e-06\\
0.245124716553288	1.7049898418664e-06\\
0.249659863945578	1.79321125415299e-06\\
0.254195011337868	1.80461463505921e-06\\
0.258730158730159	1.73551493432011e-06\\
0.267800453514739	1.49398410909107e-06\\
0.272335600907029	1.44975070293096e-06\\
0.27687074829932	1.52391825292142e-06\\
0.28140589569161	1.67956542229675e-06\\
0.2859410430839	1.84332363650785e-06\\
0.290476190476191	1.95561435989422e-06\\
0.295011337868481	1.9729557868553e-06\\
0.299546485260771	1.89225519253711e-06\\
0.308616780045352	1.60552634120847e-06\\
0.313151927437642	1.54052946910177e-06\\
0.317687074829932	1.60618928612966e-06\\
0.322222222222222	1.77276593238407e-06\\
0.326757369614513	1.96295817611603e-06\\
0.331292517006803	2.09934836548076e-06\\
0.335827664399093	2.14057555136815e-06\\
0.340362811791383	2.08435446502619e-06\\
0.344897959183673	1.96095840510058e-06\\
0.349433106575964	1.83198281411556e-06\\
0.353968253968254	1.77200418738366e-06\\
0.358503401360544	1.8276067231815e-06\\
0.363038548752834	1.97792128320273e-06\\
0.367573696145125	2.15669188619039e-06\\
0.372108843537415	2.30357217867461e-06\\
0.376643990929705	2.37442944073241e-06\\
0.381179138321995	2.35099967538233e-06\\
0.385714285714286	2.25860075375728e-06\\
0.390249433106576	2.15011671833576e-06\\
0.394784580498866	2.09106141646487e-06\\
0.399319727891156	2.12071402316827e-06\\
0.403854875283447	2.24530766856912e-06\\
0.408390022675737	2.41244456422596e-06\\
0.412925170068027	2.56375562335265e-06\\
0.417460317460318	2.66419266459206e-06\\
0.421995464852608	2.69632761518146e-06\\
0.426530612244898	2.67931545738787e-06\\
0.431065759637188	2.65088675521774e-06\\
0.435600907029479	2.65326686776923e-06\\
0.440136054421769	2.72151083435812e-06\\
0.444671201814059	2.85348152007092e-06\\
0.449206349206349	3.01379661704626e-06\\
0.453741496598639	3.1651754687223e-06\\
0.45827664399093	3.28568196317214e-06\\
0.46281179138322	3.3608126174256e-06\\
0.46734693877551	3.40409588484025e-06\\
0.4718820861678	3.44031751664531e-06\\
0.476417233560091	3.49443338665267e-06\\
0.480952380952381	3.58251772902829e-06\\
0.485487528344671	3.70129267057706e-06\\
0.490022675736961	3.83330569009299e-06\\
0.494557823129252	3.95496389848926e-06\\
0.499092970521542	4.04982853230429e-06\\
0.503628117913832	4.11513072506864e-06\\
0.512698412698413	4.21641330838884e-06\\
0.517233560090703	4.28518757056295e-06\\
0.521768707482993	4.37449013057645e-06\\
0.530839002267574	4.57595208391837e-06\\
0.535374149659864	4.65087178014301e-06\\
0.539909297052154	4.70223831638972e-06\\
0.544444444444444	4.74375182413223e-06\\
0.548979591836735	4.79634552186248e-06\\
0.553514739229025	4.87401663096243e-06\\
0.558049886621315	4.96506042180211e-06\\
0.562585034013605	5.0441568218984e-06\\
0.567120181405896	5.08689583289534e-06\\
0.571655328798186	5.08932867160361e-06\\
0.576190476190476	5.07919483661023e-06\\
0.580725623582766	5.09199776144165e-06\\
0.585260770975057	5.15271084833291e-06\\
0.589795918367347	5.26195556060927e-06\\
0.594331065759637	5.38626796275921e-06\\
0.598866213151927	5.47583751309144e-06\\
0.603401360544218	5.50112809960535e-06\\
0.607936507936508	5.45916821104164e-06\\
0.612471655328798	5.38075658074177e-06\\
0.617006802721088	5.32794259111569e-06\\
0.621541950113379	5.35648728028359e-06\\
0.626077097505669	5.48424412250004e-06\\
0.635147392290249	5.87746270497291e-06\\
0.63968253968254	5.98224556302984e-06\\
0.64421768707483	5.94404807630067e-06\\
0.64875283446712	5.76882847291849e-06\\
0.65328798185941	5.51685138250324e-06\\
0.657823129251701	5.31732793607187e-06\\
0.662358276643991	5.30931952191782e-06\\
0.666893424036281	5.53873318398512e-06\\
0.675963718820862	6.31261560579688e-06\\
0.680498866213152	6.54311440095867e-06\\
0.685034013605442	6.51681510477342e-06\\
0.689569160997732	6.23287968459785e-06\\
0.694104308390023	5.78879770585239e-06\\
0.698639455782313	5.37846891868199e-06\\
0.703174603174603	5.23897749393765e-06\\
0.707709750566893	5.48070310313342e-06\\
0.712244897959184	5.98553500108512e-06\\
0.716780045351474	6.52373173248844e-06\\
0.721315192743764	6.88475839039082e-06\\
0.725850340136054	6.93490166470134e-06\\
0.730385487528345	6.64428004537587e-06\\
0.734920634920635	6.09200108014907e-06\\
0.739455782312925	5.47441234755031e-06\\
0.743990929705215	5.09105337298532e-06\\
0.748526077097506	5.18738606913871e-06\\
0.753061224489796	5.71508814661836e-06\\
0.757596371882086	6.39559148962743e-06\\
0.762131519274376	6.94480504248585e-06\\
0.766666666666667	7.17258127325149e-06\\
0.771201814058957	7.00313903045924e-06\\
0.775736961451247	6.48203028319472e-06\\
0.780272108843537	5.77417825897708e-06\\
0.784807256235828	5.17067183014586e-06\\
0.789342403628118	5.0155370695526e-06\\
0.793877551020408	5.41764562678182e-06\\
0.798412698412698	6.13264446069905e-06\\
0.802947845804989	6.82397689555359e-06\\
0.807482993197279	7.24076158431508e-06\\
0.812018140589569	7.246561560037e-06\\
0.816553287981859	6.83259572287209e-06\\
0.82108843537415	6.11601044697707e-06\\
0.82562358276644	5.35233693410575e-06\\
0.83015873015873	4.92015795805489e-06\\
0.83469387755102	5.10439745160518e-06\\
0.839229024943311	5.77767716902766e-06\\
0.843764172335601	6.56650423969118e-06\\
0.848299319727891	7.1597376153365e-06\\
0.852834467120181	7.37502855885364e-06\\
0.857369614512472	7.15920609362097e-06\\
0.861904761904762	6.58015568932857e-06\\
0.866439909297052	5.8203018874492e-06\\
0.870975056689342	5.18336182042667e-06\\
0.875510204081633	5.02350908160309e-06\\
0.880045351473923	5.4485057687792e-06\\
0.884580498866213	6.19858704642429e-06\\
0.889115646258503	6.92649209843577e-06\\
0.893650793650794	7.37829045860394e-06\\
0.898185941043084	7.42071726644866e-06\\
0.902721088435374	7.04719510868923e-06\\
0.907256235827664	6.36676233078948e-06\\
0.911791383219955	5.60603576482688e-06\\
0.916326530612245	5.11779377276866e-06\\
0.920861678004535	5.21000248221437e-06\\
0.925396825396825	5.82164231122118e-06\\
0.929931972789116	6.61428678228049e-06\\
0.934467120181406	7.26517361082646e-06\\
0.939002267573696	7.56601876195084e-06\\
0.943537414965986	7.43461727625849e-06\\
0.948072562358277	6.91613437342475e-06\\
0.952607709750567	6.16967639244237e-06\\
0.957142857142857	5.47362555702868e-06\\
0.961678004535147	5.19458560199748e-06\\
0.966213151927438	5.52504540118903e-06\\
0.970748299319728	6.26377231316422e-06\\
0.975283446712018	7.0464452207597e-06\\
0.979818594104308	7.5868284477032e-06\\
0.984353741496599	7.72100690193955e-06\\
0.988888888888889	7.41424434003157e-06\\
0.993424036281179	6.75665231297806e-06\\
0.997959183673469	5.96377649286772e-06\\
};
\addlegendentry{ABC$_W^{1/2}$ adaptive}

\addplot [color=mycolor1, dashed]
  table[row sep=crcr]{%
0.000226757369614528	0\\
0.00929705215419496	2.77857203911935e-09\\
0.0274376417233561	1.45164456011315e-08\\
0.0319727891156463	1.4862199582133e-08\\
0.0365079365079365	2.07625842074677e-08\\
0.045578231292517	4.98861908448944e-08\\
0.0501133786848073	6.44927802273543e-08\\
0.0546485260770975	7.47726257577241e-08\\
0.0591836734693878	7.85579905571154e-08\\
0.0682539682539682	7.77235572568813e-08\\
0.0727891156462586	9.42000372150531e-08\\
0.0773242630385488	1.28235129071008e-07\\
0.0863945578231292	2.1032581742908e-07\\
0.0909297052154195	2.38932429175165e-07\\
0.0954648526077098	2.50722057337782e-07\\
0.10453514739229	2.38312900746607e-07\\
0.109070294784581	2.48327996277276e-07\\
0.113605442176871	2.94605337569642e-07\\
0.118140589569161	3.70081384870602e-07\\
0.122675736961451	4.53218788498688e-07\\
0.127210884353742	5.24059345585748e-07\\
0.131746031746032	5.68468046524551e-07\\
0.136281179138322	5.80551164586218e-07\\
0.145351473922902	5.49135144800594e-07\\
0.149886621315193	5.64248462064931e-07\\
0.154421768707483	6.38197102942328e-07\\
0.158956916099773	7.6105095458967e-07\\
0.163492063492063	8.98726536191496e-07\\
0.168027210884354	1.01636340332778e-06\\
0.172562358276644	1.08876842463967e-06\\
0.177097505668934	1.1006133716629e-06\\
0.186167800453515	1.01033811883422e-06\\
0.190702947845805	1.04483511309539e-06\\
0.195238095238095	1.16830400809231e-06\\
0.199773242630386	1.3180132405477e-06\\
0.204308390022676	1.45459500844325e-06\\
0.208843537414966	1.5492133235151e-06\\
0.213378684807256	1.57283085266968e-06\\
0.217913832199547	1.52578984324858e-06\\
0.222448979591837	1.42591275908988e-06\\
0.226984126984127	1.3384842617592e-06\\
0.231519274376417	1.33881434161154e-06\\
0.236054421768708	1.43392288343147e-06\\
0.240589569160998	1.57287582758237e-06\\
0.245124716553288	1.7049898418664e-06\\
0.249659863945578	1.79321125415299e-06\\
0.254195011337868	1.80461463505921e-06\\
0.258730158730159	1.73551493432011e-06\\
0.267800453514739	1.49398410909107e-06\\
0.272335600907029	1.44975070293096e-06\\
0.27687074829932	1.52391825292142e-06\\
0.28140589569161	1.67956542229675e-06\\
0.2859410430839	1.84332363650785e-06\\
0.290476190476191	1.95561435989422e-06\\
0.295011337868481	1.9729557868553e-06\\
0.299546485260771	1.89225519253711e-06\\
0.308616780045352	1.60552634120847e-06\\
0.313151927437642	1.54052946910177e-06\\
0.317687074829932	1.60618928612966e-06\\
0.322222222222222	1.77276593238407e-06\\
0.326757369614513	1.96295817611603e-06\\
0.331292517006803	2.09934836548076e-06\\
0.335827664399093	2.14057555136815e-06\\
0.340362811791383	2.08435446502619e-06\\
0.344897959183673	1.96095840510058e-06\\
0.349433106575964	1.83198281411556e-06\\
0.353968253968254	1.77200418738366e-06\\
0.358503401360544	1.8276067231815e-06\\
0.363038548752834	1.97792128320273e-06\\
0.367573696145125	2.15669188619039e-06\\
0.372108843537415	2.30357217867461e-06\\
0.376643990929705	2.37442996275927e-06\\
0.381179138321995	2.35093736355996e-06\\
0.385714285714286	2.25790645991264e-06\\
0.390249433106576	2.14741102733118e-06\\
0.394784580498866	2.08450051930598e-06\\
0.399319727891156	2.10922485754939e-06\\
0.403854875283447	2.22971728469634e-06\\
0.408390022675737	2.39530041212799e-06\\
0.412925170068027	2.54551341793974e-06\\
0.417460317460318	2.64432113505908e-06\\
0.421995464852608	2.67313949575687e-06\\
0.426530612244898	2.64938117722213e-06\\
0.431065759637188	2.60870323887108e-06\\
0.435600907029479	2.59406945946061e-06\\
0.440136054421769	2.64431209351379e-06\\
0.444671201814059	2.76295739420895e-06\\
0.453741496598639	3.07279136368965e-06\\
0.45827664399093	3.20262700803653e-06\\
0.46281179138322	3.28589069786833e-06\\
0.46734693877551	3.32981507777941e-06\\
0.4718820861678	3.35612004498387e-06\\
0.476417233560091	3.39468567489565e-06\\
0.480952380952381	3.46926586658469e-06\\
0.485487528344671	3.58582455539658e-06\\
0.490022675736961	3.73363150985018e-06\\
0.494557823129252	3.88609956625263e-06\\
0.499092970521542	4.02055615189578e-06\\
0.503628117913832	4.12471661370706e-06\\
0.508163265306122	4.20271591350385e-06\\
0.512698412698413	4.271535396283e-06\\
0.517233560090703	4.34913307567619e-06\\
0.521768707482993	4.44676752142659e-06\\
0.526303854875283	4.5657080797179e-06\\
0.530839002267574	4.69536953318084e-06\\
0.535374149659864	4.81376392402755e-06\\
0.539909297052154	4.91047640194253e-06\\
0.544444444444444	4.98739877585752e-06\\
0.548979591836735	5.0560086564122e-06\\
0.553514739229025	5.1399255860396e-06\\
0.558049886621315	5.2498466628359e-06\\
0.562585034013605	5.37817976087318e-06\\
0.567120181405896	5.49793950532784e-06\\
0.571655328798186	5.57508083198766e-06\\
0.576190476190476	5.59145302614361e-06\\
0.585260770975057	5.53148749515842e-06\\
0.589795918367347	5.57209411911064e-06\\
0.594331065759637	5.72646247543052e-06\\
0.598866213151927	5.97752586251055e-06\\
0.603401360544218	6.24712991592613e-06\\
0.607936507936508	6.43271187605876e-06\\
0.612471655328798	6.45293751322384e-06\\
0.617006802721088	6.28161629756541e-06\\
0.621541950113379	5.97434649063366e-06\\
0.626077097505669	5.68005078394052e-06\\
0.630612244897959	5.60134745919427e-06\\
0.635147392290249	5.86467484298137e-06\\
0.63968253968254	6.39756737641495e-06\\
0.64421768707483	6.98096113593838e-06\\
0.64875283446712	7.38920389387054e-06\\
0.65328798185941	7.46343130453297e-06\\
0.657823129251701	7.15005730178575e-06\\
0.662358276643991	6.51949645125161e-06\\
0.666893424036281	5.79723657512243e-06\\
0.671428571428571	5.36271793316612e-06\\
0.675963718820862	5.55608313568268e-06\\
0.680498866213152	6.32590393312427e-06\\
0.685034013605442	7.28679802342036e-06\\
0.689569160997732	8.06117405949269e-06\\
0.694104308390023	8.40346330399644e-06\\
0.698639455782313	8.20176240534654e-06\\
0.703174603174603	7.48019038621095e-06\\
0.707709750566893	6.42908527348496e-06\\
0.712244897959184	5.46618734698079e-06\\
0.716780045351474	5.20596175523913e-06\\
0.721315192743764	5.90955420032291e-06\\
0.725850340136054	7.11332834257838e-06\\
0.730385487528345	8.22796223043376e-06\\
0.734920634920635	8.8814386891567e-06\\
0.739455782312925	8.89930282244045e-06\\
0.743990929705215	8.26572833267925e-06\\
0.748526077097506	7.12568192906549e-06\\
0.753061224489796	5.83996620773242e-06\\
0.757596371882086	5.06953074974614e-06\\
0.762131519274376	5.43411861353427e-06\\
0.766666666666667	6.66138683447226e-06\\
0.771201814058957	8.01539188133304e-06\\
0.775736961451247	8.9949536036249e-06\\
0.780272108843537	9.34634882732155e-06\\
0.784807256235828	8.9889163462642e-06\\
0.789342403628118	7.99874272161549e-06\\
0.793877551020408	6.63670235012681e-06\\
0.798412698412698	5.43847330081171e-06\\
0.802947845804989	5.18988326236869e-06\\
0.807482993197279	6.11827026519851e-06\\
0.812018140589569	7.52997561015434e-06\\
0.816553287981859	8.74988381061037e-06\\
0.82108843537415	9.41032037471423e-06\\
0.82562358276644	9.36057266742729e-06\\
0.83015873015873	8.61161728327176e-06\\
0.83469387755102	7.34179328476525e-06\\
0.839229024943311	5.95729141594781e-06\\
0.843764172335601	5.17771883246088e-06\\
0.848299319727891	5.62715348839671e-06\\
0.852834467120181	6.94064699990893e-06\\
0.857369614512472	8.34010710326627e-06\\
0.861904761904762	9.32174851353551e-06\\
0.866439909297052	9.63670871345634e-06\\
0.870975056689342	9.21174832613225e-06\\
0.875510204081633	8.1342922059191e-06\\
0.880045351473923	6.69022569699251e-06\\
0.884580498866213	5.4671765077341e-06\\
0.889115646258503	5.29989414876919e-06\\
0.893650793650794	6.34757122053298e-06\\
0.898185941043084	7.83872560727872e-06\\
0.902721088435374	9.08963259760309e-06\\
0.907256235827664	9.73684108118711e-06\\
0.911791383219955	9.63259082775725e-06\\
0.916326530612245	8.7984557074483e-06\\
0.920861678004535	7.43415110915358e-06\\
0.925396825396825	5.98911974514227e-06\\
0.929931972789116	5.26062537598015e-06\\
0.934467120181406	5.85947020825817e-06\\
0.939002267573696	7.29201909144273e-06\\
0.943537414965986	8.73593621886215e-06\\
0.948072562358277	9.69622512436707e-06\\
0.952607709750567	9.93933600856867e-06\\
0.957142857142857	9.41093949857752e-06\\
0.961678004535147	8.22661682253756e-06\\
0.966213151927438	6.71839992849588e-06\\
0.970748299319728	5.54558450838183e-06\\
0.975283446712018	5.5564440180067e-06\\
0.979818594104308	6.75778800818261e-06\\
0.984353741496599	8.30852182298614e-06\\
0.988888888888889	9.54487870441501e-06\\
0.993424036281179	1.01250195326807e-05\\
0.997959183673469	9.92186406567619e-06\\
};
\addlegendentry{ABC$_W^{1/2}$}

\end{axis}
\end{tikzpicture}%